\documentclass[reqno]{amsbook}

\usepackage{fixltx2e}

\usepackage{graphicx}
\usepackage{amsmath}
\usepackage{amsfonts}
\usepackage{amssymb}
\usepackage{amsthm}
\usepackage{mathrsfs}

\usepackage[all]{xy}

\input diagxy

\newcommand{\pfbreak}{\medskip\begin{center}
*\qquad *\qquad *
\end{center}}

\newcommand{\slashfrac}[2]{{#1}/{#2}}

\newcommand{\fref}[1]{Fig.~\ref{#1}}

\newenvironment{morale}{\noindent\textit{Sketch of proof.}}{}

\numberwithin{section}{chapter}
\numberwithin{equation}{chapter}
\theoremstyle{plain}

\newtheorem{theorem}{Theorem}[section]
\newtheorem*{thm}{Theorem}
\newtheorem{lemma}[theorem]{Lemma}
\newtheorem{proposition}[theorem]{Proposition}
\newtheorem{corollary}[theorem]{Corollary}

\theoremstyle{definition}

\newtheorem{definition}[theorem]{Definition}
\newtheorem{notation}[theorem]{Notation}
\newtheorem{assumption}[theorem]{Assumption}
\newtheorem{remark}[theorem]{Remark}
\newtheorem{example}[theorem]{Example}
\newtheorem{examples}[theorem]{Examples}

\makeindex

\renewcommand{\labelenumi}{\arabic{enumi}.}
\renewcommand{\emph}[1]{\textsc{#1}}
\newcommand{\head}[1]{\medskip\noindent \textbf{#1.}}
\newcommand{\headd}[1]{\noindent \textbf{#1.}}
\DeclareSymbolFont{Euler}{U}{eur}{m}{n}
\DeclareSymbolFontAlphabet\matheu{Euler}

\newcommand{\versus}{vs.}

\newcommand{\noproof}{\qed}
\newcommand{\Ob}[1]{#1}

\newcommand{\monr}{\ar@{{^ (}->}[r]}
\newcommand{\monl}{\ar@{{^ (}->}[l]}
\newcommand{\monu}{\ar@{{^ (}->}[u]}
\newcommand{\mond}{\ar@{{^ (}->}[d]}
\newcommand{\monrr}{\ar@{{^ (}->}[rr]}
\newcommand{\monll}{\ar@{{^ (}->}[ll]}
\newcommand{\monuu}{\ar@{{^ (}->}[uu]}
\newcommand{\mondd}{\ar@{{^ (}->}[dd]}
\newcommand{\mondr}{\ar@{{^ (}->}[dr]}
\newcommand{\monur}{\ar@{{^ (}->}[ur]}
\newcommand{\vierkant}[4]{\vcenter{\xymatrix{{\cdot} \ar[r]^-{{#3}} \ar[d]_-{{#1}} & {\cdot} \ar[d]^-{{#4}}\\
{\cdot} \ar[r]_-{{#2}} & {\cdot}}}}

\newbox\pullbackbox
\setbox\pullbackbox=\hbox{\xy 0;<1mm,0mm>: \POS(4,0)\ar@{-} (0,0) \ar@{-} (4,4)
\endxy}
\def\pullback{\copy\pullbackbox}

\newbox\pushoutbox
\setbox\pushoutbox=\hbox{\xy 0;<1mm,0mm>: \POS(0,4)\ar@{-} (0,0) \ar@{-} (4,4)
\endxy}
\def\pushout{\copy\pushoutbox}

\newbox\skewpullbackbox
\setbox\skewpullbackbox=\hbox{\xy 0;<1mm,0mm>: \POS(4,0)\ar@{-} (0,4) \ar@{-} (4,4)
\endxy}
\def\skewpullback{\copy\skewpullbackbox}

\def\relto{\mathrel{\xy\ar|-{\circ} (200,0) \endxy}}
\def\regepito{\mathrel{\xy\ar@{->>} (200,0) \endxy}}
\def\rightharpoonup{\mathrel{\xy\ar@{-{ ^}} (200,0) \endxy}}
\def\cokernelto{\mathrel{\xy\ar@{-{ >>}} (200,0) \endxy}}
\def\monoto{\mathrel{\xy\ar@{{^ (}->} (200,0) \endxy}}
\def\kernelto{\mathrel{\xy\ar@{{ >}->} (200,0) \endxy}}
\def\To{\mathrel{\xy \ar@{=>} (200,0) \endxy}}
\def\comp{\raisebox{0.2mm}{\ensuremath{\scriptstyle\circ}}}
\def\vcomp{\raisebox{0.2mm}{\ensuremath{\scriptstyle\bullet}}}
\renewcommand{\mapsto}{\mathrel{\xy \ar@{|->} (200,0) \endxy}}
\newcommand{\del}{\partial}
\newcommand{\regepi}{\ensuremath{\mathsf{RegEpi}}}

\newcommand{\tensor}{\otimes}
\newcommand{\im}{\ensuremath{\mathsf{im\,}}}

\newcommand{\coker}{\ensuremath{\mathsf{coker\,}}}
\renewcommand{\ker}{\ensuremath{\mathsf{ker\,}}}
\newcommand{\Coker}{\ensuremath{\mathsf{Coker\,}}}
\newcommand{\Cok}{\ensuremath{Q}}
\newcommand{\Ker}{\ensuremath{\mathsf{Ker\,}}}
\newcommand{\coeq}{\ensuremath{\mathsf{Coeq\,}}}
\newcommand{\Sub}{\ensuremath{\mathsf{Sub}}}
\newcommand{\NSub}{\ensuremath{\mathsf{PSub}}}

\newcommand{\Eq}{\ensuremath{\mathsf{Eq}}}
\newcommand{\KP}{\ensuremath{\mathit{R}}}
\newcommand{\fib}{\ensuremath{\mathsf{fib}}}
\newcommand{\cof}{\ensuremath{\mathsf{cof}}}
\newcommand{\we}{\ensuremath{\mathsf{we}}}
\newcommand{\iso}{\ensuremath{\mathrm{iso}}}
\newcommand{\ner}{\ensuremath{\mathrm{ner}}}
\newcommand{\pr}{\ensuremath{\mathrm{pr}}}
\newcommand{\op}{\ensuremath{\mathrm{op}}}
\newcommand{\lali}{\ensuremath{\mathrm{lali}\,}}
\newcommand{\rari}{\ensuremath{\mathrm{rari}\,}}
\newcommand{\inj}{\ensuremath{\mathrm{in}}}
\newcommand{\id}{\ensuremath{\mathrm{id}}}
\newcommand{\tw}{\ensuremath{\mathrm{tw}}}
\newcommand{\dom}{\ensuremath{\mathrm{dom}}}
\newcommand{\cod}{\ensuremath{\mathrm{cod}}}
\newcommand{\compo}{\ensuremath{\mathrm{comp}}}
\newcommand{\ab}{\mathrm{Ab}}
\newcommand{\centr}{\ensuremath{\mathrm{Centr}}}
\newcommand{\ext}{\ensuremath{\mathrm{Ext}}}
\renewcommand{\hom}{\mathrm{Hom}}

\newcommand{\Ac}{\ensuremath{\mathcal{A}}}
\newcommand{\Bc}{\ensuremath{\mathcal{B}}}
\newcommand{\Cc}{\ensuremath{\mathcal{C}}}
\newcommand{\Dc}{\ensuremath{\mathcal{D}}}
\newcommand{\Ec}{\ensuremath{\mathcal{E}}}
\newcommand{\Fc}{\ensuremath{\mathcal{F}}}

\newcommand{\Hc}{\ensuremath{\mathcal{H}}}
\newcommand{\Ic}{\ensuremath{\mathcal{I}}}
\newcommand{\Lc}{\ensuremath{\mathcal{L}}}
\newcommand{\Pc}{\ensuremath{\mathcal{P}}}
\newcommand{\Rc}{\ensuremath{\mathcal{R}}}
\newcommand{\Sc}{\ensuremath{\mathcal{S}}}
\newcommand{\Tc}{\ensuremath{\mathcal{T}}}
\newcommand{\Uc}{\ensuremath{\mathcal{U}}}
\newcommand{\Vc}{\ensuremath{\mathcal{V}}}
\newcommand{\Wc}{\ensuremath{\mathcal{W}}}

\newcommand{\Wproj}{\ensuremath{\Wc_{\mathsf{proj}}}}
\newcommand{\Wsplit}{\ensuremath{\Wc_{\mathsf{split}}}}

\newcommand{\Ib}{\ensuremath{\mathbf{I}}}

\newcommand{\Kb}{\ensuremath{\mathbf{K}}}

\newcommand{\Ab}{\mathsf{Ab}}
\newcommand{\Alg}{\mathsf{Alg}}
\newcommand{\Act}{\mathsf{Act}}
\newcommand{\Cat}{\ensuremath{\mathsf{Cat}}}
\newcommand{\CentrExt}{\mathsf{Centr}}
\newcommand{\Ch}{\ensuremath{\mathsf{Ch}}}
\newcommand{\Chplus}{\ensuremath{\mathsf{Ch}^{+}\!}}
\newcommand{\CgHaus}{\mathsf{CGHaus}}

\newcommand{\CompHaus}{\mathsf{CompHaus}}
\newcommand{\CrossMod}{\mathsf{CrossMod}}
\newcommand{\Elts}{\mathsf{Elts}}
\newcommand{\Ext}{\mathsf{Ext}}
\newcommand{\Five}{\ensuremath{\mathsf{5}}}
\newcommand{\FinSet}{\ensuremath{\mathsf{FinSet}}}
\newcommand{\Fun}{\ensuremath{\mathsf{Fun}}}
\newcommand{\Gd}{\ensuremath{\mathsf{Grpd}}}
\newcommand{\Gp}{\ensuremath{\mathsf{Gp}}}
\newcommand{\Lie}{\ensuremath{\mathsf{Lie}}}
\newcommand{\Mag}{\ensuremath{\mathsf{Mag}}}
\newcommand{\Mod}{\ensuremath{\mathsf{Mod}}}
\newcommand{\Mon}{\ensuremath{\mathsf{Mon}}}
\newcommand{\omegaGp}{\ensuremath{\text{$\omega$-$\mathsf{Gp}$}}}
\newcommand{\PAr}{\mathsf{pAr}}
\newcommand{\PCh}{\ensuremath{\mathsf{PCh}}}
\newcommand{\PChplus}{\ensuremath{\mathsf{PCh}^{+}\!}}
\renewcommand{\Pr}{\mathsf{Pr}}
\newcommand{\prA}{\mathsf{Pr} \Ac}
\newcommand{\PXMod}{\mathsf{PreCrossMod}}
\newcommand{\pres}{\mathsf{Pr}}
\newcommand{\Pt}{\mathsf{Pt}}
\newcommand{\Psh}{\mathsf{PSh}}
\newcommand{\Rng}{\mathsf{Rng}}
\newcommand{\Rgrph}{\ensuremath{\mathsf{RG}}}
\newcommand{\Set}{\ensuremath{\mathsf{Set}}}
\newcommand{\Sh}{\mathsf{Sh}}
\newcommand{\simp}[1]{\Sc #1}
\newcommand{\simpA}{\simp{\Ac}}

\newcommand{\SplitEpi}{\mathsf{SplitEpi}}
\newcommand{\STE}{\mathsf{6tE}}
\newcommand{\Top}{\mathsf{Top}}
\newcommand{\Two}{\ensuremath{\mathsf{2}}}
\newcommand{\XMod}{\mathsf{XMod}}

\newcommand{\CCf}{\ensuremath{\mathscr{C}}}

\newcommand{\Af}{\ensuremath{\matheu{A}}}
\newcommand{\If}{\ensuremath{\matheu{I}}}

\newcommand{\Bf}{\ensuremath{\matheu{B}}}
\newcommand{\Cf}{\ensuremath{\matheu{C}}}

\newcommand{\Ef}{\ensuremath{\matheu{E}}}

\newcommand{\Xf}{\ensuremath{\matheu{X}}}

\newcommand{\Pf}{\ensuremath{\matheu{P}}}
\newcommand{\veth}{h}
\newcommand{\vetH}{H}
\newcommand{\vetp}{p}
\newcommand{\vetj}{j}
\newcommand{\vetk}{k}
\newcommand{\vetK}{K}
\newcommand{\vetL}{L}
\newcommand{\vetM}{M}
\newcommand{\vetr}{r}
\newcommand{\vets}{s}
\newcommand{\vetf}{f}
\newcommand{\vetg}{g}

\newcommand{\vetsf}{\ensuremath{\text{\textsl{f}}}}
\newcommand{\vetsg}{\ensuremath{\text{\textsl{g}}}}

\newcommand{\Lief}{\mathfrak{f}}
\newcommand{\Lieg}{\mathfrak{g}}

\newcommand{\Lien}{\mathfrak{n}}
\newcommand{\Lier}{\mathfrak{r}}

\newcommand{\G}{\ensuremath{\mathbb{G}}}
\newcommand{\K}{\ensuremath{\mathbb{K}}}
\newcommand{\N}{\ensuremath{\mathbb{N}}}
\newcommand{\Z}{\ensuremath{\mathbb{Z}}}
\newcommand{\R}{\ensuremath{\mathbb{R}}}
\newcommand{\Q}{\ensuremath{\mathbb{Q}}}
\newcommand{\T}{\ensuremath{\mathbb{T}}}

\newcommand{\cycle}{\nabla}
\newcommand{\El}{\Lambda}

\DeclareMathOperator{\Coeq}{Coeq}

\author{Tim Van der Linden}
\address{Vrije Universiteit Brussel\\
Pleinlaan 2\\
1050 Brussel\\
Belgium}
\email{tvdlinde@vub.ac.be}
\title[Homology and homotopy in semi-abelian categories]{Homology and homotopy\\
in semi-abelian categories}

\begin{document}

% XY stuff

\newdir{>>}{{}*!/4.5pt/:(1,-.2)@^{>}*!/4.5pt/:(1,+.2)@_{>}*!/8pt/:(1,-.2)@^{>}*!/8pt/:(1,+.2)@_{>}}
\newdir{ >>}{{}*!/9pt/@{|}*!/4.5pt/:(1,-.2)@^{>}*!/4.5pt/:(1,+.2)@_{>}}
\newdir{ >}{{}*!/-5.5pt/@{|}*!/-10pt/:(1,-.2)@^{>}*!/-10pt/:(1,+.2)@_{>}}
\newdir{ ^}{{}*:(1,-.2)@^{>}}
\newdir2{ >}{{}*-<1.8pt,0pt>@{ >}}
\newdir{^ (}{{}*!/-10pt/@{>}}
\newdir{_ (}{{}*!/-10pt/@{>}}
\newdir{>}{{}*:(1,-.2)@^{>}*:(1,+.2)@_{>}}
\newdir{<}{{}*:(1,+.2)@^{<}*:(1,-.2)@_{<}}
\newdir{|>}{{}*!/-5.5pt/@{|}*!<-10.5pt,-.9pt>@^{>}*!<-10.5pt,.9pt>@_{>}}
\CompileMatrices

% Text starts here

\maketitle 

\frontmatter
\chapter*{Introduction}\label{Chapter-Summary}

The theory of \textit{abelian} categories proved very useful, providing an axiomatic framework for homology and cohomology of modules over a ring and, in particular, of abelian groups~\cite{Cartan-Eilenberg, Maclane:Homology, Freyd, Weibel}. Out of it (and out of Algebraic Geometry) grew Category Theory. But for many years, a similar categorical framework has been lacking for non-abelian (co)homology, the subject of which includes the categories of groups, rings, Lie algebras etc. The point of this dissertation is that \textit{semi-abelian} categories (in the sense of Janelidze, M\'arki and Tholen~\cite{Janelidze-Marki-Tholen}) provide a suitable context for the study of non-abelian (co)homology and the corresponding homotopy theory.

A \textit{semi-abelian} category is pointed, Barr exact and Bourn protomodular with binary coproducts~\cite{Janelidze-Marki-Tholen, Borceux-Bourn}. \textit{Pointed} means that it has a zero object: an initial object that is also terminal. An \textit{exact} category is regular (finitely complete with pullback-stable regular epimorphisms) and such that any internal equivalence relation is a kernel pair~\cite{Barr-Grillet-vanOsdol}. A pointed and regular category is \textit{protomodular} when the Short Five Lemma holds: For any commutative diagram
\[
\xymatrix{K[p'] \ar@{{ >}->}[r]^-{\ker p'} \ar[d]_-u & E' \ar@{->>}[r]^-{p'} \ar[d]^-v & B' \ar[d]^-w \\ K[p] \ar@{{ >}->}[r]_-{\ker p} & E \ar@{->>}[r]_-{p} & B}
\]
such that $p$ and $p'$ are regular epimorphisms, $u$ and $w$ being isomorphisms implies that $v$ is an isomorphism~\cite{Bourn1991}. Examples include all abelian categories; all varieties of $\Omega$-groups (i.e., varieties of universal algebras with a unique constant and an underlying group structure~\cite{Higgins}), in particular the categories of groups, non-unitary rings, Lie algebras, commutative algebras~\cite{Borceux-Bourn}, crossed modules and precrossed modules~\cite{Janelidze}; internal versions of such varieties in an exact category~\cite{Janelidze-Marki-Tholen}; Heyting semilattices~\cite{Jo}; compact Hausdorff (profinite) groups (or semi-abelian algebras)~\cite{Borceux-Clementino, Borceux-Bourn}, non-unital $C^{*}$ algebras~\cite{Gran-Rosicky:Monadic} and Heyting algebras~\cite{Bourn1996, Janelidze-Marki-Tholen}; the dual of the category of pointed objects in any topos, in particular the dual of the category of pointed sets~\cite{Janelidze-Marki-Tholen}.

In any semi-abelian category, the fundamental diagram lemmas, such as the Short Five Lemma, the $3 \times 3$ Lemma, the Snake Lemma~\cite{Bourn2001} and Noether's Isomorphism Theorems~\cite{Borceux-Bourn} hold. Moreover, there is a satisfactory notion of homology of chain complexes: Any short exact sequence of proper chain complexes induces a long exact homology sequence~\cite{EverVdL2}.

It is hardly surprising that such a context would prove suitable for further development of non-abelian (co)homology. With this dissertation we provide some additional evidence. We consider homology of simplicial objects and internal categories, and cotriple (co)homology in the sense of Barr and Beck~\cite{Barr-Beck}. Using techniques from commutator theory~\cite{Smith, Huq, Pedicchio, Janelidze-Kelly, Bourn2002, Bourn-Gran, Bourn-Gran-Maltsev} and the theory of Baer invariants~\cite{Froehlich, Lue, Furtado-Coelho, EverVdL1}, we obtain a general version of Hopf's Formula~\cite{Hopf} and the Stallings-Stammbach Sequence~\cite{Stallings, Stammbach} in homology---and their cohomological counterparts, in particular the Hochschild-Serre sequence~\cite{Gran-VdL}. We recover results, well-known in the case of groups and Lie algebras: e.g., the fact that the second cohomology group classifies central extensions. And although homotopy theory of chain complexes seems to be problematic, Quillen model category structures for simplicial objects and internal categories exist that are compatible with these notions of homology~\cite{EKVdL, VdLinden:Simp}.\pfbreak

\headd{Chapter 1} Throughout the text, we shall attempt to formulate all results in their ``optimal'' categorical context. Sometimes this optimal context is the most general one; sometimes it is where most technical problems disappear. (It is, however, \textit{always} possible to replace the conditions on a category $\Ac$ by ``semi-abelian with enough projectives'', unless it is explicitly mentioned that we need ``semi-abelian and monadic over $\Set$'', and most of the time ``semi-abelian'' is sufficient.)

In the first chapter we give an overview of the relevant categorical structures: \textit{quasi-pointed} and \textit{pointed}, \textit{unital} and \textit{strongly unital}, \textit{regular} and \textit{Barr exact}, \textit{Mal'tsev}, \textit{Bourn protomodular}, \textit{sequentiable} and \textit{homological}, and \textit{semi-abelian} categories. This is what we shall loosely refer to as the ``semi-abelian'' context. We give examples and counter\-examples, fix notations and conventions. We add those (sometimes rather technical) properties that are essential to making the theory work.

\head{Chapter 2} In Chapter~\ref{Chapter-Homology} we extend the notion of chain complex to the semi-abelian context. As soon as the ambient category is quasi-pointed, exact and protomodular, homology of \textit{proper} chain complexes---those with boundary operators of which the image is a kernel---is well-behaved: It characterizes exactness. Also, in this case, the two dual definitions of homology (see Definition~\ref{Definition-Homology}) coincide. A short exact sequence of proper chain complexes induces a long exact homology sequence. 

We show that the \textit{Moore complex} of a simplicial object $A$ in a quasi-pointed exact Mal'tsev category $\Ac$ is always proper; this is the positively graded chain complex $N (A)\in \PChplus\Ac$ defined by $N_{0}A=A_{0}$,
\[
N_{n} A=\bigcap_{i=0}^{n-1}K[\del_{i}\colon A_{n}\to A_{n-1}]
\]
and with boundary operators $d_{n}=\del_{n}\comp \bigcap_{i}\ker \del_{i}\colon N_{n} A\to N_{n-1} A$, for $n\geq 1$. When, moreover, $\Ac$ is protomodular, the \textit{normalization functor} 
\[
N\colon \Sc \Ac\to\PChplus \Ac
\]
is exact; hence a short exact sequence of simplicial objects induces a long exact homology sequence. To show this, we use the result that in a regular Mal'tsev category, every simplicial object is Kan~\cite{Carboni-Kelly-Pedicchio}, and every regular epimorphism of simplicial objects is a Kan fibration. These facts also lead to a proof of Dominique Bourn's conjecture that for $n\geq 1$, the homology objects $H_{n}A$ of a simplicial object $A$ are abelian.

\head{Chapter 3} Extending the work of Fr\"ohlich, Lue and Furtado-Coelho, in Chapter~\ref{Chapter-Baer-Invariants}, we study the theory of Baer invariants in the semi-abelian context. Briefly, this theory concerns expressions in terms of a presentation of an object that are independent of the chosen presentation. For instance, presenting a group $G$ as a quotient $F/R$ of a free group $F$ and a ``group of relations'' $R$,
\[
\frac{[F,F]}{[R,F]}  \qquad \text{and}\qquad  \frac{R\cap [F,F]}{[R,F]}
\]
are such expressions.

Formally, a \textit{presentation} $p\colon {A_{0}\regepito A}$ of an object $A$ in a category $\Ac$ is a regular epimorphism. (In a sequentiable category, the kernel $\ker p\colon A_{1}\kernelto A_{0}$ of $p$ exists, and then $A$ is equal to the quotient $A_{0}/A_{1}$.) The category $\Pr \Ac$ is the full subcategory, determined by the presentations in $\Ac$, of the category $\Fun (\Two ,\Ac)$ of arrows in $\Ac$. A \textit{Baer invariant} is a functor $B\colon {\Pr\Ac \to \Ac}$ that makes \textit{homotopic} morphisms of presentations equal: Such are $(f_{0},f)$ and $(g_{0},g)\colon {p\to q}$ 
\[
\xymatrix{A_{0} \ar@{->>}[d]_-{p} \ar@<-0.5 ex>[r]_-{g_{0}}\ar@<0.5 ex>[r]^-{f_{0}} & B_{0} \ar@{->>}[d]^-{q}\\
A \ar@<-0.5 ex>[r]_-{g}\ar@<0.5 ex>[r]^-{f} & B}
\]
satisfying $f=g$.

Baer invariants may be constructed from subfunctors of the identity functor of $\Ac$. For instance, under the right circumstances, a subfunctor $V$ of $1_{\Ac}$ induces functors $\Pr \Ac \to \Ac$ that map a presentation $p\colon A_{0}\to A$ to 
\[
\frac{K[p]\cap V_{0}p}{V_{1}p} \qquad \text{or}\qquad \frac{V_{0}p}{K[p]\cap V_{0}p}.
\]
Choosing, for every object $A$, a \textit{projective} presentation $p$ (i.e., with $A_{0}$ projective), these Baer invariants induce functors $\Ac\to\Ac$, which are respectively denoted by $\Delta V$ and $\nabla V$. 

Ideally, $\Ac$ is an exact and sequentiable category and $V$ is a proper subfunctor of $1_{\Ac}$ that preserves regular epimorphisms: i.e., $V$ is a \textit{Birkhoff subfunctor} of $\Ac$. Birkhoff subfunctors correspond bijectively to \textit{Birkhoff subcategories} in the sense of Janelidze and Kelly~\cite{Janelidze-Kelly}: full and reflective, and closed in $\Ac$ under subobjects and regular quotients. Given a Birkhoff subcategory $\Bc$, the functor $V$ is the kernel of the unit of the adjunction. As a leading example, there is always the subcategory $\Ab \Ac$ of all abelian objects of $\Ac$. In this case, one interprets a value $V (A)$ of $V$ as a commutator $[A,A]$.

Baer invariants give rise to exact sequences. For instance, Theorem~\ref{Theorem-V-Exact-Sequence} states that, given an exact and sequentiable category with enough projectives $\Ac$, and a Birkhoff subfunctor $V$ of $\Ac$, any short exact sequence
\[
\xymatrix{0 \ar[r] & K \ar@{{ >}->}[r] & A \ar@{-{ >>}}[r]^-{f} & B \ar[r] & 1}
\]
in $\Ac$ induces, in a natural way, an exact sequence
\[
%\resizebox{\textwidth}{!}{
\xymatrix{\Delta V(A)  \ar[r]^-{\Delta Vf} & \Delta V(B)  \ar[r] & \frac{K}{V_{1}f}  \ar[r] & U(A) \ar@{-{ >>}}[r]^-{Uf} & U(B)  \ar[r] & 1.}
\]
The symbol $V_{1}f$ that occurs in this sequence can be interpreted as a commutator $[K,A]$ of $K$ with $A$, \textit{relative to $V$}. We relate it to existing notions in the field of categorical commutator theory and the theory of central extensions in the following manner. In the case of abelianization, the commutator $V_{1}$ may be described entirely in terms of Smith's commutator of equivalence relations~\cite{Smith}, and the corresponding notion of central extension coincides with the ones considered by Smith~\cite{Smith} and Huq~\cite{Huq}. Relative to an arbitrary Birkhoff subcategory---not just the full subcategory of all abelian objects---the central extensions are the same as those introduced by Janelidze and Kelly~\cite{Janelidze-Kelly}.

Finally we propose a notion of \textit{nilpotency}, relative to any Birkhoff subcategory of a semi-abelian category. An object is nilpotent if and only if is lower central series reaches zero (Corollary~\ref{Corollary-Nilpotent-iff-LCC-zero}), and the objects of a given nilpotency class form a Birkhoff subcategory (Proposition~\ref{Proposition-Nilpotency-Class}).

\head{Chapter 4} In this chapter we establish a connection between the theory of Baer invariants from Chapter~\ref{Chapter-Baer-Invariants} and the homology theory from Chapter~\ref{Chapter-Homology}. We do so by considering the problem of deriving the reflector of a semi-abelian category $\Ac$ onto a Birkhoff subcategory $\Bc$ of $\Ac$.

We extend Barr and Beck's definition of cotriple homology to the semi-abelian context as follows: Let $\G$ be a comonad on a category $\Ac$; let $\Bc$ be an exact and sequentiable category and $U\colon \Ac \to \Bc$ a functor. We say that the object $H_{n} (X,U)_{\G}=H_{n-1} NU (\G X)$ is the \textit{$n$-th homology object of $X$ with coefficients in $U$ relative to the comonad $\G$}. (Here ${\G X\to X}$ denotes the canonical simplicial resolution of $X$, induced by $\G$.) This defines a functor $H_{n} (\cdot,U)_{\G}\colon {\Ac \to \Bc}$, for any $n\in \N_{0}$.

We shall be interested mainly in the situation where $U$ is the reflector of an exact and sequentiable category $\Ac$ onto a Birkhoff subcategory $\Bc$ of $\Ac$, and $\G$ is a regular comonad (i.e., its counit induces projective presentations, and the functor $G$ preserves regular epis). (As an important special case, we may take $\Ac$ any semi-abelian category, monadic over $\Set$, $\G$ the induced comonad and $\Bc=\Ab \Ac$.)

Basing ourselves on Carrasco, Cegarra and Grandje\'an's homology theory of crossed modules~\cite{Carrasco-Homology}, we show that, in general, the lower homology objects admit an interpretation in terms of commutators (Baer invariants). In particular, this gives rise to a semi-abelian version of Hopf's Formula, a result of which well-known versions for groups and Lie algebras exist. We also obtain a version of the Stallings-Stammbach Sequence~\cite{Stallings, Stammbach}: Theorem~\ref{Theorem-Cotriple-Exact-Sequence} states that any short exact sequence
\[
\xymatrix{ 0 \ar[r] & K \ar@{{ >}->}[r] & A \ar@{-{ >>}}[r]^-{f} & B \ar[r] & 1 }
\]
in $\Ac$, in a natural way, induces an exact sequence
\[
\resizebox{\textwidth}{!}{\mbox{$H_2(A,U)_{\G} \to<400>^{H_2(f,U)_{\G}} H_2(B,U)_{\G} \to<150>
\textstyle{\frac{K}{V_1f}} \to<150> H_{1}(A,U)_{\G}
\to/{-{ >>}}/<400>^{H_{1}
(f,U)_{\G}} H_{1}(B,U)_{\G} \to<150> 1$}}
\]
in $\Bc$. 

\head{Chapter 5} In Chapter \ref{Chapter-Cohomology} we develop some new aspects of cohomology in the semi-abelian context. We simplify several recent investigations in this direction, in the context of crossed modules~\cite{Carrasco-Homology} or precrossed modules~\cite{AL}, and unify them with the classical theory that exists for groups and Lie algebras.

The definition of cohomology we use is an instance of Barr and Beck's cotriple cohomology: $H^{n} (X,A)$ is the $n$-th cohomology group of $X$, with coefficients in the functor 
\[
\hom (\cdot,A)\colon \Ac^{\op}\to \Ab,
\]
relative to the cotriple $\G$. Explicitly, we demand that $\Ac$ is a semi-abelian category, monadic over $\Set$, and $\G$ is the induced comonad; if $A$ is an abelian group object, $n\geq 1$ is a natural number, and $X$ is an object of $\Ac$, then 
\[
H^{n} (X,A)=H^{n-1}\hom (\ab (\G X),A).
\]
Note that we restrict ourselves to the case of abelianization.

We establish a Hochschild-Serre 5-term exact sequence extending the classical one for groups and Lie algebras: Any short exact sequence
\[
\xymatrix{ 0 \ar[r] & K \ar@{{ >}->}[r]^-{k} & X \ar@{-{ >>}}[r]^-{f} & Y \ar[r] & 0 }
\]
in $\Ac$ induces, in a natural way, an exact sequence
\[
\resizebox{\textwidth}{!}{\mbox{$0 \to<250> H^{1} (Y,A) \to/{{ >}->}/<250>^{H^{1}f} H^{1} (X,A) \to<250> \hom \bigl(\textstyle{\frac{K}{[K,X]}},A\bigr) \to<250> H^{2} (Y,A) \to<250>^{H^{2}f} H^{2} (X,A)$}}
\]
in $\Ab$ (Theorem~\ref{Theorem-Cohomology-Sequence}). Theorem~\ref{Theorem-H^2-Central-Extensions} states that the second cohomology group $H^2(Y,A)$ is isomorphic to the group $\CentrExt (Y,A)$ of isomorphism classes of central extensions of $Y$ by an abelian object $A$, equipped with (a generalization of) the Baer sum. We prove Proposition~\ref{Proposition-Universal-Central-Extension} that an object $Y$ is perfect (i.e., its abelianization $\ab (Y)$ is zero) if and only if it admits a universal central extension (i.e., the category $\CentrExt (Y)$ of central extensions of $Y$ has an initial object).

Finally, we consider a universal coefficient theorem to explain the relation between homology and cohomology (Theorem~\ref{Theorem-Universal-Coefficients}): If $Y$ is an object of $\Ac$ and $A$ is abelian then the sequence
\[
0 \to {\ext (H_{1}Y,A)}  \to/{{ >}->}/ H^{2} (Y, A) \to {\hom (H_{2}Y,A)}
\]
is exact. (This gives an indication why one uses homology \textit{objects} but cohomology \textit{groups}.)

\head{Chapter 6} The aim of the sixth chapter is to describe Quillen model category structures on the category $\Cat\Cc$ of internal categories and functors in a given finitely complete category~$\Cc$. Several non-equivalent notions of internal equivalence exist; to capture these notions, the model structures are defined relative to a given Grothendieck topology on $\Cc$.

Under mild conditions on $\Cc$, the regular epimorphism topology determines a model structure where $\we$ is the class of weak equivalences of internal categories (in the sense of Bunge and Par\'e~\cite{Bunge-Pare}). For a Grothendieck topos $\Cc$ we get a structure that, though different from Joyal and Tierney's~\cite{Joyal-Tierney}, has an equivalent homotopy category. In case $\Cc$ is semi-abelian, these weak equivalences turn out to be homology isomorphisms, and the model structure on $\Cat \Cc$ induces a notion of homotopy of internal crossed modules. In case $\Cc$ is the category $\Gp$ of groups and homomorphisms, it reduces to the well-known case of crossed modules of groups~\cite{Garzon-Miranda}.

The trivial topology on a category $\Cc$ determines a model structure on $\Cat \Cc$ where $\we$ is the class of strong equivalences (homotopy equivalences), $\fib$ the class of internal functors with the homotopy lifting property, and $\cof$ the class of functors with the homotopy extension property. As a special case, the ``folk'' Quillen model category structure on the category $\Cat =\Cat \Set$ of small categories is recovered.

\head{Chapter 7} In the final chapter~\ref{Chapter-Simplicial-Objects} we study the relationship between homology of simplicial objects and the model category structure on $\Sc \Ac$ introduced by Quillen~\cite{Quillen}. This model structure exists as soon as $\Ac$ is Mal'tsev and has enough projectives. Then the fibrations are just the Kan fibrations of $\Sc\Ac$. When $\Ac$ is moreover semi-abelian, a weak equivalence is nothing but a homology isomorphism: Hence, this model structure is compatible with the notion of homology from Chapter~\ref{Chapter-Homology} and with the model structure on $\Cat \Ac$ (induced by the regular epimorphism topology) from Chapter~\ref{Chapter-Internal-Categories}. We give a description of the trivial fibrations and of the acyclic objects, and investigate the relation with simplicial homotopy. \pfbreak

\headd{Acknowledgements} This is a revised version of my Ph.D.\ thesis, written and defended at the Vrije Universiteit Brussel under supervision of Prof.~Rudger W.~Kie\-boom. I wish to thank him, as well as my co-authors Tomas Everaert and Marino Gran. Many thanks also to Francis Borceux and George Janelidze for their constant support.

\begin{figure}[!b]
\[
\xymatrix{&&1\ar@{-}[ld] \ar@{-}[rd]\\
&2 \ar@{-}[ld] \ar@{-}[rdd] \ar@{-}[rrrdd] && 3 \ar@{-}[llld] \\
4 \ar@{-}[d]\\
5 && 6 && 7}
\]
\end{figure}
 % Introduction
\tableofcontents{}

\mainmatter
\setcounter{chapter}{0}
\chapter{The ``semi-abelian'' context}\label{Chapter-Preliminaries}

\setcounter{section}{-1}
\section{Introduction}\label{Section-Semiab-Introduction}
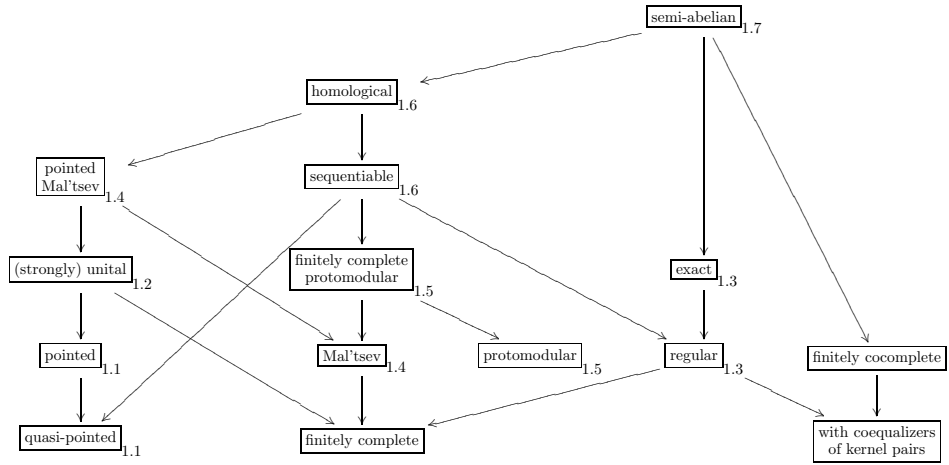
\begin{figure}[!b]
\resizebox{\textwidth}{!}{
\xymatrix{&&&&& {\fbox{semi-abelian}}_{\txt{\ref{Section-Semi-Abelian-Categories}}} \ar[lld] \ar[ddd] \ar[rdddd]\\
&&& {\fbox{homological}}_{\txt{\ref{Section-Homological-Sequentiable-Categories}}} \ar[llld] \ar[d]  \\
{\fbox{\txt{pointed\\
Mal'tsev}}}_{\txt{\ref{Section-Mal'tsev-Categories}}} \ar[d] \ar[rrrdd] &&& {\fbox{sequentiable}}_{\txt{\ref{Section-Homological-Sequentiable-Categories}}} \ar[d] \ar[rrdd] \ar[dddlll]\\
{\fbox{(strongly) unital}}_{\txt{\ref{Section-Unital-Categories}}} \ar[rrrdd] \ar[d] &&& {\fbox{\txt{finitely complete\\
protomodular}}}_{\txt{\ref{Section-Bourn-Protomodular-Categories}}} \ar[d] \ar[rd] && {\fbox{exact}}_{\txt{\ref{Section-Regular-Exact-Categories}}} \ar[d]\\
{\fbox{pointed}}_{\txt{\ref{Section-Pointed-Categories}}} \ar[d] &&& {\fbox{Mal'tsev}}_{\txt{\ref{Section-Mal'tsev-Categories}}} \ar[d] & {\fbox{protomodular}}_{\txt{\ref{Section-Bourn-Protomodular-Categories}}} & {\fbox{regular}}_{\txt{\ref{Section-Regular-Exact-Categories}}} \ar[lld] \ar[rd] & {\fbox{finitely cocomplete}} \ar[d]\\
{\fbox{quasi-pointed}}_{\txt{\ref{Section-Pointed-Categories}}} &&& {\fbox{finitely complete}} &&& {\fbox{\txt{with coequalizers\\
of kernel pairs}}}}}\caption{The ``semi-abelian'' context}\label{Figure-Semi-Abelian-Context}
\end{figure}

In this chapter we give an overview of the categorical structures, relevant to our work, as summarized in~\fref{Figure-Semi-Abelian-Context}. We consider their most important properties and examples, and the connections between them. We shall loosely refer to such categories as belonging to the ``semi-abelian'' context.

This chapter is not meant to be self-contained. Borceux and Bourn's book~\cite{Borceux-Bourn} provides an excellent introduction into this matter (and much more), and there are always the papers~\cite{Borceux-Semiab, Bourn-Gran-CategoricalFoundations, Janelidze-Marki-Tholen}, besides the literature cited below. Our main works of reference are \cite{Borceux:Cats, MacLane, AHS:Cats} for categories, \cite{Johnstone:Topos-Theory, Maclane-Moerdijk, Johnstone:Elephant} for topoi, \cite{Weibel, Maclane:Homology} for homology, \cite{Hovey, Quillen, KP:Abstact-homotopy-theory} for homotopy, and \cite{Robinson} for group theory.

Throughout the next chapters, it is \textit{always} possible to replace whatever conditions on a category $\Ac$ by ``semi-abelian with enough projectives'', except for some results in Chapter~\ref{Chapter-Cohomology}, where ``semi-abelian and monadic over $\Set$'' is needed. In fact, almost always ``semi-abelian'' is sufficient. Hence the hasty reader may imagine himself in such a nice situation and (against all advice) skim (or even skip) this chapter, using it for reference only.

\section{Quasi-pointed and pointed categories}\label{Section-Pointed-Categories}

In an abelian category one frequently considers kernels and cokernels of morphisms; but in an arbitrary category, no such concepts exist. Pointed and quasi-pointed categories form the context where a categorical definition of kernels and cokernels is possible. The former notion is apparently better known and certainly less subtle: A \emph{pointed}\index{pointed category}\index{category!pointed} category has an initial object $0$ that is also terminal (i.e., $0\cong 1$). A category is called \emph{quasi-pointed}\index{quasi-pointed category}\index{category!quasi-pointed} if it has an initial object $0$ and a terminal object $1$, and is such that the unique arrow ${0\to 1}$ is a monomorphism.

Virtually any book on category theory features pointed categories; but as a reference on quasi-pointed categories, we only know of Bourn's paper~\cite{Bourn2001}. We chose to incorporate them because we think the added generality is worth the few additional technical difficulties. For instance, thus we show that some of the constructions, carried out in the non-pointed context of the paper~\cite{EG}, actually may be seen as a consequence of the theory presented in~\cite{EverVdL1, EverVdL2}. Whoever finds this benefit doubtful may always choose to ignore quasi-pointed categories, thinking in terms of pointed categories instead.

\begin{definition}\label{Definition-Kernel}
Let $\Ac$ be a quasi-pointed category. Given a pullback square
\[
\xymatrix{K[f] \ar@{{ >}->}[r]^-{\ker f} \ar[d] \ar@{}[dr]|<{\pullback} & A \ar[d]^-{f}\\
0 \ar[r] & B}
\]
we call $\Ker f=(K[f],\ker f)$\index{K@$K[f]$}\index{ker@$\ker f$}\index{ker@$\Ker f$} a \emph{kernel}\index{kernel} of $f$. Given a pushout square
\[
\xymatrix{A \ar@{}[rd]|>{\pushout} \ar[r]^-{f} \ar[d] & B\ar@{-{ >>}}[d]^{\coker f} \\
0 \ar[r] & Q[f]}
\]
$\Coker f= (Q[f],\coker f)$\index{coker@$\Coker f$}\index{Q@$Q[f]$}\index{coker@$\coker f$} is called a \emph{cokernel}\index{cokernel} of $f$.
\end{definition}

\begin{definition}\label{Definition-Zero}
A morphism $A\to B$ that factors over $0$ is denoted $0\colon {A\to B}$ and called a \emph{zero}\index{zero} morphism. At most one such morphism exists.
\end{definition}

\begin{notation}\label{Notation-Kernel}
Most of the time we shall refer to a kernel $\Ker f$ or a cokernel $\Coker f$ by just naming the object part $K[f]$, $Q[f]$ or the morphism part $\ker f$, $\coker f$.

In a diagram, the forms $A\monoto B$\index{arrow!$\monoto$}, $A\kernelto B$\index{arrow!$\kernelto$}, $A \regepito B$\index{arrow!$\regepito$} and $A \cokernelto B$\index{arrow!$\cokernelto$} signify that the arrow is, respectively, a monomorphism, a kernel, a \emph{regular epimorphism}\index{regular epimorphism}\index{epimorphism!regular} (i.e., a coequalizer of some parallel pair of morphisms) and a cokernel (or \emph{normal epimorphism}\index{normal epimorphism}\index{epimorphism!normal}).
\end{notation}

\begin{remark}\label{Remark-Kernel-Mono}
Note that a kernel is always a monomorphism, being a pullback of a monomorphism; as a coequalizer of $f\colon A\to B$ and $0\colon A\to B$, $\coker f$ is a regular epimorphism.
\end{remark}

We recall the following important definition due to Lawvere~\cite{Lawvere}, which will allow us to describe many examples (namely, varieties in the sense of Universal Algebra) in a straightforward manner (see also~\cite{Borceux:Cats, Borceux-Bourn}).

\begin{definition}\label{Definition-Variety}
An \emph{(algebraic) theory}\index{algebraic theory}\index{theory}\index{category!theory} is a category $\T$ with objects $T^{0}, T^{1}, \dots , T^{n} , \dots$ for $n\in \N$, such that $T^{n}$ is the $n$-th power of $T=T^{1}$. (In particular, $T^{0}$ is a terminal object $1$.)

A \emph{$\T$-algebra}\index{T-algebra@$\T$-algebra} $X$ in a category $\Ac$ (or \emph{model}\index{model of a theory} $X$ of $\T$ in $\Ac$) is a functor $X\colon {\T \to \Ac}$ that preserves finite products, and a \emph{morphism of $\T$-algebras} is a natural transformation between two such functors. The category of $\T$-algebras in $\Ac$\index{category!of $\T$-algebras} is denoted by $\Alg_{\T}\Ac$\index{category!$\Alg_{\T}\Ac$} or $\Mod_{\T}\Ac$\index{category!$\Mod_{\T}\Ac$}; if $\Set$\index{category!$\Set$} denotes the category of sets and functions between them, we write $\Alg (\T)$\index{category!$\Alg (\T)$} for $\Alg_{\T}\Set$.

A \emph{variety (of universal algebras)}\index{variety}\index{variety!of universal algebras}\index{category!variety} is a category $\Vc$ that is equivalent to $\Alg (\T)$ for some theory $\T$.
\end{definition}

\begin{examples}\label{Examples-Varieties}
The category $\Mag$ of magmas\index{magma}\index{category!of magmas}\index{category!$\Mag$} is a variety: Here, $\T$ is generated by a binary operation $m:{T\times T\to T}$ and a constant $e:{1\to T}$, left and right unit for this operation. Such is also the category $\Gp$ of groups\index{category!$\Gp$} and, more generally, any \emph{variety of $\Omega$-groups}\index{variety!of $\Omega$-groups}\index{ZZ@\bigskip$\Omega$-group}\index{category!variety of $\Omega$-groups}~\cite{Higgins}: Here $\T$ contains a group operation and a unique constant (the unit of the group operation). In particular, the category $\Lie_{\K}$\index{category!$\Lie_{\K}$} of Lie algebras over a field $\K$ is a variety of $\Omega$-groups.
\end{examples}

\begin{examples}\label{Examples-Pointed-Categories}
A variety\index{variety!pointed} of universal algebras $\Vc \simeq\Alg (\T)$ is pointed if and only if its theory $\T$ contains a unique constant $e:1\to T$. Hence $\Mag$, $\Gp$, any variety of $\Omega$-groups etc.\ are examples, but not the category of unitary\index{category!of unitary rings} rings with unit-preserving morphisms (its theory contains two constants). In this text, $\Rng$\index{category!of rings}\index{category!$\Rng$} will denote the pointed variety of non-unitary rings and ring homomorphisms.

There are two dual ways of turning a category $\Ac$ into a pointed category. If it has a terminal object $1$, the coslice category $\slashfrac{1}{\Ac}$ is pointed; it is usually called the \emph{category of pointed objects in $\Ac$}\index{category!of pointed objects}. For instance, the category of pointed sets $\Set_{*}$\index{category!$\Set_{*}$} arises this way: If $1=*$ is a one-point set then a function ${1\to X}$ chooses a basepoint in $X$, and a morphism in $\slashfrac{1}{\Set}$ is a basepoint-preserving function. Dually, if $\Ac$ has an initial object $0$, we shall call \emph{category of copointed objects in $\Ac$}\index{category!of copointed objects} the slice category $\slashfrac{\Ac}{0}$.
\end{examples}

\begin{examples}\label{Examples-Quasi-Pointed-Categories}
$\Set$ is quasi-pointed. This example shows that a quasi-pointed category need not have non-trivial kernels or cokernels: The only kernels are of the form $\varnothing \to A$, and only these morphisms have a cokernel, namely $1_{A}\colon {A\to A}$. 

For us, an interesting class of non-trivial examples is formed by categories of internal groupoids over a fixed base object~\cite{Bourn2001, EG}. Let $\Ac$ be a finitely complete category and $\Gd \Ac$ the category of internal groupoids in $\Ac$. Consider the forgetful functor $(\cdot)_{0}\colon \Gd \Ac \to \Ac$, which maps a groupoid in $\Ac$ to is object of objects. This functor has a left adjoint $L$ and a right adjoint $R$. Its fibre $\Gd_{X}\Ac$\index{category!$\Gd_{X}\Ac$} over an object $X$ of $\Ac$ is the category of all groupoids in $\Ac$ with $X$ as object of objects. Any category $\Gd_{X}\Ac$ is quasi-pointed: Its initial object is the \emph{discrete}\index{discrete internal groupoid} groupoid $L (X)$ on $X$, its terminal object the \emph{indiscrete}\index{indiscrete internal groupoid} groupoid $R (X)$ on $X$ (see Section~\ref{Section-Preliminaries}).
\end{examples}

\begin{remark}\label{Remark-Kernels-Are-Special}
The example of $\Set$ also shows that in a quasi-pointed category, not every object occurs as a kernel. Only those objects occur that admit a (necessarily unique) morphism to the initial object. Indeed, for an object $A$, a morphism $A\to 0$ exists, exactly when $1_{A}\colon A\to A$ is a kernel. We call such an object \emph{copointed}\index{copointed object}. This name is well-chosen, because in a quasi-pointed category $\Ac$, the category of copointed objects $\slashfrac{\Ac}{0}$ may be considered as a subcategory of $\Ac$. Thus every quasi-pointed category canonically contains a pointed subcategory. (Note that its zero is $0$.) Hence if, in a quasi-pointed category, every object occurs as a kernel, then the category is pointed. (One may of course also see this directly by noting that the existence of a morphism $1\to 0$ implies that $0\cong 1$.)
\end{remark}

\begin{remark}\label{Remark-Cokernels-And-Cocompleteness}
A finitely complete quasi-pointed category has kernels, because it has pullbacks. But, in view of~\ref{Remark-Kernels-Are-Special}, only those morphisms can have a cokernel whose domain is copointed. In particular, a quasi-pointed category with cokernels is pointed. Nevertheless, a finitely cocomplete quasi-pointed category has all cokernels of kernels. (There is nothing contradictory about being cocomplete while not having cokernels of all morphisms: A cocomplete category has all colimits of diagrams that exist; when a morphism $f:{A\to B}$ has no cokernel, this is not because the diagram $f,0:{A\to B}$ has no coequalizer, but because there is no such diagram.)
\end{remark}

\begin{remark}\label{Remark-Copointed-Objects-Closed-Under-Strong-Epis}
Recall that a morphism $p\colon E\to B$ in a finitely complete category is a \emph{strong epimorphism}\index{strong epimorphism}\index{epimorphism!strong} when for every commutative square on the left
\[
\vcenter{\xymatrix{E \ar[d]_-{p} \ar[r] & A \mond^-{i} \\
B \ar@{.>}[ru] \ar[r] & X}}
\qquad \qquad 
\vcenter{\xymatrix{E \ar[d]_-{p} \ar[r] & 0 \mond \\
B \ar@{.>}[ru] \ar[r] & 1}}
\]
with $i$ mono there exists a (unique) diagonal making the whole diagram commute.

If $\Ac$ is finitely complete and quasi-pointed then $\slashfrac{\Ac}{0}$ is closed under strong epimorphisms: Suppose that $p\colon E\to B$ is a strong epimorphism and $E$ is in $\slashfrac{\Ac}{0}$. Then a lifting exists in the right hand side diagram above, showing that $B$ is in $\slashfrac{\Ac}{0}$.
\end{remark}

\section{Unital and strongly unital categories}\label{Section-Unital-Categories}

A well-known result due to Eckmann and Hilton~\cite{Eckmann-Hilton} states that any two internal group structures on a given group $G$ coincide, and that moreover the existence of such a structure implies that $G$ is abelian. In the context of (strongly) unital categories a general categorical version of this Eckmann-Hilton Theorem holds, giving rise to an intrinsic notion of abelian object.

\begin{definition}\label{Definition-Unital-Category}\cite{Bourn1996, Borceux-Bourn}
A pointed finitely complete category $\Ac$ is called \emph{unital}\index{unital category}\index{category!unital} when for all objects $K, K'$, the pair $(l_{K},r_{K'})$, where 
\[
l_K=(1_K,0) \colon  {K \to K \times K'}\index{l@$l_{K}$} 
\]
and $r_{K'} = (0,1_{K'}) \colon  {K' \to K \times K'}$\index{r@$r_{K}$}, is strongly epimorphic. $\Ac$ is called \emph{strongly unital}\index{category!strongly unital}\index{strongly unital category} when for every object $K$, the pair $(\Delta_{K},r_{K})$, where 
\[
\Delta_{K}= (1_{K},1_{K})\colon {K\to K\times K}\index{Delta K@$\Delta_{K}$},
\]
is a strong epimorphism.

We say that two coterminal morphisms $k\colon K\to X$ and $k'\colon K'\to X$ in a unital category \emph{cooperate}\index{cooperating morphisms} when a morphism $\varphi_{k,k'}\colon K\times K'\to X$ exists satisfying $\varphi_{k,k'}\comp l_{K}=k$ and $\varphi_{k,k'}\comp r_{K'}=k'$. The morphism $\varphi_{k,k'}$ is called a \emph{cooperator}\index{cooperator of coterminal morphisms} of $k$ and $k'$.

An object $A$ in a strongly unital category $\Ac$ is called \emph{abelian}\index{abelian!object} when $1_{A}$ cooperates with~$1_{A}$.
\end{definition}

Every strongly unital category is unital~\cite{Borceux-Bourn}. It is easily seen that, if two morphisms cooperate, their cooperator is necessarily unique.

\begin{examples}\label{Examples-Unital-Category}
Two examples of unital categories are $\Mag$\index{category!$\Mag$} magmas and $\Mon$\index{category!$\Mon$} monoids; but neither is strongly unital~\cite{Borceux-Bourn}.

Strongly unital varieties\index{variety!strongly unital}\index{strongly unital variety} may be characterized as categories $\Alg (\T)$ of algebras over a theory $\T$ that contains a unique constant $0:1\to T$ and a ternary operation 
\[
p\colon T\times T\times T\to T\colon (x,y,z)\mapsto p (x,y,z)
\]
satisfying the axioms $p (x,x,z)=z$ and $p (x,0,0)=x$ (cf.\ the case of Mal'tsev varieties, Examples~\ref{Examples-Mal'tsev}). Thus we see that groups, Heyting algebras etc.\ are strongly unital. Other examples include the category of internal (abelian) groups in a strongly unital category and finitely complete additive categories~\cite{Borceux-Bourn}. It is quite easily shown that also $\Set_{*}^{\op}$\index{category!$\Set_{*}^{\op}$}, the dual of the category of pointed sets, is strongly unital; on the other hand, $\Set_{*}$ is not even unital.
\end{examples}

\begin{theorem}[Eckmann-Hilton Theorem]\cite{Borceux-Bourn}\index{Eckmann-Hilton Theorem}\label{Theorem-Eckmann-Hilton}
Let $A$ be an object in a strongly unital category $\Ac$. The following are equivalent:
\begin{enumerate}
\item $A$ is abelian;
\item $A$ admits a structure of internal magma;
\item $A$ admits a structure of internal abelian group.
\end{enumerate}
Moreover, $A$ admits at most one structure of internal magma.\noproof 
\end{theorem}

The equivalence of 1.\ with 2.\ should be fairly clear from the definitions. Using Theorem~\ref{Theorem-Eckmann-Hilton}, it may be shown that every morphism between abelian objects induces a unique morphism of internal groups between the corresponding abelian group objects. In view of this fact, we use $\Ab \Ac$\index{category!$\Ab\Ac$} to denote both the full subcategory of $\Ac$ determined by the abelian objects and the category of abelian groups in $\Ac$. It is closed in $\Ac$ under subobjects and quotients; moreover, it is reflective. (As soon as $\Ac$ is exact---see below---such a subcategory is called a \emph{Birkhoff subcategory}\index{Birkhoff subcategory}\index{subcategory!Birkhoff ---} of $\Ac$~\cite{Janelidze-Kelly}.) This is an application of the following construction (due to Bourn~\cite{Bourn-Huq}) of the universal pair of cooperating morphisms associated with an arbitrary coterminal pair.

\begin{proposition}\cite{Bourn-Huq, Borceux-Semiab}\label{Proposition-Definition-Commutator}
In a finitely cocomplete unital category $\Ac$, let $k\colon K\to X$ and $k'\colon K'\to X$ be two morphisms with the same codomain. Then the arrow $\psi$, obtained by taking the colimit of the diagram of solid arrows
\begin{equation}\label{Diagram-Commutator}
\vcenter{\xy
\Atrianglepair(0,500)/->`.>`->``/[K`{K\times K'}`Y`X,;l_{K}``k``]
\Vtrianglepair/.>`<.`<-`<.`<-/[{K\times K'}`Y`X,`K';\varphi`\psi`r_{K'}``k']
\endxy}
\end{equation}
is a strong epimorphism. $\psi$ is the universal morphism that, by composition, makes the pair $(k,k')$ cooperate.

If $\Ac$ is strongly unital and $X$ is an object of $\Ac$, then the object $Y$ induced by the pair $(1_{X},1_{X})$ is the reflection of $X$ into $\Ab\Ac$.\noproof 
\end{proposition}

\begin{examples}\label{Examples-Abelian-Objects}
A group is abelian if and only if it is an abelian group\index{abelian!group}; a ring is abelian\index{abelian!ring} if and only if it is degenerate: Its multiplication is zero (here, more means less)~\cite{Borceux-Bourn}. A crossed module $(T,G,\del)$ is abelian if and only if both $T$ and $G$ are abelian groups and $G$ acts trivially on $T$ (see~\cite{Carrasco-Homology} and Example~\ref{Example-Crossed-Modules}).

Recall that a \emph{Lie algebra}\index{Lie algebra} $\Lieg$ is a vector space over a field $\K$ equipped with a bilinear operation $[\cdot , \cdot]\colon \Lieg \times \Lieg \to \Lieg$ (the \emph{Lie bracket}\index{Lie bracket}) that satisfies $[x,x]=0$ and $[x,[y,z]]+[y,[z,x]]+[z,[x,y]]=0$ (the \emph{Jacobi identity}\index{Jacobi identity}) for all $x,y,z\in \Lieg$. $\Lie_{\K}$\index{category!$\Lie_{\K}$} denotes the category of all Lie algebras over $\K$; morphisms are linear maps that preserve the bracket $[\cdot , \cdot]$. A~Lie algebra is called \emph{abelian}\index{abelian!Lie algebra}\index{Lie algebra!abelian} when its bracket is zero. This is the case, exactly when it is abelian: Suppose indeed that $\Lieg$ admits a structure of internal magma, $m\colon {\Lieg \times \Lieg\to \Lieg}$, $e\colon {\mathfrak{1}\to \Lieg\colon} {0\mapsto 0}$. Then, for all $x,y\in \Lieg$, 
\[
[x,y]=[m (0,x),m (y,0)]=m ([0,y],[x,0])=m (0,0)=0.
\]

In an abelian category, every object is abelian; in $\Set_{*}^{\op}$, the only abelian object is zero~\cite{Borceux-Bourn}.
\end{examples}

\section{Regular and Barr exact categories}\label{Section-Regular-Exact-Categories}

Another useful feature of abelian categories is the existence of image factorizations: Any morphism may be factored as a cokernel followed by a kernel. But in the category $\Gp$ this is no longer possible, as shows the example of the injection $H\to G$ of a subgroup $H$ that is not normal in the group $G$. However, every group homomorphism has, up to isomorphism, a unique factorization as a regular epimorphism (a coequalizer of some parallel pair) followed by a monomorphism. A category where such a factorization always exists is called \emph{regular}\index{regular epimorphism}\index{epimorphism!regular}:

\begin{definition}\label{Definition-Regular-Category}
Given a pullback diagram
\[
\xymatrix{R[f] \ar@{}[rd]|<{\pullback} \ar[r]^-{k_{1}} \ar[d]_-{k_{0}} & A \ar[d]^-{f}\\
A \ar[r]_-{f} & B}
\]
$(R[f],k_{0},k_{1})$\index{R@$R[f]$} is called a \emph{kernel pair}\index{kernel pair} of $f$.

A finitely complete category $\Ac$ with coequalizers of kernel pairs is said to be \emph{regular}\index{regular category}\index{category!regular} when in $\Ac$, pullbacks preserve regular epimorphisms.
\end{definition}

In a regular category, the \emph{image factorization}\index{image factorization} $\im f\comp p$ of a map $f\colon {A\to B}$ is obtained as follows: $p$ is a coequalizer of a kernel pair $k_{0},k_{1}\colon R[f]\to A$ and the \emph{image}\index{image} $\im f\colon I[f]\to B$\index{I@$I[f]$}\index{im@$\im f$} the universally induced arrow~\cite{Barr-Grillet-vanOsdol, Borceux:Cats, Borceux-Bourn}. Image factorizations are unique up to isomorphism, and may be chosen in a functorial way. Using the image factorization it is easily seen that in a regular category, regular and strong epimorphisms coincide\index{strong epimorphism}\index{epimorphism!strong}. (We may now forget about strong epimorphisms.) A related concept is that of direct images.

\begin{definition}\label{Definition-Direct-Image}
In a regular category $\Ac$, consider a mono $m\colon {D\to E}$ and a regular epimorphism $p\colon E\to B$. Taking the image factorization
\[
\xymatrix{D \mond_-{m} \ar@{.>>}[r] & pD \ar@{{_ (}.>}[d]^-{p (m)=\im (p\comp m)} \\
E \ar@{->>}[r]_-{p} & B}
\]
of $p\comp m$ yields a monomorphism $p (m)\colon pD\to B$ called the \emph{direct image of $m$ along $p$}\index{direct image}.
\end{definition}

We mentioned before that, even in a pointed and regular category like $\Gp$, not every morphism factors as a regular epimorphism followed by a kernel. But some do:

\begin{definition}\label{Definition-Proper-Morphism}\cite{Bourn2001}
A morphism $f$ in a quasi-pointed and regular category is called \emph{proper}\index{proper!morphism}\index{proper} when its image $\im f$ is a kernel. We call a subobject \emph{proper}\index{proper!subobject} when any representing monomorphism is proper (i.e., a kernel)\index{monomorphism!proper}\index{proper!monomorphism}. 
\end{definition}

\begin{remark}\label{Remark-Proper-Morphism-Has-Cokernel}
Following Remark~\ref{Remark-Cokernels-And-Cocompleteness}, in a quasi-pointed finitely cocomplete regular category, the image $\im f$ of a proper morphism $f$ has a cokernel. It is easily seen that this cokernel is also a cokernel of $f$. Moreover, $\im f$ is a kernel of $\coker f$.

A regular epimorphism\index{proper!regular epimorphism} is proper if and only if its codomain is copointed. For this to be the case, by Remark~\ref{Remark-Copointed-Objects-Closed-Under-Strong-Epis}, it suffices that such is its domain.
\end{remark}

Sometimes, only when this makes things easier, we shall give definitions of objects or morphisms in a regular category in terms of ``elements'' satisfying certain relations. We will only do so when it is clear how such a definition may be rewritten in terms of limits and images. It is, however, interesting to note that under the right circumstances, using elements in categories poses no problems at all, and may be justified via a Yoneda Lemma\index{Yoneda Lemma} argument. A very nice account on such considerations is given in~\cite{Borceux-Bourn}.

Another important aspect of regular categories is the behaviour of internal relations: They compose associatively. Recall that in a category with finite limits, a \emph{relation}\index{relation}\index{arrow!$\relto$} 
\[
R=(R, r_{0}, r_{1})\colon {A\relto B}
\]
from $A$ to $B$ is a subobject $(r_{0},r_{1})\colon {R\to A\times B}$. For instance, a kernel pair of some morphism $f\colon {A\to B}$ may be construed as a relation $(R[f],k_{0},k_{1})$ from $A$ to $A$ (\emph{on} $A$) called the \emph{kernel relation} of $f$\index{kernel relation}. (Most of the time, we shall simply ignore the difference between a relation as subobject of a product (a class of equivalent monomorphisms) or as jointly monic pair of arrows.) If a map $(x_{0},x_{1})\colon {X\to A\times B}$ factors through $(r_{0},r_{1})$, then the map $h\colon {X\to R}$ with $(x_{0},x_{1})= (r_{0},r_{1})\comp h$ is necessarily unique; we will denote the situation by $x_{0} (R)x_{1}$. In a regular category, $SR\colon {A\relto C}$ denotes the composition of $R\colon {A\relto B}$ with $S\colon {B\relto C}$. We recall from~\cite{Carboni-Kelly-Pedicchio} the following proposition; it follows from 2.\ that in a regular category, the composition of relations is associative.

\begin{proposition}\label{Proposition-Regularity}
Let $\Ac$ be a regular category.
\begin{enumerate}
\item A map $b\colon {X\to B}$ factorizes through the image of a map $f\colon {A\to B}$ if and only if there is a regular epimorphism $p\colon {Y\to X}$ and a map $a\colon {Y\to A}$ with $b\comp p=f\comp a$;
\item given relations $R\colon {A\relto B}$ and $S\colon {B\relto C}$ and maps $a\colon {X\to A}$ and $c\colon {X\to C}$, $c (SR) a$ iff is a regular epimorphism $p\colon {Y\to X}$ and a map $b\colon {Y\to B}$ with $b (R) a\comp p$ and $c\comp p (S)b$.\noproof
\end{enumerate}
\end{proposition}

Recall that an object $P$ in a category $\Ac$ is called \emph{(regular) projective}\index{regular projective object}\index{projective object} when for every regular epimorphism $p\colon E\to B$ in $\Ac$, the function 
\[
p\comp (\cdot)=\hom (P,p)\colon \hom (P,E)\to \hom (P,B)
\]
is a surjection. A category $\Ac$ is said to have \emph{enough (regular) projectives}\index{enough projectives} when for every object $A$ in $\Ac$ there exists a projective object $P$ and a regular epimorphism $p\colon {P\to A}$. For some categorical arguments we shall need that the category under consideration has enough projectives (see e.g., Theorem~\ref{Theorem-L-Exact-Sequence} or Chapter~\ref{Chapter-Internal-Categories} for a number of examples). 

\begin{remark}\label{Remark-Projectives}
In a regular category some arguments, which would otherwise involve projectives, may be avoided\index{projective object!avoiding}\index{avoiding projectives}, and thus the requirement that sufficiently many projective objects exist. This is how the first statement of Proposition~\ref{Proposition-Regularity} is used. (See also Section~\ref{Section-Kan}, in particular Remark~\ref{Remark-Kan-Projectives}, and Example~\ref{Example-Homotopy-not-an-Equivalence-Relation}. Using Grothendieck topologies, even the regularity requirement may sometimes be avoided: One works e.g., with the topology of pullback-stable regular epimorphisms. This is one of the possible interpretations of Proposition~\ref{Proposition-Characterization-T-Epimorphism}.)

For instance, given the assumptions of 2.\ in~\ref{Proposition-Regularity}, in case $X$ is a projective object, one easily sees that a map $b\colon {X\to B}$ exists such that $b (R)a$ and $c (S)b$. If now $X$ is arbitrary, but $\Ac$ has enough projectives, one may take a projective object $Y$ and a regular epimorphism $p\colon {Y\to X}$ to get the conclusion of 2.\ in~\ref{Proposition-Regularity}. In case $\Ac$ lacks projectives, this property needs an alternative proof; the regularity of $\Ac$ allows us to prove it using the first statement of Proposition~\ref{Proposition-Regularity}.
\end{remark}

As for ordinary relations, the notion of equivalence relation has an internal categorical counterpart. Note that for every object $X$ in a finitely complete category $\Ac$, the $\hom$ functor $\hom (X,\cdot)\colon \Ac \to \Set$ maps internal relations in $\Ac$ to relations in $\Set$.

\begin{definition}\label{Definition-Equivalence-Relation}
A relation $R$ on an object $A$ of a finitely complete category $\Ac$ is called an \emph{equivalence relation}\index{equivalence relation}\index{relation!equivalence ---} if, for every $X$ in $\Ac$, its image through $\hom (X,\cdot)$ is an equivalence relation on the set $\hom (X,A)$. $\Eq \Ac$\index{category!$\Eq \Ac$} denotes the category of internal equivalence relations in $\Ac$, considered as a full subcategory of the category $\Rgrph \Ac$\index{category!$\Rgrph \Ac$}\index{reflexive graph} of internal reflexive graphs\index{internal!reflexive graph} in $\Ac$, the functor category $\Fun (\Five,\Ac)$. Here $\Five$\index{category!$\Five$} is the Platonic Idea of a reflexive graph, the category generated by the arrows
\[
\xymatrix@1{{1} \ar@(ul,dl)[]_{{1}_1} \ar@<-1 ex>[r]_-{\varpi_{0}} \ar@<1 ex>[r]^-{\varpi_{1}} & \ar[l]|-{\varsigma} {0} \ar@(dr,ur)[]_{{1}_0}}
\]
satisfying $\varpi_{0}\comp \varsigma=\varpi_{1}\comp \varsigma=1_{0}$ and $1_{0}\comp 1_{0}=1_{0}$, $1_{1}\comp 1_{1}=1_{1}$.
\end{definition}

Modifying this definition in the obvious way, one acquires internal notions of reflexive, symmetric, etc.\ relations. It is easily seen that a kernel relation is always an equivalence relation. Such an equivalence relation is called \emph{effective}\index{effective!equivalence relation}\index{relation!effective equivalence ---}\index{equivalence relation!effective}. But in general, the converse is not necessarily true, whence the following definition.

\begin{definition}\label{Definition-Exact-Category}\cite{Barr}
A regular category $\Ac$ is called \emph{(Barr) exact}\index{exact category}\index{Barr exact category}\index{category!exact}\index{category!Barr exact} when in $\Ac$, every equivalence relation is a kernel pair.
\end{definition}

\begin{examples}\label{Examples-Barr-Exact-Category}
Any variety of universal algebras\index{variety} is exact~\cite{Barr}, as is (the dual of) any topos\index{category!topos}\index{topos}~\cite{Johnstone:Topos-Theory, Borceux:Cats, Maclane-Moerdijk, Johnstone:Elephant} (in particular, $\Set$). In a variety, the regular epimorphisms are just surjective algebra homomorphisms. (The inclusion ${\Z \to \Q}$ in the variety of unitary rings\index{category!of unitary rings} is an example of an epimorphism that is not regular.) Also any category of algebras and any category of slices or coslices in an exact category is again exact (in particular, $\Set_{*}^{\op}=\slashfrac{\hbox{$\Set^{\op}$}}{\hbox{$*$}}$). 

The category $\Cat$\index{category!$\Cat$} of small categories and functors between them is not regular: Regular epimorphisms are not pullback-stable. (The different notions of epimorphism in $\Cat$ were first explored by Giraud~\cite{Giraud}, see e.g.,~\cite{Janelidze-Sobral-Tholen}.) This is also the case for the category $\Top$\index{category!$\Top$} of topological spaces and continuous maps~\cite{Borceux:Cats}. However, both the category $\CompHaus$\index{category!$\CompHaus$}\index{category!of compact Hausdorff spaces} of compact Hausdorff spaces and its dual (the category of commutative $C^{*}$ algebras\index{category!of commutative $C^{*}$ algebras}) are exact~\cite{Barr}. The category $\CgHaus$\index{category!$\CgHaus$}\index{category!of compactly generated Hausdorff spaces} of compactly generated Hausdorff spaces (or Kelley spaces) is regular~\cite{Hardie:Regular}. 

Two examples of categories that are regular but not exact are torsion-free abelian groups\index{category!of torsion-free abelian groups}~\cite{Bourn-Gran-CategoricalFoundations} and $\Gp \Top$, groups in $\Top$ or ``topological groups''\index{category!of topological groups}~\cite{Carboni-Kelly-Pedicchio}.
\end{examples}

\section{Mal'tsev categories}\label{Section-Mal'tsev-Categories}

As regular categories constitute a natural context to work with relations, regular Mal'tsev categories constitute a natural context to work with equivalence relations.

\begin{definition}\cite{Carboni-Lambek-Pedicchio}\label{Definition-Mal'tsev-Category}
A finitely complete category $\Ac$ is called \emph{Mal'tsev}\index{Mal'tsev category} when in $\Ac$, every reflexive relation\index{reflexive relation}\index{relation!reflexive} is an equivalence relation.
\end{definition}

\begin{examples}\label{Examples-Mal'tsev}
A variety is Mal'tsev if and only if it is a \emph{Mal'tsev variety}\index{Mal'tsev variety}\index{variety!Mal'tsev ---}~\cite{Smith}. Such is the category $\Alg (\T)$ of algebras over a theory $\T$ that contains a \emph{Mal'tsev operation}\index{Mal'tsev operation}, i.e., a ternary operation 
\[
p\colon T\times T\times T\to T\colon (x,y,z)\mapsto p (x,y,z)
\]
that satisfies the axioms $p (x,x,z)=z$ and $p (x,z,z)=x$. Examples include the category of groups, where $p (x,y,z)=xy^{-1}z$, hence any variety of $\Omega$-groups, and Heyting algebras, where~\cite{Borceux-Bourn}
\[
p (x,y,z)= ((x\Rightarrow y)\Rightarrow z)\wedge ((z\Rightarrow y)\Rightarrow x)
\]
or 
\[
p (x,y,z)= (y\Rightarrow (x\wedge z))\wedge (x\vee z).
\]

A non-varietal example is the dual of any elementary topos (or, more generally, of any finitely complete pre-topos, which, in the case of compact Hausdorff spaces, yields the category of commutative $C^{*}$ algebras as an example~\cite{Carboni-Kelly-Pedicchio}). The Mal'tsev property is inherited by e.g., slice categories $\slashfrac{\Ac}{X}$ and coslice categories $\slashfrac{X}{\Ac}$  (e.g., $\Set_{*}^{\op}$ is Mal'tsev), and the domain $\Bc$ of any pullback-preserving and \emph{conservative}\index{conservative functor} (i.e., isomorphism-reflecting) functor $\Bc \to \Ac$ with a Mal'tsev codomain $\Ac$ is a Mal'tsev category. As shown by M.~Gran in~\cite{Gran:Internal}, the case for $\Cat \Ac$, categories in a regular Mal'tsev category $\Ac$, is entirely different from $\Cat$ (Examples~\ref{Examples-Barr-Exact-Category}):  $\Cat \Ac$\index{category!$\Cat \Ac$}\index{internal!category} is again regular Mal'tsev, and if $\Ac$ is exact, so is $\Cat \Ac$.
\end{examples}

A regular category is Mal'tsev if and only if the composition of equivalence relations is commutative~\cite{Carboni-Kelly-Pedicchio}, i.e., when for any two equivalence relations $S$ and $R$ on an object $A$ the equality $SR=RS$ holds. In a regular Mal'tsev category, the equivalence relations on a given object $A$ constitute a lattice, the join of two equivalence relations being their composition.

It is well known that when, in a regular category, a commutative square of regular epimorphisms
\[
\xymatrix{ A' \ar@{->>}[r]^-{f'} \ar@{->>}[d]_-{v} & B' \ar@{->>}[d]^-w \\
A \ar@{->>}[r]_-f & B}
\]
is a pullback, it is a pushout. In a regular category, a commutative square of regular epimorphisms is said to be a \emph{regular pushout}\index{regular pushout}\index{pushout!regular} when the comparison map $r\colon {A'\to P}$ to a pullback of $f$ along $w$ is a regular epimorphism (see Carboni, Kelly and Pedicchio~\cite{Carboni-Kelly-Pedicchio}).
\begin{equation}\label{Diagram-Regular-Pushout}
\vcenter{\xymatrix{ A' \ar@{->>}[drr]^-{f'} \ar@{->>}[drd]_-{v} \ar@{.>}[rd]|r \\
& P \ar@{}[rd]|<{\pullback} \ar@{.>>}[r] \ar@{.>>}[d] & B' \ar@{->>}[d]^-w \\
& A \ar@{->>}[r]_-f & B}}
\end{equation}

\begin{proposition}\cite{Bourn2003}\label{Rotlemma-kernel-pairs}
Consider, in a regular Mal'tsev category, a commutative diagram of augmented kernel pairs, such that $p$, $p'$, $q$ and $r$ are regular epimorphisms:
\[
\xymatrix{R[p'] \ar[d]_-{s} \ar@<-0.5 ex>[r]_-{k'_{1}}\ar@<0.5 ex>[r]^-{k'_{0}} & A' \ar@{->>}[d]^-{q} \ar@{->>}[r]^-{p'} & B' \ar@{->>}[d]^-{r}\\
R[p] \ar@<-0.5 ex>[r]_-{k_{1}}\ar@<0.5 ex>[r]^-{k_{0}} & A  \ar@{->>}[r]_-{p} & B.  }
\]
The right hand square is a regular pushout if and only if $s$ is a regular epimorphism.\noproof
\end{proposition}

The following characterizes exact Mal'tsev categories among the regular ones.

\begin{proposition}\cite{Carboni-Kelly-Pedicchio}\label{Proposition-Regular-Pushout}
A regular category is exact Mal'tsev if and only if, given regular epimorphisms $v\colon A'\to A$ and $f'\colon A'\to B'$ such as in~\ref{Diagram-Regular-Pushout}, their pushout (the diagram of solid arrows in~\ref{Diagram-Regular-Pushout}) exists, and moreover the comparison map $r$---where the square is a pullback---is a regular epimorphism.\noproof
\end{proposition}

It follows that in any exact Mal'tsev category, a square of regular epimorphisms is a regular pushout if and only if it is a pushout. As explained to me by Marino Gran, as a consequence, we get the following very useful result. It is well-known to hold in the pointed exact protomodular case---see below.

\begin{proposition}\label{Proposition-Image-of-Kernel}
In a quasi-pointed exact Mal'tsev category, the direct image of a kernel along a regular epimorphism is a kernel.
\end{proposition}
\begin{proof}
Suppose that $k=ker f\colon {K\to A}$ is a kernel of $f\colon {A\to C}$, and that $p\colon {A\to B}$ is a regular epimorphism. Let $m\comp r$ be the image factorization of $f$ ($r\colon {A\to I}$ regular epi, $m\colon {I\to C}$ mono). Then $k$ is still a kernel of $r$. Consider pushing out $r$ along $p$, then taking kernels:
\[
\xymatrix{K \ar@{.>}[r]^-{\overline{p}} \ar@{{ >}->}[d]_-{k} & K[\rho] \ar@{{ >}->}[d] \\
A \ar@{->>}[r]^-{p} \ar@{}[dr]|>>{\pushout} \ar@{->>}[d]_-{r} & B \ar@{->>}[d]^-{\rho} \\
I  \ar@{->>}[r]_-{\pi} & J.}
\]
We are to show that the induced arrow $\overline{p}$ is a regular epimorphism. Now according to Proposition~\ref{Proposition-Regular-Pushout}, in an exact Mal'tsev category, a pushout of a regular epimorphism along a regular epimorphism is always a regular pushout. Hence the diagram above has the following factorization.
\[
\xymatrix{K \ar@{.>}[r] \ar@{}[rd]|-{\texttt{(i)}} \ar@{{ >}->}[d]_-{k} & K[\overline{\rho}] \ar@{{ >}->}[d] \ar@{.>}[r]^{\widetilde{\pi}} & K[\rho] \ar@{{ >}->}[d] \\
A \ar@{->>}[r] \ar@{->>}[d]_-{r} & P \ar@{}[dr]|<<{\pullback} \ar@{->>}[d]_-{\overline{\rho}} \ar[r]^{\overline{\pi}} & B \ar@{->>}[d]^-{\rho} \\
I \ar@{=}[r] & I  \ar@{->>}[r]_-{\pi} & J}
\]
The morphism $\widetilde{\pi}$ is an isomorphism, because in the cube
\[
\xymatrix{{K[\overline{\rho} ]} \ar@{}[rddd]|<<<<<{\skewpullback} \ar[rr]^{\widetilde{\pi }} \ar[dd] \ar@{{ >}->}[rd] && {K[\rho ]} \ar@{}[rddd]|<<<<<{\skewpullback} \ar[dd]|{\hole} \ar@{{ >}->}[rd] \\
& {P} \ar@{}[rd]|<{\pullback} \ar[rr]^<<<<<<<{\overline{\pi}} \ar[dd]_<<<<<<<{\overline{\rho}} && {B} \ar[dd]^{\rho} \\
{0} \ar@{=}[rr]|{\hole} \ar[rd] && {0} \ar[rd] \\
& {I} \ar[rr]_-{\pi } && {J}}
\]
the front, left and right squares are pullbacks, hence so is the back square.

For $\overline{p}$ to be regular epi, it now suffices that the upper left square \texttt{(i)} is a pullback; and it is easily seen that this is indeed the case.
\end{proof}

The notion of abelian object may be extended to the context of Mal'tsev categories as follows~\cite{Carboni-Pedicchio-Pirovano, Borceux-Bourn}. An object of a Mal'tsev category admits at most one internal Mal'tsev operation\index{internal!Mal'tsev operation}; if it does, we say that the object is \emph{abelian}\index{abelian!object}. Every pointed Mal'tsev category is strongly unital~\cite{Borceux-Bourn}. Moreover, in a pointed Mal'tsev category, an object admits an internal Mal'tsev operation if and only if it admits an internal abelian group structure; hence this notion of abelian object coincides with the one from Definition~\ref{Definition-Unital-Category}. In an exact Mal'tsev category with coequalizers, abelian objects may be characterized using Smith's commutator of equivalence relations~\cite{Smith, Pedicchio, Gran:Central-Extensions}; see Section~\ref{Section-Commutators}.

Another important aspect of Mal'tsev categories concerns the Kan property; this will be the subject of Section~\ref{Section-Kan} and parts of Chapter~\ref{Chapter-Simplicial-Objects}.

\section{Bourn protomodular categories}\label{Section-Bourn-Protomodular-Categories}

The last important ingredient for the semi-abelian context is \textit{protomodularity}. This crucial notion due to Dominique Bourn is perhaps the most difficult and elusive ingredient. It is stronger than the Mal'tsev axiom; for a quasi-pointed category, it implies that every regular epimorphism is a cokernel; for a quasi-pointed regular category, it is equivalent to the validity of the Short Five Lemma. In this context, the notion of protomodularity is strong enough to imply the basic lemma's of homological algebra, such as the $3\times 3$ Lemma and the Snake Lemma. In a quasi-pointed, exact and protomodular category also Noether's Isomorphism Theorems hold.

\begin{definition}\label{Definition-Protomodularity}\cite{Bourn1991}
Recall that a \emph{split epimorphism}\index{split epimorphism}\index{epimorphism!split} or \emph{retraction}\index{retraction} is a morphism $p\colon {A\to B}$ such that there is a morphism $s\colon {B\to A}$ (called a \emph{splitting of $p$}) satisfying $p\comp s=1_{B}$. Dually, then $s$ is a \emph{split monomorphism}\index{split monomorphism}\index{monomorphism!split} or \emph{section}\index{section}.

A category is called \emph{(Bourn) protomodular}\index{protomodular category}\index{Bourn protomodular category}\index{category!protomodular}\index{category!Bourn protomodular} when, given a commutative diagram 
\[
\xymatrix{{\cdot} \ar@{}[rd]|-{\texttt{(i)}} \ar[r] \ar[d] & {\cdot} \ar@{}[rd]|-{\texttt{(ii)}} \ar@{.>}[d]\ar[r] & {\cdot} \ar[d]\\ 
{\cdot} \ar[r] & {\cdot} \ar[r] & {\cdot}}
\]
with the dotted vertical arrow a split epimorphism, if the outer rectangle \texttt{(i)} + \texttt{(ii)} and the left hand side square \texttt{(i)} are pullbacks, then so is the right hand side square \texttt{(ii)}.
\end{definition}

A split epimorphism is always a regular; hence one implication of the next statement is clear.

\begin{proposition}\label{Proposition-Protomodularity-I}
A regular category is protomodular if and only if given the commutative diagram above, where the dotted arrow is a regular epimorphism, the left hand square and the whole rectangle are pullbacks, the right hand square is also a pullback.\noproof 
\end{proposition}

\begin{proposition}\label{Proposition-Protomodularity-II}
A quasi-pointed and regular category is protomodular if and only if the \emph{Short Five Lemma}\index{Short Five Lemma} holds: This means that for any commutative diagram
\[
\xymatrix{K[p'] \ar@{{ >}->}[r]^-{\ker p'} \ar[d]_-u & E' \ar@{->>}[r]^-{p'} \ar[d]^-v & B' \ar[d]^-w \\ K[p] \ar@{{ >}->}[r]_-{\ker p} & E \ar@{->>}[r]_-{p} & B}
\]
such that $p$ and $p'$ are regular epimorphisms, $u$ and $w$ being isomorphisms implies that $v$ is an isomorphism.\noproof 
\end{proposition}

\begin{examples}\label{Examples-Protomodular-Categories}
A first large class of examples is characterized by the following theorem due to Bourn and Janelidze:

\begin{theorem}\cite{Bourn-Janelidze}\index{variety!protomodular}
A variety of universal algebras $\Vc$ is protomodular if and only if its theory $\T$ has nullary terms  $e_{1},\dots , e_{n}$, binary terms $t_1, \dots , t_n$ and one ${(n+1)}$-ary term $\tau$ satisfying the identities $\tau (x,t_1 (x,y), \dots, t_n (x,y))=y$ and $t_i(x,x)=e_{i}$ for each $i=1,\dots, n$.\noproof 
\end{theorem}

Thus we see that e.g., all varieties of $\Omega$-groups, in particular (abelian) groups, non-unitary rings, Lie algebras, crossed modules are examples of protomodular categories.

Next, internal versions of these varieties\index{internal}: Given a finitely complete category $\Ac$, one may consider groups, rings etc.\ in $\Ac$, and thus construct new protomodular categories. Other constructions also inherit protomodularity: e.g., slice categories $\slashfrac{\Ac}{X}$ and the domain $\Bc$ of any pullback-preserving and conservative functor ${\Bc \to \Ac}$ with a protomodular codomain $\Ac$.

Finally, the category of Heyting algebras~\cite{Bourn1996} and the dual of any elementary topos~\cite{Bourn:Dual-topos} are protomodular (hence so is $\Set_{*}^{\op}=\slashfrac{\hbox{$\Set^{\op}$}}{\hbox{$*$}}$), as is any additive category.
\end{examples}

\begin{proposition}\cite{Bourn1996}\label{Proposition-Protomodular-Category-is-Mal'tsev}
Every finitely complete protomodular category is a Mal'tsev category.\noproof 
\end{proposition}

\begin{remark}\label{Remark-Protomodular-Cokernels}
In a finitely complete and quasi-pointed protomodular category, every regular epimorphism is a cokernel (of its kernel). A slight modification of Remark~\ref{Remark-Cokernels-And-Cocompleteness} and Remark~\ref{Remark-Proper-Morphism-Has-Cokernel} shows that a quasi-pointed regular protomodular category has cokernels of proper morphisms; moreover, a proper morphism may be factored as a cokernel followed by a kernel. The only difference with the abelian case is that not every morphism has such a factorization.
\end{remark}

In a protomodular category $\Ac $, an intrinsic notion of normal monomorphism exists.

\begin{definition}\label{Definition-Normal-Monomorphism}~\cite{Bourn2000}
An arrow $k\colon  K \to A$ in a finitely complete category $\Ac$ is \textit{normal to an equivalence relation $R$ on $A$}\index{normal monomorphism}\index{monomorphism!normal} when: 1.\ $k^{-1}(R)$ is the largest equivalence relation $\nabla_K= (K\times K,\pr_{1},\pr_{2})$\index{pr@$\pr_1\colon K\times K\to K$}\index{nabla@$\nabla_{X}$} on $K$; 2.\ the induced map $\nabla_K \to R$ in the category $\Eq \Ac$\index{category!$\Eq \Ac$} of internal equivalence relations in $\Ac$ is a \emph{discrete fibration}\index{discrete fibration}.
\end{definition}

This means that
\begin{enumerate}
\item there is a map $\tilde k\colon K \times K \to R$ in $\Ac$ such that the diagram
\[
\xymatrix{ K \times K \ar@{.>}[r]^-{\tilde{k}} \ar@<.5ex>[d]^-{\pr_{2}} \ar@<-.5ex>[d]_-{\pr_{1}} & R \ar@<.5ex>[d]^-{d_{1}} \ar@<-.5ex>[d]_-{d_{0}}\\
K \ar[r]_-{k} & A}
\]
commutes;
\item any of the commutative squares in the diagram above is a pullback.
\end{enumerate}

One may prove that the arrow $k$ is then necessarily a monomorphism; furthermore, when the category $\Ac$ is protomodular, a monomorphism can be normal to at most one equivalence relation, so that the fact of being normal becomes a property~\cite{Bourn2000}. The notion of normal monomorphism gives an intrinsic way to express the fact that $K$ is an equivalence class of $R$.

In a quasi-pointed finitely complete category every kernel is normal. (In particular, $\ker f$ is normal with respect to~$R[f]$. The converse is not true: If $B$ is not copointed then the morphism $1_{B}$ is not a kernel, although it is normal w.r.t.\ $\nabla_{B}$. \cite[Counterexample~3.2.19]{Borceux-Bourn} gives an example in the pointed and protomodular category $\Top \Gp$.) In the pointed exact protomodular case, normal monomorphisms and kernels coincide.

There is a natural way to associate, with any equivalence relation $(R,d_{0},d_{1})$ in a quasi-pointed finitely complete category, a normal monomorphism $k_{R}$, called the \emph{normalization}\index{normalization of an equivalence relation} of $R$, or \emph{the normal subobject associated with $R$}: It is defined as the composite $k_{R}=d_{1}\comp \ker d_{0}$
\[
K[d_{0}] \to/{ >}->/^{\ker d_{0}} R \to^{d_{1}} A.
\]
In the exact Mal'tsev case, by Proposition~\ref{Proposition-Image-of-Kernel}, this normal monomorphism is a kernel. In the exact protomodular case, this construction determines a bijection between the (effective) equivalence relations on $A$ and the proper subobjects of $A$ (a proper monomorphism corresponds to the kernel pair of its cokernel). %The coequalizer of \nabla_X is the image of X\to 1; the kernel of this map is the associated proper subobject.
In the pointed protomodular case, it determines a bijection between the equivalence relations on $A$ and the normal subobjects of $A$~\cite{Bourn2000}. 

The following important result (the Non-Effective Track of the $3\times 3$ Lemma, Theorem 4.1 in~\cite{Bourn2003}) is a non-pointed version of Proposition~\ref{Proposition-Image-of-Kernel}.

\begin{proposition}\label{Proposition-non-eff-3x3-lemma}
Consider, in a regular and protomodular category, a commutative square with horizontal regular epimorphisms
\[
\xymatrix{A' \ar[d]_-{v}
\ar@{->>}[r]^-{f'} & B'\ar[d]^-{w} \\ A \ar@{->>}[r]_-{f} & B.}
\]
When $w$ is a monomorphism and $v$ a normal monomorphism, $w$ is normal.\noproof
\end{proposition}

In the pointed and exact case, it implies that the direct image of a kernel is a kernel (cf.\ Proposition~\ref{Proposition-Image-of-Kernel}).

\section{Sequentiable and homological categories}\label{Section-Homological-Sequentiable-Categories}

This is the context where exact sequences are well-behaved and the major homological lemma's hold. \emph{Sequentiable}\index{sequentiable category}\index{category!sequentiable} categories~\cite{Bourn2001} are quasi-pointed, regular and protomodular; a pointed and sequentiable category is called \emph{homological}\index{homological category}\index{category!homological}~\cite{Borceux-Bourn}.

\begin{definition}\label{Definition-Exact-Sequence}
We shall call a sequence
\begin{equation}\label{Short-Exact-Sequence}
\xymatrix{K \ar[r]^-{k} & A \ar[r]^-{f} & B}
\end{equation}
in a quasi-pointed category a \emph{short exact sequence}\index{short exact sequence}\index{exact sequence!short} if $k=\ker  f$ and $f=\coker k$. We denote this situation
\begin{equation}\label{Exact-Sequence-Short}
\xymatrix{0 \ar[r] & K \ar@{{ >}->}[r]^-{k} & A \ar@{-{ >>}}[r]^-{f} & B \ar[r] & 1.}
\end{equation}
If we wish to emphasize the object $K$ instead of the arrow $k$, we denote the cokernel $f$ by $\eta_{K}\colon A\to A/K$\index{H@$\eta_{K}$}. In a sequentiable category the exactness of Sequence~\ref{Short-Exact-Sequence} is equivalent to demanding that $k=\ker  f$ and $f$ is a regular epimorphism.

A sequence of morphisms $(f_{i})_{i\in I}$
\[
\xymatrix{\dots \ar[r] & A_{i+1} \ar[r]^-{f_{i+1}} & A_{i} \ar[r]^-{f_{i}} & A_{i-1}\ar[r] & \dots}
\]
in a sequentiable category is called \emph{exact at $A_{i}$}\index{exact sequence} if $\im f_{i+1}=\ker f_{i}$. It is called \emph{exact} when it is exact at $A_{i}$, for all $i\in I$. (Here $I$ is not supposed to be infinite.) A functor between sequentiable categories is called \emph{exact}\index{exact functor} if it preserves exact sequences.
\end{definition}

\begin{remark}\label{Remark-Exact-Sequence}
One would imagine that Sequence~\ref{Exact-Sequence-Short} is exact if and only if it represents~\ref{Short-Exact-Sequence} as a short exact sequence. If $\Ac$ is homological then this is indeed the case; moreover the notation becomes less awkward, as we may then (according to our taste) replace the $0$ and $1$ by two $0$s or two $1$s. But if $\Ac$ is only quasi-pointed this is, in general, \textit{not} true. 

In a sequentiable category, a morphism $k\colon {K\to A}$ is mono if and only if its kernel is ${0\to K}$; hence $k=\ker f$ if and only if~\ref{Exact-Sequence-Short} is exact at $K$ and at $A$. A morphism $f\colon {A\to B}$ is a cokernel (of its kernel $k$) if and only if it is a regular epimorphism if and only if its image is $1_{B}\colon {B\to B}$. But $1_{B}$ is a kernel of ${B\to 1}$ \textit{only when $B$ is copointed}. (A regular epimorphism is proper if and only if its codomain is copointed.) Hence, in general, asking that~\ref{Exact-Sequence-Short} is exact at $B$ is stronger than asking that $f$ is a cokernel. At least, we can say that the notation used in~\ref{Exact-Sequence-Short} \textit{sometimes} makes sense; and this is the best we can do. Note that, when $B$ is copointed, we can also write a $0$ at the end of the sequence; but otherwise this would be wrong.

To ``solve'' this problem, we adopt the following inconsequent convention: If we say that a sequence
\[
\xymatrix{\dots \ar[r] & A_{2} \ar[r]^-{f_{2}} & A_{1} \ar[r]^-{f_{1}} & A_{0}\ar[r] & 1}
\]
is exact, we mean that it is exact at all $A_{i}$ for $i>0$, and moreover $f_{1}$ is a regular epimorphism. \textit{It need not be exact at $A_{0}$.}
\end{remark}

Using these notations, we can write the Short Five Lemma in a slightly more familiar way.

\begin{proposition}[Short Five Lemma]\label{Proposition-Short-Five-Lemma}
Given a commutative diagram with exact rows
\[
\xymatrix{0 \ar[r] & K' \ar@{{ >}->}[r]^-{k'} \ar[d]_-u & X' \ar@{-{ >>}}[r]^-{f'} \ar[d]^-v & Y' \ar[d]^-w \ar[r] & 1\\
0 \ar[r] & K \ar@{{ >}->}[r]_-{k} & X \ar@{-{ >>}}[r]_-{f} & Y \ar[r] & 1}
\]
in a sequentiable category, $u$ and $w$ being isomorphisms implies that $v$ is an isomorphism.\noproof 
\end{proposition}

\begin{proposition}[Snake Lemma]\cite{Bourn2001, Borceux-Bourn}\label{snake-lemma}\index{Snake Lemma}
Let $\Ac$ be a sequentiable category. Any commutative diagram with exact rows as below such that $u$, $v$ and $w$ are proper, can be completed to the following diagram, where all squares commute,
\[
\xymatrix{
 & K[u] \ar@{.>}[r] \ar@{{ >}.>}[d]_-{\ker u} & K[v] \ar@{.>}[r] \ar@{{ >}.>}[d]_-{\ker v} & K[w] \ar@{{ >}.>}[d]^-{\ker w} \ar@{.>}[llddd]^>>>>>>>>>{\delta} \\
 & K' \ar[r]_-{k'} \ar[d]_-u & X' \ar@{-{ >>}}[r]_-{f'} \ar[d]^v & Y' \ar[d]^w \ar[r] & 1 \\
0 \ar[r] & K \ar@{{ >}->}[r]^-k \ar@{.{ >>}}[d]_-{\coker u} & X \ar[r]^-f \ar@{.{ >>}}[d]^-{\coker v} & Y \ar@{.{ >>}}[d]^-{\coker w}  &  \\
 & \Cok[u] \ar@{.>}[r] & \Cok[v] \ar@{.>}[r] & \Cok[w] & }
\]
in such a way that
\[
\xymatrix{
K[u] \ar[r] & K[v] \ar[r] & K[w] \ar[r]^-{\delta} & Q[u] \ar[r] & Q[v] \ar[r] & Q[w] }
\]
is exact. Moreover, this can be done in a natural way, i.e., defining a functor ${\PAr\Ac\to \STE\Ac}$\index{category!$\PAr \Ac$}\index{category!$\STE \Ac$}, where $\PAr\Ac$ is the category of proper arrows of $\Ac$ and $\STE\Ac$ the category of six-term exact sequences in $\Ac$.\noproof
\end{proposition}

In fact, in~\cite{Bourn2001} only the exactness of the sequence is proven. However, it is quite clear from the construction of the connecting morphism $\delta$ that the sequence is, moreover, natural.

\begin{proposition}\cite{Bourn1991, Bourn2001}\label{Proposition-LeftRightPullbacks}
In a sequentiable category $\Ac$, let us consider the commutative diagram with exact rows 
\begin{equation}\label{Diagram-LeftRightPullbacks}
\vcenter{\xymatrix{0 \ar[r] & K' \ar[d]_{u} \ar@{{ >}->}[r]^{k'} & X' \ar[d]^{v} \ar@{-{ >>}}[r]^{f'}   &  Y'\ar[d]^{w} \ar[r] & 1\\
0 \ar[r] & K \ar@{{ >}->}[r]_{k} & X  \ar@{-{ >>}}[r]_{f}   &  Y \ar[r] & 1. }}
\end{equation}
Then:
\begin{enumerate}
\item $u$ is an isomorphism if and only if the right-hand square is a pullback;
\item $w$ is a monomorphism if and only if the left-hand square is a pullback.\noproof 
\end{enumerate}
\end{proposition}

This may almost be dualized:

\begin{proposition}\cite{Gran-VdL, Bourn-Gran}\label{Proposition-pushout}\label{Rotlemma}
In a sequentiable category $\Ac$, let us consider Diagram~\ref{Diagram-LeftRightPullbacks}. Then:
\begin{enumerate}
\item if the left-hand square is a pushout, then $w$ is an isomorphism; conversely, when $\Ac$ is exact and $u$ is a regular epimorphism, if $w$ is an isomorphism, then the left-hand square is a pushout;
\item if $u$ is an epimorphism then the right hand square is a pushout; conversely, when $\Ac$ is exact and $v$ and $w$ are regular epi, if the right-hand square is a pushout, then $u$ is a regular epimorphism.
\end{enumerate}
\end{proposition}
\begin{proof}
1.\ If the left hand side square is a pushout then $f'$ and $0\colon K\to Y'$ induce an arrow $f''\colon X\to Y'$ satisfying $f''\comp v=f'$ and $f''\comp k=0$. Then also $w\comp f''=f$. Moreover, like $f$, $f''$ is a cokernel of $k$; hence the unique comparison map $w$ is an isomorphism.

Now again consider the diagram above; pushing out $u$ along $k'$, then taking a cokernel of $\overline{k'}$ induces the dotted arrows in the diagram below.
\[
\xymatrix{0 \ar[r] & K' \ar[d]_{u} \ar@{{ >}->}[r]^{k'} \ar@{}[rd]|>{\pushout}& X' \ar@{.>}[d]^{\overline{u}} \ar@{-{ >>}}[r]^{f'}   &  Y'\ar@{:}[d] \ar[r] & 1\\
0 \ar@{.>}[r] & K \ar@{:}[d] \ar@{.>}[r]^{\overline{k'}} & P \ar@{.>}[d]^{v'} \ar@{.{ >>}}[r]^{f''}   &  Y'\ar[d]^{w} \ar@{.>}[r] & 1\\
0 \ar[r] & K \ar@{{ >}->}[r]_{k} & X  \ar@{-{ >>}}[r]_{f}   &  Y \ar[r] & 1}
\]
Note that $\overline{k'}$ is a monomorphism because $k$ is one; being the regular image of the kernel $k'$, by Proposition~\ref{Proposition-Image-of-Kernel}, $\overline{k'}$ is a kernel as well. The Short Five Lemma~\ref{Proposition-Short-Five-Lemma} implies that if $w$ is an isomorphism then so is $v'$.

2.\ For a proof, see Lemma 1.1 in~\cite{Bourn-Gran}.
\end{proof}

\begin{proposition}[Noether's First Isomorphism Theorem]\cite{Borceux-Bourn}\label{generalexact}\index{Noether's First Isomorphism Theorem}
Let $A\subseteq B\subseteq C$ be objects of a sequentiable category $\Ac$, such that $A$ and $B$ are proper subobjects of $C$ (i.e., the inclusions are kernels). Then
\[
\xymatrix{
0 \ar[r] & \frac{B}{A} \ar@{{ >}->}[r] & \frac{C}{A} \ar@{-{ >>}}[r] & \frac{C}{B} \ar[r] & 1}
\]
is a short exact sequence in $\Ac$.\noproof
\end{proposition}

\section{Semi-abelian categories}\label{Section-Semi-Abelian-Categories}

For the results in this thesis, the ideal context is that of \emph{semi-abelian}\index{category!semi-abelian} categories~\cite{Janelidze-Marki-Tholen}: pointed, exact, protomodular with binary coproducts. Not only does it encompass all the mentioned categorical environments, so that all results in the foregoing sections are valid in a semi-abelian category; there is also a historical reason for the importance of this notion. Indeed, introducing semi-abelian categories, Janelidze, M\'arki and Tholen solved Mac\,Lane's long standing problem~\cite{MacLane:Duality} of finding a framework that reflects the categorical properties of non-abelian groups as nicely as abelian categories do for abelian groups. But over the years, many different people came up with partial solutions to this problem, proving theorems starting from various sets of axioms, which all require ``good behaviour'' of normal mono- and epimorphisms. In the paper~\cite{Janelidze-Marki-Tholen}, the relationship between these ``old-style'' axioms and the semi-abelian context is explained, and thus the old results are incorporated into the new theory.

Next to being suitable for homological algebra of non-abelian structures, semi-abelian categories provide a good foundation for the treatment of isomorphism and decomposition theorems (e.g., Proposition~\ref{generalexact} and Borceux and Grandis~\cite{Borceux-Grandis}), radical theory~\cite{BG:Torsion} and commutator theory. There is also an intrinsic notion of semi-direct product~\cite{Bourn-Janelidze:Semidirect}, internal action~\cite{BJK} and internal crossed module~\cite{Janelidze}; cf.\ Section~\ref{Section-Regular-Epimorphisms}. 

The following theorem due to Bourn and Janelidze characterizes semi-abelian varieties.

\begin{theorem}\cite{Bourn-Janelidze}\index{variety!semi-abelian}
A variety of universal algebras $\Vc\simeq \Alg (\T)$ is semi-abelian if and only if its theory $\T$ has a unique constant $0$, binary terms $t_1, \dots , t_n$ and a $(n+1)$-ary term $\tau$ satisfying the identities 
\[
\tau (x,t_1 (x,y), \dots, t_n (x,y))=y
\]
and $t_i(x,x)=0$ for each $i=1,\dots, n$.\noproof 
\end{theorem}

Since any semi-abelian category is strongly unital with finite colimits, we may revisit the results from Section~\ref{Section-Unital-Categories} and use them as a starting point to commutator theory. Here we briefly introduce one approach, due to Huq~\cite{Huq}, generalized by Borceux and Bourn to semi-abelian categories~\cite{Borceux-Bourn, Borceux-Semiab}. Different approaches, e.g., using equivalence relations, valid in more general situations, exist and are treated throughout the text; see Section~\ref{Section-Commutators}.

\begin{definition}\label{Definition-Commutator}
In a semi-abelian category, consider coterminal morphisms $k\colon K\to X$ and $k'\colon K'\to X$ and the resulting colimit diagram~\ref{Diagram-Commutator}. The kernel of $\psi$ is denoted $[k,k']\colon {[k,k']\to X}$\index{$[k,k']$} and called the \emph{commutator of $k$ and $k'$ (in the sense of Huq)}\index{commutator!of coterminal morphisms}\index{commutator!in the sense of Huq}.

A morphism $k\colon K\to X$ is called \emph{central}\index{central morphism}\index{Huq centrality} when $[k,1_{X}]=0$. 
\end{definition}

\begin{remark}\label{Remark-Commutator}
Recall that, in a semi-abelian category, all regular epimorphisms are cokernels of their kernels, hence $Y=X/[k,k']$. It follows that $k$ and $k'$ cooperate if and only if $[k,k']=0$. In particular, an object $A$ of $\Ac$ is abelian\index{abelian!object} exactly when $[1_{A},1_{A}]=0$.
\end{remark}

\begin{examples}\label{Examples-Huq-Central}
For groups $N,N'\triangleleft G$ we denote by $[N,N']$ the (normal) subgroup of $G$ generated by the elements $nn'n^{-1}{n'}^{-1}$, with $n\in N$ and $n'\in N'$. It is well-known that a group $G$ is abelian exactly when $[G,G]=0$; moreover, $[G,G]=[1_{G},1_{G}]$. More generally, if $i_{N}$ and $i_{N'}$ denote the inclusions of $N$ and $N'$ into $G$, then $[N,N']=[i_{N},i_{N'}]$. An inclusion $i_{N}$ is central if and only if $N$ lies in the \emph{centre}\index{z@$\zeta G$} 
\[
\zeta G=\{x\in G\,|\,\forall g\in G:\,xg=gx \}
\]
of $G$.

For Lie algebras, the situation is very similar. First recall that a kernel (normal monomorphism) in $\Lie_{\K}$ is an inclusion $i_{\Lien}\colon \Lien \to \Lieg$ of an ideal $\Lien$; such is a subspace $\Lien$ of a Lie algebra $\Lieg$ that satisfies $[g,n]\in \Lien$ for all $g\in \Lieg$, $n\in \Lien$. Clearly, for any two ideals $\Lien$ and $\Lien'$ of $\Lieg$, the subspace $[\Lien ,\Lien']$ of $\Lieg$ generated by $[n,n']$ for $n\in \Lien$ and $n'\in \Lien'$ is an ideal. It is almost tautological that a Lie algebra $\Lieg$ is abelianized by dividing out the ideal $[\Lieg ,\Lieg]=[1_{\Lieg},1_{\Lieg}]$. Again, $[\Lien ,\Lien ']=[1_{\Lien},1_{\Lien '}]$, and an inclusion $i_{\Lien}$ is central if and only if $\Lien$ lies in the \emph{centre}\index{z@$\zeta \Lieg$}
\[
\zeta \Lieg=\{x\in \Lieg\,|\,\forall g\in \Lieg:\,[g,x]=0 \}
\]
of $\Lieg$.
\end{examples}

We shall need the following properties of the commutator $[k,k']$.

\begin{proposition}\cite{EverVdL1}\label{Proposition-Commutator}
In a semi-abelian category $\Ac$, let $k\colon K\to X$ and $k'\colon K'\to X$ be coterminal morphisms.
\begin{enumerate}
\item If $k$ or $k'$ is zero then $[k,k']=0$.
\item For any regular epimorphism $p\colon X\to X'$, $p[k,k']=[p\comp k,p\comp k']$. (See Definition~\ref{Definition-Direct-Image}; this means that there exists a regular epimorphism $\overline{p}$ such that the square
\[
\xymatrix{[k,k'] \ar@{{ >}->}[d] \ar@{.{ >>}}[r]^-{\overline{p}} & [p\comp k,p\comp k'] \ar@{{ >}->}[d]\\
X \ar@{-{ >>}}[r]_-{p} & X'}
\]
commutes.)
\item If $k\colon K\to X$ is a kernel, then $[k,k]$ factors over $K$, and $K/[k,k]$ is abelian in $\Ac$.
\end{enumerate}
\end{proposition}
\begin{proof}
1.\ is obvious. The rest of the proof is based on Huq~\cite[Proposition 4.1.4]{Huq}. As in Proposition~\ref{Proposition-Definition-Commutator}, let $\psi\colon  X\to Y$ and $\varphi\colon K\times K'\to Y$, resp. $\psi'\colon X'\to Y'$ and $\varphi'\colon {K\times K'\to Y'}$, denote the couniversal arrows obtained from the construction of $[k,k']$ and $[p\comp k,p\comp k']$. Then
\[
\bfig
\Atrianglepair(0,500)/->`.>`->``/[K`{K\times K'}`Y'`X,;l_{K}`\psi'\comp p\comp k`k``]
\Vtrianglepair/.>`<.`<-`<.`<-/[{K\times K'}`Y'`X,`K';\varphi'`\psi'\comp p`r_{K'}`\psi'\comp p\comp k'`k']
\efig\]
is a cocone on the diagram of solid arrows~\ref{Diagram-Commutator}. The couniversal property of colimits yields a unique map $y\colon Y\to Y'$. In the commutative diagram of solid arrows
\[
\xymatrix{0 & Y \ar[l] \ar[d]_-{y} & X \ar@{-{ >>}}[l]_-{\psi} \ar[d]^-{p} & [k,k'] \ar@{.>}[d]^-{\overline{p}} \ar@{{ >}->}[l]_-{\ker \psi} & 0 \ar[l]\\
0 & Y' \ar[l] & X' \ar@{-{ >>}}[l]^-{\psi'} & [p\comp k,p\comp k'] \ar@{{ >}->}[l]^-{\ker \psi'} & 0, \ar[l]}
\]
there exits a unique map $\overline{p}\colon [k,k']\to [p\comp k,p\comp k']$ such that the right hand square commutes.

For 2.\ we must show that if $p$ is a regular epimorphism, then so is $\overline{p}$. To do so, we prove that $\kappa =\ker \psi'\comp \im \overline{p}$ is a kernel of $\psi'$. By Proposition~\ref{Proposition-Image-of-Kernel}, $\kappa $ is a kernel; hence, it is sufficient that $\psi'$ be a cokernel of $\kappa $.

Let $z\colon X'\to Z$ be a map such that $z\comp \kappa =0$. Then $z\comp p\comp \ker \psi=0$, which yields a map $y\colon Y\to Z$ with $y\comp \psi=z\comp p$. We get the following cocone.
\[
\bfig
\Atrianglepair(0,500)/->`.>`->``/[K`{K\times K'}`Z`X';l_{K}`z\comp p\comp k`p\comp k``]
\Vtrianglepair/.>`<.`<-`<.`<-/[{K\times K'}`Z`X'`K';y\comp \varphi`z`r_{K'}`z\comp p\comp k'`p\comp k']
\efig
\]
Thus we acquire a unique arrow $x\colon Y'\to Z$ such that $x\comp \psi'=z$.

$[1_{K},1_{K}]$ is a subobject of $[k,k]$: Take $p=k$ and $k=k'=1_{K}$ in the discussion above. Hence, using that $K/[1_{K},1_{K}]$ is abelian and that $\Ab\Ac$ is closed under quotients, the first statement of 3.\ implies the second one. This first statement follows from the fact that
\[
\bfig
\Atrianglepair(0,500)|aaaaa|/->`.>`->``/[K`{K\times K}`{\Cok[k]}`X,;l_{K}`0`k``]
\Vtrianglepair|aabab|/.>`<.`<-`<.`<-/[{K\times K}`{\Cok[k]}`X,`K;0`\coker k`r_{K}`0`k]
\efig
\]
is a cocone. Thus a map may be found such that the right hand square in
\[
\xymatrix{0 \ar[r] & [k,k] \ar@{.>}[d]_-{i} \ar@{{ >}->}[r] & X \ar@{=}[d] \ar@{-{ >>}}[r]^-{\psi} & \tfrac{X}{[k,k]} \ar@{.>}[d] \ar[r] & 0\\
0 \ar[r] & K \ar@{{ >}->}[r]_-{k} & X \ar@{-{ >>}}[r]_-{\coker k} & \Cok [k] \ar[r] & 0}
\]
commutes; this yields the needed map $i$.
\end{proof}

\section{Abelian categories}\label{Section-Abelian-Categories}

Of course, abelian categories do not really belong to the semi-abelian context. But, as it is important to know what not to do, this section gives a quick overview of those techniques (one uses in the abelian context) that, when valid in a semi-abelian category, make it abelian. For a more profound account on such considerations, see~\cite{Janelidze-Marki-Tholen, Bourn-Gran-CategoricalFoundations, Borceux-Bourn, Borceux-Semiab}.

Recall that a category is \emph{abelian}\index{abelian!category}\index{category!abelian}~\cite{Freyd} when it is pointed, has binary products and coproducts, has kernels and cokernels, and is such that every monomorphism is a kernel and every epimorphism is a cokernel. Examples of abelian categories include all categories of modules over a ring and (pre)sheaves of modules. The category $\Ab \Ac$ of abelian group objects in an exact category $\Ac$ is always abelian. This definition is clearly self-dual, whereas a semi-abelian category with a semi-abelian dual is an abelian category.

A category is \emph{additive}\index{additive category}\index{category!additive} when it is finitely complete and enriched in $\Ab$. The \emph{Tierney equation}\index{Tierney equation} ``abelian = additive + exact'' states that a category is abelian if and only if it is additive and Barr exact. The category $\Ab \Top$\index{category!$\Ab\Top$} of abelian topological groups is an example of a category that is additive and regular, though not exact. A category is additive if and only if it is pointed and protomodular, and such that every monomorphism is normal. Hence, a semi-abelian category where every monomorphism is a kernel is abelian. 

A pointed category with finite limits is \emph{linear}\index{linear category}\index{category!linear} when its binary products are binary coproducts (i.e., it has biproducts). Any additive category is linear. A semi-abelian and linear category is abelian, because all its objects are.
 % Context
\chapter{Homology}\label{Chapter-Homology}

\setcounter{section}{-1}
\section{Introduction}\label{Section-Homology-Introduction}
In this chapter we study homology in the semi-abelian context. We start by considering the concept of chain complex, the basic building block of any homology theory. As soon as the ambient category is quasi-pointed, exact and protomodular, homology of \textit{proper} chain complexes---those with boundary operators of which the image is a kernel---is well-behaved: It characterizes exactness. As in the abelian case, the $n$-th homology object of a proper chain complex $C$ with boundary operators $d_n$ is said to be 
\[
H_{n}C=\Cok[C_{n+1}\to K[d_n]].
\]
We prove that this equals the dual $K_{n}C={K[\Cok[d_{n+1}]\to C_{n-1}]}$. Moreover, any short exact sequence of proper chain complexes gives rise to a long exact sequence of homology objects.

In Section~\ref{Section-Simplicial-Objects} we extend the homology theory of Section~\ref{Section-Chain-Homology} to simplicial objects. To do so, we consider the \textit{normalization functor} $N\colon \simpA \to \Chplus\Ac$. Suppose that $\Ac$ is a pointed category with pullbacks. Let us write $\del_{i}$ for the face operators of a simplicial object $A$ in $\Ac$. The \textit{Moore complex} $N (A)$ of $A$ is the chain complex with $N_{0}A=A_{0}$,
\[
N_{n} A=\bigcap_{i=0}^{n-1}K[\del_{i}\colon A_{n}\to A_{n-1}]
\]
and boundary operators $d_{n}=\del_{n}\comp \bigcap_{i}\ker \del_{i}\colon N_{n} A\to N_{n-1} A$, for $n\geq 1$. We show that the Moore complex $N (A)$ of a simplicial object $A$ in a quasi-pointed exact Mal'tsev category $\Ac$ is always proper. This allows us to define the \textit{$n$-th homology object} of an object $A$ in an exact and sequentiable category $\Ac$ as $H_{n}A=H_{n}N (A)$. We prove that, if $\epsilon \colon A\to A_{-1}$ is a contractible augmented simplicial object, then $H_{0}A=A_{-1}$ and, for $n\geq 1$, $H_{n}A=0$.

The validity of the generalized Hopf Formula studied in Chapter~\ref{Chapter-Cotriples} strongly depends on the normalization functor $N\colon {\simpA \to \Chplus \Ac}$ being exact. This essentially amounts to the fact that, in a regular Mal'tsev category, any regular epimorphism of simplicial objects is a Kan fibration. We prove this in Section~\ref{Section-Kan}, thus generalizing Carboni, Kelly and Pedicchio's result~\cite{Carboni-Kelly-Pedicchio} that in a regular Mal'tsev category, every simplicial object is Kan. The exactness of $N$ is proved in Section~\ref{Section-Implications-Kan}. We get that any short exact sequence of simplicial objects induces a long exact sequence of homology objects. These facts also lead to a proof of Dominique Bourn's conjecture that for $n\geq 1$, the homology objects $H_{n}A$ of a simplicial object $A$ are abelian.\pfbreak

\noindent Some of the results in this chapter may be found in Borceux and Bourn's book~\cite{Borceux-Bourn}; they were first proved in my paper with Tomas Everaert~\cite{EverVdL2}. Of course, the possibility of considering homology of chain complexes was implicit in the very definition of semi-abelian categories, and it was known for long that chain complexes of (non-abelian) groups are well-behaved. On the other hand, the fact that the Moore complex of a simplicial object is proper (Theorem~\ref{Theorem-N-Proper}) seems to be new, and this is a crucial ingredient for the further development of homology theory along semi-abelian lines.

Around the time we proved his conjecture, Dominique Bourn found an alternative proof~\cite{BournBodeux}. His proof is more conceptual and less technical, because it avoids the use of the Kan property; and as it does not need the existence of coproducts, it is also slightly more general.

\section{Chain complexes}\label{Section-Chain-Homology}
Although usually considered in an abelian context, chain complexes of course make sense in any pointed or even quasi-pointed category $\Ac$:

\begin{definition}\label{Definition-Chain-Complex}
Let $\Ac$ be a quasi-pointed category. An \emph{unbounded chain complex}\index{unbounded chain complex} $C$ is a collection of morphisms $(d_{n}\colon C_{n}\to C_{n-1})_{n\in \Z}$ in $\Ac$ such that $d_{n}\comp d_{n+1}=0$, for all $n\in\Z$. The category of unbounded chain complexes in $\Ac$ (with, as morphisms, commutative ladders) is denoted by $\Ch\Ac$\index{category!$\Ch\Ac$}.

A \emph{bounded chain complex} or a \emph{positively graded chain complex}\index{chain complex!positively graded}\index{bounded chain complex}\index{positively graded chain complex}\index{chain complex!bounded} $C$ is a collection $(d_{n}\colon C_{n}\to C_{n-1})_{n\in \N_{0}}$ satisfying $d_{n}\comp d_{n+1}=0$, for all $n\in\N_{0}$. The category of positively graded chain complexes in $\Ac$ is denoted by $\Chplus\Ac$\index{category!$\Chplus\Ac$}.

When we use the term \emph{chain complex}\index{chain complex}, we mean the chain complexes of either one variety or of both, depending on the context.
\end{definition}

\begin{remark}\label{Remark-Positively-Graded}
In the pointed case, as usual, the category $\Chplus\Ac$ can be considered as a full subcategory of $\Ch\Ac$ by extending a positively graded chain complex with $0$s in the negative degrees. Note however that, in a quasi-pointed category, this is in general not possible: A complex $C$ has such an extension if and only if $C_{0}$ is copointed. Also, extending a chain complex with $1$s instead of $0$s would not work, because then $d_{-2}\comp d_{-1}$ is not zero.
\end{remark}

Nevertheless, obtaining a good notion of homology object $H_nC$ of a chain complex $C$ demands stronger assumptions on $\Ac$ and on $C$. When $\Ac$ is an abelian category, $H_nC$ is ${K[d_n]}/{I[d_{n+1}]}$ (see e.g.,~\cite{Weibel}). Since (in the abelian case) this is just ${\Cok[C_{n+1}\to K[d_n]]}$, it seems reasonable to define $H_nC$ this way, provided that the considered kernels and cokernels exist in $\Ac$. Yet, one could also suggest the dual $K_{n}C={K[\Cok[d_{n+1}]\to C_{n-1}]}$, because, in the abelian case, this equals $H_nC$. This is not true in general (see Example~\ref{Example-Different-H-and-K}). Now in principle, the existence of two distinct notions of homology poses no problems. But it \textit{is} a problem that, as the next example shows, neither of these two notions characterizes the exactness of the sequence under consideration, not even when the ambient category is semi-abelian. (Needless to say, detecting exactness is \textit{the} main point of homology.)

\begin{example}\label{Example-Homology-no-Exactness}
Consider, in the category $\Gp$, the chain complex $C$ defined by choosing $d_{1}$ the inclusion of $A_{4}$ into $A_{5}$, and $C_{n}=0$ for all $n\not \in \{0,1 \}$. (Any other example of an inclusion of a subgroup, maximal while not normal, would do.) Recall that for $n= 3$ or $n\geq 5$, the \emph{alternating group}\index{alternating group} $A_{n}$---the group of all even permutations of the set $\{0,\dots ,n-1 \}$---is simple, i.e., its only normal subgroups are $0$ and itself~\cite{Robinson}. Hence all $H_{n}C$ and $K_{n}C$ are zero, but of course $C$ is not exact at $C_{0}$. 
\end{example}

To cope with this problem, we shall only consider those chain complexes that allow their exactness properties to be detected using homology. Luckily, the class of such complexes is easy to describe, because we know that a morphism can only occur in an exact sequence when it is proper (or occurs at the end of the sequence). Moreover, as we shall see later---this is what makes Theorem~\ref{Theorem-N-Proper} important---such a requirement is not too restrictive. 

\begin{definition}\label{Definition-Homology}
Let $\Ac$ be quasi-pointed and regular. A chain complex $C$ in $\Ac$ is called \emph{proper}\index{chain complex!proper}\index{proper!chain complex} when all its boundary operators $d_{n}\colon {C_{n}\to C_{n-1}}$ are proper morphisms. $\PCh\Ac$\index{category!$\PCh\Ac$} (resp.\ $\PChplus\Ac$\index{category!$\PChplus\Ac$}) denotes the full subcategory of $\Ch\Ac$ (resp.\ $\Chplus\Ac$) determined by the proper complexes. 

Now suppose that $\Ac$ is, moreover, protomodular. Consider a proper complex $C\in\Ob{\PCh\Ac}$ and $n\in \Z$ (or $C\in\Ob{\PChplus\Ac}$ and $n\in \N_{0}$). $H_nC$\index{H_n@$H_{n}C$} is the \emph{$n$-th homology object}\index{homology!of proper chain complexes} of $C$, and $K_nC$\index{K_n@$K_{n}C$} its dual, defined as follows:
\begin{equation}\label{Diagram-K-H}
\vcenter{\xymatrix{
C_{n+1} \ar[rd]_{d'_{n+1}} \ar[rr]^-{d_{n+1}} & & C_n \ar[rr]^-{d_n} \ar@{-{ >>}}[rd]^{\coker d_{n+1}} & & C_{n-1}   \\
& K[d_n] \ar@{{ >}->}[ru]^-{\ker d_n} \ar@{-{ >>}}[d]_{\coker d'_{n+1}} \ar@{.>}[rrd]^<<<<<<{p_n} & & \Cok[d_{n+1}] \ar[ru]_{d''_{n}} \\
& H_nC \ar@{.>}[rru]^>>>>>>{j_n}  \ar@{.>}[rr]_{(\lambda_n)_C} & & \ar@{{ >}->}[u]_{\ker d''_{n}} K_nC. &  }}
\end{equation}
$H_{n}C$ is a cokernel of the induced morphism $d'_{n+1}$ and $K_{n}C$ is a kernel of $d''_{n}$. Thus we get functors $H_{n}$, $K_{n}\colon {\PCh \Ac \to \Ac}$, for any $n\in \Z$, and $H_{n}$, $K_{n}\colon {\PChplus \Ac \to \Ac}$, for any $n\in \N_{0}$.
\end{definition}

\begin{remark}\label{Remark-Cokernels-of-Kernels}
$H_{n}C$ and $K_{n}C$ exist, as $\Ac $, being sequentiable, has all cokernels of proper morphisms (see Remark~\ref{Remark-Protomodular-Cokernels}).
\end{remark}

\begin{proposition}\label{Proposition-Homology-Characterizes-Exactness}
Let $\Ac$ be a sequentiable category. A proper chain complex $C$ in $\Ac$ is exact at $C_{n}$ if and only if $H_{n}C=0$.
\end{proposition}
\begin{proof}
The proper morphism $d_{n+1}'$ is regular epi if and only if its cokernel $H_{n}C$ is $0$.
\end{proof}

\begin{remark}\label{Remark-Homology-for-Bounded-Below-Complexes}
In the quasi-pointed case, the definition of homology (Definition~\ref{Definition-Homology}) is not complete: A definition of $H_{0}$ (and $K_{0}$) must be added. For $C\in \PChplus\Ac$, we put 
\[
H_{0}C=K_{0}C=Q[d_{1}\colon {C_{1}\to C_{0}}],
\]
i.e., a cokernel of $d_{1}$, because this is what these homology objects are equal to in case $C_{0}$ is copointed (and hence $C$ may be extended to an unbounded proper chain complex).

To be entirely honest, this is in fact the very reason for introducing positively graded chain complexes. In the next section, we shall have to turn a simplicial object into a proper chain complex; but in the quasi-pointed case, such a chain complex hardly ever is unbounded.  
\end{remark}

\begin{examples}\label{Examples-Proper-Chain-Complexes}
Of course, any chain complex in an abelian category is proper, and using the fact that a subgroup need not be normal, as in Example~\ref{Example-Homology-no-Exactness}, one easily constructs chain complexes of groups that are not proper.

A proper morphism can always be considered as a positively graded proper chain complex (just add zeros). It is, for instance, well-known that a morphism $\del \colon T\to G$, part of a crossed module of groups $(T,G,\del)$, is proper. This leads to homology of (internal) crossed modules and internal categories: See Chapter~\ref{Chapter-Internal-Categories}, in particular Section~\ref{Section-Regular-Epimorphisms}.

In the next section we shall see that the Moore complex of a simplicial object is always proper, which provides us with an additional class of examples. This insight forms the basis for our study of cotriple homology in Chapter~\ref{Chapter-Cotriples}. In Chapter~\ref{Chapter-Simplicial-Objects} we show that homology of simplicial objects admits a compatible homotopy theory. 
\end{examples}

\begin{example}\label{Example-Different-H-and-K}
Now the promised example of a chain complex where the two notions of homology do not coincide. As the next proposition shows, such a complex $C$ can not be proper; but it can be a chain complex of groups:
\[
\xymatrix{{\cdots} \ar[r] & 0 \ar[r] & {\langle x \rangle} \monr^-{d_{2}} & {\langle x,y \rangle} \ar@{-{ >>}}[r]^-{d_{1}} & {\langle y \rangle} \ar[r] & 0 \ar[r] & {\cdots.}}
\]
Here the inclusion $d_{2}$ maps $x$ to $x$, and the quotient $d_{1}$ maps $x$ to $1$ and $y$ to $y$. It is easily checked that $K_{1}C=0$, but 
\[
H_{1}C\cong\langle y^{n}xy^{-n},\,n\in \Z_{0}\rangle.
\]
(Incidentally, this gives an example of a free group with an infinite number of generators as a normal subgroup of a group with just two generators, cf.~\cite{BJK}. These kinds of subobjects are also crucial in the definition of internal object action~\cite{BJK}: In the notation of Section~\ref{Section-Regular-Epimorphisms}, $K[d_{1}]=\langle y\rangle \flat \langle x \rangle$.)
\end{example}

\begin{proposition}\label{Proposition-H-K}
Let $\Ac$ be a sequentiable category. For any $n$, $H_n$ and $K_n$ are naturally isomorphic functors.
\end{proposition}
\begin{proof}
Consider the commutative diagram of solid arrows~\ref{Diagram-K-H}. Since $H_nC$ is a cokernel and $K_nC$ is a kernel, unique morphisms $j_n$ and $p_n$ exist that keep the diagram commutative. A~natural transformation $\lambda \colon H_{n}\To K_{n}$ is defined by the resulting unique $(\lambda_n)_C$. To prove it an isomorphism, first consider the following diagram with exact rows.
\[
\xymatrix{
0 \ar[r] & K[d_n] \ar@{{ >}->}[r] \ar[d]_{p_n} & C_{n} \ar[r]^-{d_{n}} \ar@{-{ >>}}[d]^{\coker d_{n+1}} & C_{n-1} \ar@{=}[d] \\
0 \ar[r] & K_{n}C \ar@{{ >}->}[r]_-{\ker d''_{n}} & \Cok[d_{n+1}] \ar[r]_-{d''_{n}} & C_{n-1}}
\]
$1_{C_{n-1}}$ being a monomorphism, the left hand square is a pullback, and $p_n$ is a regular epimorphism.

Considering the image factorizations of $d'_{n+1}$ and $d_{n+1}$, there is a morphism $i$ such that the diagram
\[
\bfig
 \square/-{ >>}`=``-{ >>}/[C_{n+1}`I{[d'_{n+1}]}`C_{n+1}`I{[d_{n+1}]};```]
 \square(500,0)|blra|/{^ (}->`.>` >->` >->/<700,500>[I{[d'_{n+1}]}`
K{[d_n]}` I{[d_{n+1}]}`C_n;\im d'_{n+1}`i`\ker d_n`\im d_{n+1}]
 \square(1200,0)/-{ >>}``>`-{ >>}/<800,500>[K{[ d_n]}`H_nC`C_n`\Cok{[d_{
n+1}]};\coker d'_{n+1}``j_n`\coker d_{n+1}]
 \morphism(0,500)/{@{->}@/^12pt/}/<1200,0>[C_{n+1}`K{[d_n]};d'_{n+1}]
 \morphism|b|/{@{->}@/_12pt/}/<1200,0>[C_{n+1}`C_n;d_{n+1}]
 \efig
\]
commutes. Clearly, it is both a monomorphism and a regular epimorphism, hence an isomorphism. Because $d_{n+1}$ is proper, $\im d_{n+1}$ is a kernel. We get that also $\im d'_{n+1}$ is a kernel. Now the middle square is a pullback, so Proposition~\ref{Proposition-LeftRightPullbacks} implies that $j_n$ is a monomorphism.

Accordingly, since Diagram~\ref{Diagram-K-H} commutes, $(\lambda_n)_C$ is both regular epi and mono, hence it is an isomorphism.
\end{proof}

\begin{remark}\label{Remark-Exact-Sequence-of-Proper-Chain-Complexes}
Note that, like $\Ac$, the categories $\Ch\Ac$ and $\Chplus\Ac$ are sequentiable. This is not the case for $\PCh\Ac$ or $\PChplus\Ac$ since, as the next example shows, these categories need not have kernels. By an \emph{exact sequence of proper chain complexes}\index{exact sequence!of proper chain complexes}, we mean an exact sequence in $\Ch\Ac$ or in $\Chplus\Ac$ such that the objects are proper chain complexes.
\end{remark}

\begin{example}\label{Example-Proper-Chain-Complexes-no-Kernels}
The \emph{dihedral group}\index{dihedral group} $D_{2n}$ of order $2n$ may be presented as
\[
\langle x,a\,|\,x^{2}=a^{n}=1,\,x^{-1}ax=a^{-1}\rangle.
\]
An interesting feature of these groups is that, although both inclusions ${D_{4}\to D_{8}}$ and ${D_{8}\to D_{16}}$ that send $a$ to $a^{2}$ and $x$ to $x$ are proper (as embeddings of a subgroup of index~$2$), their composition is not, because e.g., $axa^{-1}=a^{2}x$ is not in its image. Hence considering the vertical arrows in the diagram with exact rows
\[
\xymatrix{0 \ar[r] & D_{4} \mond \ar@{{ >}->}[r] & D_{8} \ar@{{ >}->}[d] \ar@{-{ >>}}[r]   &  \Z_{2} \ar@{-{ >>}}[d] \ar[r] & 0\\
0 \ar[r] & D_{16} \ar@{=}[r] & D_{16} \ar@{-{ >>}}[r] & 0}
\]
as chain complexes in $\Gp$ shows that a kernel in $\Ch \Gp$ of a morphism in $\PCh \Gp$ is not necessarily proper. It remains to show that no other chain complex can play the r\^ole of kernel in $\PCh \Gp$.

Let us denote the above morphism of proper chain complexes as $\overline{D}\to Z$. Suppose that $k\colon {K\to \overline{D}}$ is its kernel in $\PCh \Gp$. Then certainly all $k_{i}$ are injections (a kernel is always mono, and monomorphisms between bounded below proper chain complexes are degreewise). $k_{0}\colon {K_{0}\to \overline{D}_{0}=D_{16}}$ is also a split epimorphism, because by the universal property of kernels, the morphism of proper chain complexes induced by the left hand side square
\[
\vcenter{\xymatrix{0 \ar@{{ >}->}[d] \ar@{{ >}->}[r] & D_{8} \ar@{{ >}->}[d]\\
D_{16} \ar@{=}[r] & D_{16}}}
\qquad \qquad 
\vcenter{\xymatrix{D_{4} \ar@{{ >}->}[d] \ar@{{ >}->}[r] & D_{8} \ar@{{ >}->}[d]\\
D_{8} \ar@{{ >}->}[r] & D_{16}}}
\]
factors over $k$. Now consider the morphism in $\PCh \Gp$ induced by the right hand side square. Using that it factors over $k$, one shows that $d_{1}\colon {K_{1}\to K_{0}}$ is not proper. 
\end{example}

\begin{notation}\label{Notation-overline-d_{n}}
For any $n$, let $\overline{d}_{n}$ denote the unique map such that the diagram
\[
\xymatrix{C_{n} \ar@{-{ >>}}[d]_-{\coker d_{n+1}} \ar[r]^-{d_{n}} & C_{n-1}\\
\Cok[d_{n+1}] \ar@{.>}[r]_-{\overline{d}_{n}} & K[d_{n-1}] \ar@{{ >}->}[u]_-{\ker d_{n-1}} }
\]
commutes. Since $d_{n}$ is proper, so is $\overline{d}_n$.
\end{notation}

The following is a straightforward generalization of the abelian case---see, for instance, Theorem 1.3.1 in Weibel~\cite{Weibel} or Theorem 4.5.7 in Borceux and Bourn~\cite{Borceux-Bourn}.

\begin{proposition}\label{Proposition-Long-Exact-Homology-Sequence}
Let $\Ac$ be a sequentiable category. Any short exact sequence of proper chain complexes
\[
\xymatrix{
0 \ar[r] & C'' \ar@{{ >}->}[r] & C' \ar@{-{ >>}}[r] & C \ar[r] & 1 }
\]
gives rise to a long exact sequence of homology objects
\begin{equation}\label{Lange-Exacte-Homologierij}\index{homology!sequence}
%\resizebox{\textwidth}{!}
{\xymatrix{
{\dots} \ar[r] &H_{n+1}C \ar[r]^-{\delta_{n+1}} & H_nC'' \ar[r] & H_nC' \ar[r] & H_nC \ar[r]^-{\delta_{n}} & H_{n-1}C'' \ar[r] & {\dots}},}
\end{equation}
which depends naturally on the given short exact sequence. (When the chain complexes are positively graded, this sequence is bounded on the right.)
\end{proposition}
\begin{proof}
Since the $d_{n}$ and $\overline{d}_{n}$ are proper, mimicking the abelian proof---using the Snake Lem\-ma~\ref{snake-lemma} twice---we get an exact sequence
\[
\xymatrix{
K_nC'' \ar[r] & K_nC' \ar[r] & K_nC \ar[r] & H_{n-1}C'' \ar[r] & H_{n-1}C' \ar[r] & H_{n-1}C}
\]
for every $n$. By Proposition~\ref{Proposition-H-K} we can paste these together to Sequence~\ref{Lange-Exacte-Homologierij}. The naturality follows from the naturality of the Snake Lemma.
\end{proof}

\section{Simplicial objects}\label{Section-Simplicial-Objects}

In this section we extend the homology theory of Section~\ref{Section-Chain-Homology} to simplicial objects. We start by considering the \textit{normalization functor} $N\colon \simpA \to \Chplus\Ac$, which maps a simplicial object $A$ in a quasi-pointed category with pullbacks $\Ac$ to the \textit{Moore complex} $N (A)$.  We prove that, when $\Ac $ is, moreover, exact Mal'tsev, $N (A)$ is always proper; hence if $\Ac$ is quasi-pointed, exact and protomodular, we may define the \textit{$n$-th homology object} of $A$ as $H_{n}A=H_{n}N (A)$. Furthermore, if $\epsilon \colon {A\to A_{-1}}$ is a contractible augmented simplicial object, then $H_{0}A=A_{-1}$ and, for $n\geq 1$, $H_{n}A=0$.

When working with simplicial objects in a category $\Ac$, we shall use the notations of~\cite{Weibel}; see also~\cite{May, Maclane:Homology, KP:Abstact-homotopy-theory, Quillen}. The \emph{simplicial category}\index{category!simplicial}\index{simplicial category} $\Delta$\index{category!$\Delta$} has, as objects, finite ordinals $[n]=\{0,\dots,n \}$, for $n\in \N$ and, as morphisms, monotone functions.  The category $\simpA$\index{category!$\Sc \Ac$} of \emph{simplicial objects}\index{simplicial object} and \emph{simplicial morphisms} of $\Ac$ is the functor category $\Fun (\Delta^{\op },\Ac)$. Thus a simplicial object $A\colon {\Delta^{\op}\to \Ac}$ corresponds to the following data: a sequence of objects $(A_{n})_{n\in \N}$, \emph{face operators}\index{face operator} $\del_{i}\colon {A_{n}\to A_{n-1}}$ for $i\in [n]$ and $n\in \N_{0}$, and \emph{degeneracy operators}\index{degeneracy operator} $\sigma_{i}\colon {A_{n}\to A_{n+1}}$, for $i\in [n]$ and $n\in \N$, subject to the \emph{simplicial identities}
\[
\begin{aligned}
\del_{i}\comp \del_{j} &=\del_{j-1}\comp \del_{i}\quad \text{if $i<j$}\\
\sigma_{i}\comp \sigma_{j} &= \sigma_{j+1}\comp \sigma_{i}\quad \text{if $i\leq j$}
\end{aligned}
\qquad\qquad 
\del_{i}\comp\sigma_{j}=\begin{cases}\sigma_{j-1}\comp \del_{i} & \text{if $i<j$} \\
1 & \text{if $i=j$ or $i=j+1$}\\
\sigma_{j}\comp \del_{i-1} & \text{if $i>j+1$.}\end{cases}
\]
Dually, a \emph{cosimplicial object}\index{cosimplicial object} in $\Ac$ is a functor $\Delta \to \Ac$. An \emph{augmented}\index{augmented simplicial object}\index{simplicial object!augmented} simplicial object $\epsilon \colon A\to A_{-1}$ consists of a simplicial object $A$ and a map $\epsilon \colon {A_{0}\to A_{-1}}$ with $\epsilon \comp \del_{0}= \epsilon \comp \del_{1}$. It is \emph{contractible}\index{contractible augm.\ simplicial object}\index{simplicial object!augmented!contractible} when there exist morphisms $f_{n}\colon {A_{n}\to A_{n+1}}$, $n\geq -1$, with $\epsilon \comp f_{-1}=1_{A_{-1}}$, $\del_{0}\comp f_{0}=f_{-1}\comp \epsilon$, $\del_{n+1}\comp f_{n}=1_{A_{n}}$ and $\del_{i}\comp f_{n}=f_{n-1}\comp \del_{i}$, for $0\leq i \leq n$ and $n\in \N$.

\begin{figure}
\begin{center}
\resizebox{\textwidth}{!}{\includegraphics{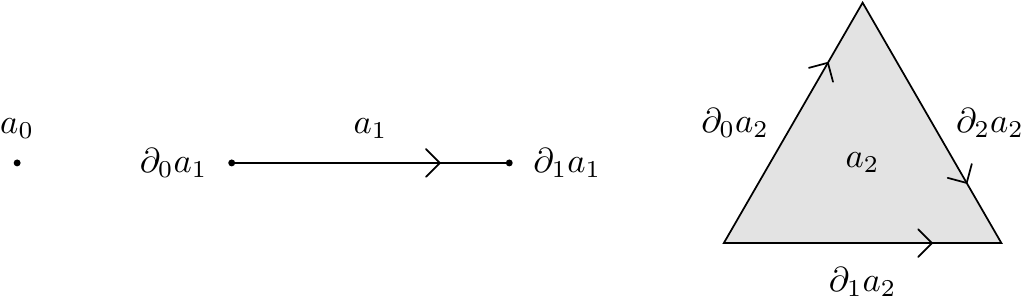}\qquad \includegraphics{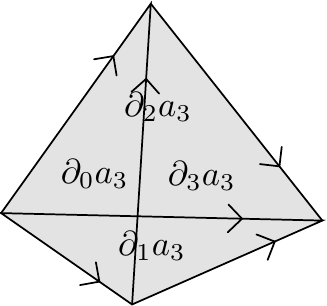}}
\end{center}
\caption{A $0$-simplex $a_{0}$, a $1$-simplex $a_{1}$, a $2$-simplex $a_{2}$ and a $3$-simplex $a_{3}$.}\label{Figure-What-is-Simplicial-Set}
\end{figure}

\begin{figure}
\begin{center}
\resizebox{\textwidth}{!}{\includegraphics{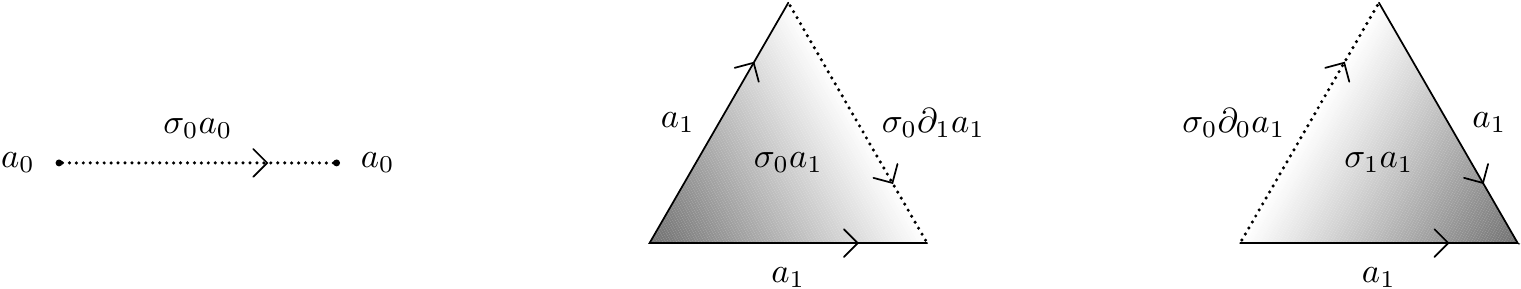}}
\end{center}
\caption{A degenerate $1$-simplex and two degenerate $2$-simplices.}\label{Figure-What-is-Simplicial-Set-2}
\end{figure}

\begin{remark}\label{Remark-What-is-Simplicial-Set}
When working with simplicial objects the following mental picture of a simplicial set $A$ might be helpful. A $0$-simplex $a_{0}$ in $A$ is a point, a $1$-simplex $a_{1}$ a segment between two points $\del_{0}a_{1}$ and $\del_{1}a_{1}$, a $2$-simplex $a_{2}$ is a triangle with three faces $\del_{0}a_{2}$, $\del_{1}a_{2}$ and $\del_{2}a_{2}$, a $3$-simplex is a tetrahedron, etc.\ (see~\fref{Figure-What-is-Simplicial-Set}). The degeneracies turn a simplex into a degenerate simplex of a higher degree (see~\fref{Figure-What-is-Simplicial-Set-2}).
\end{remark}

\begin{examples}\label{Examples-Simplicial-Objects}
A first classical example of a simplicial set is the \emph{total singular complex}\index{total singular complex} $S (X)$ of a topological space $X$. For $n\in \N$, let $\Delta_{n}$ denote the subspace 
\[
\{ (x_{0},\dots ,x_{n})\,|\,x_{i}\in [0,1],\, x_{0}+\cdots +x_{n}=1\}
\]
of $\R^{n+1}$. If $X$ is a topological space, $S (X)_{n}$ consists of all continuous functions $f\colon {\Delta_{n}\to X}$. The faces and degeneracies of $S (X)$ are defined by 
\[
(\del_{i}f) (x_{0}, \dots, x_{n-1})=f (x_{0},\dots ,x_{i-1}, 0,x_{i},\dots ,x_{n-1})
\]
and $(\sigma_{i}f) (x_{0},\dots ,x_{n+1})=f (x_{0},\dots ,x_{i-1},x_{i}+x_{i+1},x_{i+2},\dots,x_{n+1})$. A~topological space gives rise to a simplicial abelian group (or module over a ring $R$) by applying the free abelian group functor (or the free $R$-module functor) to $S (X)$. One then associates a chain complex to this simplicial object, and defines the homology of the space $X$\index{homology!of a topological space} as the homology of this complex (cf.\ Definition~\ref{Definition-Moore-Complex}, Definition~\ref{Definition-Simplicial-Homology} and Remark~\ref{Remark-Simplicial-Homology-in-Abelian-Case}).

Any small category $\Af$ determines a simplicial set $\ner \Af$\index{ner@$\ner \Af$} called its \emph{nerve}\index{nerve!of a small category}, by putting $\ner \Af_{n}=\Cat ([n],\Af)$, for $n\in \N$ (here we consider an ordinal $[n]$ as a small category). A~$0$-simplex of $\ner \Af$ is an object of $\Af$, a $1$-simplex is a morphism of $\Af$ (with its codomain and domain as faces), a $2$-simplex is a pair of composable arrows, an $n$-simplex a string of $n$ composable arrows. The faces are given by composition of two consecutive arrows in a string, and the degeneracies by insertion of an identity (see Section~\ref{Section-Regular-Epimorphisms} for an internal version of this construction).
\end{examples}

\begin{definition}\label{Definition-Moore-Complex}
Let $A$ be a simplicial object in a quasi-pointed category $\Ac$ with pullbacks. The \emph{Moore complex}\index{Moore complex}\index{chain complex!Moore ---} $N (A)$\index{N@$N (A)$} is the chain complex with $N_{0}A=A_{0}$,
\[
N_{n} A=\bigcap_{i=0}^{n-1}K[\del_{i}\colon A_{n}\to A_{n-1}]=K[(\del_{i})_{i\in [n-1]}\colon A_{n}\to A_{n-1}^{n}]
\]
and boundary operators $d_{n}=\del_{n}\comp \bigcap_{i}\ker \del_{i}\colon N_{n} A\to N_{n-1} A$, for $n\geq 1$. This gives rise to the \emph{normalization functor}\index{normalization functor} $N\colon \simpA\to \Chplus \Ac$.
\end{definition}

\begin{remark}\label{Remark-N-0}
Note that, in the above definition, 
\[
\del_{n}\comp \bigcap_{i}\ker \del_{i}\colon N_{n} A\to A_{n-1}
\]
may indeed be considered as an arrow $d_{n}\colon N_{n}A\to N_{n-1} A$: The map clearly factors over 
\[
\bigcap_{i}\ker \del_{i}\colon {N_{n-1} A\to A_{n-1}}.
\]
\end{remark}

\begin{remark}\label{Remark-N-1}
Obviously, for $n\geq 1$, the object of \emph{$n$-cycles}\index{cycle} $Z_{n}A=K[d_{n}]$\index{Z_n@$Z_{n}A$} of a simplicial object $A$ of $\Ac$ is equal to $\bigcap_{i=0}^{n}K[\del_{i}\colon A_{n}\to A_{n-1}]$. Write $Z_{0}A=A_{0}$.
\end{remark}

\begin{remark}\label{Remark-N-2}
The functor $N\colon \simpA\to \Chplus \Ac$ preserves limits. Indeed, limits in $\simpA$ and $\Chplus \Ac$ are computed degreewise, and taking kernels and intersections (pulling back), as occurs in the construction of $N$, commutes with taking arbitrary limits in $\Ac$. In Section~\ref{Section-Implications-Kan} we shall prove that $N$, moreover, preserves regular epimorphism, and hence is exact.
\end{remark}

If we recall from Section~\ref{Section-Bourn-Protomodular-Categories} how to normalize an equivalence relation, the following comes as no surprise.

\begin{theorem}\label{Theorem-N-Proper}
Let $\Ac $ be a quasi-pointed exact Mal'tsev category and $A$ a simplicial object in~$\Ac $. Then $N (A)$ is a proper chain complex of $\Ac$.
\end{theorem} 
\begin{proof}
Any boundary operator $d_{n}$, when viewed as an arrow to $A_{n-1}$, is a composition of a proper monomorphism (an intersection of kernels), and a regular epimorphism (a split epimorphism, by the simplicial identities). Hence Proposition~\ref{Proposition-Image-of-Kernel} implies that $d_{n}\colon N_{n} A\to A_{n-1}$ is proper. This clearly remains true when, following Remark~\ref{Remark-N-0}, we consider $d_n$ as an arrow to $N_{n-1} A$.
\end{proof}

\begin{definition}\label{Definition-Simplicial-Homology}
Suppose that $\Ac$ is quasi-pointed, exact and protomodular, and $A$ is a simplicial object in $\Ac$. The object $H_{n} A=H_{n} N (A)$\index{H_n@$H_{n}A$} will be called the \emph{$n$-th homology object of $A$}\index{homology!of simplicial objects}, and the resulting functor 
\[
H_{n}\colon {\simpA\to \Ac}
\]
the \emph{$n$-th homology functor}, for $n\in \N$.
\end{definition}

\begin{figure}
\begin{center}
\resizebox{\textwidth}{!}{\includegraphics{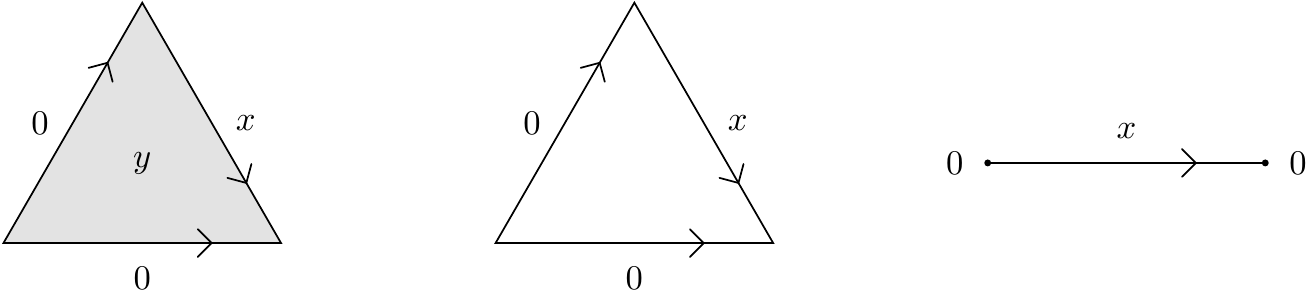}}
\end{center}
\caption{An element $y$ of $N_{2}A$, its \emph{boundary}\index{boundary of a simplex} $(\del_{0} y, \del_{1}y, \del_{2}y)$, and a $1$-cycle $x=\del_{2}y$.}\label{Figure-Meaning-Homology}
\end{figure}

\begin{remark}\label{Remark-Meaning-Homology}
A homology object $H_{n}A$ is zero if and only if the morphism $d_{n+1}\colon {N_{n+1}A\to Z_{n}A}$ is a regular epimorphism. For a simplicial group $A$ and $n=1$, this may be pictured as in~\fref{Figure-Meaning-Homology}: An element $y$ of $N_{2} A$ is mapped by its face operators on $0$, $0$ and $x=\del_{2}y$, respectively. Hence $\del_{0}x=\del_{1}x=0$, which means that $x\in Z_{1}A$. If now $H_{1}A=0$ then for every $x$ in $Z_{1}A$ there exists a $y\in N_{2}A$ such that $\del_{2}y=x$. Thus one may think of $H_{1}A$ as an obstruction to certain triangles being filled. (See also Section~\ref{Section-Weak-Equiv-Trivial-Fibrations}.)
\end{remark}

\begin{remark}\label{Remark-Simplicial-Homology-in-Abelian-Case}
In the abelian case there is an other way of constructing a chain complex out of a simplicial object. Given a simplicial object~$A$, the \emph{unnormalized chain complex}\index{unnormalized chain complex}\index{chain complex!unnormalized} $C (A)$\index{C@$C (A)$} of $A$ is defined by $C (A)_{n}=A_{n}$ and $d_{n}=\del_{0}-\del_{1}+\dots + (-1)^{n}\del_{n}$. In general, $C (A)$ and $N (A)$ do not coincide, but their homology objects are isomorphic (see e.g., \cite{Weibel}).
\end{remark}

In order to compute, in Proposition~\ref{Proposition-Contractible-Object}, the homology of a contractible augmented simplicial object, we first make the following, purely categorical, observations. We start by recalling a result due to Dominique Bourn.

\begin{lemma}\cite{Bourn1991}\label{Lemma-Pullback-Along-Regular-Epi}
If, in a protomodular category, a square with vertical regular epimorphisms
\[
\xymatrix{ {\cdot} \ar[r] \ar@{->>}[d] & {\cdot} \ar@{->>}[d] \\
{\cdot} \ar[r] & {\cdot} }
\]
is a pullback, it is also a pushout.\noproof
\end{lemma}

A \emph{fork}\index{fork} is a diagram such as~\ref{Diagram-Fork} below where $e\comp \del_{0}=e\comp \del_{1}$.

\begin{proposition}\label{Proposition-Forks}
In a quasi-pointed and protomodular category $\Ac$, let
\begin{equation}\label{Diagram-Fork}
\xymatrix{A \ar@<.5ex>[r]^-{\del_{1}} \ar@<-.5ex>[r]_-{\del_{0}} & B \ar[r]^-{e} & C}
\end{equation}
be a fork and $t\colon B\to A$ a map with $\del_{0}\comp t=\del_{1}\comp t=1_{B}$. Then the following are equivalent:
\begin{enumerate}
\item $e$ is a coequalizer of $\del_{0}$ and $\del_{1}$;
\item the square
\[
\xymatrix{A \ar[d]_-{\del_{0}} \ar[r]^-{\del_{1}} & B \ar[d]^-{e}\\
B \ar[r]_-{e} & C}
\]
is a pushout;
\item $e$ is a cokernel of $\del_{1}\comp \ker \del_{0}$.
\end{enumerate}
\end{proposition}
\begin{proof}
The equivalence of 1.\ and 2.\ is obvious. In the diagram
\[
\xymatrix{K[\del_{0}] \ar@{}[rd]|<<{\pullback}\ar[d] \ar@{{ >}->}[r]^-{\ker \del_{0}} & A \ar@{-{ >>}}[d]_-{\del_{0}} \ar[r]^-{\del_{1}} & B \ar[d]^-{e}\\
0 \ar[r] & B \ar[r]_-{e} & C,}
\]
the left hand side square is a pushout: Indeed, Lemma~\ref{Lemma-Pullback-Along-Regular-Epi} applies, since it is a pullback along the split, hence regular, epimorphism $\del_{0}$. Consequently, the outer rectangle is a pushout if and only if the right square is, which means that 2.\ and 3.\ are equivalent.
\end{proof}

\begin{corollary}\label{Corollary-H_0}
If $\Ac$ is a quasi-pointed, exact and protomodular category and $A$ a simplicial object of $\Ac$ (with face operators $\del_{0},\del_{1}\colon A_{1}\to A_{0}$), then $H_{0}A=\Coeq [\del_{0},\del_{1}]$.\noproof
\end{corollary}

\begin{remark}\label{Remark-Forks-Normalization}
This is related to the process of normalizing an equivalence relation in the following way. If $\Ac$ is exact and $(\del_{0},\del_{1})$ is an equivalence relation, $\del_{1}\comp \ker \del_{0}$ is its normalization, and Proposition~\ref{Proposition-Forks} just shows that a coequalizer of an equivalence relation is a cokernel of the normalization. Conversely, $\del_{1}\comp \ker \del_{0}$ is a kernel (of its cokernel $e$) if and only if it is a monomorphism (Proposition~\ref{Proposition-Image-of-Kernel}), but by Proposition~\ref{Proposition-Protomodularity-I} this is the case exactly when $(\del_{0},\del_{1})$ is a kernel pair (of~$e$).
\end{remark}

A fork is \emph{split}\index{fork!split}\index{split fork} when there are two more arrows $s\colon {C\to B}$, $t\colon {B\to A}$ such that $e\comp s=1_{C}$, $\del_{1}\comp t=1_{B}$ and $\del_{0}\comp t=s\comp e$. Every split fork is a coequalizer diagram.

\begin{proposition}\label{Proposition-Contractible-Object}\index{simplicial object!augmented!contractible}
If $\epsilon \colon A\to A_{-1}$ is a contractible augmented simplicial object in a quasi-pointed, exact and protomodular category $\Ac$, then $H_{0}A=A_{-1}$ and, for $n\geq 1$, $H_{n}A=0$.
\end{proposition}
\begin{proof}
The contractibility of $\epsilon \colon A\to A_{-1}$ implies that the fork
\[
\xymatrix{A_{1} \ar@<.5ex>[r]^-{\del_{1}} \ar@<-.5ex>[r]_-{\del_{0}} & A_{0} \ar[r]^-{\epsilon } & A_{-1}}
\]
is split (by the arrows $f_{-1}\colon A_{-1}\to A_{0}$ and $f_{0}\colon A_{0}\to A_{1}$). We get that it is a coequalizer diagram. The first equality now follows from Corollary~\ref{Corollary-H_0}.

In order to prove the other equalities, first recall Remark~\ref{Remark-N-1} that, for $n\geq 1$,
\[
Z_{n}A=K[d_{n}\colon N_{n}A\to N_{n-1}A]=\bigcap_{i=0}^{n}K[\del_{n}\colon A_{n}\to A_{n-1}].
\]

We are to show that the image of the morphism $d_{n+1}\colon N_{n+1}A\to N_{n}A$ is a kernel $\ker d_{n}\colon {K[d_{n}]\to N_{n}A}$ of $d_{n}$. But, for any $i\leq n$, the left hand downward-pointing arrow in the diagram with exact rows
\[
\xymatrix{0 \ar[r] & K[\del_{i}] \ar@{{ >}->}[r]^-{\ker \del_{i}} \ar@{.>}@<-0.5 ex>[d] & A_{n+1} \ar@<-0.5 ex>@{-{ >>}}[d]_-{\del_{n+1}} \ar[r]^-{\del_{i}} & A_{n} \ar[r] \ar@<-0.5 ex>@{-{ >>}}[d]_-{\del_{n}} & 1\\
0 \ar[r] & K[\del_{i}] \ar@{{ >}->}[r]_-{\ker \del_{i}} \ar@{.>}@<-0.5 ex>[u] & A_{n} \ar[r]_-{\del_{i}} \ar@<-0.5 ex>[u]_{f_{n}} & A_{n-1} \ar@<-0.5 ex>[u]_{f_{n-1}} \ar[r] & 1}
\]
is a split epimorphism, because both its upward and downward pointing squares commute. It follows that the intersection 
\[
N_{n+1}A=\bigcap_{i\in [n]}K[\del_{i}]\to \bigcap_{i\in [n]}K[\del_{i}]=K[d_{n}]
\]
is a split, hence a regular, epimorphism, and $\im d_{n+1}=\ker d_{n}$.
\end{proof}

\section{The Kan property}\label{Section-Kan}

\textit{Kan} simplicial objects and \textit{Kan fibrations} are very important in the homotopy theory of simplicial sets (or simplicial groups, algebras, etc.). In their article~\cite{Carboni-Kelly-Pedicchio}, Carboni, Kelly and Pedicchio extend the notion of Kan simplicial set to an arbitrary category $\Ac$. When $\Ac$ is regular, their definition amounts to the one stated in Definition~\ref{Definition-Kan} below. (They only consider horns with $n\geq 2$; indeed, any simplicial object fulfils the Kan property for $n=1$.) In the same spirit, we propose an extension of the notion of Kan fibration to a regular category $\Ac$.

The results in this section will allow us to prove, in Section~\ref{Section-Implications-Kan}, two important facts concerning $H_{n}\colon \PChplus \Ac \to \Ac $ and $N\colon {\Sc \Ac \to \PChplus\Ac }$: Proposition~\ref{Proposition-N-Exact} and Theorem~\ref{Theorem-Homology-Abelian}. As such, the Kan property may seem a very technical tool; in Chapter~\ref{Chapter-Simplicial-Objects}, we shall try to clarify its meaning for homotopy theory of simplicial objects. 

\begin{definition}\label{Definition-Kan}
Consider a simplicial object $K$ in a regular category $\Ac$. For $n\geq 1$ and $k\in [n]$, a family
\[
x= (x_{i}\colon X\to K_{n-1})_{i\in [n],i\neq k}
\]
is called an \emph{$(n,k)$-horn}\index{horn}\index{n,k-horn@$(n,k)$-horn} of $K$ if it satisfies $\del_{i}\comp x_{j}=\del_{j-1}\comp x_{i}$, for $i<j$ and $i,j\neq k$.

We say that \emph{$K$ is Kan}\index{Kan simplicial object}\index{simplicial object!Kan ---} if, for every $(n,k)$-horn $x= (x_{i}\colon X\to K_{n-1})$ of $K$, there is a regular epimorphism $p\colon {Y\to X}$ and a map $y\colon {Y\to K_{n}}$ such that $\del_{i}\comp y=x_{i}\comp p$ for $i\neq k$.

A map $f \colon A\to B$ of simplicial objects is said to be a \emph{Kan fibration}\index{Kan fibration}\index{fibration!Kan ---} if, for every $(n,k)$-horn $x= (x_{i}\colon X\to A_{n-1})$ of $A$ and every $b\colon X\to B_{n}$ with $\del_{i}\comp b=f_{n-1}\comp x_{i}$ for all $i\neq k$, there is a regular epimorphism $p\colon {Y\to X}$ and a map $a\colon {Y\to A_{n}}$ such that $f_{n}\comp a=b\comp p$ and $\del_{i}\comp a=x_{i}\comp p$ for all $i\neq k$.
\end{definition}

\begin{figure}
\begin{center}
\resizebox{\textwidth}{!}{\includegraphics{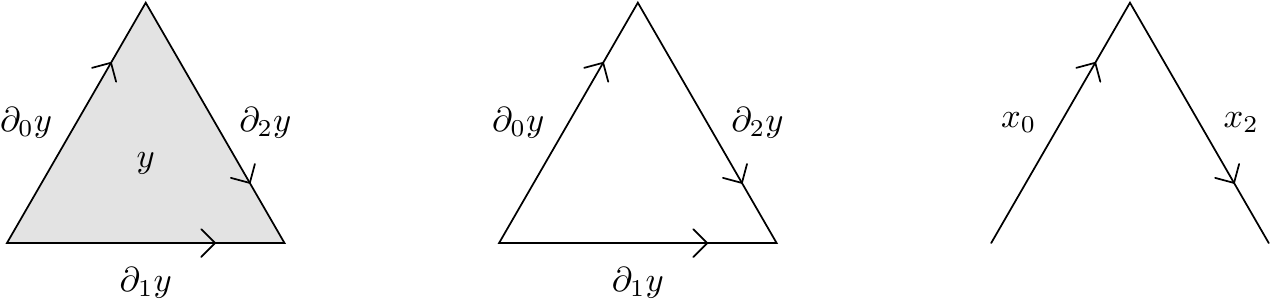}}
\end{center}
\caption{A $2$-simplex $y$, its boundary $\del y$ and a $(2,1)$-horn $x$.}\label{Figure-Kan}
\end{figure}

\begin{remark}\label{Remark-Kan-Simplicial-Sets}
In case $\Ac$ is $\Set$, these notions have an equivalent formulation in which $X$ and $Y$ are both equal to a terminal object $1$. These equivalent formulations are the classical definitions of Kan simplicial set and Kan fibration---see, for instance, Weibel~\cite{Weibel}. The Kan property is usually pictured as follows. The \emph{boundary}\index{boundary of a simplex} of an $n$-simplex $y$ of a simplicial set $K$ is the set of its faces $\del y=\{\del_{i}y\,|\,i\in [n] \}$. If we delete the $k$-th face we get an $(n,k)$-horn. (See~\fref{Figure-Kan}.) So a horn may be the boundary of a simplex, but with one face missing. If, conversely, (up to change of domain) every horn comes from the boundary of some simplex, then $K$ is Kan. Also, thus filling a horn gives a way of replacing a missing face (but since there is no uniqueness requirement in the Kan property, in general, several replacements are possible).
\end{remark}

\begin{remark}\label{Remark-Kan-vs.-Kan-Fibration}
Obviously, a simplicial object $K$ is Kan if and only if the unique map $K\to 1$ from $K$ to a terminal object $1$ of $\Sc\Ac$ is a Kan fibration.
\end{remark}

\begin{remark}\label{Remark-Kan-Projectives}
If the category $\Ac$ has enough regular projectives, the change of domain in this definition (the regular epimorphism $p$) may be interpreted as follows: $K$ is Kan if and only if for every projective object $P$ and for every $(n,k)$-horn $x= (x_{i}\colon P\to K_{n-1})$ of $K$, there is a morphism $y\colon {P\to K_{n}}$ such that $\del_{i}\comp y=x_{i}$, $i\neq k$.
\end{remark}

According to M. Barr~\cite{Barr}, a Mal'tsev category is a regular one in which every simplicial object is Kan. The equivalence with the definition from Section~\ref{Section-Mal'tsev-Categories} was shown by Carboni, Kelly and Pedicchio in their paper~\cite{Carboni-Kelly-Pedicchio}. Thus they extend Moore's classical result that every simplicial group is Kan to the categorical context. With Proposition~\ref{Proposition-Kan} below, we add that every regular epimorphism between simplicial objects in a regular Mal'tsev category is a Kan fibration. In its proof, we make heavy use of the techniques recalled with Proposition~\ref{Proposition-Regularity}.

\begin{proposition}\label{Proposition-Kan}
Let $\Ac$ be a regular Mal'tsev category. Then
\begin{enumerate}
\item every simplicial object $K$ of $\Ac$ is Kan;
\item if $f \colon A\to B$ is a regular epimorphism between simplicial objects of $\Ac$, then it is a Kan fibration.
\end{enumerate}
\end{proposition}
\begin{proof}
The first statement is part of Theorem 4.2 in~\cite{Carboni-Kelly-Pedicchio}; we prove the second one. For $n\geq 1$ and $k\in [n]$, let
\[
x= (x_{i}\colon X\to A_{n-1})_{i\in [n],i\neq k}
\]
be an $(n,k)$-horn of $A$ and let $b\colon X\to B_{n}$ be a map with $\del_{i}\comp b=f_{n-1}\comp x_{i}$, for $i\in [n]$ and $i\neq k$. Because $f_{n}$ is  regular epi, there is a regular epimorphism $p_{1}\colon {Y_{1}\to X}$ and a map $c\colon {Y_{1}\to A_{n}}$ with $f_{n}\comp c=b\comp p_{1}$. For $i\in [n]$ and $i\neq k$, put $c_{i}=\del_{i}\comp c\colon {Y_{1}\to A_{n-1}}$. Let $\KP [f]$ denote the kernel relation of $f$. Now
\[
f_{n-1}\comp c_{i}=f_{n-1}\comp \del_{i}\comp c=\del_{i}\comp f_{n}\comp c=\del_{i}\comp b\comp p_{1}=f_{n-1}\comp x_{i}\comp p_{1},
\]
and, consequently, $c_{i} (\KP [f_{n-1}])x_{i}\comp p_{1}$. By the simplicial identities, this defines an $(n,k)$-horn
\[
((c_{i},x_{i}\comp p_{1})\colon Y_{1}\to \KP [f]_{n-1})_{i\in [n],i\neq k}
\]
of $\KP [f]$. This simplicial object being Kan yields a regular epi $p_{2}\colon {Y_{2}\to Y_{1}}$ and a $(d,e)\colon Y_{2}\to \KP [f]_{n}$ such that $\del_{i}\comp d=c_{i}\comp p_{2}$ and $\del_{i}\comp e=x_{i}\comp p_{1}\comp p_{2}$, for $i\in [n]$ and $i\neq k$, and $f_{n}\comp d=f_{n}\comp e$. Hence, $c\comp p_{2} (D)d$, where $D$ is the equivalence relation $\bigwedge_{i\in [n],i\neq k}D_{i}$ and $D_{i}$ is the kernel relation of $\del_{i}\colon A_{n}\to A_{n-1}$. It follows that $c\comp p_{2} (D\KP [f_{n}])e$. By the Mal'tsev property, $D\KP [f_{n}]$ is equal to $\KP [f_{n}]D$. This, in turn, implies that there exists a regular epimorphism $p_{3}\colon {Y\to Y_{2}}$ and a map $a\colon Y\to A_{n}$ such that $a (D)e\comp p_{3}$ and $c\comp p_{2}\comp p_{3} (\KP [f_{n}])a$. The required maps are now $a$ and $p=p_{1}\comp p_{2}\comp p_{3}\colon Y\to X$: Indeed, $f_{n}\comp a=f_{n}\comp c\comp p_{2}\comp p_{3}=b\comp p$ and $\del_{i}\comp a=\del_{i}\comp e\comp p_{3}=x_{i}\comp p$.
\end{proof}

\section{Two applications of the Kan property}\label{Section-Implications-Kan}
In this section we consider two important applications of Proposition~\ref{Proposition-Kan}. We show that for a simplicial object $A$ of $\Ac$, being Kan implies that $H_{n}A$ is abelian ($n\geq 1$). But first we prove that $N\colon \Sc \Ac \to \PChplus \Ac $ is an exact functor. This is a consequence of the fact that every regular epimorphism between simplicial objects is a Kan fibration. Finally, using the exactness of $N$, we prove that every short exact sequence of simplicial objects induces a long exact homology sequence (Proposition~\ref{Proposition-Long-Exact-Homology-Sequence-Simp}). We obtain it as an immediate consequence of Proposition~\ref{Proposition-Long-Exact-Homology-Sequence}.

\subsection*{The normalization functor is exact}\label{Subsection-Normalization-Functor-is-Exact}
Recall from Remark~\ref{Remark-N-2} that $N$ preserves kernels. If the category $\Ac$ is regular Mal'tsev, Proposition~\ref{Proposition-Kan} implies that $N$ preserves regular epimorphisms. If, moreover, $\Ac $ is protomodular, all regular epimorphisms are cokernels, and we obtain

\begin{proposition}\label{Proposition-N-Exact}\index{exact functor}
If $\Ac$ is an exact and sequentiable category, then the functor 
\[
N\colon \simpA\to \PChplus \Ac
\]
is exact.
\end{proposition}
\begin{proof}
We prove that $N$ preserves regular epimorphisms. To do so, let us consider a regular epimorphism $f\colon {A\to B}$ between simplicial objects of $\Ac$. Then $N_{0}f=f_{0}$ is regular epi. For $n\geq 1$, let $b\colon {X\to N_{n}B}$ be a map. We are to show---see Proposition~\ref{Proposition-Regularity}---that there is a regular epimorphism $p\colon {Y\to X}$ and a map $a\colon {Y\to N_{n}A}$ with $b\comp p=N_{n}f\comp a$. Now 
\[
(x_{i}=0\colon {X\to A_{n-1}})_{i\in [n-1]}
\]
is an $(n,n)$-horn of $A$ with $\del_{i}\comp \bigcap_{i\in [n-1]}\ker \del_{i}\comp b=0=f_{n-1}\comp x_{i}$, for $i\in {[n-1]}$. By Proposition~\ref{Proposition-Kan}, $f$ is a Kan fibration, which implies that there is a regular epimorphism $p\colon {Y\to X}$ and a map $a'\colon {Y\to A_{n}}$ such that $f_{n}\comp a'=\bigcap_{i\in [n-1]}\ker \del_{i}\comp b\comp p$ and $\del_{i}\comp a'=x_{i}\comp p$. Then the unique factorization $a\colon {Y\to N_{n}A}$ of $a'$ over 
\[
\bigcap_{i\in [n-1]}\ker \del_{i}\colon {N_{n}A\to A_{n}}
\]
is the required map.
\end{proof}

Proposition~\ref{Proposition-Long-Exact-Homology-Sequence} may now immediately be extended to simplicial objects.

\begin{proposition}\label{Proposition-Long-Exact-Homology-Sequence-Simp}\index{homology!sequence}
Let $\Ac$ be an exact and sequentiable category. Any short exact sequence of simplicial objects
\[
\xymatrix{
0 \ar[r] & A'' \ar@{{ >}->}[r] & A' \ar@{-{ >>}}[r] & A \ar[r] & 1 }
\]
gives rise to a long exact sequence of homology objects
\[
\resizebox{\textwidth}{!}{\xymatrix@R=1pc{
{\dots} \ar[r] &H_{n+1}A \ar[r]^-{\delta_{n+1}} & H_nA'' \ar[r] & H_nA' \ar[r] & H_nA \ar[r]^-{\delta_{n}} & H_{n-1}A'' \ar[r] & {\dots}\\
{\dots} \ar[r] &H_{1}A \ar[r]_-{\delta_{1}} & H_0A'' \ar[r] & H_0A' \ar@{-{ >>}}[r] & H_0A \ar[r]_-{\delta_{0}} & 1,}}
\]
which depends naturally on the given short exact sequence of simplicial objects.\noproof
\end{proposition}

Although we shall not use it, we think it worth mentioning that this result can be formulated in terms of homological $\delta$-functors. By a \emph{(universal) homological $\delta $-functor}\index{homological $\delta$-functor}\index{universal homological $\delta$-functor} between pointed, exact and protomodular categories $\Ac$ and $\Bc$ we mean a collection of functors 
\[
(T_{n}\colon \Ac \to \Bc)_{n\in \N}
\]
that preserve binary products, together with a collection of connecting morphisms $(\delta_{n})_{n\in \N}$ as in~\cite[Definition~2.1.1 and~2.1.4]{Weibel}.

\begin{notation}\label{Notation-A^-}\index{A@$A^{-}$}
For any simplicial object $A$ in $\Ac$, let us denote by $A^-$ the simplicial object defined by $A^-_n=A_{n+1}$,  $\del^-_i=\del_{i+1}\colon {A_{n+1}\to A_n}$, and $\sigma^-_i=\sigma_{i+1}\colon {A_{n+1}\to A_{n+2}}$, for $i\in [n]$,  $n\in \N$. This is the simplicial object obtained from $A$ by leaving out $A_{0}$ and, for $n\in \N$, all $\del_0\colon {A_n\to A_{n-1}}$ and $\sigma_0\colon {A_n\to A_{n+1}}$. Observe that $\del=(\del_{0})_n$ defines a simplicial morphism from $A^-$ to $A$.
\end{notation}

\begin{proposition}\label{Proposition-Hom-Delta-Functor}
Let $\Ac$ be a pointed, exact and protomodular category. The sequence of functors $(H_{n}\colon \simpA\to \Ac)_{n\in \N}$, together with the connecting morphisms $(\delta_{n})_{n\in \N}$, form a universal homological $\delta$-functor.
\end{proposition}
\begin{proof}
To prove that a functor $H_{n}\colon \simpA\to \Ac$ preserves binary products, it suffices that, for proper complexes $C$ and $C'$ in $\Ac$, $H_{n} (C\times C')=H_{n}C\times H_{n}C'$. One shows this by using that for $f$ and $f'$ proper, $\coker (f\times f')=\coker f\times \coker f'$. This follows from the fact that in any regular category, a product of two regular epimorphisms is regular epi. The universality is proven by modifying the proof of Theorem 2.4.7 in~\cite{Weibel}, replacing, for a simplicial object $A$, the projective object $P$ by $A^{-}$.
\end{proof}

\subsection*{Homology objects are abelian}\label{Subsection-Homology-Objects-Abelian}

This result was suggested to us by Dominique Bourn. Our proof is based on the classical proof for groups and makes heavy use of the Kan property. A more conceptual and slightly more general version of Theorem~\ref{Theorem-Homology-Abelian} may be found in Bourn's paper~\cite{BournBodeux}.

Let $A$ be a Kan simplicial object in a pointed and regular category $\Ac$ and $n\geq 1$. Recall from Remark~\ref{Remark-N-1} the notation $Z_{n}A=K[d_{n}]=\bigcap_{i=0}^{n} K[\del_{i}]$. We write $z_{n}\colon Z_{n}A\to A_{n}$ for the inclusion $\bigcap_{i=0}^{n} \ker\del_{i}$. Basing ourselves on Weibel~\cite{Weibel} (see also Definition~\ref{Definition-Homotopy-Groups-Kan}), we say that morphisms $x\colon {X\to Z_{n}A}$ and $x'\colon {X'\to Z_{n}A}$ are \emph{homotopic}, and write $x\sim x'$, if there is an arrow $y\colon {Y\to A_{n+1}}$ (called a \emph{homotopy from $x$ to $x'$}\index{homotopy!of cycles}) and regular epimorphisms $p\colon {Y\to X}$ and $p'\colon {Y\to X'}$ such that
\[
\del_{i}\comp y = \begin{cases}
0,  & \text{if $i < n$;}\\
z_{n}\comp x\comp p,& \text{if $i = n$;}\\
z_{n}\comp x'\comp p',& \text{if $i = n+1$}.
\end{cases}
\]

\begin{proposition}\label{Proposition-Pseudo-Homotopy}
Let $A$ be a Kan simplicial object in a pointed and regular category $\Ac$. For $n\geq 1$, $\sim$ defines an equivalence relation on the class of arrows $x\colon X\to  Z_{n}A$ in $\Ac $ with codomain $Z_{n}A$.
\end{proposition}
\begin{proof}
Due to to fact that $A$ is Kan, the proof of 8.3.1 in Weibel~\cite{Weibel} may be copied: One just considers arrows with codomain $A_{n}$ instead of elements of $A_{n}$ and reads ``$0$'' instead of ``$*$''.
\end{proof}

\begin{theorem}\label{Theorem-Homology-Abelian}
Let $A$ be a simplicial object in a semi-abelian category $\Ac$. For all $n\geq 1$, $H_{n}A$ is an abelian object of $\Ac$.
\end{theorem}
\begin{proof}
We use the commutator from Proposition~\ref{Proposition-Commutator}. Let $n\geq 1$. Consider the subobjects $k_{n}\colon K_{n}=[z_{n},z_{n}]\to A_{n}$ of $A_{n}$ and
\[
s_{n}\colon S_{n}=[\sigma_{n-1}\comp z_{n},\sigma_{n}\comp z_{n}]\to A_{n+1}
\]
of $A_{n+1}$. By the second statement of Proposition~\ref{Proposition-Commutator},
\[
\del_{i}S_{n}=[\del_{i}\comp \sigma_{n-1}\comp z_{n}, \del_{i}\comp \sigma_{n}\comp z_{n}],
\]
for $0\leq i\leq n+1$. Hence, by the simplicial identities and the first statement of Proposition~\ref{Proposition-Commutator}, $\del_{i}S_{n}=0$, for $i\neq n$, and $\del_{n}S_{n}=[z_{n},z_{n}]=K_{n}$. This last equality means that there exists a regular epimorphism filling the square
\[
\xymatrix{S_{n} \ar@{{ >}->}[d]_-{s_{n}} \ar@{.{ >>}}[r] & K_{n} \ar@{{ >}->}[d]^-{k_{n}}\\
A_{n+1} \ar@{-{ >>}}[r]_-{\del_{n}} & A_{n}. }
\]
By the third statement of Proposition~\ref{Proposition-Commutator}, there is a map $l_{n}\colon K_{n}\to Z_{n}A$ such that $z_{n}\comp l_{n}=k_{n}$. It follows that $l_{n}\sim  0$. Now, by Proposition~\ref{Proposition-Kan}, $A$ is Kan; hence, Proposition~\ref{Proposition-Pseudo-Homotopy} implies that $0\sim l_{n}$. So there exists a morphism $x\colon {X\to A_{n+1}}$ and a regular epimorphism $p\colon {X\to  K_{n}}$ such that $\del_{n+1}\comp x = z_{n}\comp l_{n}\comp  p$, i.e., the outer rectangle in
\[
\xymatrix{ X \ar@{.>}[rr]^{x} \ar@{.{ >>}}[dd]_-{p} \ar@{.>}[rd]_-{q} && A_{n+1} \ar[dd]^{\del_{n+1}} \\
{}\ar@{}[rd]|{\texttt{(i)}} & {N_{n+1}A} \ar@{{ >}->}[ru]^-{\bigcap_{i=1}^{n}\ker \del_{i}} \ar[d]^-{d'_{n+1}} & \\
K_{n} \ar@{{ >}->}[r]_-{l_{n}} & {Z_{n}A} \ar@{{ >}->}[r]_-{z_{n}} & A_{n}}
\]
commutes, and such that, moreover, $\del_{i}\comp x = 0$, for $i\neq n+1$.  It follows that an arrow $q$ exists such that, in the diagram above, the triangle commutes. As $z_{n}$ is a monomorphism, we get the commutativity of the trapezium \texttt{(i)}. Now, in the diagram
\[
\xymatrix{ X \ar[d]_{q} \ar@{-{ >>}}[r]^-{p} & K_{n} \ar@{{ >}->}[r]^-{l_{n}} & {Z_{n}A} \ar@{=}[d] \ar@{-{ >>}}[r] & \frac{Z_{n}A}{K_{n}} \ar@{.>}[d]^{r}\\
{N_{n+1}A} \ar[rr]_{d'_{n+1}} && {Z_{n}A} \ar@{-{ >>}}[r] & H_{n}A,}
\]
an arrow $r$ exists such that the right hand square commutes. Indeed, as $p$ is an epimorphism, $\coker (l_{n}\comp p)=\coker l_n$. This map $r$ is a regular epimorphism. Since, by Proposition~\ref{Proposition-Commutator}, ${Z_{n}A}/{K_{n}}=Z_{n}A/[z_{n},z_{n}]$ is an abelian object, and since a quotient of an abelian object of $\Ac$ is abelian, so is the object $H_{n}A$.
\end{proof}
 % Homology
\chapter{Baer invariants}\label{Chapter-Baer-Invariants}

\setcounter{section}{-1}
\section{Introduction}\label{Section-Baer-Introduction}
It is classical to present a group $G$ as a quotient $F/R$ of a free group $F$ and a ``group of relations''~$R$ (see e.g., Example~\ref{Example-Proper-Chain-Complexes-no-Kernels}). The philosophy is that $F$ might be easier to understand than $G$. Working with these presentations of $G$, it is relevant to ask which expressions of the datum $R\triangleleft F$ are independent of the chosen presentation. A first answer to the question was given by Hopf in~\cite{Hopf}, where he showed that
\[
\frac{[F,F]}{[R,F]}  \qquad \text{and}\qquad  \frac{R\cap [F,F]}{[R,F]}
\]
are such expressions. In~\cite{Baer}, R.~Baer further investigated this matter, constructing several of these invariants: He constructed expressions in terms of presentations such that ``similar presentations'' induce isomorphic groups---different presentations of a group by a free group and a group of relations always being similar. Whence the name \textit{Baer invariant} to denote such an expression. As in the two examples above, he constructed his invariants using commutator subgroups.

The work of Baer was followed up by Fr\"{o}hlich~\cite{Froehlich}, Lue~\cite{Lue} and Furtado-Coelho~\cite{Furtado-Coelho}, who generalized the theory to the case of Higgins's $\Omega$-groups~\cite{Higgins}. Whereas Baer constructs invariants using commutator subgroups, these authors obtain, in a similar fashion, generalized Baer invariants from certain subfunctors of the identity functor of the variety of $\Omega$-groups considered. Fr\"ohlich and Lue use subfunctors associated with subvarieties of the given variety, and Furtado-Coelho extends this to arbitrary subfunctors of the identity functor. By considering the variety of groups and its subvariety of abelian groups, one recovers the invariants obtained by Baer.

That the context of $\Omega$-groups, however, could still be further enlarged, was already hinted at by Furtado-Coelho, when he pointed out that
\begin{quote}
\dots all one needs, besides such fundamental concepts as those of kernel, image, etc., is the basic lemma on connecting homomorphisms.
\end{quote}
\pfbreak

\noindent We will work in the context of quasi-pointed exact protomodular categories. As we have seen in Chapter~\ref{Chapter-Preliminaries}, an important consequence of these axioms is the validity of the basic diagram lemmas of homological algebra, such as the $3\times 3$ Lemma and the Snake Lemma.  Thus, keeping in mind Furtado-Coelho's remark, one could expect this context to be suitable for a general description of the theory of Baer invariants due to Fr\"ohlich, Lue and Furtado-Coelho. Sections~\ref{Section-Exact-Sequences} and~\ref{Section-Birkhoff-Subfunctors} of this text are a confirmation of that thesis.

In Section~\ref{Section-General-Theory} we give a definition of Baer invariants.  We call \textit{presentation} of an object $A$ of a category $\Ac$ any regular epimorphism $p\colon A_{0}\to A$.  $\prA$ will denote the category of presentations in $\Ac$ and commutative squares between them, and $\pr\colon \prA \to \Ac $ the forgetful functor, which maps a presentation to the object presented.  Two morphisms of presentations $\vetsf,\vetsg \colon p\to q$ are called \textit{homotopic}, notation $\vetsf \simeq\vetsg $, if $\pr\,\vetsf =\pr\,\vetsg $. Homotopically equivalent presentations are said to be \textit{similar}. A functor $B\colon \prA \to \Ac$ is called a \textit{Baer invariant} if $\vetsf\simeq\vetsg$ implies that $B\vetsf=B\vetsg$; it maps similar presentations to isomorphic objects.

We prove that any functor $L_{0}\colon \prA \to \Ac$ can be turned into a Baer invariant by dividing out a subfunctor $S$ of $L_{0}$ that is ``large enough''. Given an appropriate subcategory of $\prA$, a Baer invariant can then be turned into a functor $\Ac \to \Ac $.

In Section~\ref{Section-Exact-Sequences}, the context is reduced to quasi-pointed exact protomodular categories. We show that in case $L_{0}$ arises from a functor $L\colon \Ac \to \Ac$, the class $\Fc_{L_{0}}$ of such ``large enough'' subfunctors $S$ of $L_{0}$ has a minimum $L_{1}\colon {\prA \to \Ac}$.  (A different interpretation of the functor $L_{1}$ is given in Section~\ref{Section-Commutators}.)  Next to $L_{1}$, the functor $K[\cdot]\cap L_0$ is also a member of $\Fc_{L_{0}}$.  Using Noether's First Isomorphism Theorem, out of these functors, we construct two exact sequences of Baer invariants. Finally, as an application of the Snake Lemma, we find a six-term exact sequence of functors $\Ac \to \Ac$.

We call \textit{Birkhoff subfunctor of $\Ac$} any proper subfunctor $V$ of $1_{\Ac}$ (i.e., a kernel ${V\To 1_{\Ac}}$) that preserves regular epimorphisms. In this way, we capture Fr\"ohlich's notion of \textit{variety subfunctor}.  As introduced by Janelidze and Kelly in~\cite{Janelidze-Kelly}, a full and reflective subcategory $\Bc$ of an exact category $\Ac$ is called \textit{Birkhoff} when it is closed in $\Ac$ under subobjects and quotient objects.  A Birkhoff subcategory of a variety of universal algebra~\cite{Cohn} (thus, in particular, any variety of $\Omega$-groups) is nothing but a subvariety. In Section~\ref{Section-Birkhoff-Subfunctors}, we see that Birkhoff subfunctors correspond bijectively to the Birkhoff subcategories of $\Ac$: Assuming that the sequence
\[
\xymatrix{0 \ar@{=>}[r] & V \ar@{{|>}=>}[r]^-{\mu} & 1_{\Ac} \ar@{={ >>}}[r]^-{\eta} & U \ar@{=>}[r] & 1}
\]
is exact, $V$ is a Birkhoff subfunctor if and only if $U$ reflects $\Ac$ onto a Birkhoff
subcategory.  It follows that Baer invariants can be obtained by considering suitable subcategories of a quasi-pointed, exact and protomodular category.  This allows us to refine the six-term exact sequence from Section~\ref{Section-Exact-Sequences}.

In Section~\ref{Section-Commutators}, we show that the functor $V_{1}\colon \prA \to \Ac$ associated with a Birkhoff subcategory $\Bc$ of $\Ac$ may be interpreted as a commutator. We relate it to existing notions in the field of categorical commutator theory and the theory of central extensions. The resulting concept of centrality fits into Janelidze and Kelly's theory of central extensions~\cite{Janelidze-Kelly}. In case $\Bc$ is the Birkhoff subcategory $\Ab\Ac$ of all abelian objects in $\Ac$, this notion of centrality coincides with the ones introduced by Smith~\cite{Smith} and Huq~\cite{Huq}; our commutator corresponds to Smith's commutator of equivalence relations, which was generalized to the exact Mal'tsev context by Pedicchio~\cite{Pedicchio}.

Finally, in Section~\ref{Section-Nilpotency}, we propose a notion of \textit{nilpotency}, relative to a Birkhoff subfunctor $V$ of $\Ac$.  We prove that an object is $V$-nilpotent if and only if its \textit{$V$-lower central series} reaches zero.  The nilpotent objects of class $n$ form a Birkhoff subcategory of $\Ac$.

This theory of Baer invariants will be used in Chapter~\ref{Chapter-Cotriples} to give an interpretation to a semi-abelian version of Barr-Beck homology.

This chapter more or less contains my first paper with Tomas Everaert~\cite{EverVdL1}. Here part of the theory is considered in a slightly more general context (sequentiable instead of semi-abelian), and some mistakes are corrected. Some results in Section~\ref{Section-Commutators} (e.g., Proposition~\ref{Proposition-Huq-Centrality-is-Smith-Centrality}) come from my paper with Marino Gran~\cite{Gran-VdL}.

\section{General theory}\label{Section-General-Theory}
We start this section by giving a definition of Baer invariants. In order to turn a functor $L_{0}\colon {\prA \to \Ac}$ into a Baer invariant, we consider a class $\Fc_{L_{0}}$ of subfunctors of $L_{0}$. Proposition~\ref{balg} shows that for any $L_{0}$ and any $S\in \Fc_{L_{0}}$, $L_{0}/S\colon \prA \to \Ac$ is a Baer invariant. We give equivalent descriptions of the class $\Fc_{L_{0}}$ in case $L_{0}$ arises from a functor $\Ac \to \Ac$. Finally, in Proposition~\ref{Proposition-Bear-Inv-is-Relative-Baer-Inv}, we show that a Baer invariant can be turned into a functor ${\Ac \to \Ac} $, given an appropriate subcategory of $\prA $.

\begin{definition}\label{Definition-Presentation}
Let $\Ac$ be a category. By a \emph{presentation}\index{presentation} of an object $A$ of $\Ac$ we mean a regular epimorphism $\xymatrix{p\colon A_{0} \ar@{->>}[r] &A}$. (In a sequentiable category, the kernel $\ker p\colon {A_{1}\kernelto A_{0}}$ of $p$ exists, and then $A$ is equal to the quotient $A_{0}/A_{1}$; cf.\ Remark~\ref{Remark-Presentation-is-Short-Exact-Sequence}.) We denote by $\prA$\index{category!$\Pr \Ac$} the category of presentations of objects of $\Ac$, a morphism $\vetsf=(f_{0},f)\colon p\to q$\index{morphism of presentations} \emph{over $f$} being a commutative square
\[
\xymatrix{A_{0} \ar@{->>}[d]_-{p} \ar[r]^-{f_{0}} & B_{0} \ar@{->>}[d]^-{q} \\
A \ar[r]_-{f} & B.}
\]
($\Pr \Ac$ is a full subcategory of the category of arrows\index{category!of arrows}\index{category!$\Two$} $\Fun (\Two,\Ac)$.)

Let $\pr\colon \prA \to \Ac$\index{pr@$\pr\colon \prA\to \Ac$} denote the forgetful functor that maps a presentation to the object presented, sending a morphism of presentations $\vetsf =(f_{0},f)$ to $f$.  Two morphisms of presentations $\vetsf,\vetsg \colon p\to q$ are called \emph{homotopic}\index{homotopy!of morphisms of presentations}, notation $\vetsf\simeq\vetsg $, if $\pr\,\vetsf =\pr\,\vetsg $ (or $f=g$). Homotopically equivalent presentations are called \emph{similar}\index{similar presentations}\index{homotopically equivalent presentations}.
\end{definition}

\begin{remark}\label{Remark-Homotopy}
Note that in any category $\Ac$, a kernel pair $(\KP [p],k_{0},k_{1})$ of a morphism $p\colon {A_{0}\to A}$ is an equivalence relation, hence an internal category\index{internal!category}. In case $\Ac$ has kernel pairs of regular epimorphisms, a morphism of presentations $\vetsf\colon {p\to q}$ gives rise to an internal functor
\[
\xymatrix{
\KP [p] \ar[r]^-{\KP \vetsf} \ar@<-0.5 ex>[d]_-{k_{0}} \ar@<0.5 ex>[d]^-{k_{1}} & \KP [q] \ar@<-0.5 ex>[d]_-{l_{0}} \ar@<0.5 ex>[d]^-{l_{1}} \\
A_{0} \ar[r]_-{f_{0}} & B_{0}.}
\]
Then $\vetsf \simeq\vetsg\colon p\to q$ if and only if the corresponding internal functors\index{internal!functor} are naturally isomorphic or \emph{homotopic}\index{homotopy!of internal functors}. More on homotopy of internal categories may be found in Chapter~\ref{Chapter-Internal-Categories}.
\end{remark}

\begin{definition}\label{Definition-Baer-Invariant}
A functor $B\colon \prA \to \Ac$ is called a \emph{Baer invariant}\index{Baer invariant} if $\vetsf\simeq\vetsg$ implies that $B\vetsf=B\vetsg$.
\end{definition}

The following shows that a Baer invariant maps similar presentations to isomorphic objects.

\begin{proposition}\label{Proposition-Baer-on-Objects}
In $\Ac$ let $p$ and $p'$ be presentations and $f$ and $g$ maps
\[
\xymatrix{A_{0} \ar@{->>}[rd]_-{p} \ar@<+0.5 ex>[rr]^-{f} && A'_{0} \ar@<+0.5 ex>[ll]^-{g} \ar@{->>}[ld]^-{p'}\\
& A}
\]
such that $p'\comp f=p$ and $p\comp g=p'$. If $B\colon \prA \to \Ac$ is a Baer invariant, then $Bp\cong Bp'$.
\end{proposition}
\begin{proof}
The existence of $f$ and $g$ amounts to $(g,1_{A})\comp (f,1_{A})\simeq 1_{p}$ and 
\[
(f,1_{A})\comp (g,1_{A})\simeq 1_{p'}.
\]
Obviously now, the map $B (f,1_{A})\colon {Bp\to Bp'}$ is an isomorphism with inverse $B (g,1_{A})$.
\end{proof}

Recall that a \emph{subfunctor} of a functor $F\colon \Cc \to \Dc$ is a subobject of $F$ in the functor category $\Fun (\Cc ,\Dc)$\index{category!$\Fun (\Cc ,\Dc)$}. We shall denote such a subfunctor $\mu \in \Sub F$\index{sub@$\Sub F$} by a representing monic natural transformation $\mu \colon G\To F$, or simply by $G\subseteq F$ or $G$. Now let $\Dc$ be a quasi-pointed category. A subfunctor $\mu \colon {G\To F}$ is called \emph{proper}\index{proper!subfunctor}\index{subfunctor!proper} if every component $\mu_{C}\colon {G (C)\to F(C)}$ is a kernel. This situation will be denoted by $G\triangleleft F$\index{$G\triangleleft F$}. In case $\Dc$ is sequentiable, a proper subfunctor gives rise to a short exact sequence
\[
\xymatrix{
0 \ar@{=>}[r] & F \ar@{{|>}=>}[r]^-{\mu} & G \ar@{={ >>}}[r]^-{\eta_{F}} & \frac{G}{F} \ar@{=>}[r] & 1}
\]
of functors $\Cc \to \Dc$.  Like here, in what follows, exactness of sequences of functors will always be pointwise. If confusion is unlikely we shall omit the index $F$, and write $\eta$ for $\eta_{F}$ and $\eta_{C}$ for $(\eta_{F})_{C}$.

The following, at first, quite surprised us:

\begin{example}\label{Example-Proper-Subfunctor}
Every subfunctor $L$ of $1_{\Gp}\colon \Gp \to \Gp$ is proper. Indeed, by the naturality of $\mu\colon L\To 1_{\Gp}$, every $\mu_{G}\colon L (G)\to G$ is the inclusion of a fully-invariant\index{fully-invariant subgroup}, hence normal, subgroup $L (G)$ into $G$---and a kernel in $\Gp$ is the same thing as a normal subgroup.

This is, of course, not true in general: Consider the category $\omegaGp$\index{category!$\omegaGp$} of groups with an operator~$\omega$. An object of $\omegaGp$ is a pair $(G,\omega)$, with $G$ a group and $\omega\colon {G\to G}$ an endomorphism of $G$, and an arrow ${(G,\omega)\to (G',\omega')}$ is a group homomorphism $f\colon {G\to G'}$ satisfying $f\comp \omega=\omega'\comp f$. Then putting $L(G,\omega)=(\omega(G),\omega\vert_{\omega(G)})$ defines a subfunctor of the identity functor $1_{\omegaGp}$. $L$ is, however, not proper: For any group endomorphism $\omega\colon {G\to G}$ of which the image is not normal in $G$, the inclusion of $L(G,\omega)$ in $(G,\omega)$ is not a kernel.
\end{example}

If $\Ac$ is quasi-pointed and has cokernels of kernels, any functor $\prA\to \Ac$ can be turned into a Baer invariant by dividing out a ``large enough'' subfunctor. In order to make this precise, we make the following

\begin{definition}\cite{Furtado-Coelho}\label{Definition-Class-Fc}
Consider a functor $L_{0}\colon \prA\to \Ac$ and a presentation $q\colon B_{0}\to B$. Then $\Fc_{L_{0}}^{q}$\index{F@$\Fc_{L_{0}}^{q}$} denotes the class of proper subobjects $S\triangleleft L_{0}q$ for which $\vetsf\simeq\vetsg\colon p\to q$ implies that $\eta_{S}\comp L_{0}\vetsf =\eta_{S}\comp L_{0}\vetsg$, i.e., the two compositions in the diagram
\[
\xymatrix{
L_{0}p \ar@<0.5 ex>[r]^-{L_{0}\vetsf} \ar@<-0.5 ex>[r]_-{L_{0}\vetsg } & L_{0}q \ar@{->>}[r]^-{\eta_{S}} & \frac{L_{0}q}{S}}
\]
are equal.

$\Fc_{L_{0}}$\index{F@$\Fc_{L_{0}}$} is the class of functors $S\colon \prA\to\Ac$ with $Sq\in \Fc_{L_{0}}^{q}$, for every $q\in \Ob{\prA}$. Hence, $S\in \Fc_{L_{0}}$ if and only if the following conditions are satisfied:
\renewcommand{\labelenumi}{(\roman{enumi})}
\begin{enumerate}
\item $S\triangleleft L_{0}$;
\item $\vetsf\simeq\vetsg\colon p\to q$ implies that $\eta_{q}\comp L_{0}\vetsf =\eta_{q}\comp L_{0}\vetsg$, i.e., the two compositions in the diagram
\[
\xymatrix{
L_{0}p \ar@<0.5 ex>[r]^-{L_{0}\vetsf} \ar@<-0.5 ex>[r]_-{L_{0}\vetsg } & L_{0}q \ar@{->>}[r]^-{\eta_{q}} & \frac{L_{0}}{S}q}
\]
are equal.
\end{enumerate}
\renewcommand{\labelenumi}{(\arabic{enumi})}
The class $\Fc_{L_{0}}^{q}$ (respectively $\Fc_{L_{0}}$) may be considered a subclass of $\Sub (L_{0}q)$\index{sub@$\Sub F$} (respectively $\Sub (L_{0})$), and as such carries the inclusion order.
\end{definition}

\begin{proposition}\label{balg}\cite[Proposition 1]{Furtado-Coelho}
Suppose that $\Ac$ is quasi-pointed and has cokernels of kernels. Let $L_{0}\colon \prA\to\Ac$ be a functor and $S$ an element of $\Fc_{L_{0}}$. Then $L_{0}/S\colon \prA\to\Ac$ is a Baer invariant.
\end{proposition}
\begin{proof}
Condition (i) ensures that $L_{0}/S$ exists. $\eta\colon {L_{0}\To L_{0}/S}$ being a natural transformation implies that, for morphisms of presentations $\vetsf,\vetsg\colon p\to q$, both left and right hand downward pointing squares
\[
\xymatrix{
L_{0}p \ar@{->>}[r]^-{\eta_{p}} \ar@<-0.5 ex>[d]_-{L_{0}\vetsf} \ar@<0.5 ex>[d]^-{L_{0}\vetsg} & \frac{L_{0}}{S}p \ar@<-0.5 ex>[d]_-{\frac{L_{0}}{S}\vetsf} \ar@<0.5 ex>[d]^-{\frac{L_{0}}{S}\vetsg} \\
L_{0}q \ar@{->>}[r]_-{\eta_{q}} & \frac{L_{0}}{S}q }
\]
commute. If $\vetsf \simeq \vetsg$ then, by (ii), $\eta_{q}\comp L_{0}\vetsf =\eta_{q}\comp L_{0}\vetsg$, so 
\[
\tfrac{L_{0}}{S}\vetsf \comp \eta_{p}=\tfrac{L_{0}}{S}\vetsg \comp \eta_{p},
\]
whence the desired equality $({L_{0}}/{S}) (\vetsf)= ({L_{0}}/{S}) (\vetsg)$.
\end{proof}

Proposition~\ref{balg} is particularly useful for functors $L_{0}\colon \prA\to\Ac$ induced by a functor $L\colon \Ac \to \Ac$ by putting\index{L_0@$L_{0}$}
\[
L_0 (p\colon A_{0}\to A)=L(A_0)\quad \text{and}\quad L_0 ((f_{0},f)\colon p\to q) =Lf_{0}.
\]
In this case, the class $\Fc_{L_{0}}^{q}$ in Definition~\ref{Definition-Class-Fc} has some equivalent descriptions.

\begin{proposition}\label{Proposition-(ii')}
Suppose that $\Ac$ is quasi-pointed with cokernels of kernels. Given a functor $L\colon \Ac \to \Ac$ and a presentation $q\colon B_{0}\to B$, the following are equivalent:
\begin{enumerate}
\item $S\in \Fc_{L_{0}}^{q}$;
\item $\vetsf\simeq\vetsg\colon p\to q$ implies that $\eta_{S}\comp L_{0}\vetsf =\eta_{S}\comp L_{0}\vetsg$;
\item for morphisms $f_{0},g_{0}\colon  A_{0} \to B_{0}$ in $\Ac$ with $q\comp f_{0}=q\comp g_{0}$, $\eta_{S}\comp Lf_{0}=\eta_{S}\comp Lg_{0}$, i.e., the two compositions in the diagram
\[
\xymatrix{
L (A_{0}) \ar@<0.5 ex>[r]^-{Lf_{0}} \ar@<-0.5 ex>[r]_-{Lg_{0}} & L (B_{0}) \ar@{->>}[r]^-{\eta_{S}} & \frac{L (B_{0})}{S}}
\]
are equal.
\end{enumerate}
If, moreover, $\Ac$ has kernel pairs of regular epimorphisms, these conditions  are equivalent to
\begin{enumerate}
\setcounter{enumi}{3}
\item $L (\KP [q])\subseteq \KP \bigl[\eta_{S}\colon L (B_{0})\to {L (B_{0})}/{S}\bigr]$.
\end{enumerate}
\end{proposition}
\begin{proof}
We show that 2.\ implies 3.: A presentation $p$ and maps of presentations $\vetsf\simeq \vetsg \colon p\to q$ such as needed in condition 2.\ are given by
\[
\xymatrix{A_0  \ar@<0.5 ex>[r]^-{f_{0}} \ar@<-0.5 ex>[r]_-{g_{0}} \ar@{=}[d] & B_{0} \ar@{->>}[d]^-{q}\\
A_{0} \ar[r]_-{q\comp f_{0}} & B}
\]
$p=1_{A_{0}}$, $\vetsf = (f_{0},q\comp f_{0})$, $\vetsg = (g_{0},q\comp f_{0})$.
\end{proof}

We now show how a Baer invariant gives rise to a functor $\Ac \to \Ac $, given the following additional datum:

\begin{definition}\label{Definition-Web}
Let $\Ac$ be a category. We call a subcategory $\Wc$ of $\prA$ a \emph{web on \Ac}\index{web}\index{category!web} if
\begin{enumerate}
\item for every object $A$ of $\Ac$, a presentation $p\colon A_{0}\to A$ in $\Wc$ exists;
\item given presentations $p\colon A_{0}\to A$ and $q\colon B_{0}\to B$ in $\Wc$, for every map $f\colon {A\to B}$ in $\Ac$ there exists a morphism $\vetsf\colon p\to q$ in $\Wc$ such that $\pr\, \vetsf =f$.
\end{enumerate}
\end{definition}

\begin{example}\label{Example-Split-Presentation}
A \emph{split}\index{split presentation}\index{presentation!split} presentation is a split epimorphism $p$ of $\Ac$. (Note that an epimorphism $p\colon {A_{0}\to A}$ with splitting $s\colon {A\to A_{0}}$ is a coequalizer of the maps $1_{A_{0}}$ and $s\comp p\colon {A_{0}\to A_{0}}$.) The full subcategory of $\prA$ determined by the split presentations of $\Ac$ is a web $\Wsplit$\index{category!$\Wsplit$}, the \emph{web of split presentations}\index{web!of split presentations}. Indeed, for any $A\in \Ac$, $1_{A}\colon {A\to A}$ is a split presentation of $A$; given split presentations $p$ and $q$ and a map $f$ such as indicated by the diagram
\[
\xymatrix{A_{0} \ar@{.>}[r]^-{t\comp f\comp p} \ar@<-0.5 ex>@{->>}[d]_-{p} & B_{0} \ar@<-0.5 ex>@{->>}[d]_-{q}\\
A \ar@<-0.5 ex>[u]_-{s} \ar[r]_-{f} & B \ar@<-0.5 ex>[u]_-{t},}
\]
$(t\comp f\comp p,f)$ is the needed morphism of $\Wc$.
\end{example}

\begin{example}\label{Example-Functorial-Web}
Let $F\colon \Ac \to \Ac$ be a functor and let $\pi\colon F\To 1_{\Ac}$ be a natural transformation of which all components are regular epimorphisms. Then the presentations 
\[
\pi_A\colon {F (A)\to A},
\]
\index{pi@$\pi_{A}$}for $A\in \Ob{\Ac}$, together with the morphisms of presentations 
\[
(F (f),f)\colon {\pi_{A}\to \pi_{B}},
\]
for $f\colon {A\to B}$ in $\Ac$, constitute a web $\Wc_{F}$\index{category!$\Wc_{F}$} on $\Ac$ called the \emph{functorial web determined by $F$}\index{web!functorial}.
\end{example}

\begin{example}\label{Example-Projective-Web}
A presentation $p\colon A_{0}\to A$ is called \emph{projective}\index{projective presentation}\index{presentation!projective} when $A_{0}$ is a projective object of $\Ac$. If $\Ac$ has enough projectives, then by definition every object has a projective presentation. In this case the full subcategory $\Wproj$\index{category!$\Wproj $} of all projective presentations of objects of $\Ac$ is a web, called the \emph{web of projective presentations}\index{web!of projective presentations}.
\end{example}

Recall~\cite{Borceux:Cats} that a \emph{graph morphism}\index{graph morphism} $F\colon \Cc \to \Dc$ between categories $\Cc$ and $\Dc$ has the structure, but not all the properties, of a functor from $\Cc$ to $\Dc$: It need not preserve identities or compositions.

\begin{definition}\label{Definition-Choice-and-gBi}
Let $\Wc$ be a web on a category $\Ac$ and $i\colon \Wc \to \prA$ the inclusion.  By a \emph{choice $c$ of presentations in $\Wc$}\index{choice of presentations in a web}, we mean a graph morphism $c\colon \Ac\to \Wc$ such that $\pr\comp i\comp c=1_{\Ac}$.

A functor $B\colon \prA\to \Ac$ is called a \emph{Baer invariant relative to $\Wc$}\index{Baer invariant!relative to a web} when
\begin{enumerate}
\item for every choice of presentations $c\colon \Ac \to \Wc$, the induced graph morphism $B\comp i\comp c\colon {\Ac \to \Ac}$ is a functor;
\item for every two choices of presentations $c,c'\colon \Ac \to \Wc$, the functors $B\comp i\comp c$ and $B\comp i\comp c'$ are naturally isomorphic.
\end{enumerate}
\end{definition}

\begin{example}\label{Example-Functorial-gbi}
Any functor $B\colon \prA\to \Ac$ is a Baer invariant relative to any functorial web $\Wc_{F}$: There is only one choice $c\colon \Ac \to \Wc_{F}$, and it is a functor.
\end{example}

\begin{proposition}\label{Proposition-Bear-Inv-is-Relative-Baer-Inv}
If $B\colon \prA \to \Ac$ is a Baer invariant, it is a Baer invariant relative to any web $\Wc$ on $\Ac$.
\end{proposition}
\begin{proof}
Let $\Wc$ be a web and $c,c'\colon \Ac \to \Wc$ two choices of presentations in $\Wc$. For the proof of 1., let $f\colon A\to A'$ and $g\colon A'\to A''$ be morphisms of $\Ac$. Then $1_{ic (A)}\simeq ic1_{A}$ and $ic (f\comp g)\simeq icf\comp icg$; consequently, $1_{Bic (A)}=Bic1_{A}$ and $Bic (f\comp g)=Bicf\comp Bicg$.

The second statement is proven by choosing, for every object $A$ of $\Ac$, a morphism 
\[
\tau_{A}\colon {c (A)\to c' (A)}
\]
in the web $\Wc$. Then $\nu_{A}=Bi\tau_{A}\colon Bic (A)\to Bic' (A)$ is independent of the choice of $\tau_{A}$, and the collection $(\nu_{A})_{A\in \Ob{\Ac}}$ is a natural isomorphism ${B\comp i\comp c\To B\comp i\comp c'}$.
\end{proof}

This proposition implies that a Baer invariant $B\colon \prA \to \Ac$ gives rise to functors $B\comp i\comp c\colon {\Ac \to \Ac}$, independent of the choice $c$ in a web $\Wc$. But it is important to note that such functors \textit{do} depend on the chosen web: Indeed, for two choices of presentations $c$ and $c'$ in two different webs $\Wc$ and $\Wc '$, the functors $B\comp i\comp c$ and $B\comp i'\comp c'$ need not be naturally isomorphic---see Remark~\ref{Remark-Dependence-on-Web} for an example.

Suppose that $\Ac$ is a quasi-pointed category with cokernels of kernels and let $c$ be a choice in a web $\Wc$ on $\Ac$. If $L_{0}\colon \prA \to A$ is a functor and $S\in \Fc_{L_{0}}$, then by Proposition~\ref{balg}
\[
\frac{L_{0}}{S}\comp i\comp c
\]
is a functor $\Ac \to \Ac$. It is this kind of functor that we shall study in the next sections.

\section{The exact \& sequentiable case}\label{Section-Exact-Sequences}
In this section, $\Ac$ shall be an exact and sequentiable category. We show that the class $\Fc_{L_{0}}$ contains two canonical functors: its minimum $L_{1}$ and the functor $K[\cdot]\cap L_0$. Applying Noether's First Isomorphism Theorem, we construct two exact sequences, and thus find new Baer invariants. Finally, applying the Snake Lemma, we construct a six-term exact sequence of functors $\Ac\to \Ac$.

\begin{remark}\label{Remark-Presentation-is-Short-Exact-Sequence}
If $\Ac$ is sequentiable, giving a presentation $p$ is equivalent to giving a short exact sequence
\[
\xymatrix{
0 \ar[r] & K[p] \ar@{{ >}->}[r]^-{\ker p} & A_{0} \ar@{-{ >>}}[r]^-p & A \ar[r] & 1. }
\]
$K[\cdot]\colon \Fun (\Two , \Ac)\to \Ac$, and its restriction $K[\cdot]\colon \prA \to \Ac$, denote the kernel functor.
\end{remark}

\begin{proposition}\label{Proposition-Pullback-gBi}\cite[Proposition 2]{Furtado-Coelho} Let $\Ac $ be a sequentiable category.  If $L$ is a subfunctor of $1_{\Ac}\colon \Ac\to \Ac$ then $K[\cdot]\cap L_0\colon \prA \to \Ac$ is in $\Fc_{L_{0}}$.  Hence
 \[
\frac{L_0}{K[\cdot]\cap L_0}\colon \prA \to \Ac
\]
is a Baer invariant.
\end{proposition}
\begin{proof}
Let $\mu\colon L\To 1_{\Ac}$ be the subfunctor $L$, and $L_0\colon \prA\to\Ac$ the induced functor (see the definition above Proposition~\ref{Proposition-(ii')}). By Proposition~\ref{balg} it suffices to prove the first statement. We show that, for every presentation $q\colon {B_{0}\to B}$, $K[q]\cap L(B_0)\in \Fc_{L_{0}}^{q}$, by checking the third condition of Proposition~\ref{Proposition-(ii')}. Since pullbacks preserve kernels---see the diagram below---${K[q]\cap L(B_0)}$ is normal in ${L (B_{0})}$. Moreover, applying Proposition~\ref{Proposition-LeftRightPullbacks}, the map $m_{q}$, defined by first pulling back $\ker q$ along $\mu_{B_{0}}$ and then taking cokernels
\[
\xymatrix{
0 \ar[r] & K[q]\cap L(B_0) \ar@{}[dr]|<<<{\pullback} \ar@{{ >}->}[r] \mond & L(B_0) \ar@{-{ >>}}[r] \mond_-{\mu_{B_{0}}} & \frac{L(B_0)}{K[q]\cap L(B_0)} \ar@{.>}[d]^-{m_{q}} \ar[r] & 1\\
0 \ar[r] & K[q] \ar@{{ >}->}[r]_-{\ker q} & B_0 \ar@{-{ >>}}[r]_-{q} & B \ar[r] & 1, }
\]
is a monomorphism. Consider morphisms $f_{0},g_{0}\colon  A_{0} \to B_{0}$ in $\Ac$ with $q\comp f_{0}=q\comp g_{0}$. We are to prove that in the following diagram, the two compositions on top are equal:
\[
\xymatrix{
L(A_0)\ar@<0.5 ex>[r]^-{Lf_0} \ar@<-0.5 ex>[r]_-{Lg_0} & L(B_0) \ar@{-{ >>}}[r] \mond_-{\mu_{B_{0}}} & \frac{L(B_0)}{K[q]\cap L(B_0)} \mond^-{m_{q}} \\
 & B_0 \ar@{-{ >>}}[r]_-{q} &B.}
\]
Since $m_{q}$ is a monomorphism, it suffices to prove that $q\comp \mu_{B_{0}}\comp Lf_0=q\comp \mu_{B_{0}}\comp Lg_0$. But, by naturality of $\mu$, $q\comp\mu_{B_{0}} =\mu_{B}\comp Lq$, and so
\[
q\comp \mu_{B_{0}} \comp Lf_0=\mu_{B}\comp L(q\comp f_0)=\mu_{B}\comp L(q\comp g_0)=q\comp \mu_{B_{0}}\comp Lg_0.\qedhere
\]
\end{proof}

\begin{proposition}\cite[Proposition 4]{Furtado-Coelho}\label{L_1}
Let $\Ac$ be exact and sequentiable. For any $L\colon \Ac \to \Ac $, the ordered class $\Fc_{L_0}$ has a minimum $L_1\colon {\prA \to \Ac}$.\index{L_1@$L_{1}$}
\end{proposition}
\begin{proof}
For a morphism of presentations $\vetsf \colon p\to q$, $L_{1}\vetsf$ is defined by first taking kernel pairs
\begin{equation}\label{Diagram-L1-Step1}
\vcenter{\xymatrix{\KP [p] \ar@{}[rd]|{\texttt{(i)}} \ar@{.>}[d]_-{\KP \vetsf} \ar@<0.5 ex>[r]^-{k_{0}} \ar@<-0.5 ex>[r]_-{k_{1}} & A_{0} \ar@{}[rd]|{\texttt{(ii)}} \ar[d]^-{f_{0}} \ar@{-{ >>}}[r]^-{p} & A \ar[d]^-{f} \\
\KP [q] \ar@<0.5 ex>[r]^-{l_{0}} \ar@<-0.5 ex>[r]_-{l_{1}} & B_{0} \ar@{-{ >>}}[r]_-{q} & B,}}
\end{equation}
next applying $L$ and taking coequalizers (which exist because $\Ac$ is exact Mal'tsev)
\begin{equation}\label{Diagram-L1-Step2}
\vcenter{\xymatrix{L (\KP [p]) \ar[d]_-{L\KP \vetsf} \ar@<0.5 ex>[r]^-{Lk_{0}} \ar@<-0.5 ex>[r]_-{Lk_{1}} & L (A_{0}) \ar@{}[rd]|{\texttt{(iii)}} \ar[d]^-{Lf_{0}} \ar@{-{ >>}}[r]^-{c} & \Coeq [Lk_{0},Lk_{1}] \ar@{.>}[d]^-{\Coeq (L\KP \vetsf, Lf_{0})} \\
L (\KP [q]) \ar@<0.5 ex>[r]^-{Ll_{0}} \ar@<-0.5 ex>[r]_-{Ll_{1}} & L (B_{0}) \ar@{-{ >>}}[r]_-{d} & \Coeq [Ll_{0},Ll_{1}]}}
\end{equation}
and finally taking kernels
\begin{equation}\label{Diagram-L1-Step3}
\vcenter{\xymatrix{L_{1}p \ar@{.>}[d]_-{L_{1}\vetsf} \ar@{{ >}->}[r]^-{\ker c} & L (A_{0}) \ar@{}[rd]|{\texttt{(iii)}} \ar[d]^-{Lf_{0}} \ar@{-{ >>}}[r]^-{c} & \Coeq [Lk_{0},Lk_{1}] \ar@{.>}[d]^-{\Coeq (L\KP \vetsf, Lf_{0})} \\
L_{1}q \ar@{{ >}->}[r]_-{\ker d} & L (B_{0}) \ar@{-{ >>}}[r]_-{d} & \Coeq [Ll_{0},Ll_{1}].}}
\end{equation}
It is easily seen that this defines a functor $L_1\colon \prA \to \Ac$, a minimum in $\Fc_{L_0}$ for the inclusion order.
\end{proof}

\begin{remark}\label{Remark-L_{1}-Normalization}
One may, alternatively, compute $L_{1}p$ as the normalization of the reflexive graph\index{normalization of a reflexive graph} $L (\KP [p])$: The image of $Lk_{1}\comp \ker Lk_{0}$ is a kernel of $c$ (cf.\ Definition~\ref{Definition-Normal-Monomorphism}, Remark~\ref{Remark-Forks-Normalization}). Thus the two-step construction involving coequalizers may be avoided.
\end{remark}

\begin{remark}\label{Remark-Galois-Groupoid}
Since the category $\Gd \Ac$ of internal groupoids in $\Ac$ is closed under subobjects in the category $\Rgrph \Ac$ of reflexive graphs in $\Ac$~\cite{Pedicchio}, and since $R[p]$, being an equivalence relation, is an internal groupoid in $\Ac$, its image 
\[
(Lk_{0},Lk_{1})\colon {L(R[p])\to L (A_{0})}
\]
through $L$ is an internal groupoid. It is related to Janelidze's \emph{Galois groupoid}\index{Galois groupoid} of $p$~\cite{Janelidze:Como, Janelidze-Kelly}.
\end{remark}

\begin{remark}\label{Remark-L_1-Locally-Minimum}
Observe that the construction of $L_1$ above is such that, for every $p\in\Ob{\prA}$, $L_1p$ is a minimum in the class $\Fc^p_{L_0}$.
\end{remark}

\begin{corollary}\label{Corollary-L_1}
Let $\Ac$ be an exact and sequentiable category, $L\colon {\Ac \to \Ac} $ a functor. Then
\begin{enumerate}
\item $S\in \Fc_{L_{0}}$ if and only if $L_{1}\leq S \triangleleft L_{0}$;
\item $S\in \Fc_{L_{0}}^{p}$ if and only if $L_{1}p\leq S \triangleleft L_{0}p$.\noproof
\end{enumerate}
\end{corollary}

\begin{notation}\label{Notation-DL}\index{DL@$DL_{\Wc}$}\index{Baer invariant!$DL_{\Wc}$}
If $L\colon \Ac \to \Ac $ is a functor and $\Wc$ a web on $\Ac$ then, by Proposition~\ref{L_1}, we have the canonical functor
\[
DL_{\Wc}=\frac{L_{0}}{L_{1}}\comp i\comp c\colon \Ac \to \Ac.
\]
\end{notation}

\begin{remark}\label{Remark-L-in-L_1}
For $L\subseteq 1_{\Ac}$ we now have the following inclusions of functors $\prA \to \Ac$:
\[
L\comp K[\cdot]\quad \triangleleft \quad L_{1} \quad \triangleleft \quad K[\cdot]\cap L_{0}\quad  \triangleleft \quad L_{0}.
\]
Only the left hand side inclusion is not entirely obvious. For a presentation $p$, let $r$ denote the cokernel of the inclusion $L_{1}p\to L (A_{0})$. Then $p\comp \ker p=p\comp 0$, hence $r\comp L\ker p=r\comp L0=0$. This yields the required map.
\end{remark}

\begin{remark}\label{Remark-L-uit-L_1}
Suppose that $\Ac$ is pointed. Since then $L (A)={L_{1} (A\to 0)}$, the functor $L$ may be regained from $L_{1}$: Indeed, evaluating the above inclusions in $A\to 0$, we get
\[
L (A) = (L\comp K[\cdot]) (A\to 0)\subseteq L_{1} (A\to 0)\subseteq L_{0} (A\to 0)=L (A).
\]
\end{remark}

\begin{proposition}\label{exS}\cite[Proposition 6]{Furtado-Coelho}
Let $\Ac$ be exact and sequentiable; consider a subfunctor $L$ of $1_{\Ac}$. For $S\in \Fc_{L_{0}}$ with $S\subseteq K[\cdot]\cap L_0$, the following sequence of functors $\prA \to \Ac$ is exact, and all its terms are Baer invariants.
\begin{equation}\label{3exS}
\xymatrix{
0 \ar@{=>}[r] & \frac{K[\cdot]\cap L_0}{S} \ar@{{|>}=>}[r] & \frac{L_0}{S} \ar@{={ >>}}[r] & \frac{L_0}{K[\cdot]\cap L_0} \ar@{=>}[r] & 1}
\end{equation}
If, moreover, $L\triangleleft 1_{\Ac}$, the sequence
\begin{equation}\label{4exS}
\xymatrix{
0 \ar@{=>}[r] & \frac{K[\cdot]\cap L_0}{S} \ar@{{|>}=>}[r] & \frac{L_0}{S} \ar@{=>}[r] & \pr \ar@{={ >>}}[r] & \frac{\pr}{L_0/(K[\cdot]\cap L_0)} \ar@{=>}[r] & 1}
\end{equation}
is exact, and again all terms are Baer invariants.
\end{proposition}
\begin{proof}
The exactness of Sequence~\ref{3exS} follows from Proposition~\ref{generalexact}, since $S$ and $K[\cdot]\cap L_0$ are proper in $L_{0}$, both being in $\Fc_{L_{0}}$. The naturality is rather obvious.

Now suppose that $L\triangleleft 1_{\Ac}$. To prove the exactness of Sequence~\ref{4exS}, it suffices to show that ${L_0}/{(K[\cdot]\cap L_0)}$ is proper in $\pr$. Indeed, if so,
\[
\xymatrix{
0 \ar@{=>}[r] & \frac{L_0}{K[\cdot]\cap L_0} \ar@{{|>}=>}[r] & \pr \ar@{={ >>}}[r] & \frac{\pr}{L_0/(K[\cdot]\cap L_0)} \ar@{=>}[r] & 1 }
\]
is a short exact sequence, and by pasting it together with~\ref{3exS}, we get~\ref{4exS}. Reconsider, therefore, the following commutative diagram from the proof of Proposition~\ref{Proposition-Pullback-gBi}:
\[
\xymatrix{
L(A_0) \ar@{-{ >>}}[r] \ar@{{ >}->}[d]_-{\mu_{A_{0}}} & \frac{L(A_0)}{K[p]\cap L(A_0)} \mond^-{m_{p}}\\
A_0 \ar@{-{ >>}}[r]_-{p} & A. }
\]
Proposition~\ref{Proposition-Image-of-Kernel} implies that $m_{p}$ is a kernel, because $\mu_{A_{0}}$ is, and the category $\Ac$ is exact.

The terms of~\ref{3exS} and~\ref{4exS} being Baer invariants follows from the exactness of the sequences, and from the fact that ${L_0}/{S}$, ${L_0}/{(K[\cdot]\cap L_0)}$ and $\pr$ are Baer invariants.
\end{proof}

\begin{corollary}\label{Corollary-Split-Presentation}
Let $L$ be a subfunctor of $1_{\Ac}\colon \Ac \to \Ac$ and $S\in \Fc_{L_{0}}$ with $S\subseteq K[\cdot]\cap L_0$. Then
\begin{enumerate}
\item for any $A\in \Ob{\Ac}$, $S1_{A}=0$;
\item for any split presentation $p\colon A_{0}\to A$ in $\Ac$, $Sp\cong K[p]\cap L (A_{0})$;
\item for any presentation $p\colon A_{0}\to A$ of a projective object $A$ of $\Ac$, $Sp\cong K[p]\cap L (A_{0})$.
\end{enumerate}
\end{corollary}
\begin{proof}
The first statement holds because 
\[
S1_{A}\subseteq (K[\cdot]\cap L_0) (1_{A})= 0\cap L (A)=0.
\]

For a proof of the second statement, let $s\colon A\to A_{0}$ be a splitting of $p$. Then $(s,1_{A})\colon {1_{A}\to p}$ and $(p,1_{A})\colon {p\to 1_{A}}$ are morphisms of presentations. By Proposition~\ref{exS},
\[
\tfrac{K[\cdot]\cap L_0}{S}\colon \prA \to \Ac
\]
is a Baer invariant. Hence, Proposition~\ref{Proposition-Baer-on-Objects} implies that
\[
0=\tfrac{K[\cdot]\cap L_0}{S} (1_{A})\cong \tfrac{K[\cdot]\cap L_0}{S} (p)=\tfrac{K[p]\cap L (A_{0})}{Sp},
\]
whence the isomorphism $Sp\cong K[p]\cap L (A_{0})$.

The third statement follows immediately from 2.
\end{proof}

\begin{notation}\label{Notation-Baer-Invariants}\index{Baer invariant!$\Delta L_{\Wc}$}\index{Delta L@$\Delta L_{\Wc}$}\index{Baer invariant!$M_{\Wc}$}\index{M@$M_{\Wc}$}\index{Baer invariant!$\nabla L_{\Wc}$}\index{nabla@$\nabla L_{\Wc}$}
Let $c$ be a choice in a web $\Wc$ on $\Ac$. If $L\subseteq 1_{\Ac}$, Proposition~\ref{exS} and~\ref{L_1} tell us that
\[
\Delta L_{\Wc}=\frac{K[\cdot]\cap L_0}{L_1}\comp i\comp c\qquad \text{and}\qquad \nabla L_{\Wc}=\frac{L_0}{K[\cdot]\cap L_0}\comp i\comp c
\]
are functors $\Ac \to \Ac$. If, moreover, $L\triangleleft 1_{\Ac}$, then
\[
M_{\Wc}=\frac{\pr}{L_0/(K[\cdot]\cap L_0)}\comp i\comp c\colon \Ac \to \Ac
\]
is a functor as well. We omit the references to $c$, since any other choice $c'$ gives naturally isomorphic functors. In the following we shall always assume that a particular $c$ has been chosen. When the index $\Wc$ is omitted, it is understood that $\Wc$ is the web $\Wproj$ of projective presentations.
\end{notation}

\begin{remark}\label{Remark-Dependence-on-Web}
Note that these functors \textit{do} depend on the chosen web; for instance, by Corollary~\ref{Corollary-Split-Presentation}, $\Delta L_{\Wc_{1_{\Ac}}}=0$, for any $L\subseteq 1_{\Ac}$. But, as will be shown in Example~\ref{Examples-V-Exact-Sequence}, if $\Ac =\Gp$ and $L$ is the subfunctor associated with the Birkhoff subcategory $\Ab$ of abelian groups, then for any $R\triangleleft F$ with $F$ projective,
\[
\Delta (\eta_{R}\colon F\to \tfrac{F}{R})=\frac{R\cap [F,F]}{[R,F]}.
\]
By Hopf's Formula~\cite{Hopf}, this is the second integral homology group of $G=F/R$---and it need not be $0$.
\end{remark}

 From Proposition~\ref{exS} we deduce the following.

\begin{proposition}\label{exL}\cite[Proposition 6, 7]{Furtado-Coelho}
Let $\Wc$ be a web on an exact sequentiable category $\Ac$. If $L\subseteq 1_{\Ac}$, then
\begin{equation}\label{3exL}
\xymatrix{0 \ar@{=>}[r] & \Delta L_{\Wc} \ar@{{|>}=>}[r] & DL_{\Wc} \ar@{={ >>}}[r] & \nabla L_{\Wc} \ar@{=>}[r] & 1}
\end{equation}
is exact. If, moreover, $L\triangleleft 1_{\Ac}$, then also
\begin{equation}\label{4exL}
\xymatrix{0 \ar@{=>}[r] & \Delta L_{\Wc} \ar@{{|>}=>}[r] & DL_{\Wc} \ar@{=>}[r] & 1_{\Ac} \ar@{={ >>}}[r] & M_{\Wc} \ar@{=>}[r] & 1}
\end{equation}
is exact.\noproof
\end{proposition}

\begin{remark}\label{Remark-Birkhoff-Nabla}\cite[Proposition 8]{Furtado-Coelho}
In the important case that the functor $L\colon \Ac \to \Ac$ preserves regular epimorphisms, we get that $L$ and $\nabla L_{\Wc}$ represent the same subfunctor of $1_{\Ac}$. (Accordingly, in this case, $\nabla L_{\Wc}$ does not depend on the chosen web.) Of course, the exact sequence~\ref{3exL} then simplifies to
\begin{equation}\label{3exL-simplified}
\xymatrix{0 \ar@{=>}[r] & \Delta L_{\Wc} \ar@{{|>}=>}[r] & DL_{\Wc} \ar@{={ >>}}[r] & L \ar@{=>}[r] & 1.}
\end{equation}
Let, indeed, $A$ be an object of $\Ac$ and $p\colon A_{0}\to A$ a presentation of $A$. Then Proposition~\ref{Proposition-LeftRightPullbacks} implies that the upper sequence in the diagram
\[
\xymatrix{
0 \ar@{.>}[r] & K[p]\cap L(A_0) \ar@{}[dr]|<<<{\pullback} \ar@{{ >}->}[r] \mond & L(A_0) \ar@{-{ >>}}[r]^-{Lp} \mond_{\mu_{A_{0}}} & L(A) \mond^{\mu_{A}} \ar@{.>}[r] & 1 \\
0 \ar[r] & K[p] \ar@{{ >}->}[r]_-{\ker p} & A_0 \ar@{-{ >>}}[r]_-p & A \ar[r] & 1. }
\]
is short exact, since $\mu_{A}$ is a monomorphism and the left hand square a pullback.
\end{remark}

To prove Theorem~\ref{Theorem-L-Exact-Sequence}, we need the following two technical lemmas concerning the functor~$L_{1}$. The second one essentially shows that Fr\"ohlich's $V_{1}$~\cite{Froehlich} and Furtado-Coelho's $L_{1}$~\cite{Furtado-Coelho} coincide, although Fr\"ohlich demands $V_{1}p$ to be proper in $A_{0}$, and Furtado-Coelho only demands that $L_{1}p$ is proper in $L_{0} p\subset A_{0}$.

\begin{lemma}\label{Lemma-L1-Regepis}
Let $L\colon \Ac \to \Ac$ be a functor. If a morphism of presentations $\vetsf \colon p\to q$ is a pullback square
\[
\xymatrix{A_{0} \ar@{}[rd]|<<{\pullback} \ar@{-{ >>}}[d]_-{p} \ar@<-0.5 ex>@{-{ >>}}[r]_-{f_{0}} & B_{0} \ar@<-0.5 ex>[l]_-{s_{0}} \ar@{-{ >>}}[d]^-{q} \\
A  \ar@{-{ >>}}[r]_-{f} & B}
\]
in $\Ac$ with $f_{0}$ split epi, then $L_{1}\vetsf$ is a regular epimorphism.
\end{lemma}
\begin{proof}
First note that, if square \texttt{(ii)} of Diagram~\ref{Diagram-L1-Step1} is a pullback, then so are both squares \texttt{(i)}. Consequently, if, moreover, $f_{0}$ is split epi, $LR\vetsf$ is a split, hence regular, epimorphism. Because, in this case, also $Lf_{0}$ is split epi, Proposition~\ref{Rotlemma-kernel-pairs} implies that square \texttt{(iii)} of Diagram~\ref{Diagram-L1-Step2} is a regular pushout, so $L_{1}\vetsf$ is regular epi by Proposition~\ref{Rotlemma}.
\end{proof}

\begin{lemma}\label{Lemma-L_1}
Suppose that $L\triangleleft 1_{\Ac}$, and let $p\colon A_{0}\to A$ be a presentation in $\Ac$. Then $L_{1}p$ is proper in $A_{0}$.
\end{lemma}
\begin{proof}
Taking the kernel pair of $p$
\[
\xymatrix{\KP [p] \ar@{-{ >>}}[r]^-{k_{0}} \ar@{-{ >>}}[d]_-{k_{1}} \ar@{}[rd]|<{\pullback}& A_{0} \ar@{-{ >>}}[d]^-{p}\\
A_{0} \ar@{-{ >>}}[r]_-{p} & A}
\]
yields split epis $k_{0}$ and $k_{1}$ and defines a map of presentations ${(k_{0},p)\colon k_{1}\to p}$. Now, by Lemma~\ref{Lemma-L1-Regepis}, the arrow $L_{1} (k_{0},p)$ in the square
\[
\xymatrix{L_{1} k_{1} \ar[d]_-{L_{1} (k_{0},p)} \ar@{{ >}->}[r] & R[p] \ar@{-{ >>}}[d]^-{k_{0}} \\
L_{1}p \monr & A_{0}}
\]
 is regular epi.  (In fact it is an isomorphism.)  Moreover, $L_{1} k_{1}$ is proper in $R[p]$, since, according to Corollary~\ref{Corollary-Split-Presentation}, it is the intersection $K[k_{1}]\cap L (R[p])$ of two proper subobjects of $R[p]$.  Proposition~\ref{Proposition-Image-of-Kernel} now implies that $L_{1}p$ is proper in $A_{0}$.
\end{proof}

We shall now apply the Snake Lemma to obtain a six term exact sequence of Baer invariants. In what follows, $\Ac$ is an exact and sequentiable category with enough projectives. Let
\[
\xymatrix{0 \ar[r] & K \ar@{{ >}->}[r] & A \ar@{-{ >>}}[r]^-{f} & B \ar[r] & 1}
\]
be a short exact sequence in $\Ac$. By the naturality of Sequence~\ref{4exL} we have a commutative square~\texttt{(ii)} and thus we get a factorization $\gamma$, such that \texttt{(i)} commutes:
\begin{equation}
\label{gamma-normaal?}
\vcenter{\xymatrix{
0 \ar[r] & K[DLf] \ar@{}[rd]|{\texttt{(i)}} \ar@{{ >}->}[r] \ar@{.>}[d]_-{\gamma} & DL(A) \ar@{}[rd]|{\texttt{(ii)}} \ar[r]^-{DLf} \ar[d]_-{\alpha_A} & DL(B) \ar[d]^-{\alpha_B} \ar@{.>}[r] & 1 \\
0 \ar[r] & K \ar@{{ >}->}[r] & A \ar@{-{ >>}}[r]_-{f} & B \ar[r] & 1.}}
\end{equation}
A priori, only the left exactness of the upper row is clear. Nevertheless, it is easily shown that $DLf$ is a regular epimorphism by choosing a projective presentation $p\colon A_{0}\to A$ of $A$ and then using the map $(1_{A_{0}},f)\colon p\to f\comp p$ of $\Wproj$.

\begin{theorem}\label{Theorem-L-Exact-Sequence}\cite[Theorem 9]{Furtado-Coelho}
Let $\Ac$ be an exact and sequentiable category with enough projectives. Consider a short exact sequence
\[
\xymatrix{0 \ar[r] & K \ar@{{ >}->}[r] & A \ar@{-{ >>}}[r]^-{f} & B \ar[r] & 1}
\]
in $\Ac$. If $L\triangleleft 1_{\Ac}$ there is an exact sequence
\begin{equation}\label{langexactL}
\resizebox{\textwidth}{!}{\xymatrix{0 \ar[r] & K[\gamma] \ar@{{ >}->}[r] & \Delta L(A)  \ar[r]^-{\Delta Lf} & \Delta L(B)  \ar[r] & \frac{K}{I[\gamma ]}  \ar[r] & M(A) \ar@{-{ >>}}[r]^-{Mf} & M(B)  \ar[r] & 1.}}
\end{equation}
This exact sequence depends naturally on the given short exact sequence.
\end{theorem}
\begin{proof}
It suffices to prove that $\gamma$ is proper. Indeed, we already know that $\alpha_A$ and $\alpha_B$ are proper, because $L$ is a proper subfunctor, and thus we get~\ref{langexactL} by applying the Snake Lemma~\ref{snake-lemma}. The naturality then follows from the Snake Lemma and Proposition~\ref{exL}.

Choose projective presentations $p\colon A_{0}\to A$ and $f\comp p\colon A_{0}\to B$ of $A$ and $B$ as above.  Then, by the Noether's First Isomorphism Theorem~\ref{generalexact}, Diagram~\ref{gamma-normaal?} becomes
\[
\xymatrix{
0 \ar[r] & \frac{L_1(f\comp p)}{L_1p} \ar@{{ >}->}[r] \ar[d]_-{\gamma} & \frac{L(A_0)}{L_1p} \ar@{-{ >>}}[r] \ar[d]_-{\alpha_A} & \frac{L (A_0)}{L_1(f\comp p)} \ar[d]^-{\alpha_B} \ar[r] & 1 \\
0 \ar@{{ >}->}[r] & \frac{K[f\comp p]}{K[p]} \ar@{{ >}->}[r] & \frac{A_0}{K[p]} \ar@{-{ >>}}[r] & \frac{A_0}{K[f\comp p]} \ar[r] & 1. }
\]
The map $\gamma$ is unique for the diagram with exact rows
\[
\xymatrix{
0 \ar[r] & {L_{1}p} \ar@{{ >}->}[r] \ar@{{ >}->}[d] & L_1(f\comp p) \ar@{-{ >>}}[r] \ar@{{ >}->}[d] & \frac{L_1(f\comp p)}{L_1p} \ar[d]^-{\gamma} \ar[r] & 1\\
0 \ar[r] & K[p] \ar@{{ >}->}[r] & K[f\comp p] \ar@{-{ >>}}[r] & \frac{K[f\comp p]}{K[p]} \ar[r] & 1 }
\]
 to commute.  By Lemma~\ref{Lemma-L_1}, we get that $L_1(f\comp p)$ is proper in $A_{0}$, hence in $K[f\comp p]$. Proposition~\ref{Proposition-Image-of-Kernel} now implies that $\im \gamma$ is a kernel; hence $\gamma$ is proper.
\end{proof}

\section{Birkhoff subfunctors}\label{Section-Birkhoff-Subfunctors}
We now refine the main result of Section~\ref{Section-Exact-Sequences}, Theorem~\ref{Theorem-L-Exact-Sequence}, by putting extra conditions on the subfunctor $L\subseteq 1_{\Ac }$. To this aim we introduce the notion of \textit{Birkhoff subfunctor}. We show that these subfunctors correspond bijectively to the Birkhoff subcategories of $\Ac$~\cite{Janelidze-Kelly}. But first we prove that $L_{1}\colon \prA \to \Ac$ preserves regular epimorphisms if $L\colon \Ac \to \Ac$ does so; and what are the regular epimorphisms of $\prA$?

\begin{proposition}\label{Proposition-PrA-Regepis}
If $\Ac$ is an exact Mal'tsev category, then the regular epimorphisms of $\prA$ are exactly the regular pushout squares in $\Ac$.
\end{proposition}
\begin{proof}
$\prA$ is a full subcategory of the ``category of arrows'' $\Fun (\Two , \Ac)$, which has the limits and colimits of $\Ac$, computed pointwise (see Section II.4 of Mac Lane~\cite{MacLane} or Section I.2.15 of Borceux~\cite{Borceux:Cats}). Hence, given a pair $\vetsf ,\vetsg \colon p\to q$ of parallel arrows in $\prA$, their coequalizer in $\Fun (\Two , \Ac)$ exists:
\[
\xymatrix{A_{0} \ar@{-{ >>}}[d]_-{p} \ar@<0.5 ex>[r]^-{f_{0}} \ar@<-0.5 ex>[r]_-{g_{0}} & B_{0} \ar@{}[rd]|{\texttt{(i)}} \ar@{-{ >>}}[d]^-{q} \ar@{-{ >>}}[r]^-{c_{0}} & \coeq (f_{0},g_{0}) \ar@{-{ >>}}[d]^-{r} \\
A \ar@<0.5 ex>[r]^-{f} \ar@<-0.5 ex>[r]_-{g} & B \ar@{-{ >>}}[r]_-{c} & \coeq (f,g).}
\]
Since $r$ is a regular epimorphism, $(c_{0},c)$ is an arrow in $\prA$; moreover, by Proposition~\ref{Rotlemma-kernel-pairs}, the right hand square \texttt{(i)} is a regular pushout.

Conversely, given a regular pushout square such as \texttt{(i)}, Proposition~\ref{Rotlemma-kernel-pairs} ensures that its kernel pair in $\Fun (\Two , \Ac)$ is a pair of arrows in $\prA$. Their coequalizer is again the regular pushout square \texttt{(i)}.
\end{proof}

\begin{proposition}\label{V_1 bewaart regulier epis}
Let $\Ac$ be an exact and sequentiable category, and let $L\colon {\Ac \to \Ac} $ be a functor. If $L$ preserves regular epimorphisms, then so does $L_{1}$.
\end{proposition}
\begin{proof}
Suppose that $L$ preserves regular epis. If $\vetsf$ is a regular epimorphism of $\prA$, square \texttt{(ii)} of Diagram~\ref{Diagram-L1-Step1} is a regular pushout. Proposition~\ref{Rotlemma-kernel-pairs} implies that $R\vetsf$ is a regular epimorphism. By assumption, $Lf_{0}$ and $LR\vetsf$ are regular epis; hence, again using Proposition~\ref{Rotlemma-kernel-pairs}, we get that square \texttt{(iii)} of Diagram~\ref{Diagram-L1-Step2} is a regular pushout. An application of Proposition~\ref{Rotlemma} on~\ref{Diagram-L1-Step3} shows that $L_{1}\vetsf$ is regular epi, which proves the first statement.
\end{proof}

Recall the following definition from Janelidze and Kelly~\cite{Janelidze-Kelly}.

\begin{definition}\label{Definition-Birkhoff-Subcategory}\index{U@$U$}
Let $\Ac$ be an exact category. A reflective subcategory $\Bc$ of $\Ac$ that is full and closed in $\Ac$ under subobjects and quotient objects is said to be a \emph{Birkhoff subcategory}\index{subcategory!Birkhoff ---}\index{Birkhoff subcategory} of $\Ac$. We will denote by $U\colon \Ac \to \Bc$ the left adjoint of the inclusion $I\colon \Bc \to \Ac$ and by $\eta \colon 1_{\Ac}\To I\comp U$ the unit of the adjunction. We shall omit all references to the functor $I$, writing $\eta_{A}\colon {A\to U (A)}$ for the component of $\eta$ at $A\in \Ob{\Ac}$.
\end{definition}

\begin{examples}\label{Examples-Birkhoff}
We have already seen in Chapter~\ref{Chapter-Preliminaries} that, in an exact Mal'tsev category, the subcategory of all abelian objects is Birkhoff. A subcategory of a variety of universal algebras is Birkhoff if and only if it is a subvariety. 

In Section~\ref{Section-Nilpotency} we shall see that, in a semi-abelian category, the objects of a given nilpotency class form a Birkhoff subcategory. The observation that, in an exact Mal'tsev category with coequalizers, the category of internal groupoids is Birkhoff in the category of internal reflexive graphs (see \cite{Gran:Internal} and Section~\ref{Section-Preliminaries}), forms one of the starting points to the paper~\cite{EG}.
\end{examples}

\begin{proposition}\label{Proposition-Birkhoff-Properties}\cite{Janelidze-Kelly}
Let $\Ac$ be an exact category and $U\colon \Ac \to \Bc$ a reflector (with unit $\eta$) onto a full replete subcategory $\Bc$ of $\Ac$. Then
\begin{enumerate}
\item $\Bc$ is closed in $\Ac$ under subobjects when the component $\eta_{A}\colon {A\to U (A)}$ of $\eta$ at an $A\in \Ob{\Ac}$ is a regular epimorphism;
\item if 1.\ holds, then $\Bc$ is closed in $\Ac$ under quotient objects if and only if, for any regular epimorphism $f\colon A\to B$ in $\Ac$, the naturality square
\[
\xymatrix{A \ar@{-{ >>}}[d]_-{f} \ar@{-{ >>}}[r]^-{\eta_{A}} & U (A) \ar@{-{ >>}}[d]^-{Uf}\\
B \ar@{-{ >>}}[r]_-{\eta_{B}} & U (B)}
\]
is a pushout.\noproof
\end{enumerate}
\end{proposition}

\begin{corollary}\label{Corollary-Birkhoff-Subcategory-is-Semiabelian}
Suppose that $\Ac$ is an exact category and $\Bc$ is a Birkhoff subcategory of~$\Ac$. Then $\Bc$ is exact. Moreover, if $\Ac$ is, respectively, quasi-pointed, pointed, protomodular or finitely cocomplete, then so is $\Bc$. In particular, a Birkhoff subcategory of a semi-abelian category is semi-abelian.
\end{corollary}
\begin{proof}
In~\cite{Janelidze-Kelly} it is shown that a Birkhoff subcategory $\Bc$ of an exact category $\Ac$ is necessarily exact.

If $\Ac$ has an initial object, then $U (0)$ is initial in $\Bc$, and $\eta_{0}\colon {0\to U (0)}$ is a regular epimorphism. If $\Ac$ is quasi-pointed then $\eta_{0}$ is an isomorphism, hence $0$ belongs to $\Bc$. Now the regular epimorphism $\eta_{1}\colon {1\to U (1)}$ turns out to be split monic, and $1$ is in $\Bc$. It follows that $\Bc$ has an initial and a terminal object, and one sees immediately that the monomorphism ${0\to 1}$ in $\Ac$ is still a monomorphism in $\Bc$. Obviously, if now $\Ac$ has a zero object, then so has $\Bc$. 

$\Bc$ fulfils the requirement in Definition~\ref{Definition-Protomodularity} when it holds in $\Ac$ because $I$ preserves pullbacks and $\Bc$ is full. Finally, a full and reflective subcategory $\Bc$ of $\Ac$ has all colimits that exist in~$\Ac$.
\end{proof}

If the functor $L$ from Remark~\ref{Remark-Birkhoff-Nabla} is a proper subfunctor of $1_{\Ac}$, $L$ satisfies the conditions Fr\"ohlich used in the article~\cite{Froehlich} to obtain his Baer invariants. Because of Proposition~\ref{Proposition-Birkhoff-Subfunctor}, and its importance in what follows, we think this situation merits a name and slightly different notations.

\begin{definition}\label{Definition-Birkhoff-Subfunctor}\index{V@$V$}
Let $\Ac$ be an exact sequentiable category. We call \emph{Birkhoff subfunctor}\index{Birkhoff subfunctor}\index{subfunctor!Birkhoff ---} of $\Ac$ any proper subfunctor $V$ of $1_{\Ac}$ that preserves regular epimorphisms. The corresponding functors are denoted $V_{0}$, $V_{1}$, etc.\index{V_0@$V_{0}$}\index{V_1@$V_{1}$}
\end{definition}

The following establishes a bijective correspondence between Birkhoff subfunctors and Birkhoff subcategories of a given exact and sequentiable category $\Ac$.

\begin{proposition}\label{Proposition-Birkhoff-Subfunctor}\cite[Theorem 1.2]{Froehlich}
Let $\Ac$ be an exact and sequentiable category.
\begin{enumerate}
\item If $U\colon \Ac \to \Bc$ is the reflector of $\Ac$ onto a Birkhoff subcategory $\Bc$, setting $V (A)=K[\eta_{A}]$ and $\mu_{A}=\ker \eta_{A}$ defines a Birkhoff subfunctor $\mu \colon {V\To 1_{\Ac}}$ of $\Ac$.
\item Conversely, if $\mu \colon V\To 1_{\Ac}$ is a Birkhoff subfunctor of $\Ac$, putting $U (A)=\Cok [\mu_{A}]$, $\eta_{A}=\coker \mu_{A}$ defines a full functor $U\colon \Ac \to \Ac$ and a natural transformation $\eta \colon 1_{\Ac}\To U$. Here we can, and will, always choose $\coker (0\to A)=1_A$. The image of $U$ is a Birkhoff subcategory $\Bc$ of $\Ac$. Furthermore, $U$, considered as a functor $\Ac \to \Bc$, is left adjoint to the inclusion $\Bc \to \Ac$, and $\eta$ is the unit of this adjunction.
\end{enumerate}
In both cases, for any $A\in\Ob{\Ac}$, the sequence
\[
\xymatrix{0 \ar[r] & V (A) \ar@{{ >}->}[r]^-{\mu_{A}} & A \ar@{-{ >>}}[r]^-{\eta_{A}} & U (A) \ar[r] & 1}
\]
is exact.
\end{proposition}
\begin{proof}
For 1.\ we only need to prove that $V$ preserves regular epimorphisms. But this follows immediately from Proposition~\ref{Proposition-Birkhoff-Properties} and Proposition~\ref{Rotlemma}, applied to the diagram
\begin{equation}\label{Diagram-Rot-Toepassing}
\vcenter{\xymatrix{0 \ar[r] & V (A)\ar[d]_-{Vf} \ar@{{ >}->}[r]^-{\mu_{A}} &  A \ar@{-{ >>}}[d]_-{f} \ar@{-{ >>}}[r]^-{\eta_{A}} \ar@{}[rd]|>>{\pushout} & U (A) \ar@{-{ >>}}[d]^-{Uf} \ar[r] & 1\\
0 \ar[r] & V (B) \ar@{{ >}->}[r]_-{\mu_{B}} & B  \ar@{-{ >>}}[r]_-{\eta_{B}} & U (B) \ar[r] & 1.}}
\end{equation}

To prove 2, we first show that, for any $C\in\Ob\Ac$,  $V(U(C))=0$. In the diagram
\[
\xymatrix{
0 \ar[r] & V(C) \ar@{{ >}->}[r]^-{\mu_A} \ar@{-{ >>}}[d]_-{V\eta_C} &C \ar@{-{ >>}}[r]^-{\eta_C} \ar@{-{ >>}}[d]_-{\eta_C} & U(C) \ar[r] \ar[d]^-{U\eta_C} & 1
\\ 0 \ar[r] & V(U(C)) \ar@{{ >}->}[r]_-{\mu_{U(C)}} & U(C) \ar@{-{ >>}}[r]_-{\eta_{U(C)}} & U(U(C)) \ar[r] & 1,}
\]
the right square is a pushout, because $V\eta_{C}$ is an epimorphism, as $V$ preserves regular epimorphisms. It follows easily that $\eta_{U(C)}$ is a split monomorphism (hence an isomorphism), and $V(U(C))=0$.

For the fullness of $U$, any map $h\colon U (A)\to U (B)$ should be $Ug$ for some $g$. But $h$ itself is such a $g$: Indeed $\eta_{U (C)}=1_{U (C)}$, as we just demonstrated that $\mu_{U (C)}=0$, for any $C\in \Ob{\Ac}$.

Now let $A$ be any object, and $f\colon A\to U (B)$ any map, of $\Ac$. To prove that a unique arrow $\overline{f}$ exists in $\Bc$ such that
\[
\xymatrix{
0 \ar[r] & V(A) \ar@{{ >}->}[r]^-{\mu_A} & A \ar@{-{ >>}}[r]^-{\eta_A} \ar[rd]_-f & U(A) \ar@{.>}[d]^-{\overline{f}} \ar[r] & 1\\
&&& U(B)& }
\]
commutes, it suffices to show that $f\comp \mu_A=0$. (The unique $\overline{f}$ in $\Ac$ such that the diagram commutes is a map of $\Bc$, because $U$ is full.) But, as $V$ is a subfunctor of the identity functor, $f\comp \mu_A=\mu_{U(B)}\comp Vf=0$.

This shows that for any object $A$ of $\Ac$, the morphism $\eta_{A}\colon A\to U (A)$ is a universal arrow from $A$ to $U(A)$. It follows that $U$ is a reflector with unit $\eta$. The subcategory $\Bc$ is then Birkhoff by Proposition~\ref{Proposition-Birkhoff-Properties}, as in Diagram~\ref{Diagram-Rot-Toepassing} the right square is a pushout, $Vf$ being a regular epimorphism.
\end{proof}

\begin{remark}\label{Remark-Birkhoff-Exact-Sequence}
In case $V$ is a Birkhoff subfunctor of $\Ac$, note that, by Remark~\ref{Remark-Birkhoff-Nabla} and Proposition~\ref{Proposition-Birkhoff-Subfunctor}, the exact sequence~\ref{4exL} becomes\index{Baer invariant!$\Delta V$}\index{Delta V@$\Delta V$}\index{Baer invariant!$DV$}\index{DV@$DV$}\index{Baer invariant!$U$}\index{U@$U$}
\begin{equation}\label{4exV}
\xymatrix{0 \ar@{=>}[r] & \Delta V \ar@{{|>}=>}[r] & DV \ar@{=>}[r] & 1_{\Ac} \ar@{={ >>}}[r] & U \ar@{=>}[r] & 1.}
\end{equation}
\end{remark}

Theorem~\ref{Theorem-L-Exact-Sequence} now has the following refinement.

\begin{theorem}\label{Theorem-V-Exact-Sequence}\cite[Theorem 3.2]{Froehlich}
Let $\Ac$ be an exact and sequentiable category with enough projectives. Consider a short exact sequence
\begin{equation}\label{Sequence-SE}
\xymatrix{0 \ar[r] & K \ar@{{ >}->}[r] & A \ar@{-{ >>}}[r]^-{f} & B \ar[r] & 1}
\end{equation}
in $\Ac$. If $V$ is a Birkhoff subfunctor of $\Ac$, then the sequence
\begin{equation}\label{langexactV}
\resizebox{\textwidth}{!}{\xymatrix{0 \ar[r] & K[\gamma] \ar@{{ >}->}[r] & \Delta V(A)  \ar[r]^-{\Delta Vf} & \Delta V(B)  \ar[r] & \frac{K}{V_{1}f}  \ar[r] & U(A) \ar@{-{ >>}}[r]^-{Uf} & U(B)  \ar[r] & 1}}
\end{equation}
is exact and depends naturally on the given short exact sequence.
\end{theorem}
\begin{proof}
To get~\ref{langexactV} from Sequence~\ref{langexactL}, it suffices to recall Remark~\ref{Remark-Birkhoff-Exact-Sequence}, and to prove that the image $I[\gamma ]$ of $\gamma$ is $V_1f$. If we choose presentations $p\colon {A_{0}\to A}$ and $f\comp p\colon A_{0}\to B$ as in the proof of Theorem~\ref{Theorem-L-Exact-Sequence}, $\gamma$ becomes a map
\[
\frac{V_1(f\comp p)}{V_1p}\to\frac{K[f\comp p]}{K[p]}.
\]
Hence, to prove that $I[\gamma]=V_1f$, we are to show that the arrow
\[
\frac{V_1(f\comp p)}{V_1p}\to V_1f
\]
is regular epi. But this is equivalent to $V_{1}(p,1_{B})\colon  V_1(f\comp p)\to V_1f$ being regular epi. This is the case, since by Proposition~\ref{V_1 bewaart regulier epis}, $V_1$ preserves regular epimorphisms, and
\[
\xymatrix{A_{0} \ar@{-{ >>}}[r]^-{p} \ar@{-{ >>}}[d]_-{f\comp p} & A \ar@{-{ >>}}[d]^-{f}\\
B \ar@{=}[r] & B}
\]
is a regular pushout square, which means that $(p,1_{B})\colon  (f\comp p)\to f$ is a regular epimorphism of $\prA$.
\end{proof}

\begin{examples}\label{Examples-V-Exact-Sequence}
Consider the category $\Gp$ of groups and its Birkhoff subcategory $\Ab =\Ab\Gp$ of abelian groups. The associated Birkhoff subfunctor (cf.\ Proposition~\ref{Proposition-Birkhoff-Subfunctor} and Examples~\ref{Examples-Huq-Central}) sends $G$ to the commutator subgroup $[G,G]$. It is indeed well known that the abelianization $\ab (G)$ of a group $G$, i.e., the reflection of $G$ along the inclusion ${\Ab \to \Gp}$, is just ${G}/{[G,G]}$. Proposition~19 of Furtado-Coelho~\cite{Furtado-Coelho} states that, if $V\colon \Gp\to\Gp$ is the Birkhoff subfunctor defined by $V(G)=[G,G]$, then
\[
V_1\bigl(\eta_{N}\colon G\to \tfrac{G}{N}\bigr)=[N,G].
\]
Hence the exact sequence~\ref{langexactV} induced by a short exact sequence of groups~\ref{Sequence-SE} becomes the well-known sequence
\[
\xymatrix{{\frac{R\cap [F,F]}{[R,F]}} \ar[r] & {\frac{R'\cap [F',F']}{[R',F']}} \ar[r] & {\frac{K}{[K,A]}} \ar[r] & {\ab(A)} \ar@{-{ >>}}[r]^-{\ab f} & {\ab(B)}  \ar[r] & 0,}
\]
with $A=F/R$ and $B=F'/R'$, and where $F$ and $F'$ are free groups.

For abelianization of Lie algebras we can say the same: Since $\Lie_{\K}$ is a variety of $\Omega$-groups and $\ab (\Lieg)=\Lieg /[\Lieg ,\Lieg]$, \cite[Proposition 19]{Furtado-Coelho} implies that 
\[
V_1(\eta_{\Lien}\colon \Lieg\to {\Lieg}/{\Lien})=[\Lien,\Lieg].
\]
\end{examples}

\begin{remark}\label{Remark-Splitting}
Recalling Corollary~\ref{Corollary-Split-Presentation}, note that in case $A$ (or $B$) is projective, the exact sequences~\ref{langexactL} and~\ref{langexactV} become much shorter, because then $\Delta L(A)$ and $\Delta V(A)$ (or $\Delta L(B)$ and $\Delta V(B)$) are $0$.
\end{remark}

\section{Several commutators compared}\label{Section-Commutators}\index{commutator}
In this section we compare several notions of commutator. There is Huq's commutator of coterminal arrows~\cite{Huq}\index{Huq commutator} considered in Section~\ref{Section-Semi-Abelian-Categories}, Smith's commutator\index{Smith commutator} of equivalence relations~\cite{Smith} and Fr\"ohlich's functor $V_{1}$~\cite{Froehlich}. Over the last decade, many people have been studying categorical versions of these commutators and the related topic of central extensions. The foundational work includes~\cite{Carboni-Pedicchio-Pirovano, Pedicchio} (where a categorical version of Smith's commutator is considered),~\cite{Janelidze-Kelly} (where a general theory of central extensions is introduced) and~\cite{Bourn2002, Bourn2004, Bourn-Huq, Borceux-Bourn, Borceux-Semiab, BG} (where a categorical version of Huq's commutator is studied, see Section~\ref{Section-Semi-Abelian-Categories}). Janelidze and Pedicchio~\cite{Janelidze-Pedicchio} study commutator theory in a very general context (that of \emph{finitely well-complete}\index{finitely well-complete category}\index{category!finitely well-complete} categories) by basing it on the theory of internal categorical structures. Here we consider $V_{1}$ from the point of view of categorical commutator theory and investigate when and how it relates to Huq's or Smith's commutator.

Recall from Chapter~\ref{Chapter-Preliminaries} that in the semi-abelian context an intrinsic notion of abelian object\index{abelian!object} exists. In a strongly unital category, an object is abelian if and only if it can be provided with a (necessarily unique) structure of internal abelian group; in a Mal'tsev category, if and only if it admits an internal Mal'tsev operation. The two notions coincide as soon as the category is pointed and Mal'tsev.

In~\cite{Gran:Central-Extensions} it is proven that the full subcategory of abelian objects $\Ab\Ac$ in any exact Mal'tsev category with coequalizers $\Ac$ is a Birkhoff subcategory of $\Ac$. Thus it is fairly clear that the associated Birkhoff subfunctor $V$ behaves like a commutator. Moreover, one could ask whether it makes sense to view $V_{1}$ this way. In view of Examples~\ref{Examples-V-Exact-Sequence}, this seems feasible. In what follows, we answer the question positively, and show that, with respect to abelianization, the categorical counterparts of Huq's and Smith's commutators induce the same notion of central extension.\pfbreak

In~\cite{Pedicchio}, M.\,C. Pedicchio extends Smith's notion of commutator of equivalence relations~\cite{Smith} (defined in the context of Mal'tsev varieties) to the context of exact Mal'tsev categories with coequalizers. 

\begin{definition}\label{Definition-Centralizing-Equivalence-Relations}~\cite{Smith, Carboni-Pedicchio-Pirovano, Pedicchio}
Let $\Ac$ be a finitely complete category. An equivalence relation in $\Eq \Ac$\index{category!$\Eq \Ac$} is called a \emph{double equivalence relation}\index{double equivalence relation} in $\Ac$. Let $\nabla_{X}= (X\times X,\pr_{1},\pr_{2})$\index{nabla@$\nabla_{X}$} and $\Delta_{X}= (X,1_{X},1_{X})$\index{Delta X@$\Delta_{X}$} denote, respectively, the largest and the smallest equivalence relation on $X$. 

Two equivalence relations $(R,r_{0},r_{1})$ and $(S,s_{0},s_{1})$ on an object $X$ of $\Ac$ \emph{centralize (in the sense of Smith)}\index{centralize (in the sense of Smith)} or \emph{admit a centralizing relation}\index{admit a centralizing relation} when there exists a double equivalence relation $C$ on $R$ and $S$ such that any commutative square in the diagram
\[
\xymatrix{
C \ar@<-0.5 ex>[r]_-{\sigma _1} \ar@<0.5 ex>[r]^-{\sigma _0} \ar@<-0.5 ex>[d]_-{\rho_0} \ar@<0.5 ex>[d]^-{\rho_1} & R \ar@<-0.5 ex>[d]_-{r_{0}} \ar@<0.5 ex>[d]^-{r_{1}} \\
S  \ar@<-0.5 ex>[r]_-{s_{1}} \ar@<0.5 ex>[r]^-{s_{0}} & X }
\]
is a pullback. In particular, an equivalence relation $R$ on $X$ is \emph{central}\index{central equivalence relation}\index{Smith central} when $R$ and $\nabla_{X}$ admit a centralizing relation.
\end{definition}

In an exact category with coequalizers, an object $A$ is abelian if and only if $\nabla_{A}$ and $\nabla_{A}$ centralize. Thus we see that, in the pointed case, Huq's commutator and Smith's at least give rise to the same notion of abelian object. However, it is well-known that, in general, two equivalence relations $R$ and $S$ need not centralize when their normalizations $k_{R}$ and $k_{S}$ cooperate, not even in a variety of $\Omega$-groups (see~\cite{Bourn2004} for a counterexample). Yet, we are now going to show that this \textit{is} the case in any pointed protomodular category, whenever $R$ (or $S$) is $\nabla_{X}$:

\begin{proposition}\label{Proposition-Huq-Centrality-is-Smith-Centrality}
An equivalence relation $R$ in a pointed protomodular category is central if and only if its associated normal subobject $k_{R}$ is central.
\end{proposition}
\begin{proof}
Let us first assume that the equivalence relation 
\[
\xymatrix{{R}  \ar@<-1 ex>[rr]_-{\pr_{1}} \ar@<1 ex>[rr]^-{\pr_{2}} && {X} \ar[ll]|-{\Delta}}
\]
is central, and let $C$ be the associated centralizing double relation on $R$ and $\nabla_{X}$. We can then consider the following diagram
\[
\xymatrix{C \ar[d]_{p_1}  \ar@<-.5 ex>[r]_{p_2} \ar@<.5 ex>[r]^{p_1} & R \ar@{}[rd]|-{\texttt{(i)}} \ar[d]_-{\pr_{1}} \ar@{.{ >>}}[r] & K  \ar@{.>}[d] \\
X \times X \ar@<-.5 ex>[r]_(.6){\pr_{2}} \ar@<.5 ex>[r]^(.6){\pr_{1}}  & X \ar@{.{ >>}}[r] & 0,}
\]
where the square \texttt{(i)} is obtained by taking coequalizers. Since $C$ is centralizing, both left hand side squares are pullbacks, hence so is \texttt{(i)}, and $R \cong X \times K$. This moreover induces the commutative diagram
\[
\xymatrix{K \ar[d] \ar[r]^-{k} & R \ar@{}[rd]|<{\pullback} \ar[d]_-{\pr_{1}} \ar[r]^-{p} & K \ar[d] \\
0 \ar[r]  & X \ar[r] & 0}
\]
where $p\comp k=1_{K}$, and both the outer rectangle and the right hand square are pullbacks. It follows that $K$ is the normal subobject associated with $R$. By considering also the second projection $\pr_{2}$ from $R$ to $X$ one can easily check that $R$ is then canonically isomorphic to the equivalence relation
\begin{equation}\label{Diagram-Centralizing-Relation}
\xymatrix{{K\times X} \ar@<-1 ex>[rr]_-{\pr_{X}} \ar@<1 ex>[rr]^-{\varphi_{k_{R},1_{X}}} && {X}; \ar[ll]|-{r_{X}}}
\end{equation}
the arrow $\varphi_{k_{R},1_{X}}$ is the needed cooperator.

Conversely, let us suppose that $k$ is a central monomorphism with cooperator $\varphi_{k,1_{X}}$. One can then form the reflexive graph~\ref{Diagram-Centralizing-Relation}---call it $R_{k}$; it is a relation because 
\[
(\varphi_{k,1_{X}},\pr_{2})\colon {K\times X\to X\times X}
\]
equals $k\times 1_{X}\colon K\times X\to X\times X$. Since $R_{k}$ is a reflexive relation in a Mal'tsev category, it is an equivalence relation. One may check that this equivalence relation corresponds to $k$ via the bijection between normal subobjects and equivalence relations.

Furthermore, consider the double equivalence relation determined by the kernel pair $R[\pr_{K}]$ of $\pr_{K}\colon K\times X\to K$:
\[
\xymatrix{R[\pr_{K}] \ar@<-.5 ex>[d] \ar@<.5 ex>[d] \ar@<-.5 ex>[r] \ar@<.5 ex>[r] & X\times X \ar@<-.5 ex>[d]_-{\pr_{1}} \ar@<.5 ex>[d]^-{\pr_{2}}\\
K\times X \ar@<-.5 ex>[r]_-{\pr_{X}} \ar@<.5 ex>[r]^-{\varphi_{k,1_{X}}} & X.}
\]
It is clearly a centralizing double relation on $R_{k}$ and $\nabla_{X}$, as desired.
\end{proof}

\begin{definition}\label{Definition-Central-Extension}
In a sequentiable category, an \emph{extension} of an object $Y$ (\emph{by} an object~$K$)\index{extension} is a regular epimorphism $f\colon X\to Y$ with its kernel $K$:
\[
0 \to K \to/{{ >}->}/^{k} X \to/-{ >>}/^{f} Y \to 1.
\]
The category of extensions of $Y$ (considered as a full subcategory of the slice category $\slashfrac{\Ac}{Y}$) is denoted by $\Ext (Y)$\index{category!$\Ext (Y)$}; the category of all extensions in $\Ac$ (considered as a full subcategory of the arrow category $\Fun (\mathsf{2},\Ac)$: morphisms are commutative squares) by $\Ext \Ac$\index{category!$\Ext\Ac$}.

An extension $f\colon X\to Y$ is called \emph{central}\index{central extension} when its kernel pair is central in the sense of Smith: this is the case when $R[f]$ and $\nabla_{X}$ centralize. We write $\CentrExt (Y)$\index{category!$\CentrExt (Y)$} for the full subcategory of $\Ext (Y)$ determined by the central extensions. 
\end{definition}

\begin{remark}\label{Remark-Ext-is-Pr}
Of course, an extension is the same thing as a presentation, and the categories $\Ext \Ac$ and $\Pr \Ac$ coincide. For historical reasons, we keep both notions, using the one that locally seems most appropriate. 
\end{remark}

\begin{remark}\label{Remark-Huq-Smith-Central-Extension}
If $\Ac$ is a homological category, Proposition~\ref{Proposition-Huq-Centrality-is-Smith-Centrality} implies that, in this definition, we could as well have asked that the kernel of $f$ is central in the sense of Huq, i.e., that $\ker f$ cooperates with $1_{X}$.
\end{remark}

\begin{examples}\label{Examples-Central-Extensions}
Both for groups and Lie algebras, this notion of central extension coincides with the classical one: In each case, such is an extension of which the kernel lies in the centre of the domain.
\end{examples}

The following well-known property of central extensions of groups will now be shown to hold in any homological category:

\begin{proposition}\label{Proposition-Subobject-of-Central-is-Normal}
Let $\Ac$ be a homological category and $f\colon X\to Y$ a central extension in $\Ac$. Every subobject of $k=\ker f\colon K\to X$ is proper in $X$. 
\end{proposition}
\begin{proof}
Let $i\colon M\to K$ be a monomorphism and denote $m=k\comp i$. We are to show that $m$ is a kernel. Now $k$ is central, hence so is $m$, because $m$ cooperates with $1_{X}$; in particular, there exists an arrow $\varphi_{m,1_{X}}\colon M\times X\to X$ satisfying $\varphi_{m,1_{X}}\comp l_{M}=m$. But the arrow $l_{M}$ is the kernel of $\pr_{2}$; hence $m$ is proper by Proposition~\ref{Proposition-Image-of-Kernel}, since it is the regular image of $\varphi_{m,1_{X}}\comp l_{M}$.
\end{proof}

To further investigate centrality using commutators we need some additional notations.

\begin{notation}\label{Notation-Eqrel-Image}
If $\Ac$ is exact and $(R,r_{0},r_{1})$ is an equivalence relation on an object $X$ of $\Ac$, then $\eta_{R}\colon {X\to X/R}$ denotes the coequalizer of $r_{0}$ and $r_{1}$. Now let $\Ac$ be, moreover, sequentiable. We denote the kernel $K[\eta_{R}]$ of the coequalizer $\eta_{R}\colon {X\to X/R}$ by $k_{R}\colon {K_{R}\to X}$.\index{K@$K_{R}$}\index{k@$k_{R}$}\index{H@$\eta_{R}$} Note that this is just the normalization of $R$\index{normalization of an equivalence relation} (see Section~\ref{Section-Bourn-Protomodular-Categories}). This defines a one-one correspondence between (isomorphism classes of) equivalence relations on $X$ and proper subobjects of $X$. Furthermore, this bijection $\Eq (X)\cong \NSub (X)$\index{Eq@$\Eq (X)$}\index{psub@$\NSub (X)$} is an order isomorphism, the order on both classes being (induced by) the usual order on subobjects.

For a regular epimorphism $f\colon X\to Y$ in a regular category $\Ac$ and a relation $(R,r_{0},r_{1})$ on $X$, $f (R)$ denotes the relation obtained from the image factorization
\[
\vcenter{\xymatrix{
R\ar@<-0.5 ex>[d]_-{r_0} \ar@<0.5 ex>[d]^-{r_1} \ar@{-{ >>}}[r] & f(R) \ar@<-0.5 ex>[d]_-{\rho_0} \ar@<0.5 ex>[d]^-{\rho_1}\\
X \ar@{-{ >>}}[r]_-{f} & Y}}
\]
of $(f\times f)\comp (r_{0},r_{1})\colon R\to Y\times Y$. It is called the \emph{direct image of $R$ along $f$}\index{$f (R)$}\index{direct image!of an equivalence relation}. Clearly, when $R$ is reflexive, $f (R)$ is reflexive as well. Hence, when $\Ac$ is Mal'tsev, if $R$ is an equivalence relation, then so is $f(R)$.
\end{notation}

In~\cite{Pedicchio}, a notion of commutator of equivalence relations is introduced, characterized by the following properties.

\begin{definition}\label{Definition-Pedicchio-Commutator}
Let $\Ac$ be an exact Mal'tsev category with coequalizers and $R$ and $S$ equivalence relations on an object $X\in\Ob{\Ac}$. The \emph{commutator (in the sense of Smith)}\index{commutator!of equivalence relations}\index{commutator!in the sense of Smith} $[R,S]$\index{$[R,S]$} is the equivalence relation on $X$ defined by the following properties:
\begin{enumerate}
\item
$\eta_{[R,S]}(R)$ and $\eta_{[R,S]}(S)$ admit a centralizing relation;
\item
for any equivalence relation $T$ on $X$ such that $\eta_T(R)$ and $\eta_T(S)$ admit a centralizing relation, there exists a unique $\delta\colon  X/[R,S]\to X/T$ with $\eta_T=\delta\comp \eta_{[R,S]}$.
\end{enumerate}
\end{definition}

There exist several explicit constructions of this commutator; we give one due to Bourn~\cite{Bourn2004}, because of its beauty and the formal resemblance to Proposition~\ref{Proposition-Definition-Commutator}.

\begin{proposition}\label{Proposition-Construction-Pedicchio-Commutator}\cite{Bourn-Huq, Borceux-Semiab}
In a finitely cocomplete exact Mal'tsev category $\Ac$, let $(R,r_{0},r_{1})$ and $(S,s_{0},s_{1})$ be two equivalence relations on an object $X$. Consider the left hand side pullback diagram below and the dotted colimit of the right hand side diagram of solid arrows
\[
\vcenter{\xymatrix{{R\times_{X} S} \ar@{}[dr]|<<{\pullback} \ar@<.5 ex>[r] \ar@<-.5 ex>[d] & S \ar@<.5 ex>[d]^-{s_{0}} \ar@<.5 ex>[l]^-{r_{S}}\\
R \ar@<-.5 ex>[r]_-{r_{1}} \ar@<-.5 ex>[u]_-{l_{R}} & X \ar@<.5 ex>[u] \ar@<-.5 ex>[l]}}
\qquad \qquad 
\vcenter{\xy
\Atrianglepair(0,500)/->`.>`->``/[R`{R\times_{X} S}`Y`X.;l_{R}``r_{0}``]
\Vtrianglepair/.>`<.`<-`<.`<-/[{R\times_{X} S}`Y`X.`S;`\psi`r_{S}``s_{1}]
\endxy}
\]
The kernel pair of $\psi$ is the commutator $[R,S]$ of $R$ and $S$. $R$ and $S$ centralize if and only if $[R,S]=\Delta_{X}$.\noproof 
\end{proposition}
\vspace{-\baselineskip}\pfbreak

\noindent The notions of centrality we discussed so far are special in a very specific sense: They are relative to the Birkhoff subcategory of abelian objects. Janelidze and Kelly's notion of central extension goes beyond abelianization; their definition is relative to a so-called \emph{admissible}\index{admissible subcategory}\index{subcategory!admissible} subcategory $\Bc$ of an exact category $\Ac$. If $\Ac$ is Mal'tsev or has the weaker \emph{Goursat}\index{Goursat category}\index{category!Goursat} property, then its admissible subcategories are exactly the Birkhoff subcategories of $\Ac$.

We shall not introduce this notion here. Instead, following Fr\"ohlich~\cite{Froehlich}, we use $V_{1}$ to define a notion of central extension relative to a Birkhoff subcategory, and show that it coincides with Janelidze and Kelly's.

\begin{definition}\label{Definition-Centrality}
Let $\Ac$ be an exact and sequentiable category, $V$ a Birkhoff subfunctor of $\Ac$. Then an extension $f\colon X\to Y$ is called \emph{$V$-central}\index{central extension}\index{V-central extension@$V$-central extension}\index{extension!central} when $V_{1}f =0$.
\end{definition}

\begin{examples}\label{Examples-V-Central-Extension}
As will be shown below, in the case of abelianization, this notion of central extension coincides with the one discussed above. Thus we get the classical examples of groups and Lie algebras etc.

A different kind of example is considered by Everaert and Gran in~\cite{EG}. They study the category of precrossed modules\index{category!of precrossed module} over a fixed base group with its Birkhoff subcategory of crossed modules\index{category!of crossed modules} over this base, and give a description of the resulting notion of central extension. (In fact, their characterization is much more general, and includes e.g., the case of precrossed vs.\ crossed rings over a fixed base ring.)
\end{examples}

\begin{proposition}\cite{Janelidze-Kelly, Bourn-Gran}\label{Proposition-Centrality}
Let $\Ac$ be an exact and sequentiable category, $\Bc$ a Birkhoff subcategory of $\Ac$ and $V\colon \Ac \to \Ac$ associated Birkhoff subfunctor. Let $f\colon X\to Y$ be an extension in $\Ac$. Then the following are equivalent:
\begin{enumerate}
\item $f$ is central (relative to $\Bc$) in the sense of Janelidze and Kelly;
\item for every $a,b\colon X'\to X$ with $f\comp a=f\comp b$, one has $V a=V b$;
\item $V (\KP [f])=\Delta_{V (X)}$;
\item $f$ is $V$-central.\noproof
\end{enumerate}
\end{proposition}

The correspondence between these two notions---centrality on one hand, and commutators on the other---was made clear, first by Janelidze and Kelly in~\cite{Janelidze-Kelly:Univ} and~\cite{Janelidze-Kelly:Maltsev}, next by Bourn and Gran in~\cite{Bourn-Gran}, then, most generally, by Gran, in~\cite{Gran-Alg-Cent}. In this article, he proves that, in any \emph{factor permutable}\index{factor permutable category}\index{category!factor permutable} category (hence, in any exact Mal'tsev category with coequalizers), an extension $f\colon {X\to Y}$ is central relative to the subcategory of ``abelian objects'', precisely when (a generalization of) the Smith commutator $[R[f],\nabla_{X}]$ is equal to $\Delta_{X}$. Accordingly, when we choose the subcategory $\Ab\Ac$\index{category!$\Ab\Ac$} of abelian objects as Birkhoff subcategory of $\Ac $, the list of equivalent statements of Proposition~\ref{Proposition-Centrality} can be enlarged as follows:

\begin{theorem}\label{Theorem-Abelian-Centrality}
Consider an exact sequentiable category $\Ac$, its Birkhoff subcategory of abelian objects $\Ab\Ac$ and the associated Birkhoff subfunctor $V$. If $f$ be a presentation in $\Ac$, then the following are equivalent:
\begin{enumerate}
\item $f$ is central in the sense of Janelidze and Kelly;
\item $f$ is $V$-central;
\item $R[f]$ is central in the sense of Smith;
\item (when $\Ac$ is pointed) $\ker f$ is central in the sense of Huq.
\end{enumerate}
\end{theorem}
\begin{proof} The equivalence of 1.\ and 3.\ is a combination of Theorem~6.1 of Gran~\cite{Gran-Alg-Cent} (in the context of factor permutable categories) and Proposition 3.6 of Pedicchio~\cite{Pedicchio} (in the context of exact Mal'tsev categories with coequalizers). 1.\ and 2.\ are equivalent by Proposition~\ref{Proposition-Centrality} and 3.\ and 4.\ by Proposition~\ref{Proposition-Huq-Centrality-is-Smith-Centrality}.
\end{proof}

\begin{remark}\label{Remark-Centrext-Birkhoff}
It is well-known that the category $V\CentrExt (Y)$\index{VCentr@$V\CentrExt (Y)$} of $V$-central extensions of a given object $Y$ is Birkhoff in $\Ext (Y)$. Using the commutator $V_{1}$ we can easily construct the reflection in the following way. Given an extension $f\colon {X\to Y}$ with kernel $K$, consider the diagram of solid arrows
\[
\xymatrix{& V_{1}f \ar@{=}[r] \ar@{{ >}->}[d] & V_{1}f \ar@{{ >}->}[d] \\
0 \ar[r] & K \ar@{}[rd]|{\texttt{(i)}} \ar@{{ >}->}[r]^-{k} \ar@{-{ >>}}[d] & X \ar@{-{ >>}}[d] \ar@{-{ >>}}[r]^-{f} & Y \ar[r] \ar@{:}[d] & 0\\
0 \ar@{.>}[r] & {\tfrac{K}{V_{1}f}} \ar@{{ >}.>}[r]_-{\overline{k}} & {\tfrac{X}{V_{1}f}} \ar@{.{ >>}}[r]_{\overline{f}} & Y \ar@{.>}[r] & 0.}
\]
Proposition~\ref{Proposition-pushout} implies that \texttt{(i)} is a pushout; by Proposition~\ref{Proposition-LeftRightPullbacks}, it is also a pullback. Now $\overline{k}$ is a monomorphism, and it is moreover a kernel by Proposition~\ref{Proposition-Image-of-Kernel}. Taking its cokernel yields the rest of the diagram. Using that $V_{1}$ preserves regular epimorphisms (Proposition~\ref{V_1 bewaart regulier epis}), one easily checks that $\overline{f}$ is the needed reflection of $f$ into $V\CentrExt (Y)$.

This kind of argument will be crucial to Chapter~\ref{Chapter-Cohomology}. Conversely, thus $V_{1}$ may be constructed as a Birkhoff subfunctor. 
\end{remark}

Theorem~\ref{Theorem-Abelian-Centrality} characterizes when, in the case of abelianization, the commutator $V_{1}$ is zero. We shall sharpen this result with Theorem~\ref{theorem-commutator}, which states that the commutator $V_{1}$ may be expressed entirely in terms of Smith's. To do so, we need the following technical property.

\begin{proposition}\label{Proposition-V_{1}-Commutator}
Let $\Ab\Ac$ be the Birkhoff subcategory of abelian objects of an exact and sequentiable category $\Ac$, $V$ the associated Birkhoff subfunctor and $p\colon A_0\to A\in\prA$. For a subobject $F\subseteq V (A_{0})$ such that $F\triangleleft A_{0}$, the following are equivalent:
\begin{enumerate}
\item
$F\in \Fc_{V_{0}}^p$;
\item
$V_1\bigl(\eta_{\eta_{F}R[p]}\colon \tfrac{A_0}{F}\to \tfrac{A_0}{F}/\eta _F(R[p])\bigr)=0$;
\item
$\eta _F(R[p])$ and $\eta _F(\nabla_{A_0})=\nabla_{{A_0}/{F}}$ admit a centralizing relation.
\end{enumerate}
\end{proposition}
\begin{proof}
The equivalence of 2.\ and 3.\ follows from Proposition~\ref{Proposition-Centrality} and Theorem~\ref{Theorem-Abelian-Centrality}. We will prove the equivalence of 1.\ and 2.

Consider the image $\eta_{F}R[p]$ of $R[p]$ along $\eta_{F}$.
\[
\xymatrix{
&& R[p] \ar@<-0.5 ex>[d]_-{k_0} \ar@<0.5 ex>[d]^-{k_1} \ar@{-{ >>}}[r]^-{\pi } & \eta_{F}(R[p]) \ar@<-0.5 ex>[d]_-{\kappa _0} \ar@<0.5 ex>[d]^-{\kappa _1}\\
0 \ar[r] & F \ar@{{ >}->}[r] & A_{0} \ar@{-{ >>}}[r]_-{\eta_{F}} & \tfrac{A_{0}}{F} \ar[r] & 1}
\]
Applying $V$ yields the following commutative diagram of $\Ac$, where the isomorphism exists due to the fact---see Remark~\ref{Remark-Birkhoff-Nabla}---that $\nabla_V$ and $V$ represent the same subfunctor.
\[
\xymatrix{ & R[V\eta _{F}] \ar@<-0.5 ex>[d] \ar@<0.5 ex>[d] & \\
V (R[p]) \ar@<-0.5 ex>[r]_-{Vk_{1}} \ar@<0.5 ex>[r]^-{V k_{0}} \ar@{-{ >>}}[d]_-{V \pi} & V (A_{0})\ar@{-{ >>}}[d]_-{V \eta _{F}} \ar@{-{ >>}}[rd]^-{\eta'_F} & \\
V \eta _{F} (R[p]) \ar@<-0.5 ex>[r]_-{V\kappa _{1}} \ar@<0.5 ex>[r]^-{V\kappa _{0}} & V (\frac{A_{0}}{F}) \ar[r]_-{\cong} & \frac{V(A_0)}{F}}
\]
$V\pi$ is a regular epimorphism because $V$ is Birkhoff; it follows that $V(R[p])\subseteq R[V\eta _F]=R[\eta'_F]$ precisely when $V\eta _F(R[p])=\Delta_{V ({A_0}/{F})}$.

Now, on one hand, by Proposition~\ref{Proposition-(ii')}, $F\in \Fc_{V_{0}}^{p}$ if and only if
\[
V (R[p])\subseteq R\bigl[\eta'_{F}\colon V (A_{0})\to \tfrac{V (A_{0})}{F}\bigr];
\]
on the other hand, Proposition~\ref{Proposition-Centrality} implies that 
\[
V\eta _F(R[p])=V (R[\eta_{\eta _F(R[p])}]) =\Delta_{V ({A_0}/{F})}
\]
if and only if $V_1\eta_{\eta_{F}R[p]}=0$. This shows that 1.\ and 2.\ are equivalent.
\end{proof}

\begin{theorem} \label{theorem-commutator}
Consider an exact and sequentiable category $\Ac $, its Birkhoff subcategory $\Ab\Ac$ of abelian objects, the associated Birkhoff subfunctor $V$, and the resulting functor $V_{1}$. For any presentation $p\colon A_{0}\to  A$ in $\Ac$, the following equality holds:
\[
V_{1}p=K_{[R[p],\nabla_{A_{0}}]}.
\]
\end{theorem}
\begin{proof}
By the order isomorphism mentioned in Notations~\ref{Notation-Eqrel-Image} and by Definition~\ref{Definition-Pedicchio-Commutator}, $N=K_{[R[p],\nabla_{A_0}]}$ is the smallest proper subobject of $A_{0}$ such that $\eta _{N}(R[p])$ and $\eta _{N}(\nabla_{A_0})=\nabla_{{A_0}/{N}}$ admit a centralizing relation. Thus, by Proposition~\ref{Proposition-V_{1}-Commutator}, it is the smallest element in $\Fc_{V_{0}}^{p}$; hence, it is $V_{1}p$.
\end{proof}

\section{$V_{1}$ and nilpotency}\label{Section-Nilpotency}
Interpreting $V_{1}$ as a commutator, we now focus on the resulting notion of nilpotency. We show that, as in the case of abelianization of groups, an object is nilpotent if and only if its lower central series reaches $0$. The nilpotent objects of class $n$ form a Birkhoff subcategory.

This section is based mainly on Huq~\cite{Huq}, where nilpotency is studied in a context that is essentially equivalent to that of semi-abelian categories (he uses old-style axioms~\cite{Janelidze-Marki-Tholen}, and supposes his categories to be well-powered and (co)complete). In their book~\cite{Freese-McKenzie}, Freese and McKenzie consider it in the context of \emph{congruence modular}\index{congruence modular variety}\index{variety!congruence modular} varieties.

 From now on, $\Ac$ will be an exact and homological category and $V$ a Birkhoff subfunctor of~$\Ac$.

\begin{definition}
Let $A$ be an object of $\Ac$. A \emph{$V$-central series}\index{central series}\index{V-central series@$V$-central series} of $A$ is a descending sequence
\[
A=A_0\supseteq A_1\supseteq \cdots \supseteq A_n\supseteq\cdots
\]
of proper subobjects of $A$, such that, for all $n\in\N $, the arrow ${{A}/{A_{n+1}}\to{A}/{A_{n}}}$ is $V$-central.

$A$ is called \emph{$V$-nilpotent}\index{nilpotent}\index{V-nilpotent@$V$-nilpotent} when there exists in $\Ac$ a $V$-central series of $A$ that reaches zero, i.e., such that $A_{n}=0$ for some $n\in \N$. In that case, $A$ is said to be $V$-nilpotent \emph{of class $n$}\index{nilpotent!of class $n$}.

The \emph{$V$-lower central series}\index{lower central series}\index{central series!lower} of $A$ is the descending sequence
\begin{equation}\label{Lower-Central-Series}
A=V_1^0 (A)\supseteq V_1^1 (A) \supseteq \cdots\supseteq V_{1}^{n} (A)\supseteq\cdots
\end{equation}
defined, for $n\in \N$, by putting $V_1^0 (A)=A$ and 
\[
V_1^{n+1} (A) =V_1(\eta_{V_{n}^{1} (A)} \colon  A\to\tfrac{A}{V^{n}_{1} (A)}).
\]
\end{definition}

\begin{remark}\label{Remark-Nilpotency}
Note that, by Lemma~\ref{Lemma-L_1}, $V_{1}^{i} (A)$ is a proper subobject of $A$, for all $i\geq 1$. Now we also need $V_{1}^{0}A$ to be proper in $A$, and this is true for all objects $A$ if and only if $\Ac$ is pointed. This is the reason why, in this section, we suppose that $\Ac$ is homological instead of sequentiable.
\end{remark}

\begin{example}
In case $\Ac=\Gp$ and $V$ the kernel of the abelianization functor---see Example~\ref{Examples-V-Exact-Sequence}---one has $V_{1} (\eta_{B}\colon A\to {A}/{B})=[B,A]$. Thus the sequence~\ref{Lower-Central-Series} becomes
\[
A \supseteq [A,A]\supseteq [[A,A],A]\supseteq [[[A,A],A],A]\supseteq\cdots,
\]
whence the name ``lower central series''. For Lie algebras, the situation is exactly the same.  
\end{example}

\begin{proposition}\label{Proposition-Lower-C-S}
Let $A\supseteq B\supseteq C$ be objects in $\Ac $ with $C\triangleleft A$ and $B\triangleleft A$. The following are equivalent:
\begin{enumerate}
\item $V_{1}\bigl (p\colon \frac{A}{C}\to \frac{A}{B})=0$;
\item $C\supseteq V_{1}\bigl (\eta _{B}\colon A\to \frac{A}{B}\bigr)$.
\end{enumerate}
\end{proposition}
\begin{proof}
Consider the following diagram.
\[
\xymatrix{R[\eta _{B}] \ar@{.>}[d]_-{\sigma}\ar@<-0.5 ex>[r] \ar@<0.5 ex>[r] &A \ar@{-{ >>}}[r]^-{\eta _{B}} \ar@{-{ >>}}[d]_-{\eta _{C}} & \frac{A}{B} \\
R[p] \ar@<-0.5 ex>[r] \ar@<0.5 ex>[r] & \frac{A}{C} \ar@{-{ >>}}[ru]_-{p} }
\]
By Proposition~\ref{Rotlemma-kernel-pairs}, the factorization $\sigma$ is a regular epimorphism. Now, applying $V$, we get the following diagram; the isomorphism exists due to the fact---see Remark~\ref{Remark-Birkhoff-Nabla}---that $\nabla_V$ and $V$ represent the same subfunctor.
\[
\xymatrix{ & R[V \eta _{C}] \ar@<-0.5 ex>[d] \ar@<0.5 ex>[d] & \\
V (R[\eta _{B}]) \ar@<-0.5 ex>[r] \ar@<0.5 ex>[r] \ar@{-{ >>}}[d]_-{V \sigma} & V (A)\ar@{-{ >>}}[d]_-{V \eta _{C}} \ar@{-{ >>}}[rd]^-{\eta'_{V (A)\cap C}} & \\
V (R[p]) \ar@<-0.5 ex>[r] \ar@<0.5 ex>[r] & V (\frac{A}{C}) \ar[r]_-{\cong} & \frac{V(A)}{V (A)\cap C}}
\]
$V$ being Birkhoff, $V \sigma$ is an epimorphism, hence $V(R[\eta _{B}])\subseteq R[V\eta _C]=R[\eta'_{V (A)\cap  C}]$ precisely when $V(R[p])=\Delta_{V ({A}/{C})}$.

Using Proposition~\ref{Proposition-(ii')} and Proposition~\ref{Proposition-Centrality}, we get that 1.\ holds if and only if $V (A)\cap C\in \Fc_{V_{0}}^{\eta_{B}}$. But, by Corollary~\ref{Corollary-L_1}, $V (A)\cap C$ is an element of $\Fc_{V_{0}}^{\eta_{B}}$ if and only if $V_{1}\eta_{B}\subseteq V (A)\cap C$. This last statement is equivalent to condition 2.
\end{proof}

\begin{corollary}\label{Corollary-Lower-Central=Central}
For any object $A\in\Ac$, the $V$-lower central series in $A$
\[
V_1^0 (A)\supseteq V_1^1 (A)\supseteq V_1^2 (A) \supseteq \cdots
\]
is a $V$-central series, i.e.,
\[
V_{1} \bigl(\tfrac{A}{V_{1}^{n+1} (A)}\to \tfrac{A}{V_{1}^{n} (A)}\bigr)=0
\]
for all $n\in\N$.
\end{corollary}
\begin{proof}
Take $A=A$, $B=V_{1}^{n} (A)$ and $C=V_{1}^{n+1} (A)$ in Proposition~\ref{Proposition-Lower-C-S}; then 2.\ is trivially fulfilled.
\end{proof}

\begin{corollary}\label{Corollary-Nilpotent-iff-LCC-zero}
An object $A$ of $\Ac$ is $V$-nilpotent of class $n$ if and only if $V_1^n(A)=0$.
\end{corollary}
\begin{proof}
By Corollary~\ref{Corollary-Lower-Central=Central}, one implication is obvious. For the other, suppose that there exists a descending sequence
\[
A=A_0\supseteq A_1\supseteq \cdots \supseteq A_n=0
\]
with all objects proper in $A$, and that for all $i\in\{0,\dots,n-1\}$, 
\[
V_1\bigl(\tfrac{A}{A_{i+1}}\to\tfrac{A}{A_{i}}\bigr)=0.
\]
We will show by induction that, for all $i\in\{0,\dots,n\}$, $V_1^i(A)\subseteq A_i$. If so, $V_1^n(A)\subseteq A_n=0$, which proves our claim.

The case $i=0$ is clear. Now suppose $V_{1}^{i} (A)\subseteq A_{i}$ for a certain $i\in\{0,\dots,n-1\}$. Recall from Remark~\ref{Remark-L-in-L_1} that there is an inclusion of functors $V_{1}\subseteq K[\cdot]\cap V_{0} \subseteq K[\cdot]$. Hence we get a commutative diagram of inclusions (all proper, since the four objects are proper subobjects of $A$):
\[
\xymatrix{
V_{1}\bigl(A\to \frac{A}{V_{1}^{i} (A)}\bigr) \ar@{{ >}->}[r]  \ar@{{ >}->}[d] & V_{1}^{i} (A) \ar@{{ >}->}[d]\\
V_{1}\bigl(A\to \frac{A}{A_{i}}\bigr)  \ar@{{ >}->}[r] & A_{i}. }
\]
In particular, we have $V_{1}(A\to {A}/{V_{1}^{i} (A)})\subseteq V_{1}(A\to {A}/{A_{i}})$ and thus
\[
V_{1}^{i+1} (A)=V_{1}\bigl(A\to  \tfrac{A}{V_{1}^{i} (A)}\bigr) \subseteq  V_{1}\bigl (A\to \tfrac{A}{A_{i}}\bigr) \subseteq A_{i+1},
\]
where the last inclusion follows from Proposition~\ref{Proposition-Lower-C-S}, as 
\[
{V_{1}({A}/{A_{i+1}}\to {A}/{A_{i}})}=0.\qedhere
\]
\end{proof}

\begin{remark}
Since, by Remark~\ref{Remark-L-uit-L_1}, $V_{1}^{1}=V$, a $1$-nilpotent object in an exact homological category is nothing but an object in the Birkhoff subcategory associated with $V$.
\end{remark}

The following was inspired by Section 4.3 in Huq~\cite{Huq}.

\begin{proposition}\label{Proposition-Nilpotency-Class}
For $n\in \N$, $V_{1}^{n}\colon \Ac \to \Ac$ is a Birkhoff subfunctor of $\Ac$. The corresponding Birkhoff subcategory is the full subcategory of all objects of $V$-nilpotency class $n$.
\end{proposition}
\begin{proof}
Let $p\colon A\to B$ be a regular epimorphism. The first statement is clear in case $n=0$, so suppose that $V_{1}^{n}$ preserves regular epimorphisms. Then, by Proposition~\ref{Rotlemma}, the induction hypothesis implies that the right hand square in the diagram with exact rows
\[
\xymatrix{0 \ar[r] & V_{1}^{n} (A) \ar@{{ >}->}[r] \ar@{-{ >>}}[d]_-{V_{1}^{n}p} & A \ar@{-{ >>}}[r] \ar@{-{ >>}}[d]^-p & \tfrac{A}{V_{1}^{n} (A)} \ar@{-{ >>}}[d] \ar[r] & 1\\
0 \ar[r] & V_{1}^{n} (B) \ar@{{ >}->}[r] & B \ar@{-{ >>}}[r] & \frac{B}{V_{1}^{n} (B)} \ar[r] & 1}
\]
is a pushout. Hence, the category $\Ac$ being exact Mal'tsev, it is a regular pushout; by Proposition~\ref{Proposition-Regular-Pushout} we get that it is a regular epimorphism $\vetsf$ of $\prA$. We conclude with Proposition~\ref{V_1 bewaart regulier epis} that $V_{1}^{n}p=V_{1}\vetsf$ is regular epi.

The second statement immediately follows from Proposition~\ref{Proposition-Birkhoff-Subfunctor}.
\end{proof}
 % Baer invariants
\chapter{Cotriple homology \protect\versus\ Baer invariants}\label{Chapter-Cotriples}

\setcounter{section}{-1}
\section{Introduction}\label{Section-Cotriples-Introduction}
The main result of Chapter~\ref{Chapter-Baer-Invariants} is Theorem~\ref{Theorem-V-Exact-Sequence}:

\bigskip
\begin{thm}
Let $\Ac$ be an exact and sequentiable category with enough projectives. Let $V$ be a Birkhoff subfunctor on $\Ac$. Then every short exact sequence
\[
\xymatrix{0 \ar[r] & K \ar@{{ >}->}[r] & A \ar@{-{ >>}}[r]^-{f} & B \ar[r] & 1}
\]
in $\Ac$ induces an exact sequence
\begin{equation}\label{Sequence-langexactV-Local}
\resizebox{\textwidth}{!}{\xymatrix{0 \ar[r] & K[\gamma] \ar@{{ >}->}[r] & \Delta V(A)  \ar[r]^-{\Delta Vf} & \Delta V(B)  \ar[r] & \frac{K}{V_{1}f}  \ar[r] & U(A) \ar@{-{ >>}}[r]^-{Uf} & U(B)  \ar[r] & 1,}}
\end{equation}
which depends naturally on the given short exact sequence.\noproof
\end{thm}
\bigskip

The aim of this chapter is to interpret Sequence~\ref{Sequence-langexactV-Local} as a generalization of the Stallings-Stammbach Sequence from integral homology of groups~\cite{Stallings, Stammbach}. We do so by proving a semi-abelian version of Hopf's Formula~\cite{Hopf}:
\begin{equation}\label{Our-Hopf-Formula}
H_{2} (X,U)_{\G}\cong \Delta V (X).
\end{equation}
Thus we relate cotriple homology with the theory of Baer invariants. Alternatively, one may take the opposite point of view, and ask whether in the semi-abelian context, cotriple homology makes sense. Then this chapter gives a partial positive answer, by interpreting (one specific flavour of) cotriple homology in terms of Baer invariants.
\pfbreak

\noindent Let us start by explaining the right hand side of this formula, briefly recalling the main concepts from Chapter~\ref{Chapter-Baer-Invariants}. Let $\Ac$ be an exact and sequentiable category. A normal subfunctor $V$ of $1_{\Ac}$ is called \textit{Birkhoff subfunctor of $\Ac$} when $V$ preserves regular epimorphisms. The Birkhoff subfunctors of $\Ac$ correspond bijectively to its Birkhoff subcategories (Proposition~\ref{Proposition-Birkhoff-Subfunctor}). 

Any Birkhoff subfunctor $V$ of $\Ac$ induces a functor $V_{1}\colon {\Pr\Ac\to\Ac}$ from the category $\Pr \Ac$ of presentations in $\Ac$ to $\Ac$, constructed as follows. Let $f\colon {X\to Y}$ be a presentation in $\Ac $ and $(R[f],k_{0},k_{1})$ its kernel pair. Applying $V$, next taking a coequalizer of $Vk_{0}$ and $Vk_{1}$ and then a kernel of $\Coeq (Vk_{0},Vk_{1})$, one gets $V_{1}f$:
\[
\xymatrix{
 &&& V(\KP[f]) \ar@<-0.5 ex>[d]_{Vk_{0}} \ar@<0.5 ex>[d]^{Vk_{1}} &\\
0 \ar[r] & V_1f \ar@{{ >}->}[rr] && V(X) \ar@{-{ >>}}[rr]_-{\coeq (Vk_{0},Vk_{1})} && \Coeq[Vk_{0},Vk_{1}] \ar[r] & 1.}
\]
Now given a Birkhoff subfunctor $V$ of $\Ac$ and a projective presentation $f\colon {X\to Y}$, 
\[
\Delta V\colon {\Ac \to \Ac}
\]
is the functor induced by the Baer invariant $\prA \to \Ac$ that maps $f$ to
\[
\frac{K[f]\cap V (X)}{V_{1}f}.
\]

The meaning of this functor will probably be more clear if we consider an example: for instance, the category $\Gp$ of groups and its Birkhoff subcategory $\Ab =\Ab\Gp$ of abelian groups (cf.\ Examples~\ref{Examples-Huq-Central}). The associated Birkhoff subfunctor $V$ sends $G$ to $[G,G]$. If $R\triangleleft F$ and $f$ denotes the quotient ${F\to F/R}$, then $V_{1}f=[R,F]$. Now let
\[
\xymatrix{0 \ar[r] & R \ar@{{ >}->}[r] & F \ar@{-{ >>}}[r]^-{f} & G \ar[r] & 0}
\]
be a presentation of a group $G$ by a ``group of generators'' $F$ and a ``group of relations'' $R$, i.e., a short exact sequence with $F$ a free (or, equivalently, projective) group. Then it follows that
\[
\frac{R\cap [F,F]}{[R,F]}=\Delta V (G).
\]
\vspace{-\baselineskip}
\pfbreak

\noindent Hopf's formula~\cite{Hopf} is the isomorphism
\begin{equation}\label{Hopf}
H_2(G,\Z)\cong\frac{R\cap [F,F]}{[R,F]};
\end{equation}
here, the \emph{Schur multiplicator}\index{Schur multiplicator} $H_2(G,\Z)$ is the second integral homology group of $G$. It is clear that the left hand side of the formula~\ref{Our-Hopf-Formula} should, in a way, generalize this homology group.

That such a generalization is possible comes as no surprise, since Carrasco, Cegarra and Grandje\'an (in their article~\cite{Carrasco-Homology}) already prove a generalized Hopf formula in the category $\XMod$ of crossed modules. They also obtain a five term exact sequence that generalizes the Stallings-Stammbach Sequence~\cite{Stallings, Stammbach}. Their five term exact sequence becomes a particular case of ours (Theorem~\ref{Theorem-Cotriple-Exact-Sequence}).

Since the right hand side of Formula~\ref{Our-Hopf-Formula} is defined relative to a Birkhoff subcategory $\Bc$ of an exact and sequentiable category, so must be the left-hand side. We restrict ourselves to the following situation: $\Ac$ is exact and sequentiable with enough projectives; $U\colon \Ac \to \Bc$ is a reflector onto a Birkhoff subcategory $\Bc$ of $\Ac$; $\G= (G,\delta,\epsilon)$ is a regular comonad on $\Ac$ (i.e., any value $\epsilon_{X}\colon {GX\to X}$ of its counit is a projective presentation, and its functor $G$ preserves regular epis). Under these assumptions, the formula~\ref{Our-Hopf-Formula} holds---given the appropriate notion of homology.

Carrasco, Cegarra and Grandje\'an define their homology of crossed modules by deriving the functor $\ab\colon \XMod\to\Ab\XMod$, which sends a crossed module $(T,G,\del)$ to its abelianization $\ab (T,G,\del)$, an object of the abelian category $\Ab\XMod$. More precisely, the monadicity of the forgetful functor $\Uc \colon {\XMod\to\Set}$ yields a comonad $\G$ on $\XMod$. This, for any crossed module $(T,G,\del)$, gives a canonical simplicial object $\G(T,G,\del)$ in $\XMod$. The $n$-th homology object (an abelian crossed module) of a crossed module $(T,G,\del)$ is then $H_{n-1}C\ab\G(T,G,\del)$, the $(n-1)$-th homology object of the unnormalized chain complex associated with the simplicial object $\ab\G(T,G,\del)$. This is an application of Barr and Beck's cotriple homology theory~\cite{Barr-Beck}, which gives a way of deriving any functor $U\colon {\Ac \to \Bc}$, from an arbitrary category $\Ac$ equipped with a comonad $\G$, to an abelian category $\Bc$.

To prove our Hopf formula, we use methods similar to theirs---in fact, our proof of Theorem~\ref{Theorem-Hopf-Formula} is a modification of~\cite[Theorem 12]{Carrasco-Homology}. We use the ``semi-abelian'' notion of homology introduced in Chapter \ref{Chapter-Homology} to generalize Barr and Beck's notion of cotriple homology to the situation where $\Bc$ is an exact and sequentiable category. Their definition is modified to the following. Let $\Ac$ be a category, $\G$ a comonad on $\Ac$ and $U\colon {\Ac \to \Bc}$ a functor. The \textit{$n$-th homology object $H_{n} (X,U)_{\G}$ of $X$ with coefficients in $U$ relative to the cotriple $\G$} is the object $H_{n-1}NU (\G X)$, the $(n-1)$-th homology object of the Moore complex of the simplicial object $U (\G X)$ of $\Bc$. 

This is the subject of Section~\ref{Section-Cotriple-Homology}, where we also consider some examples and show that $H_{1} (X,U)_{\G}=U (X)$. In Section~\ref{Section-Hopf} we prove Formula~\ref{Our-Hopf-Formula} and obtain a general version of the Stallings-Stammbach Sequence.

This chapter is based on the paper~\cite{EverVdL2}, modified to accommodate the quasi-pointed and non-monadic cases.

\section{Cotriple homology}\label{Section-Cotriple-Homology}

\begin{notation}\label{Notation-Comonad}
Let $\Ac$ be an arbitrary category. A comonad or cotriple $\G$\index{G@$\G$}\index{comonad}\index{cotriple} on $\Ac$ will be denoted by
\[
\G =(G\colon \Ac \to \Ac ,\quad \delta\colon G\To G^{2} ,\quad \epsilon\colon G\To 1_{\Ac}).
\]
Recall that the axioms of comonad state that $\epsilon_{GX}\comp \delta_{X}=G\epsilon_{X}\comp \delta_{X}=1_{GX}$ and $\delta_{GX}\comp \delta_{X}=G\delta_{X}\comp \delta_{X}$, for any object $X$ of $\Ac$. Putting
\[
\del_{i}=G^{i}\epsilon_{G^{n-i}X}\colon G^{n+1}X\to G^nX
\]
and
\[
\sigma_{i}=G^{i}\delta_{G^{n-i}X}\colon G^{n+1}X\to G^{n+2}X,
\]
for $0\leq i \leq n$, makes the sequence $(G^{n+1}X)_{n\in \N}$ a simplicial object $\G X$ of $\Ac$. This induces a functor $\Ac \to \simpA$, which, when confusion is unlikely, will be denoted by $\G$.
\end{notation}

The following extends Barr and Beck's definition of cotriple homology~\cite{Barr-Beck} to the semi-abelian context.

\begin{definition}\label{Definition-Cotriple-Homology}
Let $\G$ be a comonad on a category $\Ac$. Let $\Bc$ be an exact and sequentiable category and $U\colon \Ac \to \Bc$ a functor. We say that the object
\[
H_{n} (X,U)_{\G}=H_{n-1} NU (\G X)
\]\index{H_n@$H_{n} (\cdot,U)_{\G}$}
is the \emph{$n$-th homology object of $X$ with coefficients in $U$ relative to the comonad $\G$}\index{homology!cotriple ---}\index{cotriple homology}\index{homology!w.r.t.\ a comonad}. This defines a functor $H_{n} (\cdot,U)_{\G}\colon \Ac \to \Bc $, for any $n\in \N_{0}$.
\end{definition}

We shall be interested mainly in the situation where $U$ is the reflector of an exact and sequentiable category $\Ac$ onto a Birkhoff subcategory $\Bc$ (which is itself exact and sequentiable by Corollary~\ref{Corollary-Birkhoff-Subcategory-is-Semiabelian}). Our aim is then to characterize the lowest homology objects in terms of Baer invariants. To do so, we would like to use Theorem~\ref{Theorem-V-Exact-Sequence}.

For this to work, $\G$ must induce projective presentations $\epsilon_{X} \colon {GX\to X}$, and moreover, $G$ must preserve regular epimorphisms. In other words, the comonad $\G$ must be such that for every regular epimorphism $f:X\to Y$ in $\Ac$, all arrows in the naturality square
\begin{equation}\label{Diagram-Comonad-Regepi}
\vcenter{\xymatrix{G X \ar[r]^-{Gf} \ar[d]_-{\epsilon_{X}} & G Y \ar[d]^-{\epsilon_{Y}}\\
X \ar[r]_-{f} & Y}}
\end{equation}
are regular epimorphisms, and the objects $GX$ and $GY$ are projective. (Hence $Gf$ is a split epimorphism.) Such a diagram represents a regular epimorphism of projective presentations: Indeed, a regular epimorphism $f\colon X\to Y$ with kernel pair $R[f]$ gives rise to a diagram
\[
\xymatrix{G (R[f]) \ar@<.5ex>[r]^-{Gk_{1}} \ar@<-.5ex>[r]_-{Gk_{0}} \ar@{-{ >>}}[d]_-{\epsilon _{R[f]}} & G (X) \ar@{-{ >>}}[r]^-{Gf} \ar@{-{ >>}}[d]_-{\epsilon_{X}} & G (Y) \ar@{-{ >>}}[d]^-{\epsilon_{Y}} \\
R[f] \ar@<.5ex>[r]^-{k_{1}} \ar@<-.5ex>[r]_-{k_{0}} & X \ar@{-{ >>}}[r]_-{f} &Y.}
\]
A kernel pair of $(Gf,f)\colon \epsilon_{X}\to \epsilon_{Y}$ in $\Fun (\Two,\Ac)$ is a morphism 
\[
R [(Gf,f)]\colon {R[Gf]\to R[f]}
\]
in $\Ac$, and this morphism is a regular epimorphism, because $\epsilon_{R[f]}$ is. Proposition~\ref{Rotlemma-kernel-pairs} now implies that $(Gf,f)$ is regular epi in $\prA$ (the square~\ref{Diagram-Comonad-Regepi} is a regular pushout).

When, in particular, $f$ is $\epsilon_{X}\colon G X\to X$, we get that
\begin{equation}\label{Diagram-Comonad-Pushout}
\vcenter{\xymatrix{G^{2} X \ar@{-{ >>}}[r]^-{G\epsilon_{X}} \ar@{-{ >>}}[d]_-{\epsilon_{G X}} & G X \ar@{-{ >>}}[d]^-{\epsilon_{X}}\\
G X \ar@{-{ >>}}[r]_-{\epsilon_{X}} & X}}
\end{equation}
is a pushout, and hence Proposition~\ref{Proposition-Forks} implies that 
\begin{equation}\label{Diagram-H1}
\xymatrix{G^{2}X \ar@<.5ex>[r]^-{G\epsilon_{X}} \ar@<-.5ex>[r]_-{\epsilon_{GX}} & GX \ar[r]^-{\epsilon_{X}} & X}
\end{equation}
is a coequalizer diagram.

\begin{definition}\label{Definition-Regular-Comonad}
We call a comonad $\G$ on an exact and sequentiable category $\Ac$ \emph{regular}\index{regular comonad}\index{comonad!regular} when, for every regular epimorphism $f:X\to Y$ in $\Ac$, the naturality square~\ref{Diagram-Comonad-Regepi} is a regular epimorphism of projective presentations.
\end{definition}

\begin{assumption}\label{Assumption-Cotriples}
From now on, we suppose that
\begin{enumerate}
\item $U$ is the reflector of an exact and sequentiable category $\Ac$ onto a Birkhoff subcategory $\Bc$ of $\Ac$;
\item $\G$ a regular comonad on $\Ac$.
\end{enumerate}
\end{assumption}

\begin{remark}\label{Remark-Comonad-and-Delta}
For any Birkhoff subfunctor $V$ and any $X\in \Ac$, $\Delta V G (X)=0$, because $\epsilon_{GX}\colon G^{2}X\to GX$ is a split epimorphism.
\end{remark}

\begin{example}\label{Example-Functorial-Web-Monadic}
An important special case is the following: $\Ac$ is quasi-pointed, protomodular and monadic\index{category!monadic}\index{monadic category} over $\Set$, and $\G$ is the induced comonad on $\Ac$.

Let $\Upsilon \colon \Ac \to \Set$ and $\Phi \colon \Set \to \Ac$ denote the respective right and left adjoint functors and $\epsilon \colon \Phi \comp \Upsilon \To 1_{\Ac}$ and $\zeta \colon 1_{\Set}\To \Upsilon \comp \Phi $ the counit and unit. Then $G=\Phi \comp \Upsilon$, $\epsilon$ is just the counit and $\delta$ is the natural transformation defined by $\delta_{X}=\Phi \zeta_{\Upsilon (X)}$, for $X\in \Ob{\Ac}$.

Because the forgetful functor $\Upsilon \colon \Ac \to \Set$ preserves regular epimorphisms (see e.g., Borceux~\cite[Theorem II.4.3.5]{Borceux:Cats}) and because, in $\Set$, every object is projective, $\Ac$ is easily seen to have enough projectives. In particular, any $GX$ is projective, and for any object $X$ of $\Ac$, the map $\epsilon_{X}\colon GX\to X$ is a projective presentation, the \emph{$\G$-free}\index{presentation!$\G$-free} presentation of $X$. Moreover, $G$ turns regular epimorphisms into split epimorphisms; hence $\G$ is regular.

The requirement that $\Ac$ be monadic over $\Set$ also implies that $\Ac$ is complete, cocomplete and exact (again~\cite[Theorem II.4.3.5]{Borceux:Cats}); hence, if $\Ac$ is, moreover, quasi-pointed and protomodular, it is exact sequentiable. Reciprocally, any variety of algebras over $\Set$ is monadic---see e.g., Cohn~\cite{Cohn}, Borceux~\cite{Borceux:Cats} or Mac\,Lane~\cite{MacLane}---and thus sequentiable varieties (and in particular, semi-abelian varieties) form an example of the situation considered. A characterization of semi-abelian varieties of algebras over $\Set$ is given by Bourn and Janelidze in their paper~\cite{Bourn-Janelidze}. More generally, in~\cite{Gran-Rosicky:Monadic}, Gran and Rosick\'y characterize semi-abelian categories, monadic over $\Set$. In \cite{EverVdL2} and \cite{Gran-VdL} we work in this latter context.
\end{example}

The next characterization of $H_{1} (X,U)_{\G}$ is an immediate consequence of Proposition~\ref{Proposition-Forks}. In particular, it shows that $H_{1} (\cdot , U)_{\G}$ is independent of the chosen comonad $\G$.

\begin{proposition}\label{Proposition-First-Homology}
For any object $X$ of $\Ac$, $H_{1} (X,U)_{\G}\cong U (X)$.
\end{proposition}
\begin{proof}
On one hand, $H_{1} (X,U)_{\G}=H_{0}N U (\G X)$ is a cokernel of the map $UG\epsilon_{X}\comp \ker U\epsilon_{GX}$. On the other hand, Diagram~\ref{Diagram-H1} is a coequalizer diagram, and $U$ preserves coequalizers. Since, moreover, $U\delta_{X}\colon {UGX\to UG^{2}X}$ is a splitting for both $UG\epsilon_{X}$ and $U\epsilon_{GX}$, Proposition~\ref{Proposition-Forks} applies, and $U (X)$ is a cokernel of $UG\epsilon_{X}\comp \ker U\epsilon_{GX}$ as well.
\end{proof}

\section{Hopf's Formula and the Stallings-Stammbach Sequence}\label{Section-Hopf}

Keeping in mind that a functor category $\Fun (\Cc ,\Ac)$ has the limits and colimits of $\Ac$, computed pointwise, the following follows immediately from the definitions.

\begin{lemma}\label{Lemma-Simplicial-Category}
Let $\Cc$ be a small category, $\Ac$ an exact sequentiable category and $U\colon {\Ac\to \Bc}$ a reflector onto a Birkhoff subcategory $\Bc$ of $\Ac$. Then
\begin{enumerate}
\item the functor category $\Fun (\Cc ,\Ac)$ is exact and sequentiable;
\item $\Fun (\Cc ,\prA)=\pres \Fun (\Cc ,\Ac) $;
\item $\Fun (\Cc ,\Bc)$ is a Birkhoff subcategory of $\Fun (\Cc ,\Ac)$;
\item the functor $\Fun (\Cc ,U)=U\comp (\cdot)\colon \Fun (\Cc ,\Ac)\to \Fun (\Cc ,\Bc)$ is its reflector;
\item $V_{1}^{\Fun (\Cc ,\Bc)}= \Fun \bigl(\Cc ,V_{1}^{\Bc}\bigr)\colon \Fun (\Cc ,\prA )\to \Fun (\Cc ,\Ac)$, where $V^{\Fun (\Cc ,\Bc)}$ and $V^{\Bc}$ are the Birkhoff subfunctors associated with $\Fun (\Cc ,\Bc)$ and $\Bc$, respectively.
\end{enumerate}
These properties hold, in particular, when $\Cc$ is $\Delta^{\op }$, i.e., when all functor categories are categories of simplicial objects and simplicial maps.\noproof
\end{lemma}

\begin{proposition}\label{Proposition-General-H_0-Sequence}
Any short exact sequence
\[
\xymatrix{ 0 \ar[r] & K \ar@{{ >}->}[r] & X \ar@{-{ >>}}[r]^-{f} & Y \ar[r] & 1}
\]
in $\Ac$ induces the exact sequence 
\begin{equation}\label{Sequence-simpB}
0 \to \tfrac{K[\G f]}{V_{1}\G f} \to/{{ >}->}/^{\ker U\G f} U\G X \to/{-{ >>}}/^{U\G f} U\G Y \to 1
\end{equation}
 in $\Sc \Bc$. 
\end{proposition}
\begin{proof}
Since $G$ turns regular epimorphisms into split epimorphisms (because $\G$ is a regular comonad) the simplicial morphism $U\G f$ is degreewise split epimorphic in $\Sc \Bc$.

By Remark~\ref{Remark-Comonad-and-Delta}, for any $n\geq 1$, the short exact sequence 
\[
0 \to K[G^{n} f] \to/{{ >}->}/^{\ker G^{n} f} G^{n} X \to/{-{ >>}}/^{G^{n} f} G^{n} Y \to 1,
\]
through Theorem~\ref{Theorem-V-Exact-Sequence}, induces the exact sequence
\[
0 \to \tfrac{K[G^{n} f]}{V_{1}G^{n} f} \to/{{ >}->}/^{\ker UG^{n} f} UG^{n} X \to/{-{ >>}}/^{UG^{n} f} UG^{n} Y \to 1
\]
As a consequence, 
\[
K[U\G f]\cong \tfrac{K[\G f]}{V_{1}\G f}.\qedhere
\]
\end{proof}

\begin{proposition}\label{Proposition-pre-H_0-Sequence}
Let
\[
\xymatrix{ 0 \ar[r] & K \ar@{{ >}->}[r] & X \ar@{-{ >>}}[r]^-{f} & Y \ar[r] & 1}
\]
be a short exact sequence in $\Ac$. Then
\[
H_{0}\tfrac{K[\G f]}{V_{1}\G f}\cong\tfrac{K}{V_{1}f}.
\]
\end{proposition}
\begin{proof}
Consider the following diagram, where the $\kappa_i$ are the face operators of $K[\G f]$ and $\kappa$ is the kernel object $K[\vetsf]$ of the arrow $\vetsf=(Gf,f)\colon \epsilon_{X}\to \epsilon_{Y}$ represented by Diagram~\ref{Diagram-Comonad-Regepi}.
\begin{equation}\label{Diagram-Forks}
\vcenter{\xymatrix{0 \ar[r] & V_{1} G^2 f \ar@{{ >}->}[rr] \ar@<.5ex>[d]^-{V_{1} (G\epsilon_{X},G\epsilon_{Y})} \ar@<-.5ex>[d]_-{V_{1} (\epsilon_{GX},\epsilon_{GY})} && K[G^2 f] \ar@{-{ >>}}[rr] \ar@<.5ex>[d]^-{\kappa_{1}} \ar@<-.5ex>[d]_-{\kappa_{0}} && \frac{K[G^2 f]}{V_{1} G^2 f} \ar[r] \ar@<.5ex>[d]^-{\frac{\kappa_{1}}{V_{1} (G\epsilon_{X},G\epsilon_{Y})}} \ar@<-.5ex>[d]_-{\frac{\kappa_{0}}{V_{1} (\epsilon_{GX},\epsilon_{GY})}} & 1\\
0 \ar[r] & V_{1} G f \ar@{{ >}->}[rr] \ar@{-{ >>}}[d]_{V_{1}\vetsf} && K[G f] \ar@{-{ >>}}[rr]^-{p} \ar@{}[rrd]|{\texttt{(i)}} \ar[d]_{\kappa } && \frac{K[G f]}{V_{1} G f}\ar[r] \ar[d]^{\frac{\kappa }{V_{1} \vetsf}} & 1\\
0 \ar[r] & {V_{1} f} \ar@{{ >}->}[rr]  && K \ar@{-{ >>}}[rr]_-{q} && \frac{K}{V_{1} f} \ar[r] & 1}}
\end{equation}
We have to show that the right hand side fork is a coequalizer diagram. Now, since $\G$ is regular, $\vetsf$ is a regular epimorphism in $\prA$; hence $V_{1}\vetsf$ is regular epi (Proposition~\ref{V_1 bewaart regulier epis}), and the square \texttt{(i)} is a pushout (Proposition~\ref{Proposition-pushout}). It follows that when the middle fork is a coequalizer diagram, so is the right hand side fork. This then yields the needed isomorphism.

To see that the middle fork in Diagram~\ref{Diagram-Forks} is indeed a coequalizer diagram, note that it is defined by the exactness of the rows in the diagram
\begin{equation}\label{Diagram-Kappa}
\vcenter{\xymatrix{0 \ar[r] & K[G^{2}f] \ar@{{ >}->}[r] \ar@{.>}@<.5ex>[d]^-{\kappa_{1}} \ar@{.>}@<-.5ex>[d]_-{\kappa_{0}} & G^2 X \ar@{-{ >>}}[r]^{G^{2}f} \ar@<-.5ex>[d]_-{\epsilon_{GX}} \ar@<.5ex>[d]^-{G\epsilon_{X}} & G^{2} Y \ar[r] \ar@<-.5ex>[d]_-{\epsilon_{GY}} \ar@<.5ex>[d]^-{G\epsilon_{Y}} & 1\\
0 \ar[r] & K[G f] \ar@{{ >}->}[r] \ar@{.>}[d]_{\kappa } & GX \ar@{-{ >>}}[r]^-{G f} \ar@{-{ >>}}[d]_-{\epsilon_{X}} & G Y \ar[r] \ar@{-{ >>}}[d]^-{\epsilon_{Y}} & 1\\
0 \ar[r] & K  \ar@{{ >}->}[r]  & X \ar@{-{ >>}}[r]_-{f} & Y \ar[r] & 1.}}
\end{equation}
The square~\ref{Diagram-Comonad-Regepi} being a pushout, certainly $\kappa $ is a regular epimorphism; hence we only need to show that the factorization $(\kappa_{0},\kappa_{1})\colon {K[G^{2}f]\to \KP [\kappa]}$ over the kernel pair of $\kappa$ is regular epi as well. Now the diagram
\[
\xymatrix{0 \ar[r] & K[G^{2}f] \ar@{{ >}->}[r] \ar@{.>}[d]_-{(\kappa_{0},\kappa_{1})} & G^2 X \ar@{-{ >>}}[r]^{G^{2}f} \ar@{.{ >>}}[d]|-{(\epsilon_{G X},G\epsilon_{X})} & G^{2}Y \ar[r] \ar@{.{ >>}}[d]^-{(\epsilon_{GY},G\epsilon_{Y})} & 1\\
0 \ar[r] & \KP [\kappa ] \ar@{{ >}->}[r] & \KP [\epsilon_{X}] \ar@{-{ >>}}[r] & \KP [\epsilon_{Y}] \ar[r] & 1}
\]
has exact rows, and its right hand side square is a pushout: Indeed, by Proposition~\ref{Proposition-N-Exact}, the morphism $N_{1}\G f$ in
\[
\xymatrix{0 \ar[r] & N_{1}\G X \ar@{{ >}->}[r] \ar@{.>}[d]_-{N_{1}\G f} & G^2 X \ar@{-{ >>}}[rr]^{(\epsilon_{G X},G\epsilon_{X})} \ar@{-{ >>}}[d]_-{G^{2}f} && \KP [\epsilon_{X}] \ar[r] \ar@{-{ >>}}[d] & 1\\
0 \ar[r] & N_{1}\G Y \ar@{{ >}->}[r] & GY \ar@{-{ >>}}[rr]_{(\epsilon_{GY},G\epsilon_{Y})} && \KP [\epsilon_{Y}] \ar[r] & 1}
\]
is a regular epimorphism. We get that $(\kappa_{0},\kappa_{1})$ is regular epi and the result follows.
\end{proof}

\begin{theorem}[Hopf's Formula]\label{Theorem-Hopf-Formula}\index{Hopf's Formula}
For any object $X$ of $\Ac$, 
\[
H_{2} (X,U)_{\G}\cong \Delta V (X).
\]
\end{theorem}
\begin{proof}
On one hand, by Remark~\ref{Remark-Comonad-and-Delta}, the tail of the exact sequence~\ref{Sequence-langexactV-Local}, induced by the short exact sequence
\[
0 \to K[\epsilon_{X}] \to/{{ >}->}/^{\ker \epsilon_{X}} GX \to/{-{ >>}}/^{\epsilon_{X}} X \to 1,
\]
becomes
\[
\xymatrix{0 \ar[r] & {\Delta V (X)} \ar@{{ >}->}[r] & \frac{K[\epsilon_{X}]}{V_{1} \epsilon_{X}} \ar[r]^-{\psi} & UGX.}
\]
On the other hand, by Proposition~\ref{Proposition-General-H_0-Sequence}, the sequence 
\[
0 \to \tfrac{K[\G \epsilon_{X}]}{V_{1}\G \epsilon_{X}} \to/{{ >}->}/^{\ker U\G \epsilon_{X}} U\G GX \to/{-{ >>}}/^{U\G \epsilon_{X}} U\G X \to 1
\]
is exact in $\Bc$. Via Proposition~\ref{Proposition-Long-Exact-Homology-Sequence-Simp}, it induces the exact homology sequence 
\begin{equation}\label{N-Long-Exact-Sequence}
\xymatrix{0 \ar[r] & H_2(X,U)_{\G} \ar@{{ >}->}[r] & H_0\frac{K[\G\epsilon_X]}{V_1\G\epsilon_X} \ar[r]^-{\varphi} & UGX }
\end{equation}
in $\Bc$. Indeed, it is well-known that $\epsilon_{GX}\colon {\G GX\to GX}$ is a contractible augmented simplicial object\index{contractible augm.\ simplicial object}\index{simplicial object!augmented!contractible}\index{augmented simplicial object!contractible}. Proposition~\ref{Proposition-Contractible-Object} implies that $H_1U\G GX=0$ and $H_0U\G GX=UGX$.

Accordingly, $\Delta V (X)$ is a kernel of $\psi$ and $H_2(X,U)_{\G}$ is a kernel of $\varphi$. We prove that $\psi$ and $\varphi$ are equal.

Reconsider Diagram~\ref{Diagram-Kappa}. Using Proposition~\ref{Proposition-pre-H_0-Sequence} and the functoriality of $H_{0}$, choosing the right coequalizers, we get that $\varphi$ is the unique morphism such that the right hand square \texttt{(iii)} in the diagram
\[
\vcenter{\xymatrix{
K[G\epsilon_X] \ar@{}[rrd]|{\texttt{(ii)}} \ar@{-{ >>}}[rr]^-p \ar@{{ >}->}[d]_{\ker G\epsilon_X} && \frac{K[G\epsilon_X]}{V_1G\epsilon_X} \ar@{}[rrd]|{\texttt{(iii)}} \ar@{-{ >>}}[rr]^-{\frac{\kappa}{V_{1} (\epsilon_{GX},\epsilon_{X})}}  \ar@{{ >}->}[d]|-{\ker UG\epsilon_X} && \frac{K[\epsilon_X]}{V_1\epsilon_X} \ar@{.>}[d]^{\varphi} \\
G^2X \ar@{-{ >>}}[rr]_{\eta_{G^2X}}  && UG^2X \ar@{-{ >>}}[rr]_{U\epsilon_{GX}} && UGX }}
\]
commutes. But also the square \texttt{(ii)} is commutative, and moreover $\psi$ is unique in making
\[
\xymatrix{
K[\epsilon_X] \ar@{{ >}->}[d]_-{\ker\epsilon_X}  \ar@{-{ >>}}[r]^-{q} & \frac{K[\epsilon_X]}{V_1\epsilon_X} \ar[d]^-{\psi}\\
GX \ar@{-{ >>}}[r]_-{\eta_{GX}} & UGX }
\]
commute. Now
\begin{align*}
\varphi\comp \tfrac{\kappa}{V_{1} (\epsilon_{GX},\epsilon_{X})}\comp p &= U\epsilon_{GX}\comp \eta_{G^2X}\comp \ker G\epsilon_X\\
&= \eta_{GX}\comp \epsilon_{GX}\comp \ker G\epsilon_X \\
&= \eta_{GX}\comp \ker \epsilon_X \comp \kappa\\
&= \psi \comp q \comp \kappa \\
&= \psi\comp \tfrac{\kappa}{V_{1} (\epsilon_{GX},\epsilon_{X})}\comp p,
\end{align*}
where the second equality follows from the naturality of $\eta$, and the third one holds by definition of $\kappa$. This proves our claim that $\varphi=\psi$.
\end{proof}

Combining Theorem~\ref{Theorem-V-Exact-Sequence} with Theorem~\ref{Theorem-Hopf-Formula} yields the following categorical version of the Stallings-Stammbach Sequence.

\begin{theorem}[Stallings-Stammbach Sequence]\label{Theorem-Cotriple-Exact-Sequence}\index{Stallings-Stammbach Sequence}
Let
\[
\xymatrix{ 0 \ar[r] & K \ar@{{ >}->}[r] & A \ar@{-{ >>}}[r]^-{f} & B \ar[r] & 1 }
\]
be a short exact sequence in $\Ac$. There exists an exact sequence
\[
\resizebox{\textwidth}{!}{\mbox{$ H_2(A,U)_{\G} \to<400>^{H_2(f,U)_{\G}} H_2(B,U)_{\G} \to<150>
\textstyle{\frac{K}{V_1f}} \to<150> H_{1}(A,U)_{\G}
\to/{-{ >>}}/<400>^{H_{1}
(f,U)_{\G}} H_{1}(B,U)_{\G} \to<150> 1$}}
\]
in $\Bc$, which depends naturally on the given short exact sequence.\noproof
\end{theorem}

\begin{example}[Crossed modules]\label{Example-Crossed-Modules}
Recall that a \emph{crossed module} $(T,G,\del)$\index{crossed module} is a group homomorphism $\del \colon T\to G$ together with an action of $G$ on $T$ (mapping a couple $(g,t)\in  G\times T$ to $^{g}t\in T$) satisfying
\begin{enumerate}
\item $\del (^{g}t)=g\del tg^{-1}$, for all $g\in G$, $t\in T$;
\item $^{\del t}s=tst^{-1}$, for all $s,t\in T$.
\end{enumerate}
The second requirement is called the \emph{Peiffer identity}\index{Peiffer identity}. A triple $(T,G,\del)$ that satisfies the first identity (but not necessarily the second one) is called a \emph{precrossed module}\index{precrossed module}. (We consider an internal notion of (pre)crossed module in Section~\ref{Section-Regular-Epimorphisms}.) A \emph{morphism of crossed modules}\[
(f,\phi)\colon (T,G,\del)\to (T',G',\del')
\]
is a pair of group homomorphisms $f\colon T\to T'$, $\phi \colon G\to G'$ with
\begin{enumerate}
\item $\del '\comp f=\phi \comp \del$;
\item $f (^{g}t)={^{\phi (g)}\!f (t)}$, for all $g\in G$, $t\in T$.
\end{enumerate}

It is well known that the category $\XMod$\index{category!$\XMod $} of crossed modules is equivalent to a variety of $\Omega$-groups, namely to the variety of $1$-categorical groups\index{one-categorical group@$1$-categorical group} (see Loday~\cite{Loday}). Hence, it is semi-abelian~\cite{Janelidze-Marki-Tholen}. Moreover, under this equivalence, a crossed module $(T,G,\del)$ corresponds to the semi-direct product $G\ltimes T$ equipped with the two appropriate homomorphisms. $G\ltimes T$ and $G\times T$ have the same underlying set; hence the forgetful functor $\Uc\colon \XMod\to\Set$ sends a crossed module $(T,G,\del)$ to the product $T\times G$ of its underlying sets. This determines a comonad $\G$ on $\XMod$. We get the cotriple homology of crossed modules described by Carrasco, Cegarra and Grandje\'an in~\cite{Carrasco-Homology} as a particular case of Definition~\ref{Definition-Cotriple-Homology} by putting $U$ the usual abelianization functor $\ab\colon \XMod\to \Ab\XMod$ (as defined for $\Omega$-groups, see Higgins~\cite{Higgins}).

In~\cite{Carrasco-Homology}, Carrasco, Cegarra and Grandje\'an give an explicit proof that
\[
\bigl[(T,G,\del), (T, G,\del)\bigr]=V (T,G,\del)
\]
equals $([T,G], [G,G],\del)$. For any crossed module $(T,G,\del)$, and any two normal subgroups $K\triangleleft G$, $S\triangleleft T$, let $[K,S]$ denote the (normal) subgroup of $T$ generated by the elements $(^{k}s)s^{-1}$, for $k\in K$, $s\in S$. It may be shown that, for $(N,R,\del)\triangleleft (Q,F,\del)$,
\[
\bigl[(N,R,\del),(Q,F,\del)\bigr]=V_{1}\Bigl((Q,F,\del)\to \tfrac{(Q,F,\del)}{(N,R,\del)} \Bigr)
\]
is equal to $ ([R,Q][F,N],[R,F],\del)$.

We recall the definition from~\cite{Carrasco-Homology}. Let $(T,G,\del)$ be a crossed module and $n\geq 1$. The $n$-th homology object $H_n(T,G,\del)$ of $(T,G,\del)$ is
\[
H_{n-1}C\Ab\G(T,G,\del).
\]
Here $C\ab\G(T,G,\del)$ is the \emph{unnormalized chain complex}\index{unnormalized chain complex} of $\ab\G(T,G,\del)$, defined by $(C\ab\G(T,G,\del))_n=(\ab\G(T,G,\del))_n$ and
\[
d_n=\del_0^{\ab}-\del_1^{\ab} + {\dots} + (-1)^n\del_n^{\ab}.
\]
$\ab\G(T,G,\del)$ denotes the simplicial abelian crossed module $\ab\comp \G(T,G,\del)$. The $\del_i^{\ab}$ are its face operators.

$C\ab\G(T,G,\del)$ need not be the same as $N\ab\G(T,G,\del)$. However, their homology objects are isomorphic, since $\ab\G(T,G,\del)$ is a simplicial object in the abelian category $\Ab \XMod$ (see Remark~\ref{Remark-Simplicial-Homology-in-Abelian-Case}). Hence $H_n(T,G,\del)=H_n((T,G,\del),\ab)_{\G}$.

The Hopf formula obtained in~\cite[Theorem 13 and above]{Carrasco-Homology} is a particularization of our Theorem~\ref{Theorem-Hopf-Formula}. In this particular case, the exact sequence of Theorem~\ref{Theorem-Cotriple-Exact-Sequence} becomes the exact sequence of~\cite[Theorem 12 (i)]{Carrasco-Homology}.
\end{example}

\begin{example}[Groups]\label{Example-Group-Homology}\index{category!$\Gp$}
It is shown in~\cite[Theorem 10]{Carrasco-Homology} that cotriple homology of crossed modules encompasses classical group homology; hence, so does our theory. If
\[
\xymatrix{ 0 \ar[r] & R \ar@{{ >}->}[r] & F \ar@{-{ >>}}[r] & G \ar[r] & 1 }
\]
is a presentation of a group $G$ by generators and relations, 
\[
U=\ab\colon \Gp \to \Ab\Gp=\Ab
\]
is the abelianization functor and $\G$ the ``free group on a set''-monad, the sequence in Theorem~\ref{Theorem-Cotriple-Exact-Sequence} becomes the Stallings-Stammbach sequence in integral homology of groups~\cite{Stallings, Stammbach}. The isomorphism in Theorem~\ref{Theorem-Hopf-Formula} is nothing but Hopf's formula~\cite{Hopf}
\[
H_{2} (G,\Z)\cong\frac{R\cap [F,F]}{[R,F]}.
\]
In fact, by proving this formula, Hopf showed that the right hand side expression is independent of the chosen presentation of $G$; this brings us back to the aim of Chapter~\ref{Chapter-Baer-Invariants}. 

Although, as mentioned above, it follows from \cite{Carrasco-Homology}, it is probably not immediately clear \textit{why} this type of cotriple homology coincides with integral homology of groups. We give a direct argument. On one hand, given a group $X$, the augmented simplicial group 
\[
\xymatrix{{\cdots} \ar@<1ex>[r] \ar[r] \ar@<-1ex>[r] & G^{2}X \ar@<.5ex>[r] \ar@<-.5ex>[r] & GX \ar[r] & X}
\]
is a free simplicial resolution of $X$. Applying the abelianization functor $\ab$ yields an augmented simplicial abelian group 
\begin{equation}\label{Simplicial-Group-G}
\xymatrix{{\cdots} \ar@<1ex>[r] \ar[r] \ar@<-1ex>[r] & {\ab (G^{2}X)} \ar@<.5ex>[r] \ar@<-.5ex>[r] & {\ab (GX)} \ar[r] & {\ab (X).}}
\end{equation}
Its homology groups are our $H_{n} (X,\ab)_{\G}$.

On the other hand, the integral homology groups $H_{n} (X,\Z)$ of $X$ are defined as follows. One takes, in the category ${_{\Z X}\Mod}$\index{category!${_{\Z X}\Mod}$} of left modules over the groupring $\Z X$, a free simplicial resolution of~$\Z$:
\begin{equation}\label{Simplicial-Resolution}
\xymatrix{{\cdots} \ar@<1ex>[r] \ar[r] \ar@<-1ex>[r] & B_{1} \ar@<.5ex>[r] \ar@<-.5ex>[r] & B_{0} \ar[r] & \Z.}
\end{equation}
When considering \emph{integral}\index{homology!integral --- of groups}\index{integral homology of groups} homology of the group $X$, one chooses for module of coefficients the trivial right $\Z X$-module $\Z$. Tensoring on the left yields the following simplicial abelian group:
\[
\xymatrix{{\cdots} \ar@<1ex>[r] \ar[r] \ar@<-1ex>[r] & {\Z \tensor_{\Z X}B_{1}} \ar@<.5ex>[r] \ar@<-.5ex>[r] & {\Z \tensor_{\Z X}B_{0}.}}
\]
Its objects are free abelian groups, as $\Z \tensor_{\Z X} (\cdot)\colon {{_{\Z X}\Mod}\to \Ab}$ is just the forgetful functor, and being a left adjoint, it sends a free module on a free abelian group. The homology groups of this simplicial abelian group are the integral homology groups $H_{n} (X,\Z)$ of $X$.

Now it is well-known that its $0$-th homology, the coequalizer of 
\[
\del_{0},\del_{1}\colon {\Z \tensor_{\Z X}B_{1} \to \Z \tensor_{\Z X}B_{0}},
\]
is equal to $\ab (X)$; hence
\begin{equation}\label{Simplicial-Group-Tensor}
\xymatrix{{\cdots} \ar@<1ex>[r] \ar[r] \ar@<-1ex>[r] & {\Z \tensor_{\Z X}B_{1}} \ar@<.5ex>[r] \ar@<-.5ex>[r] & {\Z \tensor_{\Z X}B_{0}} \ar@{-{ >>}}[r] & {\ab (X)}}
\end{equation}
is also an augmented simplicial abelian group. Modifying the argument in the proof of~\cite[Theorem~4.1]{Barr-Beck-LaJolla}, one sees that the augmented simplicial abelian groups~\ref{Simplicial-Group-Tensor} and~\ref{Simplicial-Group-G} induce homotopically equivalent chain complexes; hence their homology groups are isomorphic. For the sake of completeness, we give a brief outline of Barr and Back's argument. 

Let $\Cc$ be a category and $K$ and $L$ functors from $\Cc$ to the category $\Ch_{\geq -1}\Ab$\index{category!$\Ch_{\geq -1}\Ab$} of chain complexes of abelian groups with degrees bigger than or equal to $-1$ (i.e., $C_{n}=0$ if $C\in \Ch_{\geq -1}\Ab$ and $n<-1$). Let $G\colon {\Cc \to \Cc}$ be a functor and $\epsilon \colon {G\To 1_{\Cc}}$ a natural transformation. Then $K$ is called \emph{$G$-acyclic}\index{G-acyclic@$G$-acyclic chain complex} when there is a functorial contracting homotopy in $K\comp G$\index{homotopy!of chain complexes} (i.e., for every $X$ in $\Cc$, $K (G (X))$ is chain homotopic to $0$, and there is a functorial choice of such chain homotopies). $L$ is called \emph{$G$-representable}\index{G-representable@$G$-representable chain complex} when there are natural transformations $\theta_{n}\colon {L_{n}\To L_{n}\comp G}$ such that $(1_{L_{n}}\comp \epsilon )\vcomp \theta_{n}=1_{L_{n}}$, for $n\in \N$.

\begin{theorem}\cite[Theorem~3.1, Corollary~3.2]{Barr-Beck-LaJolla}\label{Theorem-Acyclic-Models}
Let $K$ be $G$-acyclic and $L$ $G$-representable. Any natural transformation $f_{-1}\colon {K_{-1}\To L_{-1}}$ induces a natural transformation $f\colon {K\To L}$, unique up to chain homotopy. 

In particular, if $K$ and $L$ are both $G$-acyclic and $G$-representable, then $K_{-1}\cong L_{-1}$ implies $K\simeq L$. Hence $K$ and $L$ have isomorphic homology groups.\noproof 
\end{theorem}

Let $\G = (G,\delta,\epsilon)$ denote the comonad on $\Gp$ from above. It is well-known and easily seen that the augmented simplicial abelian group~\ref{Simplicial-Group-G} induces a $G$-acyclic and $G$-representable chain complex of abelian groups; both a contraction and natural transformations $\theta_{n}$ are induced by the maps $G^{n}\delta_{X}\colon {G^{n+1}X\to G^{n+2}X}$. It is also possible to choose the simplicial resolutions~\ref{Simplicial-Resolution} in such a way that the augmented simplicial abelian group~\ref{Simplicial-Group-Tensor} induces a $G$-acyclic and $G$-representable chain complex of abelian groups. The result now follows from Theorem~\ref{Theorem-Acyclic-Models}.
\end{example}

\begin{example}[Lie algebras]\label{Example-Lie-Homology}
Like for groups, it is well-known that, presenting a Lie algebra $\Lieg$ as a quotient $\Lief/\Lier$ of a free Lie algebra $\Lief$, one may characterize the second homology object with coefficients in the base field $\K$ as
\[
H_{2} (\Lieg,\K)\cong\frac{\Lier\cap [\Lief,\Lief]}{[\Lier,\Lief]}.
\] 
This may be seen as a special case of Theorem~\ref{Theorem-Hopf-Formula}. Also the Stallings-Stammbach Sequence reduces to a classical exact sequence of Lie algebras.
\end{example}
 % Cotriple homology
\chapter{Cohomology}\label{Chapter-Cohomology}

\setcounter{section}{-1}
\section{Introduction}\label{Section-Cohomology-Introduction}

In this chapter we prove some new results in the study of cohomology in the semi-abelian context. More precisely, we can establish: 
\begin{itemize}
\item a natural Hochschild-Serre $5$-term exact sequence~\cite{Hochschild-Serre}
in cohomology with trivial coefficients (Theorem~\ref{Theorem-Cohomology-Sequence}); 
\item an interpretation of the second cohomology group $H^2(Y,A)$ as the group $\CentrExt (Y,A)$ of isomorphism classes of central extensions of $Y$ by an abelian object $A$, equipped with (a generalization of) the Baer sum (Theorem~\ref{Theorem-H^2-Central-Extensions});
\item a universal coefficient theorem relating homology and cohomology (Theorem~\ref{Theorem-Universal-Coefficients}).
\end{itemize}
Thus we simplify several recent investigations in this direction, in the context of crossed modules~\cite{Carrasco-Homology} or precrossed modules~\cite{AL}, and unify them with the classical theory that exists for groups and Lie algebras. Our approach is based on the work of Fr\"ohlich's school---in particular, of Lue~\cite{Lue}---and also on the work of Janelidze and Kelly~\cite{Janelidze-Kelly} (see also~\cite{EG, DKV} and Section~\ref{Section-Commutators}). And indeed, we make heavy use of the categorical tools designed for understanding the universal algebraic property of centrality~\cite{Pedicchio, Janelidze-Pedicchio, Janelidze-Kelly, Janelidze-Kelly:Maltsev}. To avoid technical complications we restrict the context to that of a semi-abelian category with its Birkhoff subcategory of abelian objects. This also allows for a simplification in the notations: To denote commutators, we may use brackets $[K,X]$ instead of the more explicit $V_{1}f$.

In the first section we recall some important facts concerning semi-abelian categories, which will be needed throughout the chapter, and prove some useful technical properties involving exact sequences and central extensions. In Section~\ref{Section-Perfect-Case}, we show that an object is perfect if and only if it admits a universal central extension. In Section~\ref{Section-Cohomology} we establish a cohomological version of Hopf's formula as well as the Hochschild-Serre $5$-term exact sequence for cohomology: An extension $f\colon {X\to Y}$ with kernel $K$ induces the exact sequence
\[
 0 \to<250> H^{1} Y \to/{{ >}->}/<250>^{H^{1}f} H^{1} X \to<250> \hom \bigl(\textstyle{\frac{K}{[K,X]}},A\bigr) \to<250> H^{2} Y \to<250>^{H^{2}f} H^{2} X.
\]
In Section~\ref{Section-Second-Cohomology} we prove that the second cohomology group $H^{2} (Y,A)$ is isomorphic to the group $\CentrExt (Y,A)$ of isomorphism classes of central extension of $Y$ by $A$, equipped with the Baer sum~\cite{Gerstenhaber}. In the last section we obtain a universal coefficient theorem for cohomology in semi-abelian categories and give some applications.

This chapter contains most of my paper \cite{Gran-VdL} with Marino Gran; some results on commutators however appear in Section \ref{Section-Commutators}, which seemed more appropriate.

\section{Abelian objects and central extensions}\label{Section-Abelianization}

In this section, we recall some important facts from the previous chapters and adapt some notations to the present situation: We find ourselves in a semi-abelian category and concentrate on abelianization. We also prove a technical result on central extensions (Corollary~\ref{Corollary-Central-Extension-Pushout}), which will be needed in the subsequent sections.

Recall from Section~\ref{Section-Commutators} that an \emph{extension of an object $Y$ (by an object $K$)}\index{extension} is a regular epimorphism $f\colon X\to Y$ with its kernel $K$:
\[
0 \to K \to/{{ >}->}/^{k} X \to/-{ >>}/^{f} Y \to 0.
\]
The category of extensions of $Y$ (considered as a full subcategory of the slice category $\slashfrac{\Ac}{Y}$) is denoted by $\Ext (Y)$\index{category!$\Ext (Y)$}; the category of all extensions in $\Ac$ (considered as a full subcategory of the arrow category $\Fun (\mathsf{2},\Ac)$: morphisms are commutative squares) by $\Ext \Ac$\index{category!$\Ext\Ac$}. An extension $f\colon {X\to Y}$ is called \emph{central}\index{central extension}\index{extension!central} if its kernel is central in the sense of Huq, i.e., if $\ker f$ cooperates with $1_{X}$. We write $\CentrExt (Y)$\index{category!$\CentrExt (Y)$} for the full subcategory of $\Ext (Y)$ determined by the central extensions. Recall Proposition~\ref{Proposition-Subobject-of-Central-is-Normal} that, given a central extension $f\colon {X\to Y}$ in $\Ac$, every subobject of $k=\ker f\colon {K\to X}$ is normal in $X$.

Recall that an object $X$ in $\Ac$ is called \emph{abelian}\index{abelian!object} when there is a centralizing relation on $\nabla_{X}$ and $\nabla_{X}$ (see Definition~\ref{Definition-Centralizing-Equivalence-Relations} and~\ref{Definition-Unital-Category}). The full subcategory of $\Ac$ determined by the abelian objects is denoted $\Ab\Ac$\index{category!$\Ab\Ac$}. Since $\Ac$ is pointed, $\Ab\Ac$ coincides with the category of internal abelian group objects in $\Ac$. $\Ab\Ac$ is an abelian category and a Birkhoff subcategory of $\Ac$.
\[
\xymatrix{\Ac \ar@<.8 ex>[r]^-{\ab} \ar@{}[r]|-{\perp} & {\,\Ab\Ac} \ar@<.8 ex>@{_(->}[l]}
\]
For any object $X$, we shall denote the $X$-component of the unit of the adjunction by 
\[
\eta_{X}\colon {X\to X/[X,X]}=\ab (X),
\]
and its kernel by $\mu_{X}\colon [X,X]\to X$, the $X$-component of a natural transformation 
\[
\mu \colon V\To 1_{\Ac }\colon \Ac \to \Ac
\]
(cf.\ Proposition~\ref{Proposition-Birkhoff-Subfunctor}).

Similarly, for any object $Y$, the category $\CentrExt (Y)$ is reflective and closed under subobjects and quotients in $\Ext (Y)$.
\[
\xymatrix{{\Ext (Y)} \ar@<.8 ex>[r]^-{\centr} \ar@{}[r]|-{\perp} & {\,\CentrExt (Y)} \ar@<.8 ex>@{_(->}[l]}
\]
The $f$-component of the unit of the adjunction is given by the horizontal arrows in the diagram
\[
\xymatrix{X \ar@{-{ >>}}[r] \ar@{-{ >>}}[d]_-{f} & \textstyle{\frac{X}{[K,X]}} \ar@{-{ >>}}[d]^-{\centr f}\\
Y \ar@{=}[r] & Y,}
\]
where $k\colon K\to X$ denotes a kernel of $f\colon X\to Y$ (Remark~\ref{Remark-Centrext-Birkhoff} or \cite[Theorem 2.8.11]{Borceux-Bourn}). The functor ${\Ext \Ac \to \Ac}$ that maps an extension $f$ to $[K,X]$ is denoted by $V_{1}$\index{V_1@$V_{1}$} (see Remark~\ref{Remark-Centrext-Birkhoff}).

\begin{proposition}\label{Proposition-Split-Mono-Central-Extension}
Consider the diagram of solid arrows
\[
\xymatrix{0 \ar[r] & K \ar@{{ >}->}[r]^-{k} \ar@{{ >}->}@<-.5 ex>[d]_-{m} & X \ar@{-{ >>}}[r]^-{f} \ar@{{ >}.>}@<-.5 ex>[d]_-{\overline{m}} & Y \ar@{:}[d] \ar[r] & 0\\
0 \ar@{.>}[r] & A \ar@{-{ >>}}@<-.5 ex>[u]_-{s} \ar@{{ >}.>}[r] & Z \ar@{.{ >>}}@<-.5 ex>[u] \ar@{.{ >>}}[r] & Y \ar@{.>}[r] & 0.}
\]
If the above sequence is exact, $m$ is a normal monomorphism split by $s$, $f$ is a central extension and $A$ is abelian, then a central extension of $Y$ by $A$ exists making the diagram commutative.
\end{proposition}
\begin{proof}
Let $q\colon A\to Q$ denote a cokernel of $m$, and let us consider the sequence
\[
\xymatrix{0 \ar[r] & K \ar@<-.5 ex>@{{ >}->}[r]_-{m} & A \ar@<-.5 ex>@{-{ >>}}[l]_-{s} \ar@{-{ >>}}[r]^-{q} & Q \ar[r] & 0}
\]
in $\Ab\Ac$, a split exact sequence. It follows that $A$ is a product of $K$ with $Q$ and that, up to isomorphism, $m$ is $l_{K}\colon K\to K\times Q$ and $s$ is the projection $\pr_{1}\colon {K\times Q\to K}$ on $K$. In the diagram
\[
\xymatrix{0 \ar[r] & K \ar@{{ >}->}[r]^-{k} \ar@{}[rd]|-{\texttt{(i)}} \ar@{{ >}->}@<-.5 ex>[d]_-{l_{K}} & X \ar@{-{ >>}}[r]^-{f} \ar@{{ >}->}@<-.5 ex>[d]_-{l_{X}} & Y \ar@{=}[d] \ar[r] & 0\\
& K\times Q \ar@{-{ >>}}@<-.5 ex>[u]_-{\pr_{1}} \ar@{{ >}->}[r]_-{k\times 1_{Q}} & {X\times Q} \ar@{-{ >>}}@<-.5 ex>[u]_-{\pr_{1}} \ar@{-{ >>}}[r]_-{f\comp \pr_{1}} & Y}
\]
the upward-pointing square \texttt{(i)} is a pullback; it follows that $k\times 1_{Q}$ is a kernel of $f\comp \pr_{1}$, and then $f\comp \pr_{1}$ is the cokernel of $k\times 1_{Q}$. From the fact that $Q$ is abelian and $k$ is central one concludes that $k\times 1_{Q}$ is central, and this completes the proof.
\end{proof}

\begin{proposition}\label{Proposition-Reg-Epi-Central-Extension}
Consider the diagram of solid arrows
\[
\xymatrix{0 \ar[r] & K \ar@{{ >}->}[r]^-{k} \ar@{-{ >>}}[d]_-{r} & X \ar@{-{ >>}}[r]^-{f} \ar@{.{ >>}}[d]^-{\overline{r}} & Y \ar@{:}[d] \ar[r] & 0\\
0 \ar@{.>}[r] & A \ar@{{ >}.>}[r] & Z \ar@{.{ >>}}[r] & Y \ar@{.>}[r] & 0.}
\]
If the above sequence is exact, $r$ is a regular epimorphism, $f$ is a central extension and $A$ is abelian, then a central extension of $Y$ by $A$ exists making the diagram commutative.
\end{proposition}
\begin{proof}
We may define $\overline{r}$ as a cokernel of $k\comp \ker r$:
\[
\xymatrix{& K[r] \ar@{=}[r] \ar@{{ >}->}[d]_-{\ker r} & K[r] \ar@{{ >}->}[d]^-{k\comp \ker r} \\
0 \ar[r] & K \ar@{}[rd]|{\texttt{(i)}} \ar@{{ >}->}[r]^-{k} \ar@{-{ >>}}[d]_-{r} & X \ar@{.{ >>}}[d]_-{\overline{r}} \ar@{-{ >>}}[r]^-{f} & Y \ar[r] \ar@{:}[d] & 0\\
0 \ar@{.>}[r] & A \ar@{{ >}.>}[r]_-{\overline{k}} & Z \ar@{.{ >>}}[r] & Y \ar@{.>}[r] & 0.}
\]
The morphism $k\comp \ker r$ is a kernel thanks to Proposition~\ref{Proposition-Subobject-of-Central-is-Normal}. Taking cokernels induces the square \texttt{(i)}; it is a pushout by Proposition~\ref{Proposition-pushout}. Remark also that, by Proposition~\ref{Proposition-LeftRightPullbacks}, $\overline{k}$ is a monomorphism, hence a kernel, since it is the regular image of $k$ along $\overline{r}$ (Proposition~\ref{Proposition-Image-of-Kernel}). Taking a cokernel of $\overline{k}$ gives rise to the rest of the diagram, thanks to Proposition~\ref{Proposition-pushout}. The induced extension is central, because $\CentrExt \Ac$ is closed under quotients in $\Ext \Ac$---and a quotient in $\Ext \Ac$ is a pushout in $\Ac$ of a regular epimorphism along a regular epimorphism.
\end{proof}

\begin{corollary}\label{Corollary-Central-Extension-Pushout}
Consider the diagram of solid arrows
\[
\xymatrix{0 \ar[r] & K \ar@{{ >}->}[r]^-{k} \ar[d]_-{a} & X \ar@{-{ >>}}[r]^-{f} \ar@{.>}[d]^-{\overline{a}} & Y \ar@{:}[d] \ar[r] & 0\\
0 \ar@{.>}[r] & A \ar@{{ >}.>}[r]_-{k_{a}} & Z_{a} \ar@{.{ >>}}[r]_-{p_{a}} & Y \ar@{.>}[r] & 0.}
\]
If the above sequence is exact, $f$ is a central extension and $A$ is abelian, then a central extension of $Y$ by $A$ exists making the diagram commutative.
\end{corollary}
\begin{proof}
Since $K$ and $A$ are abelian, the arrow $a\colon K\to A$ lives in the abelian category $\Ab\Ac$, and thus may be factorized as the composite
\[
K\to/{ >}->/^{l_{K}} K\oplus A \to/-{ >>}/^{[a,1_{A}]} A. 
\]
The result now follows from the previous propositions. 
\end{proof}

\section{The perfect case: universal central extensions}\label{Section-Perfect-Case}

Suppose that $\Ac$ is a semi-abelian category with enough projectives and $Y$ an object of $\Ac$. Then the category $\CentrExt (Y)$ always has a weakly initial object: For if $f\colon X\to Y$ is a \emph{projective presentation}\index{projective presentation}\index{presentation!projective} of $Y$, i.e., a regular epimorphism with $X$ projective, then the reflection $\centr f\colon X/[K,X]\to Y$ of $f$ into $\CentrExt (Y)$ is a central extension of $Y$. It is weakly initial, as any other central extension $g\colon Z\to Y$ induces a morphism $\centr f\to g$ in $\CentrExt (Y)$, the object $X$ being projective.

An initial object in $\CentrExt (Y)$ is called a \emph{universal}\index{universal central extension}\index{central extension!universal} central extension of $Y$. In contrast with the existence of weakly initial objects, in order that a universal central extension of $Y$ exists, the object $Y$ must be \emph{perfect}\index{perfect object}: Such is an object $Y$ of $\Ac$ that satisfies $[Y,Y]=Y$, i.e., that has $0$ as a reflection $\ab (Y)$ into $\Ab\Ac$.

To see this, recall from Chapter~\ref{Chapter-Baer-Invariants} the definition of \emph{Baer invariant}\index{Baer invariant}: Such is a functor $B\colon {\Ext \Ac \to \Ac} $ that makes \emph{homotopic}\index{homotopy!of morphisms of extensions} morphisms of extensions equal: morphisms $(f_{0},f)$ and $(g_{0},g)\colon {p\to q}$ 
\[
\xymatrix{A_{0} \ar@{-{ >>}}[d]_-{p} \ar@<-0.5 ex>[r]_-{g_{0}}\ar@<0.5 ex>[r]^-{f_{0}} & B_{0} \ar@{-{ >>}}[d]^-{q}\\
A \ar@<-0.5 ex>[r]_-{g}\ar@<0.5 ex>[r]^-{f} & B}
\]
satisfying $f=g$. Such a functor $B$ sends homotopically equivalent extensions to isomorphic objects. The functor $V_{0}/V_{1}\colon \Ext \Ac \to \Ac$ that maps an extension
\[
\xymatrix{0 \ar[r] & K \ar@{{ >}->}[r]^{k} & A_{0}  \ar@{-{ >>}}[r]^{p} & A \ar[r] & 0}
\] 
to $[A_{0},A_{0}]/[K,A_{0}]$ is an example of a Baer invariant (Proposition~\ref{L_1}).

\begin{proposition}\label{Proposition-Universal-Central-Extension}
Let $\Ac$ be a semi-abelian category with enough projectives. An object $Y$ of $\Ac$ is perfect if and only if $Y$ admits a universal central extension.
\end{proposition}
\begin{proof}
First suppose that $\CentrExt (Y)$ has an initial object $u\colon U\to Y$. Since
\[
\pr_{1}\colon {Y\times \ab (Y)\to Y}
\]
is central, a unique morphism $(u,y)\colon  U\to Y\times \ab (Y)$ exists. But then the map $0\colon {U\to \ab (Y)}$ is equal to $\eta_{Y}\comp u\colon {U\to \ab (Y)}$, and $\ab (Y)=0$.

Conversely, consider a presentation $f\colon X\to Y$ of a perfect object $Y$ (with its kernel $k\colon {K\to X}$) and the diagram of solid arrows
\[
\xymatrix{& [K,X] \ar@{=}[r] \ar@{{ >}->}[d] & [K,X] \ar@{{ >}->}[d] \\
0 \ar[r] & K\cap [X,X] \ar@{{ >}->}[r] \ar@{.{ >>}}[d] & [X,X] \ar@{.{ >>}}[d] \ar@{-{ >>}}[r]^-{f\comp \mu_{X}} & Y \ar[r] \ar@{:}[d] & 0\\
0\ar@{.>}[r] & {\tfrac{K\cap [X,X]}{[K,X]}} \ar@{{ >}.>}[r] & {\tfrac{[X,X]}{[K,X]}} \ar@{.{ >>}}[r]_-{u} & Y \ar@{.>}[r] & 0.}
\]
Its short sequence is exact because $Y$ is perfect. As in the proof of Proposition~\ref{Proposition-Reg-Epi-Central-Extension}, taking cokernels, we get a diagram with all rows exact as above.

We claim that $u$ is the needed universal central extension. Note that the construction of $u$ is independent of the choice of the presentation, because its objects are Baer invariants. Moreover, for the same reason, we could as well have constructed $u$ out of the extension $f\comp \mu_{X}\colon {[X,X]\to Y}$: Indeed, $f$ and $f\comp  \mu_{X}$ are homotopically equivalent. Thus we see that $u$ is the reflection of $f\comp \mu_{X}$ into $\CentrExt (Y)$. It follows at once that $u$ is weakly initial, as $\mu_{X}$ induces an arrow in $\CentrExt (Y)$ from $u$ to the reflection $\centr f$ of $f$.

This weakly initial object is, in fact, initial. First note that, since $u$ is central, the bottom short exact sequence in the diagram above induces the following exact sequence (by Theorem~\ref{Theorem-V-Exact-Sequence}):
\[
\xymatrix{{\tfrac{K\cap [X,X]}{[K,X]}} \ar[r] & {\tfrac{K\cap [X,X]}{[K,X]}} \ar[r] & \tfrac{{\tfrac{[X,X]}{[K,X]}}}{\Bigl[\tfrac{[X,X]}{[K,X]},\tfrac{[X,X]}{[K,X]}\Bigr]} \ar[r] & \tfrac{Y}{[Y,Y]} \ar[r] & 0.}
\]
Hence $Y$ perfect implies that $[X,X]/[K,X]$ is perfect as well. In particular, 
\[
\mu_{\tfrac{[X,X]}{[K,X]}}\colon \Bigl[\tfrac{[X,X]}{[K,X]},\tfrac{[X,X]}{[K,X]}\Bigr]\to \tfrac{[X,X]}{[K,X]}
\]
is an isomorphism. Now suppose that $g\colon Z\to Y$ is a central extension and 
\[
a,b\colon {[X,X]}/{[K,X]\to Z}
\]
satisfy $g\comp a=g\comp b=u$. Then the centrality of $u$ and $g$ and the fact that the functor $V_{0}/V_{1} \colon {\Ext \Ac \to \Ac}$ is a Baer invariant imply that $V (a)=V (b)$. Finally, by naturality of $\mu$, 
\[
a\comp \mu_{\tfrac{[X,X]}{[K,X]}}=\mu_{Z}\comp V (a)=\mu_{Z}\comp V (b)=b\comp \mu_{\tfrac{[X,X]}{[K,X]}}
\]
and $a=b$. 
\end{proof}

\section{Cohomology and the Hochschild-Serre Sequence}\label{Section-Cohomology}

From now on, $\Ac$ will be a semi-abelian category, monadic\index{monadic category}\index{category!monadic} over the category $\Set$ of sets. Let $\G$ denote the comonad on $\Ac$, induced by the monadicity requirement, as well as the resulting functor from $\Ac$ to the category $\simpA$ of simplicial objects in $\Ac$. (We ask that $\Ac$ is monadic for simplicity only; we could as well have asked that Assumption~\ref{Assumption-Cotriples} is satisfied.) $\Ac$ has enough projective objects; we shall consider all Baer invariants relative to the web $\Wproj$ of projective presentations, and never mention it. We recall some concepts, results and notations from Chapter~\ref{Chapter-Cotriples}, adapted to our present situation. (See also Example~\ref{Example-Functorial-Web-Monadic}.)

Recall that a chain complex in a semi-abelian category is called \emph{proper} when its boundary operators have normal images. As in the abelian case, the $n$-th homology object $H_{n}C$ of a proper chain complex $C$ with boundary operators $d_n$ is the cokernel of $C_{n+1}\to K[d_n]$. The \emph{normalization functor}\index{normalization functor} $N\colon \simpA \to \PCh \Ac$ turns a simplicial object $A$ into the \emph{Moore complex}\index{Moore complex} $N (A)$ of $A$, the chain complex with $N_{0}A=A_{0}$,
\[
N_{n} A=\bigcap_{i=0}^{n-1}K[\del_{i}\colon A_{n}\to A_{n-1}]
\]
and boundary operators $d_{n}=\del_{n}\comp \bigcap_{i}\ker \del_{i}\colon N_{n} A\to N_{n-1} A$, for $n\geq 1$, and $A_{n}=0$, for $n<0$. Since $N (A)$ is a proper chain complex in a semi-abelian category, one can define its homology objects in the usual way.

\begin{definition} (Definition~\ref{Definition-Cotriple-Homology}, case $\Bc =\Ab\Ac$)\label{Definition-Homology-5}
For $n\in \N $, the object
\[
H_{n} X=H_{n-1} N \ab (\G X)
\]
is the \emph{$n$-th homology object of $X$ (with coefficients in $\ab $) relative to the comonad $\G$}\index{homology!w.r.t.\ a comonad}\index{cotriple homology}. This defines a functor $H_{n} \colon \Ac \to \Ab\Ac $\index{H_n@$H_{n} X$}, for any $n\in \N_{0}$. 
\end{definition}

\begin{proposition}\label{Proposition-H_0-Sequence}
Let
\[
\xymatrix{ 0 \ar[r] & K \ar@{{ >}->}[r] & X \ar@{-{ >>}}[r]^-{f} & Y \ar[r] & 0 }
\]
be a short exact sequence in $\Ac$. Then the induced sequence in $\Ab \Sc \Ac$ 
\begin{equation}\label{Sequence-simpB-5}
 0 \to K[\ab\G f] \to/{{ >}->}/^{\ker \ab\G f} \ab\G X \to/{-{ >>}}/^{\ab\G f} \ab\G Y \to 0
\end{equation}
is degreewise split exact and is such that 
\[
H_{0}K[\ab \G f]\cong\textstyle{\frac{K}{[K,X]}}.
\]
\end{proposition}
\begin{proof}
This is a combination of Proposition~\ref{Proposition-General-H_0-Sequence} and Proposition~\ref{Proposition-pre-H_0-Sequence}.\end{proof}

The following is an application of Theorems~\ref{Theorem-Hopf-Formula} and~\ref{Theorem-Cotriple-Exact-Sequence}.

\begin{theorem}[Stallings-Stammbach sequence and Hopf formula]\label{Theorem-Homology}\index{Stallings-Stammbach Sequence}\index{Hopf's Formula}
If
\[
\xymatrix{ 0 \ar[r] & K \ar@{{ >}->}[r] & X \ar@{-{ >>}}[r]^-{f} & Y \ar[r] & 0 }
\]
is a short exact sequence in $\Ac$, then there exists an exact sequence
\begin{equation}\label{Stallings-Stammbach-Sequence}
 H_2 X \to<250>^{H_2 f} H_{2} Y \to<250> \textstyle{\frac{K}{[K,X]}} \to<250> H_{1} X \to/{-{ >>}}/<250>^{H_{1}f} H_{1} Y \to<250> 0
\end{equation}
in $\Ab\Ac$, which depends naturally on the given short exact sequence. Moreover $H_{1} Y\cong \ab (Y)$ and, when $X$ is projective, $H_{2} Y\cong {(K\cap [X,X])}/{[K,X]}$.\noproof 
\end{theorem}

Let $A$ be abelian group object in $\Ac$. Recall that the sum $a+b$ of two elements $a,b\colon {X\to A}$ of a group $\hom (X,A)$ is the composite 
\[
m\comp (a,b)\colon {X\to A}
\]
of $(a,b)\colon {X\to A\times A}$ with the multiplication $m\colon {A\times A\to A}$ of $A$. Homming into $A$ defines a functor 
\[
\hom (\cdot,A)\colon {\Ac^{\op}\to \Ab }.
\]
Given a simplicial object $S$ in $\Ac$, its image $\hom (S,A)$ is a cosimplicial object of abelian groups; as such, it has cohomology groups $H^{n}\hom (S,A)$.

\begin{definition}\label{Definition-Cotriple-Cohomology}\index{cotriple cohomology}\index{cohomology!cotriple ---}
Let $\Ac$ be a semi-abelian category, monadic over $\Set$, and let $\G$ be the induced comonad. Let $X$ be an object of $\Ac$ and $A$ an abelian object. Consider $n\in \N_{0}$. We say that
\[
H^{n} (X,A)=H^{n-1}\hom (\G X,A)
\]
is the \emph{$n$-th cohomology group of $X$ with coefficients in $A$ (relative to the cotriple $\G$)}.\index{cohomology!w.r.t\ a comonad} This defines a functor $H^{n} (\cdot,A)\colon \Ac \to \Ab $\index{H^n@$H^{n} (\cdot,A)$}, for any $n\in \N$. When it is clear which abelian group object $A$ is meant, we shall denote it just $H^{n} (\cdot)$.\index{H^n@$H^{n}X$}
\end{definition}

\begin{remark}\label{Remark-Is-Barr-Beck}
Unlike the definition of homology, this is merely an instance of Barr and Beck's general definition of cotriple cohomology~\cite{Barr-Beck}: $H^{n} (X,A)$ is the $n$-th cohomology group of $X$, with coefficients in the functor $\hom (\cdot,A)\colon {\Ac^{\op}\to \Ab}$, relative to the cotriple $\G$.
\end{remark}

\begin{proposition}\label{Proposition-First-Cohomology}
For any object $X$ and any abelian object $A$ of $\Ac$, 
\[
H^{1} (X,A)\cong\hom (H_{1} X, A)\cong\hom (\ab (X),A)\cong\hom (X,A).
\]
If $X$ is projective then $H^{n} X=0$, for any $n\geq 2$.
\end{proposition}
\begin{proof}
The first isomorphism is a consequence of the fact that $\hom (\cdot, A)$ turns coequalizers in $\Ac$ into equalizers in $\Ab$. The second isomorphism follows from Theorem~\ref{Theorem-Homology} and the third one by adjointness of the functor $\ab$.

The second statement follows because if $X$ is projective then $\G X$ is contractible (see~\cite{Barr-Beck}).
\end{proof}
The following result extends Theorem~12 in~\cite{Carrasco-Homology} and Theorem~1 in~\cite{AL}:
\begin{theorem}[Hochschild-Serre Sequence]\index{Hochschild-Serre Sequence}\label{Theorem-Cohomology-Sequence}
Let
\[
\xymatrix{ 0 \ar[r] & K \ar@{{ >}->}[r]^-{k} & X \ar@{-{ >>}}[r]^-{f} & Y \ar[r] & 0 }
\]
be a short exact sequence in $\Ac$. There exists an exact sequence
\[
 0 \to<250> H^{1} Y \to/{{ >}->}/<250>^{H^{1}f} H^{1} X \to<250> \hom \bigl(\textstyle{\frac{K}{[K,X]}},A\bigr) \to<250> H^{2} Y \to<250>^{H^{2}f} H^{2} X
\]
in $\Ab$, which depends naturally on the given short exact sequence.
\end{theorem}
\begin{proof}
The sequence~\ref{Sequence-simpB-5} is degreewise split exact; hence homming into $A$ yields an exact sequence of abelian cosimplicial groups 
\[
 0 \to \hom (\ab\G Y,A)  \to/{{ >}->}/ \hom (\ab\G X, A) \to/{-{ >>}}/ \hom (K[\ab \G f],A) \to 0.
\]
This gives rise to an exact cohomology sequence
\[
%\resizebox{\textwidth}{!}{\mbox{$
0 \to H^{1} Y \to/{{ >}->}/^{H^{1}f} H^{1} X \to H^{0}\hom (K[\ab \G f],A) \to H^{2} Y \to^{H^{2}f} H^{2} X.
%$}}
\]
By Propositions~\ref{Proposition-H_0-Sequence} and~\ref{Proposition-First-Cohomology}, 
\[
H^{0}\hom (K[\ab \G f],A)\cong\hom (H_{0}K[\ab \G f], A)\cong\hom \bigl(\textstyle{\frac{K}{[K,X]}}, A\bigr),
\]
and the result follows.
\end{proof}

As a special case we get the following cohomological version of Hopf's formula.

\begin{corollary}\label{Corollary-CoHopf}\index{Hopf's Formula}
Let
\[
\xymatrix{ 0 \ar[r] & K \ar@{{ >}->}[r]^-{k} & X \ar@{-{ >>}}[r]^-{f} & Y \ar[r] & 0 }
\]
be a short exact sequence in $\Ac$, with $X$ a projective object. Then the sequence 
\[
\hom (X,A) \to<250> \hom \bigl(\textstyle{\frac{K}{[K,X]}},A\bigr) \to/{-{ >>}}/<250> H^{2} (Y,A) \to<250> 0
\]
is exact.
\end{corollary}
\begin{proof}
This follows immediately from the sequence in Theorem~\ref{Theorem-Cohomology-Sequence}, if we use Proposition~\ref{Proposition-First-Cohomology} asserting that $H^{2} X=0$ when $X$ is projective.
\end{proof}

This means that an element of $H^{2} (Y,A)$ may be considered as an equivalence class $[a]$ of morphisms $a\colon {K}/{[K,X]}\to A$, where $[a]=[0]$ if and only if $a$ extends to $X$.

\section{The second cohomology group}\label{Section-Second-Cohomology}
In any semi-abelian category, the \emph{Baer sum}\index{Baer sum} of two central extensions may be defined in the same way as it is classically done in the category of groups: This fact essentially depends on Proposition~\ref{Proposition-Subobject-of-Central-is-Normal}, as we are now going to explain. First, in any semi-abelian category, the validity of the Short Five Lemma allows one to speak of the isomorphism class $\{ f \}$ of a central extension $f$. Let us then recall explicitly the construction of the Baer sum of (the isomorphism classes of) two central extensions
\[
0 \to K \to/{ >}->/^{k} X \to/-{ >>}/^{f} Y \to 0
\]
and
\[
0 \to K \to/{ >}->/^{l} Z \to/-{ >>}/^{g} Y \to 0.
\]
One first forms the pullback $X\times_{Y}Z$ of $g$ along $f$; call 
\[
h=f\comp \pr_{1}=g\comp \pr_{2}\colon {X \times_{Y}Z\to Y}
\]
and $k\times l\colon K\times K\to X\times_{Y}Z$ its kernel. $h$ is a central extension because central extensions are pullback-stable; more precisely, $f\times g$ is a central extension and $h=\Delta^{-1}_{Y} (f\times g)$ (with $\Delta_{Y} \colon  {Y \to Y \times Y}$ the diagonal arrow). As an abelian object, $K$ carries a multiplication $m\colon {K\times K\to K}$, so that we may from the following diagram:
\[
\xymatrix{& 0 \ar[d] & 0 \ar[d]\\
& K[m] \ar@{=}[r] \ar@{{ >}->}[d]_{i=\ker m} & K \ar@{{ >}->}[d]^{(k\times l)\circ i} \ar[r] & 0 \ar[d]\\
0 \ar[r] & K\times K \ar@{{ >}->}[r]^-{k\times l} \ar@{}[rd]|>>{\pushout} \ar@{-{ >>}}[d]_{m} & X\times_{Y}Z \ar@{.{ >>}}[d]^-{n} \ar@{-{ >>}}[r]^-{h} & Y \ar[r] \ar@{:}[d] & 0\\
0\ar@{.>}[r] & K \ar@{{ >}.>}[r]_{\overline{k\times l}} \ar[d] & W \ar@{.>}[d] \ar@{.{ >>}}[r]_-{f+g} & Y \ar@{.>}[r] & 0.\\
& 0 & 0}
\]
The arrow $(k\times l)\comp i$ is a kernel, thanks to Proposition~\ref{Proposition-Subobject-of-Central-is-Normal} and the fact that $h$ is central. The argument given in the proof of Proposition~\ref{Proposition-Reg-Epi-Central-Extension} shows that the bottom sequence is a central extension; its isomorphism class, denoted $\{f \} + \{ g \}$, will be the Baer sum of the equivalence classes $\{f\}$ and $\{g\}$. In summary, one computes the Baer sum of two classes $\{f \}$ and $\{g \}$ as the isomorphism class of the cokernel $f+g$ of the pushout $\overline{k\times l}$ of the arrow $k\times l$ along the multiplication $m$ of $K$.

\begin{theorem}\label{Theorem-H^2-Central-Extensions}
Let $A$ be an abelian object of a semi-abelian category $\Ac$, monadic over $\Set$. The group $H^{2} (Y,A)$ is isomorphic to the group of isomorphism classes of central extensions of $Y$ by $A$, equipped with the Baer sum.
\end{theorem}
\begin{proof}
Let $f\colon X\to Y$ be a presentation of $Y$ with kernel $K$; consider the reflection 
\[
0\to \tfrac{K}{[K,X]} \to/{ >}->/ \tfrac{X}{[K,X]} \to/{-{ >>}}/^{\centr f} Y \to 0
\]
of $f$ into $\CentrExt (Y)$. In view of Corollary~\ref{Corollary-CoHopf} it suffices to show two things: (1) there is a bijection $F$ from the set of equivalence classes $[a]$ of morphisms $a\colon {{K}/{[K,X]}\to A}$, where $[a]=[0]$ if and only if $a$ extends to $X$, to the set of isomorphism classes of central extensions of $Y$ by $A$; (2) $F([a]+[b])= F([a]) + F([b])$ (where the addition on the right-hand side is the Baer sum of central extensions we recalled above).

The function $F$ is defined using Corollary~\ref{Corollary-Central-Extension-Pushout}: Because $\centr f$ is a central extension, a morphism $a\colon {K}/{[K,X]}\to A$ gives rise to a central extension of $Y$ by $A$---of which the isomorphism class $F ([a])$ is the image of $[a]$ through $F$. 

$F$ is well-defined: If $[a]=[0]$ then $F ([a])=F ([0])$. Indeed, it is easily seen that $F ([0])$ is the isomorphism class of the projection $\pr_{Y}\colon {Y\times A\to Y}$, a central extension. If $a\colon {K/[K,X]\to A}$ factors over $X$ then it factors over $X/[K,X]$, and as a consequence the extension associated to $a$ has a split monic kernel. It follows that this extension is isomorphic to $\pr_{Y}\colon {Y\times A\to Y}$.

Next, $F$ is a surjection because $X$ is projective and $\centr f$ is the reflection of $f$ into $\CentrExt (Y)$, and $F$ is injective because $F ([a])=\{\pi_{Y}\colon Y\times A\to Y\}$ entails that $a$ factors over~$X$.

Finally consider $a,b\colon K/[K,X]\to A$ for a proof of (2). The left hand side diagram
\[
\vcenter{\xymatrix{\tfrac{K}{[K,X]} \ar@{{ >}->}[r]^-{k} \ar[d]_-{l_{\tfrac{K}{[K,X]}}} & \tfrac{X}{[K,X]} \ar[d]^-{l_{\tfrac{X}{[K,X]}}}\\
\tfrac{K}{[K,X]}\oplus A \ar@{}[rd]|>>>{\pushout} \ar@{{ >}->}[r]_{k\times 1_{A}} \ar[d]_-{[m\circ (a,b),1_{A}]} & \tfrac{X}{[K,X]}\times A \ar[d]\\
A \ar@{{ >}->}[r]_-{k_{m\circ (a,b)}} & P}}
\qquad \qquad 
\vcenter{\xymatrix{
A\times A \ar@{}[rd]|>>{\pushout} \ar[d]_-{m} \ar@{{ >}->}[r]^-{k_{a}\times k_{b}} & Z_{a}\times_{Y}Z_{b}\ar[d]^{n}\\
A \ar@{{ >}->}[r]_{\overline{k_{a}\times k_{b}}} & W
}}
\]
shows the construction of $F ([a]+[b])$, the isomorphism class of a cokernel $p_{m\circ (a,b)}$ of $k_{m\circ (a,b)}$; the right hand side one, the construction of $F ([a])+F ([b])$, the isomorphism class of a cokernel $p_{a}+p_{b}$ of $\overline{k_{a}\times k_{b}}$. In order to prove that these two classes are equal, it suffices to give a morphism $g\colon {P\to W}$ that satisfies $(p_{a}+p_{b})\comp g=p_{m\circ (a,b)}$ and $g\comp k_{m\circ (a,b)}=\overline{k_{a}\times k_{b}}$. One checks that such a morphism is induced by the arrows 
\[
([a,0],[b,1_{A}])\colon \tfrac{K}{[K,X]}\oplus A\to A\times A
\]
and 
\[
(\overline{a}\comp \pr_{1},\varphi_{\overline{b},k_{b}})\colon \tfrac{X}{[K,X]}\times A\to Z_{a}\times_{Y}Z_{b},
\]
where $\varphi_{\overline{b},k_{b}}$ denotes the cooperator of $\overline{b}$ and $k_{b}$.
\end{proof}

\section{The Universal Coefficient Theorem}\label{Section-Universal-Coefficient-Theorem}

This section treats the relationship between homology and cohomology.

Recall that a central extension always has an abelian kernel; but given a central extension 
\[
0 \to A  \to/{{ >}->}/ B \to/{-{ >>}}/ C \to 0
\]
of an abelian object $C$ by an abelian object $A$, the object $B$ is not necessarily abelian. If it is, we say that the extension is \index{abelian!extension}\index{extension!abelian}\emph{abelian}---and the group of isomorphism classes of such extensions of $C$ by $A$ (a subobject of $H^{2} (C,A)$) will be denoted by $\ext (C,A)$\index{ext@$\ext (C,A)$}. Note that, because an extension in $\Ab \Ac$ is always central, an abelian extension is the same thing as an extension in $\Ab \Ac$.

Since the regular epimorphisms in $\Ab\Ac $ are just the regular epimorphisms of $\Ac$ that happen to lie in $\Ab\Ac$, the reflection $\ab (X)$ of a projective object $X$ of $\Ac$ is projective in $\Ab\Ac$. It follows that $\Ab\Ac$ has enough projectives if $\Ac$ has, and one may then choose a presentation
\begin{equation}\label{Abelian-Presentation}
0 \to R  \to/{{ >}->}/ F \to/{-{ >>}}/^{p} C \to 0
\end{equation}
of an abelian object $C$ in $\Ab\Ac$ instead of in $\Ac$.

\begin{proposition}\label{Proposition-Characterization-ext}
If $A$ is an abelian object and~\ref{Abelian-Presentation} is a presentation in $\Ab \Ac$ of an abelian object $C$, then the sequence
\[
\hom (F,A) \to \hom (R,A) \to/{-{ >>}}/ \ext (C,A) \to 0
\]
is exact.
\end{proposition}
\begin{proof}
This is an application of Corollary~\ref{Corollary-Central-Extension-Pushout}. It suffices to note that the arrow $p$ is central, and that in a square
\[
\xymatrix{R \ar[d]_-{a} \ar@{{ >}->}[r] & F \ar[d] \\
A \ar@{{ >}->}[r] & Z_{a}}
\]
induced by Corollary~\ref{Corollary-Central-Extension-Pushout}, all objects are abelian. Indeed, $Z_{a}$ being an abelian object follows from the fact that $\Ab (\Ac)$ is closed under products and regular quotients in $\Ac$, and that $a$ can be decomposed as $[a,1_{A}]\comp l_{R}$.
\end{proof}

\begin{theorem}[Universal Coefficient Theorem]\label{Theorem-Universal-Coefficients}\index{Universal Coefficient Theorem}
If $Y$ is an object of $\Ac$ and $A$ is abelian then the sequence
\begin{equation}\label{Sequence-UCT}
0 \to {\ext (H_{1}Y,A)}  \to/{{ >}->}/ H^{2} (Y, A) \to {\hom (H_{2}Y,A)}
\end{equation}
is exact.
\end{theorem}
\begin{proof}
One diagram says it all:
\[
\xymatrix{& {\hom (H_{1} X,A)} \ar@{=}[r] \ar[d] & {\hom (X,A)} \ar[d] \\
0 \ar[r] & {\hom \bigl(\tfrac{K}{K\cap [X,X]},A\bigr)} \ar@{{ >}->}[r] \ar@{-{ >>}}[d] & {\hom \bigl(\tfrac{K}{[K,X]},A\bigr)} \ar@{-{ >>}}[d] \ar[r] & {\hom (H_{2}Y,A)} \ar@{=}[d]\\
0\ar@{.>}[r] & {\ext (H_{1}Y,A)} \ar@{{ >}.>}[r] \ar[d] & H^{2} (Y,A) \ar@{.>}[r] \ar[d] & {\hom (H_{2}Y,A).}\\
& 0 & 0}
\]
Here $f\colon X\to Y$ is a presentation with kernel $K$, and the vertical sequences are exact by Proposition~\ref{Proposition-Characterization-ext}---the sequence
\[
0 \to \tfrac{K}{K\cap [X,X]}  \to/{{ >}->}/ H_{1}X \to/{-{ >>}}/ H_{1}Y \to 0
\]
being a presentation of $H_{1}Y$---and Corollary~\ref{Corollary-CoHopf}, respectively. The middle horizontal sequence is exact by the Hopf formula (Theorem~\ref{Theorem-Homology}) and the fact that, by the First Noether Isomorphism Theorem~\ref{generalexact},
\begin{equation}\label{Exact-Sequence-Sometimes-Split}
0 \to \tfrac{K\cap [X,X]}{[K,X]} \to/{{ >}->}/ \tfrac{K}{[K,X]} \to/{-{ >>}}/ \tfrac{K}{K\cap [X,X]} \to 0
\end{equation}
is an exact sequence in $\Ab\Ac$.  
\end{proof}

Recalling that an object $Y$ is perfect if and only if $H_{1}Y=0$, Theorem~\ref{Theorem-Universal-Coefficients} yields the following classical result.

\begin{corollary}\label{Corollary-Perfect-Homology-vs-Cohomology}
If $Y$ is a perfect object and $A$ is abelian then $H^{2} (Y,A)\cong\hom (H_{2}Y,A)$.
\end{corollary}
\begin{proof}
Comparing the Stallings-Stammbach Sequence~\ref{Stallings-Stammbach-Sequence} with Seq.~\ref{Exact-Sequence-Sometimes-Split} and using that $Y$ is perfect we see that $K/ (K\cap [X,X])$ is isomorphic to $\ab (X)$. The latter object being projective in $\Ab\Ac$, Sequence~\ref{Exact-Sequence-Sometimes-Split} is split exact in the abelian category $\Ab\Ac$. It follows that 
\[
0 \to {\ext (H_{1}Y,A)}  \to/{{ >}->}/ H^{2} (Y, A) \to/-{ >>}/ {\hom (H_{2}Y,A)}\to 0
\]
is a split exact sequence; but because $Y$ is perfect, $\ext (H_{1}Y,A)$ is zero.
\end{proof}

Given an object $Y$ and a presentation $F\to H_{1}Y$ with kernel $R$ in $\Ab\Ac$, Proposition~\ref{Proposition-Characterization-ext} entails the exactness of the sequence 
\[
{\hom (F,A)} \to {\hom (R, A)} \to/{-{ >>}}/ {\ext (H_{1}Y,A)} \to 0.
\]
If now, for every abelian object $A$ of $\Ac$, $\ext (H_{1}Y,A)=0$, then all induced maps ${\hom (F,A)\to \hom (R,A)}$ are epi, which means that $R\to F$ is a split monomorphism. It follows that $F=R\oplus H_{1}Y$ and $H_{1}Y$ is projective. As a consequence, we get the following partial converse to Corollary~\ref{Corollary-Perfect-Homology-vs-Cohomology}.

\begin{corollary}\label{Corollary-Perfect-Homology-vs-Cohomology-2}
If, for every $A$ abelian in $\Ac$, $H^{2} (Y,A)\cong\hom (H_{2}Y,A)$, then $H_{1}Y$ is projective in $\Ab\Ac$.\noproof 
\end{corollary}

Note that, in the case of groups~\cite{Robinson, Weibel} or crossed modules~\cite{Carrasco-Homology}, Sequence~\ref{Sequence-UCT} is also exact on the right: Then ${ H^{2} (Y, A) \to \hom (H_{2}Y,A)}$ is a regular epimorphism. This is due to homological dimension\index{homological dimension} properties specific to the categories of groups and crossed modules. In the case of groups, using the Nielsen-Schreier Theorem\index{Nielsen-Schreier Theorem}---a subgroup of a free group is free---one easily modifies the proof of Theorem~\ref{Theorem-Universal-Coefficients} to get the right exactness of the sequence. A similar argument is used in~\cite{Carrasco-Homology} to get the result for crossed modules. There is no hope of generalizing these kinds of arguments to arbitrary semi-abelian categories.
 % Cohomology
\chapter{Homotopy of internal categories}\label{Chapter-Internal-Categories}

\setcounter{section}{-1}
\section{Introduction}\label{Section-Introduction-Intcats}
It is well-known that the following choices of morphisms define a Quillen model category~\cite{Quillen} structure---known as the ``folk'' structure---on the category $\Cat$ of small categories and functors between them: $\we$ is the class of equivalences of categories, $\cof$ the class of functors, injective on objects and $\fib$ the class of functors $p\colon E\to B$ such that for any object $e$ of $E$ and any isomorphism $\beta \colon b\to p (e)$ in $B$ there exists an isomorphism $\epsilon$ with codomain $e$ such that $p (\epsilon)=\beta$; this notion was introduced for groupoids by R. Brown in~\cite{Brown:Fibrations}. We are unaware of who first proved this fact; certainly, it is a special case of Joyal and Tierney's structure~\cite{Joyal-Tierney}, but it was probably known before. A very explicit proof may be found in an unpublished paper by Rezk~\cite{Rezk}.

Other approaches to model category structures on $\Cat $ exist: Golasi{\'n}ski uses the homotopy theory of cubical sets to define a model structure on the category of pro-objects in $\Cat$~\cite{Golasinski}; Thomason uses an adjunction to simplicial sets  to acquire a model structure on $\Cat$ itself~\cite{Thomason}. Both are very different from the folk structure. Related work includes folk-style model category structures on categories of $2$-categories and bicategories (Lack~\cite{Lack:two-categories},~\cite{Lack:bicategories}) and a Thomason-style model category structure for $2$-categories (Worytkiewicz, Hess, Parent and Tonks~\cite{WHPT}).

If $\Ec$ is a Grothendieck topos there are two model structures on the category $\Cat \Ec$ of internal categories in $\Ec$. On can define the cofibrations and weak equivalences ``as in $\Cat$'', and then define the fibrations via a right lifting property. This gives Joyal and Tierney's model structure~\cite{Joyal-Tierney}. Alternatively one can define the fibrations and weak equivalences ``as in $\Cat$'' and than define the cofibrations via a left lifting property. This gives the model structure in this chapter. The two structures coincide when every object is projective, as in the case $\Ec =\Set$.

More generally, if $\Cc$ is a full subcategory of $\Ec$, one gets a full embedding of $\Cat \Cc$ into $\Cat \Ec$, and one can then define the weak equivalences and fibrations in $\Cat \Cc$ ``as in $\Cat \Ec$'', and the cofibrations via a left lifting property. In particular one can do this when $\Ec =\Sh (\Cc ,\Tc)$, for a subcanonical Grothendieck topology $\Tc$ on an arbitrary category $\Cc$. Starting with such a $\Cc$, one may also view this as follows: The notions of fibration and weak equivalence in the folk structure may be internalized, provided that one specifies what is meant by essential surjectivity and the existence claim in the definition of fibration. Both of them require some notion of surjection; this will be provided by a topology $\Tc$ on $\Cc$.

There are three main obstructions on a site $(\Cc,\Tc )$ for such a model category structure to exist. First of all, by definition, a model category has finite colimits. We give some sufficient conditions on $\Cc$ for $\Cat \Cc$ to be finitely cocomplete: either $\Cc$ is a topos with natural numbers object; or it is a locally finitely presentable category; or it is a finitely cocomplete regular Mal'tsev category. Next, in a model category, the class $\we$ of weak equivalences has the \textit{two-out-of-three} property. This means that if two out of three morphisms $f$, $g$, $g\comp f$ belong to $\we$ then the third also belongs to $\we$. A sufficient condition for this to be the case is that $\Tc$ is subcanonical. Finally, we want $\Tc$ to induce a weak factorization system in the following way. Let $Y_{\Tc}\colon {\Cc \to \Sh (\Cc ,\Tc)}$ denote the composite of the Yoneda embedding with the sheafification functor. A morphism $p\colon E\to B$ in $\Cc$ will be called a \textit{$\Tc$-epimorphism} if $Y_{\Tc} (p)$ is an epimorphism in $\Sh (\Cc ,\Tc)$. The class of $\Tc$-epimorphisms is denoted by $\Ec_{\Tc}$. If $({^{\Box}\Ec_{\Tc}}, \Ec_{\Tc})$ forms a weak factorization system, we call it \textit{the weak factorization system induced by $\Tc$}. This is the case when $\Cc$ has enough $\Ec_{\Tc}$-projective objects.

Joyal and Tierney's model structure~\cite{Joyal-Tierney} is defined as follows. Let $(\Cc,\Tc)$ be a site and $\Sh (\Cc, \Tc)$ its category of sheaves. Then a weak equivalence in $\Cat \Sh (\Cc ,\Tc) $ is a \textit{weak equivalence} of internal categories in the sense of Bunge and Par\'e~\cite{Bunge-Pare}; a cofibration is a functor, monic on objects; and a fibration has the right lifting property with respect to trivial cofibrations. Using the functor $Y_{\Tc}$ we could try to transport Joyal and Tierney's model structure from $\Cat \Sh (\Cc ,\Tc)$ to $\Cc$ as follows. For a subcanonical topology $\Tc$, the Yoneda embedding, considered as a functor ${\Cc \to \Sh (\Cc ,\Tc)}$, is equal to $Y_{\Tc}$. It follows that $Y_{\Tc}$ is full and faithful and preserves and reflects limits. Hence it induces a $2$-functor $\Cat Y_{\Tc}\colon  {\Cat \Cc \to \Cat \Sh (\Cc ,\Tc)}$. Say that an internal functor $\vetf \colon {\Af \to \Bf}$ is an equivalence or cofibration, resp., if and only if so is the induced functor $\Cat Y_{\Tc} (\vetf)$ in $\Sh (\Cc ,\Tc)$, and define fibrations using the right lifting property. 

We shall, however, consider a different structure on $\Cat \Cc$, mainly because of its application in the semi-abelian context. The weak equivalences, called \textit{$\Tc$-equivalences}, are the ones described above. (As a consequence, in the case of a Grothendieck topos, we get a structure that is different from Joyal and Tierney's, but has an equivalent homotopy category.) Where Joyal and Tierney internalize the notion of cofibration, we do so for the fibrations: $\vetp\colon \Ef \to \Bf$ is called a \textit{$\Tc$-fibration} if and only if in the diagram
\[
\xymatrix{\iso (E) \ar@/_/[rdd]_-{\delta_{1}} \ar@/^/[rrd]^-{\iso (p)_{1}} \ar@{.>}[rd]|{(\vetr_{\vetp})_{0} }\\
& (P_{p})_{0} \ar@{}[rd]|<{\pullback} \ar[r]^-{\overline{p_{0}}} \ar[d]_-{\overline{\delta_{1}}} & \iso(B) \ar[d]^-{\delta_{1}}\\
& E_{0} \ar[r]_-{p_{0}} & B_{0}}
\]
where $\iso (E)$ denotes the object of invertible arrows in the category $\Ef$, the induced universal arrow $(r_{p})_{0}$ is in $\Ec_{\Tc}$. \textit{$\Tc$-cofibrations} are defined using the left lifting property.

The chapter is organized as follows. In Section~\ref{Section-Cocylinder}, we study a cocylinder on $\Cat \Cc$ that characterizes homotopy of internal categories, i.e., such that two internal functors are homotopic if and only if they are naturally isomorphic. This cocylinder is used in Section~\ref{Section-Equivalences} where we study the notion of internal equivalence, relative to the Grothendieck topology $\Tc$ on $\Cc$ defined above. For the trivial topology (the smallest one), a $\Tc$-equivalence is a \textit{strong equivalence}, i.e., a homotopy equivalence with respect to the cocylinder. We recall that the strong equivalences are exactly the adjoint equivalences in the $2$-category $\Cat \Cc$. If $\Tc$ is the regular epimorphism topology (generated by covering families consisting of a single pullback-stable regular epimorphism), $\Tc$-equivalences are the so-called \textit{weak equivalences}~\cite{Bunge-Pare}. There is no topology $\Tc$ on $\Set$ for which the $\Tc$-equivalences are the equivalences of Thomason's model structure on $\Cat$: Any adjoint is an equivalence in the latter sense, whereas a $\Tc$-equivalence is always fully faithful. 

In Section~\ref{Section-Main-Theory} we study $\Tc$-fibrations. We prove---this is Theorem~\ref{Theorem-Model-Structure}---that the $\Tc$-e\-qui\-va\-len\-ces form the class $\we (\Tc)$ and the $\Tc$-fibrations the class $\fib (\Tc)$ of a model category structure on $\Cat \Cc$, as soon as the three obstructions mentioned above are taken into account.

Two special cases are subject to a more detailed study: in Section~\ref{Section-Regular-Epimorphisms}, the model structure induced by the regular epimorphism topology; in Section~\ref{Section-Split-Epimorphisms}, the one induced by the trivial topology.
In the first case we give special attention to the situation where $\Cc$ is a semi-abelian category, because then weak equivalences turn out to be homology isomorphisms, and the fibrations, Kan fibrations. Moreover, the category of internal categories in a semi-abelian category $\Cc$ is equivalent to Janelidze's category of internal crossed modules in $\Cc$~\cite{Janelidze}. Reformulating the model structure in terms of internal crossed modules (as is done in Theorem~\ref{Theorem-Categorical-Model-Structure}) simplifies its description. If $\Cc$ is the category of groups and homomorphisms, we obtain the model structures on the category $\Cat \Gp$ of categorical groups and the category $\XMod$ of crossed modules of groups, as described by Garz\'on and Miranda in~\cite{Garzon-Miranda}. 

The second case models the situation in $\Cat$, equipped with the folk model structure, in the sense that here, weak equivalences are homotopy equivalences, fibrations have the homotopy lifting property (Proposition~\ref{Proposition-Fibration-is-Star-Surjective}) and cofibrations the homotopy extension property (Proposition~\ref{Proposition-LLP}) with respect to the cocylinder defined in Section~\ref{Section-Cocylinder}.

This chapter almost entirely coincides with my paper \cite{EKVdL} with Tomas Everaert and Rudger Kieboom; Proposition~\ref{Proposition-Fundamental-Groupoid} is new.

\section{Preliminaries}\label{Section-Preliminaries}

\subsection*{Internal categories and groupoids}\label{Subsection-Internal-Categories}
If $\Cc$ is a finitely complete category then $\Rgrph \Cc$\index{category!$\Rgrph \Ac$} (resp. $\Cat \Cc$\index{category!$\Cat \Ac$}, $\Gd \Cc$\index{category!$\Gd \Ac$}) denotes the category of internal reflexive graphs\index{internal!reflexive graph} (resp. categories\index{internal!category}, groupoids\index{internal!groupoid}) in $\Cc$. Let 
\[
\Gd \Cc \to^{J} \Cat \Cc \to^{I} \Rgrph \Cc 
\]
denote the forgetful functors. It is well-known that $J$ embeds $\Gd \Cc$ into $\Cat \Cc$ as a coreflective subcategory. Carboni, Pedicchio and Pirovano prove in~\cite{Carboni-Pedicchio-Pirovano} that, if $\Cc$ is Mal'tsev, then $I$ is full, and $J$ is an isomorphism. Moreover, an internal reflexive graph carries at most one structure of internal groupoid; hence $\Gd\Cc$ may be viewed as a subcategory of  $\Rgrph \Cc$. As soon as $\Cc$ is, moreover, finitely cocomplete and regular, this subcategory is reflective (see Borceux and Bourn~\cite[Theorem 2.8.13]{Borceux-Bourn}). This extends the result of M.\,C.~Pedicchio~\cite{Pedicchio} that, if $\Cc$ is an exact Mal'tsev category with coequalizers, then the category $\Gd \Cc$ is regular epi-reflective in $\Rgrph \Cc$. In any case, this implies that $\Gd \Cc$ is closed in $\Rgrph \Cc$ under subobjects. In~\cite{Gran:Internal}, Gran adds to this result that $\Cat \Cc$ is closed in $\Rgrph \Cc$ under quotients. It follows that $\Cat \Cc$ is Birkhoff~\cite{Janelidze-Kelly} in $\Rgrph \Cc$. This, in turn, implies that if $\Cc$ is semi-abelian, so is $\Cat \Cc$ (Corollary~\ref{Corollary-Birkhoff-Subcategory-is-Semiabelian}). Gran and Rosick\'y~\cite{Gran-Rosicky} extend these results to the context of modular varieties. For any variety $\Vc$, the category $\Rgrph \Vc$ is equivalent to a variety. They show that, if, moreover, $\Vc$ is modular, $\Vc$ is Mal'tsev if and only if $\Gd \Vc$ is a subvariety of $\Rgrph \Vc$~\cite[Proposition 2.3]{Gran-Rosicky}.

Let $\Cc$ be finitely complete. Sending an internal category
\[
\Af=\bigl(\xymatrix{{A_{1}\times_{A_{0}}A_{1}} \ar[r]^-{m} & A_{1} \ar@<-1 ex>[r]_-{d_{0}} \ar@<1 ex>[r]^-{d_{1}} & A_{0} \ar[l]|-{i}}\bigr)
\]
to its object of objects $A_{0}$ and an internal functor $\vetf= (f_{0},f_{1})\colon \Af \to \Bf$ to its object morphism $f_{0}$ defines a functor $(\cdot)_{0}\colon \Cat \Cc \to \Cc$. Here $A_{1}\times_{A_{0}}A_{1}$ denotes a pullback of $d_{1}$ along $d_{0}$; by convention, $d_{1}\comp\pr_{1}=d_{0}\comp\pr_{2}$. It is easily seen that $(\cdot)_{0}$ has both a left and a right adjoint, respectively denoted $L$ and $R\colon {\Cc \to \Cat \Cc}$. Given an object $X$ of $\Cc$ and an internal category $\Af$, the natural bijection $\psi \colon {\Cc (X,A_{0})\to \Cat\Cc (L (X),\Af)}$ maps a morphism $f_{0}\colon {X\to A_{0}}$ to the internal functor 
\[
\vetf=\psi (f_{0}) = (f_{0},i\comp f_{0})\colon \Xf=L (X) \to \Af,
\]
where $\Xf$ is the \emph{discrete}\index{discrete internal groupoid}\index{internal!groupoid!discrete} internal groupoid $d_{0}=d_{1}=i=m=1_{X}\colon X\to X$.

The right adjoint $R$ maps an object $X$ to the \emph{indiscrete}\index{indiscrete internal groupoid}\index{internal!groupoid!indiscrete} groupoid $R (X)$ on $X$, i.e., $R (X)_{0}=X$, $R (X)_{1}=X\times X$, $d_{0}$ is the first and $d_{1}$ the second projection, $i$ is the diagonal and $m\colon R (X)_{1}\times_{R (X)_{0}}R (X)_{1}\to R (X)_{1}$ is the projection on the first and third factor.

Sending an internal category $\Af$ to its object of arrows $A_{1}$ defines a functor $(\cdot)_{1}\colon {\Cat \Cc \to \Cc}$. Since limits in $\Cat \Cc$ are constructed by first taking the limit in $\Rgrph \Cc$, then equipping the resulting reflexive graph with the unique category structure such that the universal cone in $\Rgrph \Cc$ becomes a universal cone in $\Cat \Cc$, the functor $I\colon \Cat \Cc \to \Rgrph \Cc$ creates limits. Hence $(\cdot)_{1}$ is limit-preserving. 

\subsection*{When is $\Cat \Cc$ (co)complete?}\label{Subsection-Cocompleteness}
One of the requirements for a category to be a model category is that it is finitely complete and cocomplete. Certainly the completeness poses no problems since it is a pretty obvious fact that $\Cat \Cc$ has all limits $\Cc$ has (see e.g., Johnstone~\cite[Lemma 2.16]{Johnstone:Topos-Theory}); hence $\Cat \Cc$ is always finitely complete.

The case of cocompleteness is entirely different, because in general cocompleteness of $\Cc$ need not imply the existence of colimits in $\Cat \Cc$. (Conversely, $\Cc$ has all colimits $\Cat \Cc$ has, because $(\cdot)_{0}\colon \Cat \Cc \to \Cc$ has a right adjoint.) As far as we know, no characterization exists of those categories $\Cc$ that have a finitely cocomplete $\Cat \Cc$; we can only give sufficient conditions for this to be the case.

We get a first class of examples by assuming that $\Cc$ is a topos\index{topos}\index{category!topos} with a natural numbers object\index{natural numbers object} (NNO\index{NNO}) (or, in particular, a Grothendieck topos, like Joyal and Tierney do in~\cite{Joyal-Tierney}). It is well-known that, for a topos $\Cc$, the existence of a NNO is equivalent to $\Cat \Cc$ being finitely cocomplete. Certainly, if $\Cat \Cc$ has \textit{countable} coproducts, then so has $\Cc$, hence it has a NNO: Take a countable coproduct of $1$. But the situation is much worse, because $\Cat \Cc$ does not even have arbitrary \textit{coequalizers} if $\Cc$ lacks a NNO. Considering the ordinals $1$ and $2=1+1$ (equipped with the appropriate order) as internal categories $\matheu{1}$ and $\matheu{2}$, the coproduct inclusions induce two functors $\matheu{1}\to \matheu{2}$. If their coequalizer in $\Cat \Cc$ exists, it is the free internal monoid on $1$, considered as a one-object category (its object of objects is equal to $1$). But by Remark~D5.3.4 in~\cite{Johnstone:Elephant}, this implies that $\Cc$ has a NNO! Conversely, in a topos with NNO, mimicking the construction in $\Set$, the functor $I$ may be seen to have a left adjoint; using this left adjoint, we may construct arbitrary finite colimits in $\Cat \Cc$.

Locally finitely presentable categories\index{category!locally finitely presentable}\index{category!l.f.p.} form a second class of examples. Indeed, every l.f.p.\ category is cocomplete, and if a category $\Cc$ if l.f.p., then so is $\Cat \Cc$---being a category of models of a sketch with finite diagrams~\cite[Proposition 1.53]{Adamek-Rosicky}.  (Note that in particular, we again find the example of Grothendieck topoi.)

A third class is given by supposing that $\Cc$ is finitely cocomplete and regular Mal'tsev. Then $\Cat \Cc=\Gd \Cc$ is a reflective~\cite[Theorem 2.8.13]{Borceux-Bourn} subcategory of the functor category $\Rgrph \Cc$, and hence has all finite colimits. This class, in a way, dualizes the first one, because the dual of any topos is a finitely cocomplete (exact) Mal'tsev category~\cite{Carboni-Kelly-Pedicchio},~\cite[Example A.5.17]{Borceux-Bourn},~\cite{Bourn:Dual-topos}.

\subsection*{Weak factorization systems and model categories}\label{Subsection-Weak-Factorization-Systems}

In this paper we use the definition of model category as presented by Ad\'amek, Herrlich, Rosick\'y and Tholen~\cite{AHRT}. For us, next to its elegance, the advantage over Quillen's original definition~\cite{Quillen} is its explicit use of weak factorization systems. We briefly recall some important definitions.

\begin{definition}\label{Definition-Lifting}
Let $l\colon A\to B$ and $r\colon C\to D$ be two morphisms of a category $\Cc$. $l$ is said to have the \emph{left lifting property}\index{left lifting property} with respect to $r$ and $r$ is said to have the \emph{right lifting property}\index{right lifting property} with respect to $l$ if every commutative diagram
\[\xymatrix{A \ar[r] \ar[d]_l & C \ar[d]^{r} \\ B \ar[r] \ar@{.>}[ru]^{h} & D}\]
has a lifting $h\colon B\to C$. This situation is denoted $l\Box r$\index{$l\Box r$} (and has nothing to do with double equivalence relations).

If $\Hc$ is a class of morphisms then $\Hc^{\Box}$ is the class of all morphisms $r$ with $h\Box r$ for all $h\in \Hc$; dually, ${^{\Box}\Hc}$ is the class of all morphisms $l$ with $l\Box h$ for all $h\in \Hc$.
\end{definition}

\begin{definition}\label{Definition-Weak-Factorization-System}
A \emph{weak factorization system}\index{weak factorization system} in $\Cc$ is a pair $(\Lc,\Rc)$ of classes of morphisms such that 
\begin{enumerate}
\item every morphism $f$ has a factorization $f=r\comp l$ with $r\in \Rc$ and $l\in \Lc$;
\item $\Lc^{\Box} =\Rc$ and $\Lc ={^{\Box}\Rc}$.
\end{enumerate}
\end{definition}

In the presence of condition~1., 2.\ is equivalent to the conjunction of $\Lc \Box \Rc$ and the closedness in the category of arrows $\Fun (\Two ,\Cc)$ of $\Lc$ and $\Rc$ under the formation of retracts.

\begin{definition}\cite[Remark~3.6]{AHRT}\label{Definition-Model-Category}
Let $\Cc$ be a finitely complete and cocomplete category. A~\emph{model structure}\index{model category structure} on $\Cc$ is determined by three classes of morphisms, $\fib $\index{fib@$\fib$} (\emph{fibrations}\index{fibration}\index{fibration}), $\cof$\index{cof@$\cof$} (\emph{cofibrations}\index{cofibration}) and $\we$\index{we@$\we$} (\emph{weak equivalences}\index{weak equivalence}), such that
\begin{enumerate}
\item $\we$ has the \emph{2-out-of-3 property}\index{two-out-of-three property}, i.e., if two out of three morphisms $f$, $g$, $g\comp f$ belong to $\we$ then the third morphism also belongs to $\we$, and $\we$ is closed under retracts in $\Fun (\Two,\Cc)$;
\item $(\cof ,\fib \cap \we)$ and $(\cof \cap \we ,\fib)$ are weak factorization systems.
\end{enumerate}
A category equipped with a model structure is called a \emph{(Quillen) model category}\index{Quillen model category}. A morphism in $\fib \cap \we$ (resp. $\cof \cap \we$) is called a \emph{trivial fibration}\index{trivial fibration}\index{fibration!trivial} (resp. \emph{trivial cofibration}\index{trivial cofibration}\index{cofibration!trivial}). Let $0$ denote an initial and $1$ a terminal object of $\Cc$. A \emph{cofibrant}\index{cofibrant object} object $A$ is such that the unique arrow $0\to A$ is a cofibration; $A$ is called \emph{fibrant}\index{fibrant object} if $A\to 1$ is in $\fib$. 
\end{definition}

\subsection*{Grothendieck topologies}\label{Subsection-Grothendieck-Topologies}
We shall consider model category structures on $\Cat \Cc$ that are defined relative to some \emph{Grothendieck topology}\index{Grothendieck topology}\index{topology!in the sense of Grothendieck}\index{topology} $\Tc$ on $\Cc$. Recall that such is a function that assigns to each object $C$ of $\Cc$ a collection $\Tc (C)$ of sieves on $C$ (a \emph{sieve}\index{sieve} $S$ on $C$ being a class of morphisms with codomain $C$ such that $f\in S$ implies that $f\comp g\in S$, whenever this composite exists), satisfying
\begin{enumerate}
\item the maximal sieve on $C$ is in $\Tc (C)$;
\item (stability axiom) if $S\in \Tc (C)$ then its pullback $h^{*} (S)$ along any arrow $h\colon {D\to C}$ is in $\Tc (D)$;
\item (transitivity axiom) if $S\in \Tc (C)$ and $R$ is a sieve on $\Cc$ such that $h^{*} (R)\in \Tc (D)$ for all $h\colon D\to C$ in $S$, then $R\in \Tc (C)$.
\end{enumerate}
A sieve in some $\Tc (C)$ is called \emph{covering}\index{covering sieve}\index{sieve!covering}. We would like to consider sheaves over arbitrary sites $(\Cc ,\Tc)$, not just small ones (i.e., where $\Cc$ is a small category). For this to work flawlessly, a standard solution is to use the theory of universes, as introduced in~\cite{SGA4}. The idea is to extend the Zermelo-Fraenkel axioms of set theory with the axiom (U) ``every set is an element of a universe'', where a \emph{universe}\index{universe} $\Uc$ is a set satisfying
\begin{enumerate}
\item if $x\in \Uc$ and $y\in x$ then $y\in \Uc$;
\item if $x,y\in \Uc$ then $\{x,y \}\in \Uc$;
\item if $x\in \Uc$ then the power-set $\Pc (x)$ of $x$ is in $\Uc$;
\item if $I\in \Uc$ and $(x_{i})_{i\in I}$ is a family of elements of $\Uc$ then $\bigcup_{i\in I}x_{i}\in \Uc$.
\end{enumerate}
A set is called \emph{$\Uc$-small}\index{small set}\index{U-small set@$\Uc$-small set} if it has the same cardinality as an element of $\Uc$. (We sometimes, informally, use the word \emph{class}\index{class} for a set that is not $\Uc$-small.) We shall always consider universes containing the set $\N$ of natural numbers, and work in ZFCU (with the ZF axioms\index{ZF axioms} + the axiom of choice + the universe axiom). A category consists of a set of objects and a set of arrows with the usual structure; $\Uc\Set$ ($\Uc \Cat$) denotes the category whose objects are elements of $\Uc$ (categories with sets of objects and arrows in $\Uc$) and whose arrows are functions (functors) between them. Now given a site $(\Cc,\Tc )$, the category $\Cc$ is in $\Uc \Cat$ for some universe $\Uc$; hence it makes sense to consider the category of presheaves $\Uc \Psh \Cc=\Fun (\Cc^{\op},\Uc \Set)$\index{presheaf} and the associated category $\Uc \Sh (\Cc ,\Tc)$ of sheaves\index{sheaf}. In what follows, we shall not mention the universe $\Uc$ we are working with and just write $\Set$\index{category!$\Set$}, $\Cat$\index{category!$\Cat$}, $\Psh \Cc$\index{category!$\Psh \Cc$}, $\Sh (\Cc,\Tc)$\index{category!$\Sh (\Tc ,\Cc)$}, etc.

\begin{examples}\label{Examples-Topologies}
On a finitely complete category $\Cc$, the \emph{regular epimorphism topology}\index{regular epimorphism topology}\index{topology!regular epimorphism ---} is generated by the following basis: A covering family on an object $A$ consists of a single pullback-stable regular epimorphism $A'\to A$. It is easily seen that this topology is \emph{subcanonical}\index{subcanonical topology}\index{topology!subcanonical}, i.e., that every representable functor is a sheaf. Hence the Yoneda embedding $Y\colon \Cc \to \Psh \Cc $ may be considered as a functor $\Cc \to \Sh (\Cc ,\Tc)$.

The \emph{trivial topology}\index{topology!trivial}\index{trivial topology} is the smallest one: The only covering sieve on an object $A$ is the sieve of all morphisms with codomain $A$. Every presheaf if a sheaf for the trivial topology.

The largest topology is called \emph{cotrivial}\index{cotrivial topology}: Every sieve is covering. The only sheaf for this topology is the terminal presheaf.
\end{examples}

We shall consider the weak factorization system on a category $\Cc$, generated by a Grothendieck topology in the following way.

\begin{definition}\label{Definition-T-Epimorphism}
Let $\Tc$ be a topology on $\Cc$ and let $Y_{\Tc}\colon {\Cc \to \Sh (\Cc ,\Tc)}$ denote the composite of the Yoneda embedding\index{Yoneda embedding} $Y\colon {\Cc \to \Psh \Cc}$ with the sheafification functor\index{sheafification functor} ${\Psh  \Cc \to \Sh (\Cc ,\Tc)}$. A morphism $p\colon E\to B$ will be called a \emph{$\Tc$-epimorphism}\index{T-epimorphism@$\Tc$-epimorphism} if $Y_{\Tc} (p)$ is an epimorphism in $\Sh (\Cc ,\Tc)$. The class of $\Tc$-epimorphisms is denoted by $\Ec_{\Tc}$. If $({^{\Box}\Ec_{\Tc}}, \Ec_{\Tc})$ forms a weak factorization system, we call it \emph{the weak factorization system induced by $\Tc$}\index{weak factorization system!induced by a topology}.
\end{definition}

\begin{remark}\label{Remark-Subcanonical}
Note that if $\Tc$ is subcanonical, then $Y_{\Tc}$ is equal to the Yoneda embedding; hence it is a full and faithful functor.
\end{remark}

\begin{remark}\label{Remark-Enough-Projectives}
The only condition a subcanonical $\Tc$ needs to fulfil, for it to induce a model structure on $\Cat \Cc$, is that $({^{\Box}\Ec_{\Tc}}, \Ec_{\Tc})$ is a weak factorization system. When $\Cc$ has binary coproducts, this is equivalent to $\Cc$ having enough $\Ec_{\Tc}$-projectives~\cite{AHRT}.
\end{remark}

One way of avoiding universes is by avoiding sheaves: Indeed, $\Tc$-epi\-morph\-isms have a well-known characterization in terms of the topology alone.

\begin{proposition}\cite[Corollary~III.7.5 and~III.7.6]{Maclane-Moerdijk}\label{Proposition-Characterization-T-Epimorphism}
Let $\Tc$ be a topology on a category $\Cc$. Then a morphism $p\colon E\to B$ in $\Cc$ is $\Tc$-epic if and only if for every $g\colon X\to B$ there exists a covering family $(f_{i}\colon U_{i}\to X)_{i\in I}$ and a family of morphisms $(u_{i}\colon U_{i}\to E)_{i\in I}$ such that for every $i\in I$, $p\comp u_{i}=g\comp f_{i}$.\noproof 
\end{proposition}

\begin{examples}\label{Examples-Regular-and-Split-T-Epimorphisms}
If $\Tc$ is the trivial topology, it is easily seen that the $\Tc$-epimorphisms are exactly the split epimorphisms.

When $\Tc$ is the cotrivial topology, every morphism is $\Tc$-epic.

In case $\Tc$ is the regular epimorphism topology, a $\Tc$-epimorphism is nothing but a pullback-stable regular epimorphism: Certainly, every pullback-stable regular epimorphism is $\Tc$-epic; conversely, one shows that if $p\comp u=f$ is a pullback-stable regular epimorphism then so is $p$.
\end{examples}

\section{A cocylinder on $\protect\Cat \Cc$}\label{Section-Cocylinder}\index{homotopy!of internal functors}

One way of defining homotopy in a category $\Cc$ is relative to a \emph{cocylinder}\index{cocylinder} on $\Cc$. Recall (e.g., from Kamps~\cite{Kamps:Kan-Bedingungen} or Kamps and Porter~\cite{KP:Abstact-homotopy-theory}) that this is a structure
\[
((\cdot)^{I}\colon \Cc \to \Cc,\quad \epsilon_{0},\epsilon_{1}\colon (\cdot)^{I}\To 1_{\Cc},\quad s\colon 1_{\Cc}\To (\cdot)^{I})
\]
such that $\epsilon_{0}\vcomp s=\epsilon_{1}\vcomp s=1_{1_{\Cc}}$. Given a cocylinder $((\cdot)^{I},\epsilon_{0},\epsilon_{1},s)$ on $\Cc$, two morphisms 
\[
f,g\colon {X\to Y}
\]
are called \emph{homotopic}\index{homotopy!w.r.t.\ a cocylinder} (or, more precisely, \emph{right homotopic}\index{right homotopic}\index{homotopy!right- ---}, to distinguish with the notion of \emph{left} homotopy defined using a cylinder) if there exists a morphism $H\colon {X\to Y^{I}}$ such that $\epsilon_{0} (Y)\comp H=f$ and $\epsilon_{1} (Y)\comp H=g$. The morphism $H$ is called a \emph{homotopy from $f$ to $g$} and the situation is denoted $H\colon f\simeq g$.

Let $\Cc$ be a finitely complete category. In this section, we describe a cocylinder on $\Cat \Cc$ such that two internal functors are homotopic if and only if they are naturally isomorphic. We follow the situation in $\Cat$ very closely. Let $\Ic$ denote the \emph{interval groupoid}\index{interval groupoid}, i.e., the category with two objects $\{0,1 \}$ and the following four arrows.
\[
\xymatrix@1{{0} \ar@(ul,dl)[]_{{1}_0} \ar@/^/[r]^{\tau} & {1} \ar@/^/[l]^{{\tau}^{-1}}
\ar@(dr,ur)[]_{{1}_1}}
\]
Then putting $\Cc^{\Ic}=\Fun (\Ic ,\Cc)$, the category of functors from $\Ic$ to $\Cc$, defines a cocylinder on $\Cat$. It is easily seen that an object of $\Cc^{\Ic}$, being a functor ${\Ic \to \Cc}$, is determined by the choice of an isomorphism in $\Cc$; a morphism of $\Cc^{\Ic}$, being a natural transformation $\mu \colon F\To G\colon \Ic \to \Cc$ between two such functors, is determined by a commutative square
\[
\xymatrix{F (0) \ar[r]^{\mu_{0}} \ar[d]_{F (\tau)}^{\cong} & G (0)\ar[d]_{\cong}^{G (\tau )} \\
F (1) \ar[r]_{\mu_{1}} & G (1)}
\]
in $\Cc$ with invertible downward-pointing arrows.

It is well-known that the category $\Gd \Cc$ is coreflective in $\Cat \Cc$; let 
\[
\iso \colon \Cat \Cc \to \Gd \Cc
\]
denote the right adjoint of the inclusion $J\colon {\Gd \Cc \to \Cat \Cc}$. Given a category $\Af$ in $\Cc$, the functor $\iso$ may be used to describe the object $\iso (A)$ of ``isomorphisms in $\Af$'' (cf.\ Bunge and Par\'e~\cite{Bunge-Pare}) as the object of arrows of $\iso (\Af)$\index{iso@$\iso (\Af)$}, the couniversal groupoid associated with $\Af$. The counit $\epsilon_{\Af}\colon {\iso (\Af)\to \Af}$ at $\Af$ is a monomorphism, and will be denoted
\[
\xymatrix{\iso (A)\ar@<-1 ex>[d]_-{\delta_{0}} \ar@<1 ex>[d]^-{\delta_{1}} \monr^-{j} & A_{1} \ar@<-1 ex>[d]_-{d_{0}} \ar@<1 ex>[d]^-{d_{1}} \\
A_{0} \ar[u]|{\iota} \ar@{=}[r] & A_{0}. \ar[u]|{i}}
\]

The object $A^{I}_{1}$ of ``commutative squares with invertible downward-pointing arrows in $\Af$'' is given by the the pullback
\[
\xymatrix{A^{I}_{1} \ar[r]^-{\pr_{2}} \ar[d]_-{\pr_{1}} \ar@{}[rd]|<<<{\pullback} & A_{1}\times_{A_{0}}\iso (A) \ar[d]^-{m\comp (1_{A_{1}}\times_{1_{A_{0}}}j)}\\
\iso (A)\times_{A_{0}}A_{1} \ar[r]_-{m\comp (j\times_{1_{A_{0}}}1_{A_{1}})} & A_{1}.}
\]
The unique morphism induced by a cone on this diagram, represented by 
\[
(f,g,h,k)\colon X \to \iso(A)\times A_{1}\times A_{1}\times \iso (A),
\]
will be denoted by
\[
\vierkant{f}{g}{h}{k}\colon X \to A_{1}^{I}.
\] 
Put $A^{I}_{0}=\iso (A)$. Horizontal composition
\[
\compo=\vierkant{\pr_{1}\comp\pr_{1}\comp \pr_{1}}{m\comp (\pr_{2}\comp \pr_{1}\times_{\delta_{1}}\pr_{2}\comp \pr_{1})}{m\comp (\pr_{1}\comp \pr_{2}\times_{\delta_{0}}\pr_{1}\comp \pr_{2})}{\pr_{2}\comp\pr_{2}\comp \pr_{2}},\qquad \id=\vierkant{1_{\iso(A)_{1}}}{i\comp \delta_{1}}{i\comp \delta_{0}}{1_{\iso(A)_{1}}},
\]
$\dom=\pr_{1}\comp \pr_{1}$ and $\cod =\pr_{2}\comp \pr_{2}$ now define an internal category
\[
\Af^{\If}=\bigl(\xymatrix{{A^{I}_{1}\times_{A_{0}^{I}}A_{1}^{I}} \ar[r]^-{\compo} & A^{I}_{1} \ar@<-1 ex>[r]_-{\dom} \ar@<1 ex>[r]^-{\cod} & A^{I}_{0} \ar[l]|-{\id}}\bigr).
\]  
Thus we get a functor $(\cdot)^{\If}\colon \Cat \Cc \to \Cat \Cc$\index{A I@$\Af^{\If}$}. Putting 
\[
\epsilon_{0} (\Af)= (\delta_{0}, \pr_{1}\comp \pr_{2})\colon \Af^{\If}\to \Af, \qquad \epsilon_{1} (\Af)= (\delta_{1}, \pr_{2}\comp \pr_{1})\colon \Af^{\If}\to \Af
\]
and $s (\Af)= (\iota,s (\Af)_{1})$ with
\[
s (\Af)_{1}=\vierkant{\iota \comp d_{0}}{1_{A_{1}}}{1_{A_{1}}}{\iota\comp d_{1}}\colon A_{1} \to A_{1}^{I}
\]
gives rise to natural transformations $\epsilon_{0}, \epsilon_{1}\colon {(\cdot)^{\If}\To 1_{\Cat \Cc}}$, $s\colon {1_{\Cat \Cc} \To (\cdot)^{\If}}$ such that $\epsilon _{0}\vcomp s=\epsilon_{1}\vcomp s=1_{1_{\Cat \Cc}}$.

Recall that, for internal functors $\vetf,\vetg\colon \Af \to \Bf$, an \emph{internal natural transformation} $\mu \colon \vetf\To \vetg$ is a morphism $\mu \colon A_{0}\to B_{1}$ such that $d_{0}\comp \mu =f_{0}$, $d_{1}\comp \mu =g_{0}$ and $m\comp (f_{1},\mu \comp d_{1})=m\comp (\mu \comp d_{0},g_{1})$. Categories, functors and natural transformations in a given category $\Cc$ form a $2$-category $\Cat \Cc$. For two internal natural transformations $\mu \colon \vetf\To \vetg$ and $\nu \colon \vetg\To \veth$, $\nu \vcomp \mu=m\comp (\nu, \mu)$ is their (vertical) composition; for $\mu\colon \vetf\To \vetg\colon  \Af \to \Bf$ and $\mu'\colon {\vetf'\To \vetg'\colon  \Bf \to \Cf}$, the (horizontal) composition is 
\[
\mu' \comp \mu=m\comp (\mu '\comp f_{0},g'_{1}\comp \mu)=m\comp (f'_{1}\comp \mu, \mu '\comp g_{0})\colon \vetf'\comp \vetf\To \vetg'\comp \vetg\colon \Af \to \Cf.
\]
An internal natural transformation $\mu\colon \vetf\To \vetg\colon \Af \to \Bf$ is an \emph{internal natural isomorphism}\index{internal!natural isomorphism} if and only if an internal natural transformation $\mu^{-1}\colon \vetg\To \vetf$ exists such that $\mu\vcomp\mu^{-1}=1_{\vetg}=i\comp g_{0}$ and $\mu^{-1}\vcomp\mu=1_{\vetf}=i\comp f_{0}$. Hence an internal natural isomorphism is nothing but an isomorphism in a hom-category $\Cat \Cc (\Af ,\Bf)$. Moreover, this is the case, exactly when $\mu$ factors over $j\colon \iso(A)\to A_{1}$. Note that, if $\Bf$ is a groupoid, and $\tw \colon B_1 \to B_{1}$ denotes its ``twisting isomorphism'', then $\mu^{-1}=\tw \comp \mu$.

\begin{example}\label{Example-Epsilon-Homotopy}
For every internal category $\Af$ of $\Cc$, the morphism 
\[
\vierkant{\iota\comp \delta_{0}}{j}{i\comp \delta_{0}}{1_{\iso(A)}}\colon \iso(A)=A^{I}_{0}\to A^{I}_{1}
\]
is a natural isomorphism $s (\Af)\comp \epsilon_{0} (\Af)\To 1_{\Af^{\If}}\colon \Af^{\If}\to \Af^{\If}$.
\end{example}

As expected (cf.\ Exercise 2.3 in Johnstone~\cite{Johnstone:Topos-Theory}):

\begin{proposition}\label{Proposition-Homotopy-vs-Natural-Transformation}
If $\mu \colon \vetf\To \vetg\colon \Af \to \Bf$ is an internal natural isomorphism, then $\vetH= (\mu , H_{1})\colon \Af \to \Bf^{\If}$ with 
\[
H_{1}=\vierkant{\mu \comp d_{0}}{g_{1}}{f_{1}}{\mu \comp d_{1}}\colon A_{1}\to B^{I}_{1}
\]
is a homotopy $H\colon f\simeq g$. If $\vetH\colon \Af \to \Bf^{\If}$ is a homotopy $\vetf\simeq \vetg\colon {\Af \to \Bf}$ then 
\[
j\comp H_{0}\colon {A_{0}\to B_{1}}
\]
is an internal natural isomorphism ${\vetf\To \vetg}$. Hence the homotopy relation $\simeq$ is an equivalence relation on every $\Cat \Cc (\Af ,\Bf).$\noproof 
\end{proposition}

\begin{proposition}\label{Proposition-Double-Category}
For any internal category $\Af$ of $\Cc$, putting 
\[
d_{0}=\epsilon_{0} (\Af),\quad d_{1}=\epsilon_{1} (\Af)\colon {\Af^{\If}\to \Af}
\]
and $i=s (\Af)\colon {\Af \to \Af^{\If}}$ defines a reflexive graph in $\Cat \Cc$, which carries a structure of internal groupoid; hence it is a double category in $\Cc$.\noproof 
\end{proposition}

The following well-known construction will be very useful.

\begin{definition}[Mapping path space construction]\label{Definition-Mapping-Path-Space}\index{mapping path space construction}
Let $\vetf\colon \Af\to \Bf$ be an internal functor. Pulling back the split epimorphism $\epsilon_{1} (\Bf)$ along $\vetf$ yields the following diagram, where both the upward and downward pointing squares commute, and $\overline{\epsilon_{1} (\Bf)}\comp \overline{s (\Bf)}=1_{\Af}$.
\begin{equation}\label{Diagram-Mapping-Path-Space}
\vcenter{\xymatrix{{\Pf_{\vetf}}\ar@{}[rd]|<{\pullback} \ar[r]^-{\overline{\vetf}} \ar@<0.5 ex>[d]^>>{\overline{\epsilon_{1} (\Bf)}} & {\Bf^{\If}} \ar@<0.5 ex>[d]^-{\epsilon_{1} (\Bf)}\\
{\Af} \ar@<0.5 ex>[u]^-{\overline{s (\Bf)}} \ar[r]_-{\vetf} & {\Bf} \ar@<0.5 ex>[u]^>>>>{s (\Bf)}}}
\end{equation}
The object $\Pf_{\vetf}$ is called a \emph{mapping path space} of $\vetf$. We write 
\[
\vetr_{\vetf}\colon {\Af^{\If}\to \Pf_{\vetf}}
\]
for the universal arrow induced by the commutative square $\epsilon_{1} (\Bf)\comp \vetf^{\,\If}=\vetf\comp \epsilon_{1} (\Af)$.
\end{definition}

\section{$\Tc$-equivalences}\label{Section-Equivalences}

Let $\Cc$ be a finitely complete category. Recall (e.g., from Bunge and Par\'e~\cite{Bunge-Pare}) that an internal functor $\vetf\colon \Af \to \Bf$ in $\Cc $ is called \emph{full}\index{full functor} (resp. \emph{faithful}\index{faithful functor}, \emph{fully faithful}\index{fully faithful functor}) when, for any internal category $\Xf$ of $\Cc$, the functor 
\[
\Cat \Cc (\Xf ,\vetf)\colon \Cat \Cc (\Xf ,\Af )\to \Cat \Cc (\Xf,\Bf)
\]
is full (resp. faithful, fully faithful). There is the following well-known characterization of full and faithful functors.

\begin{proposition}\label{Proposition-Full-Faithful-Functor}
Let $\vetf\colon \Af \to \Bf$ be a functor in a finitely complete $\Cc$.
\begin{enumerate}
\item If $\vetf$ is full, then the square
\begin{equation}\label{Diagram-Full-Functor}
\vcenter{\xymatrix{A_{1}\ar[d]_{(d_{0},d_{1})} \ar[r]^{f_{1}} & B_{1}\ar[d]^{(d_{0},d_{1})}\\
A_{0}\times A_{0}\ar[r]_{f_{0}\times f_{0}}& B_{0}\times B_{0}}}
\end{equation}
is a weak pullback in $\Cc$.
\item $\vetf$ is faithful if and only if the morphisms $d_{0},d_{1}\colon A_{1}\to A_{0}$ together with $f_{1}\colon A_{1} \to B_{1}$ form a mono-source.
\item $\vetf$ is fully faithful if and only if~\ref{Diagram-Full-Functor} is a pullback.\noproof 
\end{enumerate}
\end{proposition}

\begin{remark}\label{Remark-Iso-Preserves-Full-Faithful-Functors}
Since fully faithful functors reflect isomorphisms, the Yoneda Lemma (e.g., in the form of Metatheorem 0.1.3 in~\cite{Borceux-Bourn}) implies that the functor $\iso \colon \Cat \Cc \to \Gd  \Cc$ preserves fully faithful internal functors. Quite obviously, they are also stable under pulling back. 
\end{remark}

The following lifting property of fully faithful functors will prove very useful (cf.\ the proof of Lemma~2.1 in Joyal and Tierney~\cite{Joyal-Tierney}).

\begin{proposition}\label{Proposition-Lifting-we}
Consider a commutative square
\begin{equation}\label{Diagram-Lifting}
\vcenter{\xymatrix{{\Af} \ar[d]_-{\vetj} \ar[r]^-{\vetf} & {\Ef }\ar[d]^-{\vetp} \\
{\Xf}\ar@{.>}[ru]^-{\veth} \ar[r]_-{\vetg} & {\Bf}}}
\end{equation}
in $\Cat \Cc$ with $\vetp$ fully faithful. This square has a lifting $\veth\colon \Xf \to \Ef$  if and only if there exists a morphism $h_{0}\colon X_{0}\to E_{0}$ such that $p_{0}\comp h_{0}=g_{0}$ and $h_{0}\comp j_{0}=f_{0}$.\noproof 
\end{proposition}

For us, the notion of essential surjectivity has several relevant internalizations, resulting in different notions of internal equivalence. Our weak equivalences in $\Cat \Cc$ will be defined relative to some class of morphisms $\Ec$ in $\Cc$, which in practice will be the class of $\Tc$-epimorphisms for a topology $\Tc$ on $\Cc$.

\begin{definition}\label{Definition-Relative-Weak-Equivalence}
Let $\Ec$ be a class of morphisms and $\vetf\colon \Af\to \Bf$ an internal functor in $\Cc$. If the morphism $\delta_{0}\comp \overline{f_{0}}$ in the diagram 
\[
\xymatrix{(P_{f})_{0}\ar@{}[rd]|<{\pullback} \ar[d]_-{\overline{\delta_{1}}} \ar[r]^-{\overline{f_{0}}}& \iso (B) \ar[r]^-{\delta_{0}} \ar[d]^-{\delta_{1}} & B_{0} \\
A_{0} \ar[r]_-{f_{0}} & B_{0}}
\]
is in $\Ec$, then $\vetf$ is called \emph{essentially $\Ec$-surjective}. An \emph{$\Ec$-equivalence} is an internal functor that is full, faithful and essentially $\Ec$-surjective. If $\Ec=\Ec_{\Tc}$ is the class of $\Tc$-epimorphisms for a Grothendieck topology $\Tc$ on $\Cc$, the respective notions become \emph{essentially $\Tc$-surjective}\index{essentially $\Tc$-surjective} and \emph{$\Tc$-equivalence}\index{T-equivalence@$\Tc$-equivalence}\index{internal!equivalence!$\Tc$- ---}. The class of $\Tc$-equivalences for a topology $\Tc$ is denoted by $\we (\Tc)$\index{we T@$\we (\Tc)$}.
\end{definition}

\begin{example}\label{Example-Weakest-Equivalence}
In case $\Tc$ is the cotrivial topology, any functor is essentially $\Tc$-surjective, and hence the $\Tc$-equivalences are exactly the fully faithful functors.
\end{example}

\begin{example}\label{Example-Strong-Equivalence}
If $\Tc$ is the trivial topology, an internal functor $\vetf\colon {\Af \to \Bf}$ is essentially $\Tc$-surjective if and only if the functor $\Cat \Cc (\Xf,\vetf)$ is essentially surjective for all $\Xf$. If $\vetf$ is moreover fully faithful, it is called a \emph{strong equivalence}\index{strong equivalence}\index{internal!equivalence!strong}. This name is justified by the obvious fact that a strong equivalence is a $\Tc$-equivalence for every topology $\Tc$. If $\vetf$ is a strong equivalence, a functor $\vetg\colon {\Bf \to \Af}$ exists and natural isomorphisms $\epsilon \colon  \vetf\comp \vetg\To 1_{\Bf} $ and $\eta \colon 1_{\Af}\To \vetg\comp \vetf $; hence $\vetf$ is a homotopy equivalence with respect to the cocylinder from Section~\ref{Section-Cocylinder}. There is even more:
\end{example}

Recall that an \emph{internal adjunction}\index{internal!adjunction} is a quadruple
\[
(\vetf\colon \Af \to \Bf , \quad \vetg\colon \Bf \to \Af, \quad \epsilon \colon \vetf\comp \vetg\To 1_{\Bf},\quad \eta \colon 1_{\Af}\To \vetg\comp \vetf )
\]
such that the \emph{triangular identities}\index{triangular identities} $(\epsilon \comp 1_{\vetf})\vcomp (1_{\vetf}\comp \eta)=1_{\vetf}$ and $(1_{\vetg}\comp \epsilon)\vcomp (\eta \comp 1_{\vetg})=1_{\vetg}$ hold. Then $\vetf$ is \emph{left adjoint} to $\vetg$, $\vetg$ \emph{right adjoint} to $\vetf$, $\epsilon$ the \emph{counit} and $\eta$ the \emph{unit} of the adjunction. Using J. W. Gray's terminology~\cite{Gray}, we shall call \emph{lali}\index{lali} a left adjoint left inverse functor, and, dually, \emph{rari}\index{rari} a right adjoint right inverse functor. In case $\vetf$ is left adjoint left inverse to $\vetg$, we denote the situation $\vetf=\lali \vetg$ or $\vetg=\rari \vetf$.

\begin{remark}\label{Remark-lali-Triangular-Identities}
Since then $\vetf\comp \vetg=1_{\Bf}$ and $\epsilon =1_{1_{\Bf}}\colon 1_{\Bf}\To 1_{\Bf}$, the triangular identities reduce to $1_{\vetf}= (1_{1_{\Bf}} \comp 1_{\vetf})\vcomp (1_{\vetf} \comp \eta) = 1_{\vetf}\comp \eta $, which means that 
\[
f_{1}\comp i=m\comp (f_{1}\comp i,f_{1}\comp \eta)=f_{1}\comp m \comp (i,\eta)=f_{1}\comp m \comp (i\comp d_{0}, 1_{\Af})\comp \eta =f_{1}\comp \eta,
\]
and $1_{\vetg}=(1_{\vetg}\comp 1_{1_{\Af}})\vcomp (\eta \comp 1_{\vetg})=\eta \comp 1_{\vetg}$, meaning that $i\comp g_{0}=\eta \comp g_{0}$.
\end{remark}
 
An \emph{adjoint equivalence}\index{adjoint equivalence} is a (left and right) adjoint functor with unit and counit natural isomorphisms. It is well known that every equivalence of categories is an adjoint equivalence; see e.g., Borceux~\cite{Borceux:Cats} or Mac\,Lane~\cite{MacLane}. It is somewhat less known that this is still the case for strong equivalences of internal categories. In fact, in any $2$-category, an equivalence between two objects is always an adjoint equivalence; see Blackwell, Kelly and Power~\cite{Blackwell}. More precisely, the following holds.

\begin{proposition}\cite{Blackwell}\label{Proposition-2-Adjoints}
Let $\CCf$ be a $2$-category and $f\colon {C\to D}$ a $1$-cell of $\CCf$. Then $f$ is an adjoint equivalence if and only if for every object $X$ of $\CCf$, the functor 
\[
\CCf (X,f)\colon {\CCf (X,C)\to \CCf (X,D)}
\]
is an equivalence of categories. \noproof 
\end{proposition}

Hence, in the $2$-category $\Cat \Cc$ of internal categories in a given finitely complete category $\Cc$, every strong equivalence is adjoint; and in the $2$-category $\Gd \Cc$ of internal groupoids in $\Cc$, the notions ``adjunction'', ``strong equivalence'' and ``adjoint equivalence'' coincide.

\begin{remark}\label{Remark-lali-we-Triangular-Identities}
If $\vetf\colon \Af \to \Bf$ is a split epimorphic fully faithful functor, it is always a strong equivalence. Denote $\vetg=\rari \vetf\colon \Bf \to \Af$ its right adjoint right inverse. Then the unit $\eta$ of the adjunction induces a homotopy $\vetH\colon {\Af \to \Af^{\If}}$ from $1_{\Af}$ to $\vetg\comp \vetf$. It is easily checked that the triangular identities now amount to $\vetf^{\If}\comp s (\Af)=\vetf^{\If}\comp \vetH$ and $s (\Af)\comp \vetg=\vetH\comp \vetg$. 
\end{remark}

\begin{example}\label{Example-Epsilon-Strong-Equivalence}
Example~\ref{Example-Epsilon-Homotopy} implies that for any internal category $\Af$, $s (\Af)$ is a right adjoint right inverse of $\epsilon_{0} (\Af)$ and $\epsilon_{1} (\Af)$. \textit{A fortiori}, the three internal functors are strong equivalences.
\end{example}

\begin{example}\label{Example-Weak-Equivalence}
If $\Tc$ is the regular epimorphism topology then an internal functor $\vetf$ is in $\we (\Tc)$ if and only if it is a \emph{weak equivalence}\index{weak equivalence} in the sense of Bunge and Par\'e~\cite{Bunge-Pare}. In case $\Cc$ is semi-abelian, weak equivalences may be characterized using homology (Proposition~\ref{Proposition-Characterization-Weak-Equivalence}).
\end{example}

In order, for a class of morphisms in a category, to be the class of weak equivalences in a model structure, it needs to satisfy the two-out-of-three property (Definition~\ref{Definition-Model-Category}). The following proposition gives a sufficient condition for this to be the case.

\begin{proposition}\label{Proposition-Two-out-of-Three}
If $\Tc$ is a subcanonical topology on a category $\Cc$ then the class of $\Tc$-equivalences has the two-out-of-three\index{two-out-of-three property} property.
\end{proposition}
\begin{proof}
For a subcanonical topology $\Tc$, the Yoneda embedding, considered as a functor 
\[
Y\colon {\Cc \to \Sh (\Cc ,\Tc)},
\]
is equal to $Y_{\Tc}$. It follows that $Y_{\Tc}$ is full and faithful and preserves and reflects limits. Hence it induces a $2$-functor $\Cat Y_{\Tc}\colon  \Cat \Cc \to \Cat \Sh (\Cc ,\Tc)$. Moreover, this $2$-functor is such that an internal functor $\vetf \colon  \Af \to \Bf$ in $\Cc$ is a $\Tc$-equivalence if and only if the functor $\Cat Y_{\Tc} (\vetf)$ in $\Cat \Sh (\Cc ,\Tc)$ is a weak equivalence. According to Joyal and Tierney~\cite{Joyal-Tierney}, weak equivalences in a Grothendieck topos have the two-out-of-three property; the result follows. 
\end{proof}

Not every topology induces a class of equivalences that satisfies the two-out-of-three property, as shows the following example.

\begin{example}\label{Example-No-Two-out-of-Three}
Let $\vetg\colon \Bf \to \Cf$ a functor between small categories that preserves terminal objects. Let $\vetf\colon {\matheu{1}\to \Bf}$ be a functor from a terminal category to $\Bf$ determined by the choice of a terminal object in $\Bf$. Then $\vetg\comp \vetf$ and $\vetf$ are fully faithful functors, whereas $\vetg$ need not be fully faithful. Hence the class of $\Tc$-equivalences induced by the cotrivial topology on $\Set$ does not satisfy the two-out-of-three property.
\end{example}

\section{The $\Tc$-model structure on $\protect\Cat\Cc$}\label{Section-Main-Theory}

In this section we suppose that $\Cc$ is a finitely complete category such that $\Cat \Cc$ is finitely complete and cocomplete.

\begin{definition}\label{Definition-Relative-Fibration}
Let $\Ec$ be a class of morphisms in $\Cc$ and $\vetp\colon \Ef\to \Bf$ an internal functor. $\vetp$~is called an \emph{$\Ec$-fibration} if and only if in the left hand side diagram
\begin{equation}\label{Diagram-T-Fibration}
\vcenter{\xymatrix@!@=0ex{\iso (E) \ar@/_/[rrdddd]_-{\delta_{1}} \ar@/^/[rrrrdd]^-{\iso (p)_{1}} \ar@{.>}[rrdd]|{(\vetr_{\vetp})_{0} } &&&&& X \ar@/_/[rrdddd]_{e} \ar@/^/[rrrrdd]^-{\beta }  \\
&&&&&& U_{i} \ar@{.>}[lu]|-{f_{i}} \ar@{.>}[rd]|-{\epsilon_{i}}\\
&& (P_{p})_{0} \ar@{}[rrdd]|<{\pullback} \ar[rr]^-{\overline{p_{0}}} \ar[dd]_-{\overline{\delta_{1}}} && \iso(B) \ar[dd]^-{\delta_{1}} &&& \iso (E) \ar@{}[rrdd]|{\texttt{(i)}} \ar[rr]^-{\iso (p)_{1}} \ar[dd]_-{\delta_{1}} && \iso(B) \ar[dd]^-{\delta_{1}}\\\\
&& E_{0} \ar[rr]_-{p_{0}} && B_{0} &&& E_{0} \ar[rr]_-{p_{0}} && B_{0}}}
\end{equation}
the induced universal arrow $(r_{p})_{0}$ is in $\Ec$. If $\Ec =\Ec_{\Tc}$ comes from a topology $\Tc$ on $\Cc$ we say that $\vetp$ is a \emph{$\Tc$-fibration}\index{T-fibration@$\Tc$-fibration}\index{fibration!$\Tc$- ---}. The functor $\vetp$ is said to be \emph{star surjective, relative to $\Tc$}\index{star surjective!relative to $\Tc$} if, given an object $X$ in $\Cc$ and arrows $e$ and $\beta$ such as in the right hand side diagram above, there exists a covering family $(f_{i}\colon U_{i}\to X)_{i\in I}$ and a family of morphisms $(\epsilon_{i}\colon U_{i}\to \iso (E))_{i\in I}$ keeping it commutative for all $i\in I$.
\end{definition}

By Proposition~\ref{Proposition-Characterization-T-Epimorphism}, an internal functor $\vetp$ is a $\Tc$-fibration if and only if it is star surjective, relative to $\Tc$.

\begin{example}\label{Example-Strong-Fibration}
If $\Tc$ is the trivial topology then an internal functor $\vetp$ is a $\Tc$-fibration if and only if the square \texttt{(i)} is a weak pullback. Such a $\vetp$ is called a \emph{strong fibration}\index{strong fibration}\index{fibration!strong}. In case $\Cc$ is $\Set$, the strong fibrations are the \emph{star surjective}\index{star surjective} functors~\cite{Brown:Fibrations}. It is easily seen that the unique arrow $\Af \to \matheu{1}$ from an arbitrary internal category $\Af$ to a terminal object $\matheu{1}$ of $\Cat \Cc$ is always a strong fibration; hence every object of $\Cat\Cc$ is strongly fibrant. 
\end{example}

\begin{example}\label{Example-Weakest-Fibration}
Obviously, if $\Tc $ is the cotrivial topology, any functor is a $\Tc$-fibration.
\end{example}

\begin{example}\label{Example-Discrete-Fibration}
An internal functor $\vetp\colon \Ef \to \Bf$ is called a \emph{discrete fibration}\index{discrete fibration}\index{fibration!discrete} if the square 
\[
\xymatrix{E_{1} \ar[d]_-{d_{1}} \ar[r]^-{p_{1}} & B_{1}\ar[d]^-{d_{1}}\\
E_{0} \ar[r]_-{p_{0}} & B_{0}}
\]
is a pullback. Every discrete fibration is a strong fibration. Note that this is obvious in case $\Ef$ is a groupoid; in general, one proves it by considering morphisms $e\colon {X\to E_{0}}$ and $\beta \colon {X\to \iso (B)}$ such that $p_{0}\comp e=\delta_{1}\comp \beta=d_{1}\comp j\comp \beta$. Then a unique morphism $\epsilon \colon X\to E_{1}$ exists such that $p_{1}\comp \epsilon =j\comp \beta$ and $d_{1}\comp \epsilon =e$. This $\epsilon$ factors over $\iso (E)$: Indeed, since $e'=d_{0}\comp \epsilon$ is such that $p_{0}\comp e'=d_{0}\comp p_{1}\comp \epsilon =d_{1}\comp j\comp  \tw \comp \beta$, there exists a unique arrow $\epsilon '\colon X\to E_{1}$ such that $d_{1}\comp \epsilon '=e'$ and $p_{1}\comp \epsilon '=j\comp \tw \comp \beta$. Using the fact that the square above is a pullback, it is easily shown that $\epsilon'$ is the inverse of $\epsilon$. 
\end{example}

Given a topology $\Tc$ on $\Cc$, we shall consider the following structure on $\Cat \Cc$:  $\we (\Tc)$\index{we T@$\we (\Tc)$} is the class of $\Tc$-weak equivalences; $\fib (\Tc)$\index{fib T@$\fib (\Tc)$} is the class of $\Tc$-fibrations; $\cof (\Tc)$\index{cof T@$\cof (\Tc)$} is the class ${^{\Box} (\fib (\Tc) \cap \we (\Tc))}$ of \emph{$\Tc$-cofibrations}\index{T-cofibration@$\Tc$-cofibration}\index{cofibration!$\Tc$- ---}, internal functors having the left lifting property with respect to all \emph{trivial $\Tc$-fibrations}\index{trivial fibration!$\Tc$- ---}. 

The aim of this section is to prove the following 

\begin{theorem}\label{Theorem-Model-Structure}
If $\we (\Tc)$ has the two-out-of-three property and $\Cc$ has enough $\Ec_{\Tc}$-projectives then $(\Cat\Cc, \fib (\Tc) , \cof (\Tc), \we (\Tc))$ is a model category.
\end{theorem}

\begin{proposition}\label{Proposition-Characterization-Trivial-Fibrations}
A functor $\vetp\colon \Ef \to \Bf$ is a trivial $\Tc$-fibration if and only if it is fully faithful, and such that $p_{0}$ is a $\Tc$-epimorphism.
\end{proposition}
\begin{proof}
If $\vetp$ is a trivial $\Tc$-fibration then it is a fully faithful functor. Suppose that $b\colon {X\to B_{0}}$ is an arbitrary morphism. Since $\vetp$ is essentially $\Tc$-surjective, a covering family $(f_{i}\colon {U_{i}\to X})_{i\in I}$ exists and families of maps $(\beta_{i} \colon {U_{i}\to B_{1}})_{i\in I}$ and $(e_{i}\colon {U_{i}\to E_{0}})_{i\in I}$ such that $p_{0}\comp e_{i}=\delta_{1}\comp \beta_{i}$ and $\delta_{0}\comp \beta_{i} =b\comp f_{i}$. Since $\vetp$ is a $\Tc$-fibration, for every $i\in I$, $e_{i}$ and $\beta_{i}$ induce a covering family $(g_{ij}\colon {V_{ij}\to U_{i}})_{j\in I_{i}}$ and a family of morphisms $(\epsilon_{ij} \colon {V_{ij}\to \iso (E)})_{j\in I_{i}}$. Put $b'_{ij}=\delta_{0}\comp \epsilon_{ij} \colon {V_{ij}\to E_{0}}$, then $p_{0}\comp b'_{ij}=b\comp f_{i}\comp g_{ij}$; because by the transitivity axiom, 
\[
(f_{i}\comp g_{ij}\colon {V_{ij}\to X})_{j\in I_{i}, i\in I}
\]
forms a covering family, this shows that $p_{0}$ is a $\Tc$-epimorphism.

Conversely, we have to prove that $\vetp$ is an essentially $\Tc$-surjective $\Tc$-fibration. Given $b\colon {X\to B_{0}}$, $p_{0}$ being in $\Ec_{\Tc}$ induces a covering family $(f_{i}\colon {U_{i}\to X})_{i\in I}$ and a family of morphisms $(e_{i}\colon {U_{i}\to E_{0}})_{i\in I}$ such that $p_{0}\comp e_{i}=b\comp f_{i}$. Using the equality $\delta_{0}\comp (\iota\comp b)=\delta_{1}\comp (\iota\comp b)=b$, this shows that $\vetp$ is essentially $\Tc$-surjective. To prove $\vetp$ a $\Tc$-fibration, consider the right hand side commutative diagram of solid arrows~\ref{Diagram-T-Fibration} above. Because $p_{0}$ is a $\Tc$-epimorphism, there is a covering family 
\[
(f_{i}\colon U_{i}\to X)_{i\in I}
\]
and a family $(e'_{i}\colon {U_{i}\to E_0})_{i\in I}$ such that $p_{0}\comp e'_{i}=\delta_{0}\comp \beta\comp f_{i} $. This gives rise to a diagram
\[
\xymatrix{U_{i} \ar@/_/[rdd]_-{(e\comp f_{i},e'_{i})} \ar@/^/[rrd]^-{\beta\comp f_{i} } \ar@{.>}[rd]|{\epsilon_{i}} \\
& \iso (E) \ar@{}[rd]|<{\pullback} \ar[r]^-{\iso (p)_{1}} \ar[d]^>>>>{(\delta_{0},\delta_{1})} & \iso (B) \ar[d]^-{(\delta_{0},\delta_{1})}\\
& E_{0}\times E_{0} \ar[r]_-{p_{0}\times p_{0}} & B_{0}\times B_{0}}
\]
for every $i\in I$. $\vetp$ being fully faithful, according to Remark~\ref{Remark-Iso-Preserves-Full-Faithful-Functors}, its square is a pullback, and induces the needed family of dotted arrows 
\[
(\epsilon_{i}\colon U_{i}\to \iso (E))_{i\in I}\qedhere 
\]
\end{proof}

\begin{example}\label{Example-epsilon-Strong-Fibrations}
For any internal category $\Af$, $\epsilon_{0} (\Af)$ and $\epsilon_{1} (\Af)$, being strong equivalences that are split epimorphic on objects, are strong fibrations.
\end{example}

\begin{corollary}\label{Corollary-Characterization-Fibrations}
A functor $\vetp\colon \Ef \to \Bf$ is a $\Tc$-fibration if and only if the universal arrow $\vetr_{\vetp}\colon \Ef^{\If}\to \Pf_{\vetp}$ is a trivial $\Tc$-fibration.
\end{corollary}
\begin{proof}
This is an immediate consequence of Proposition~\ref{Proposition-Characterization-Trivial-Fibrations}, Example~\ref{Example-Epsilon-Strong-Equivalence} and the fact that strong equivalences have the two-out-of-three property (Proposition~\ref{Proposition-Two-out-of-Three}).
\end{proof}

\begin{corollary}\label{Corollary-Characterization-Cofibrations}
An internal functor $\vetj\colon \Af \to \Xf$ is a $\Tc$-cofibration if and only if $j_{0}\in{}^{\Box}\Ec_{\Tc}$.
\end{corollary}
\begin{proof}
This follows from Proposition~\ref{Proposition-Lifting-we} and Proposition~\ref{Proposition-Characterization-Trivial-Fibrations}.
\end{proof}

\begin{proposition}\label{Proposition-Third-Factorization}
Every internal functor of $\Cc$ may be factored as a strong equivalence (right inverse to a strong trivial fibration) followed by a strong fibration.
\end{proposition}
\begin{proof}
This is an application of K.\,S. Brown's Factorization Lemma~\cite{K.S.Brown}. To use it, we must show that $(\Cat \Cc , \Fc ,\Wc)$, where $\Fc$ is the class of strong fibrations and $\Wc$ is the class of strong equivalences, forms a category of fibrant objects. Condition (A) is just Proposition~\ref{Proposition-Two-out-of-Three} and (B) follows from the fact that split epimorphisms are stable under pulling back. Proving (C) that strong fibrations are stable under pulling back is easy, as is the stability of $\Fc \cap \Wc$, the class of split epimorphic fully faithful functors. The path space needed for (D) is just the cocylinder from Section~\ref{Section-Cocylinder}. Finally, according to Example~\ref{Example-Strong-Fibration}, every internal category is strongly fibrant, which shows condition (E).
\end{proof}

\begin{proposition}\label{Proposition-Characterization-Trivial-Cofibrations}
Any trivial $\Tc$-cofibration is a split monic adjoint equivalence.
\end{proposition}
\begin{proof}
Using Proposition~\ref{Proposition-Third-Factorization}, factor the trivial $\Tc$-cofibration $\vetj\colon \Af \to \Xf$ as a strong equivalence $\vetf\colon \Af \to \Bf$ (in fact, a right adjoint right inverse) followed by a strong fibration $\vetp \colon {\Bf \to \Xf}$. By the two-out-of-three property of weak equivalences, $\vetp$ is a trivial $\Tc$-fibration; hence the commutative square
\[
\xymatrix{{\Af} \ar[d]_-{\vetj} \ar[r]^-{\vetf} & {\Bf} \ar[d]^{\vetp} \\
{\Xf} \ar@{=}[r] \ar@{.>}[ru]^{\vets} & {\Xf}}
\]
has a lifting $\vets \colon \Xf \to \Bf$. It follows that $\vetj$ is a retract of $\vetf$. The class of right adjoint right inverse functors being closed under retracts, we may conclude that $\vetj$ is a split monic adjoint equivalence. 
\end{proof}

\begin{proposition}[Covering Homotopy Extension Property]\label{Proposition-CHEP}\index{Covering Homotopy Extension Property}
Consider the commutative diagram of solid arrows
\[
\xymatrix{{\Af} \ar[d]_-{\vetj} \ar[r]^-{\vetH} & {\Ef^{\If} } \ar[d]|{\hole}^>>>{\epsilon_{1} (\Ef)} \ar[r]^-{\vetp^{\If}} & {\Bf^{\If}} \ar[d]^-{\epsilon_{1} (\Bf)} \\
{\Xf} \ar@{.>}[ru]^-{\vetL} \ar[rru]_>>>>>{\vetK} \ar[r]_-{\vetf} & {\Ef} \ar[r]_-{\vetp} & {\Bf}.}
\]
If $\vetj\in \cof (\Tc)$ and $\vetp\in \fib (\Tc)$, then a morphism $\vetL\colon \Xf \to \Ef^{\If}$ exists keeping the diagram commutative. 
\end{proposition}
\begin{proof}
Since $\vetp$ is a $\Tc$-fibration, by Corollary~\ref{Corollary-Characterization-Fibrations}, the associated universal arrow $\vetr_{\vetp}\colon \Ef^{\If}\to \Pf_{\vetp}$ is a trivial $\Tc$-fibration. Let $\vetM\colon \Xf \to \Pf_{\vetp}$ be the unique morphism such that $\overline{\epsilon_{1} (\Bf)}\comp \vetM=\vetf$ and $\overline{\vetp}\comp \vetM=\vetK$ (cf.\ Diagram~\ref{Diagram-Mapping-Path-Space}); then the square
\[
\xymatrix{{\Af} \ar[d]_-{\vetj} \ar[r]^-{\vetH} & {\Ef^{\If} }\ar[d]^-{\vetr_{\vetp}} \\
{\Xf} \ar@{.>}[ru]^-{\vetL} \ar[r]_-{\vetM} & {\Pf_{\vetp}}}
\]
commutes, and yields the needed lifting $\vetL$.
\end{proof}

\begin{proposition}\cite[Proposition I.3.11]{KP:Abstact-homotopy-theory}\label{Proposition-RLP}
Any $\Tc$-fibration has the right lifting property with respect to any trivial $\Tc$-cofibration.
\end{proposition}
\begin{proof}
Suppose that in square~\ref{Diagram-Lifting}, $\vetj\in \cof (\Tc)\cap \we (\Tc)$ and $\vetp\in \fib (\Tc)$. By Proposition~\ref{Proposition-Characterization-Trivial-Cofibrations}, $\vetj$ is a split monomorphic adjoint equivalence; denote $\vetk=\lali \vetj\colon \Xf \to \Af$. According to Remark~\ref{Remark-lali-we-Triangular-Identities}, a homotopy $\vetH\colon {\Xf\to \Xf^{\If}}$ from $1_{\Xf }$ to $\vetj\comp \vetk$ may be found such that $\vetH\comp \vetj=s (\Xf)\comp \vetj\colon {\Af \to \Xf^{\If}}$. Because the diagram of solid arrows
\[
\xymatrix{{\Af} \ar[d]_-{\vetj} \ar[r]^-{s (\Ef)\comp \vetf} & {\Ef^{\If} } \ar[d]|{\hole}^>>>{\epsilon_{1} (\Ef)} \ar[r]^-{\vetp^{\If}} & {\Bf^{\If}} \ar[d]^{\epsilon_{1} (\Bf)} \\
{\Xf} \ar@{.>}[ru]^-{\vetL} \ar[rru]|>>>>>>>{\vetg^{\If}\comp \vetH} \ar[r]_-{\vetf\comp \vetk} & {\Ef} \ar[r]_-{\vetp} & {\Bf}}
\]
commutes, the Covering Homotopy Extension Property gives rise to a morphism $\vetL\colon {\Xf \to \Ef^{\If}}$; the morphism $\veth=\epsilon_{0} (\Ef)\comp \vetL$ is the desired lifting for Diagram~\ref{Diagram-Lifting}.
\end{proof}

\begin{lemma}\cite[Lemma 2.1]{Johnstone:Herds}\label{Lemma-Pullback}
Let $\Bf$ be a category in $\Cc$ and $p_{0}\colon {E_{0}\to B_{0}}$ a morphism in $\Cc$. Form the pullback
\[
\xymatrix{E_{1} \ar@{.>}[r]^{p_{1}} \ar@{.>}[d]_-{(d_{0},d_{1})} \ar@{}[rd]|<<<{\pullback} & B_{1} \ar[d]^-{(d_{0},d_{1})} \\
E_{0}\times E_{0} \ar[r]_-{p_{0}\times p_{0}} & B_{0}\times B_{0}.}
\]
Then the left hand side graph $\Ef$ carries a unique internal category structure such that $\vetp=(p_{0},p_{1})\colon \Ef \to \Bf$ is a functor. Moreover, $\vetp$ is a fully faithful functor, and if $p_{0}$ is a $\Tc$-epimorphism, then $\vetp$ is a trivial $\Tc$-fibration.\noproof 
\end{lemma}

\begin{proposition}\label{Proposition-Second-Factorization}
Every internal functor of $\Cc$ may be factored as a $\Tc$-cofibration followed by a trivial $\Tc$-fibration.
\end{proposition}
\begin{proof}
Let $\vetf\colon \Af \to \Bf$ be an internal functor, and, using Remark~\ref{Remark-Enough-Projectives} that $({}^{\Box} \Ec_{\Tc},\Ec_{\Tc})$ forms a weak factorization system, factor $f_{0}$ as an element $j_{0}\colon {A_{0}\to E_{0}}$ of ${}^{\Box} \Ec_{\Tc}$ followed by a $\Tc$-epimorphism $p_{0}\colon {E_{0}\to B_{0}}$. Then the construction in Lemma~\ref{Lemma-Pullback} yields a trivial $\Tc$-fibration $\vetp\colon \Ef \to \Bf$. Let $j_{1}\colon {A_{1}\to E_{1}}$ be the unique morphism such that $p_{1}\comp j_{1}=f_{1}$ and $(d_{0},d_{1})\comp j_{1}= (j_{0}\times j_{0})\comp (d_{0},d_{1})$. Since $\vetp$ is faithful, $\vetj= (j_{0},j_{1})\colon {\Af \to \Ef}$ is a functor; according to Corollary~\ref{Corollary-Characterization-Cofibrations}, it is a $\Tc$-cofibration.
\end{proof}

\begin{proposition}\label{Proposition-First-Factorization}
Every internal functor of $\Cc$ may be factored as a trivial $\Tc$-cofibration followed by a $\Tc$-fibration.
\end{proposition}
\begin{proof}
For an internal functor $\vetf\colon \Af \to \Bf$, let $\vetf=\vetp\comp \vetj'$ be the factorization of $\vetf$ from Proposition~\ref{Proposition-Third-Factorization}. Then $\vetp$ is a $\Tc$-fibration and $\vetj'$ a $\Tc$-equivalence. Using Proposition~\ref{Proposition-Second-Factorization}, this $\vetj'$ may be factored as a (necessarily trivial) $\Tc$-cofibration $\vetj$ followed by a $\Tc$-fibration $\vetp'$. Thus we get a trivial $\Tc$-cofibration $\vetj$ and a $\Tc$-fibration $\vetp\comp \vetp'$ such that $f= (\vetp\comp \vetp')\comp \vetj$.
\end{proof}

\begin{proof}[Proof of Theorem~\ref{Theorem-Model-Structure}]
We only need to comment on the closedness under retracts of the classes $\fib (\Tc)$, $\cof (\Tc)$ and $\we (\Tc)$. For the $\Tc$-fibrations, $\Tc$-cofibrations and essentially $\Tc$-surjective morphisms, closedness under retracts follows from the closedness of the classes $\Ec_{\Tc}$ and ${^{\Box}\Ec_{\Tc}}$; for fully faithful functors, the property has a straightforward direct proof.
\end{proof}

\section{Case study: the regular epimorphism topology}\label{Section-Regular-Epimorphisms}\index{topology!regular epimorphism ---}\index{regular epimorphism topology}

This section treats the model category structure on $\Cat \Cc$ induced by choosing $\Tc$ the regular epimorphism topology on $\Cc$. Let us recall what we already know about it: The class $\we (\Tc)$ consists of weak equivalences of internal categories; fibrations are $\regepi$-fibrations (where $\regepi$\index{regepi@$\regepi$} denotes the class of all regular epimorphisms in $\Cc$) and cofibrations have an object morphism in ${^{\Box}\regepi}$. Hence all of its objects are fibrant, and an internal category is cofibrant if and only if its object of objects is $\regepi$-projective.

We shall be focused mainly on semi-abelian categories and internal crossed modules, but we start by explaining a connection with Grothendieck topoi and Joyal and Tierney's model structure.

\subsection*{Grothendieck topoi}\label{Subsection-Grothendieck-Topoi}\index{Grothendieck topos}\index{topos!Grothendieck ---} Let $\Cc$ be a Grothendieck topos equipped with the regular epimorphism topology $\Tc $. Since all epimorphisms are regular, it is clear that epimorphisms and $\Tc$-epimorphisms coincide. It follows that the notion of $\Tc$-equivalence coincides with the one considered by Joyal and Tierney in~\cite{Joyal-Tierney}. Consequently, the two structures have equivalent homotopy categories. It is, however, clear that in general, Joyal and Tierney's fibrations and cofibrations are different from ours.

\subsection*{Semi-abelian categories}\label{Subsection-Semi-Abelian-Categories}\index{category!semi-abelian} From now on we suppose that $\Ac$ is a semi-abelian category with enough $\regepi$-projectives, to give an alternative characterization of weak equivalences, and to describe the induced model category structure on the category of internal crossed modules.

Recall from Section~\ref{Section-Simplicial-Objects} the following notion of homology of simplicial objects in a semi-abelian category. First of all, a morphism is called \emph{proper} when its image is a kernel, and a chain complex is \emph{proper} whenever all its boundary operators are. As in the abelian case, the \emph{$n$-th homology object} of a proper chain complex $C$ with boundary operators $d_n$ is said to be $H_{n}C=\Cok[C_{n+1}\to K[d_n]]$. The category of proper chain complexes in $\Ac$ is denoted $\PCh \Ac$, and $\Sc \Ac=\Fun (\Delta^{\op},\Ac)$ is the category of simplicial objects in $\Ac$. The \emph{normalization functor} $N\colon {\Sc \Ac \to \PCh \Ac}$ maps a simplicial object $A$ in $\Ac$ with face operators $\del_{i}$ and degeneracy operators $\sigma_{i}$ to the chain complex $N (A)$ in $\Ac$ given by
\[
N_{n} A=\bigcap_{i=0}^{n-1}K[\del_{i}\colon A_{n}\to A_{n-1}],\qquad d_{n}=\del_{n}\comp \bigcap_{i}\ker \del_{i}\colon N_{n} A\to N_{n-1} A,
\]
for $n\geq 1$, and $N_{n}A=0$, for $n<0$. Because $N (A)$ is proper (Theorem~\ref{Theorem-N-Proper}), it makes sense to consider its homology. Indeed, the \emph{$n$-th homology object} $H_{n}A$ is defined to be the $n$-th homology object $H_{n}N (A)$ of the associated proper chain complex $N (A)$. This process sends a short exact sequences of simplicial objects to a long exact sequence in $\Ac$ (Proposition~\ref{Proposition-Long-Exact-Homology-Sequence-Simp}).

Since, via the nerve construction, any internal category may be considered as a simplicial object, we can apply this homology theory to internal categories. More precisely, recall (e.g., from Johnstone~\cite[Remark 2.13]{Johnstone:Topos-Theory}) that there is an embedding $\ner\colon \Cat \Ac \to \Sc \Ac$ of $\Cat \Ac$ into $\Sc \Ac$ as a full subcategory. Given a category $\Af$ in $\Ac$, its \emph{nerve}\index{nerve!of an internal category} $\ner \Af$\index{ner@$\ner \Af$} is the simplicial object defined on objects by $n$-fold pullback $\ner_{n}\Af =A_{1}\times_{A_{0}}\dots \times_{A_{0}}A_{1}$ if $n\geq 2$, $\ner_{1}\Af =A_{1}$ and $\ner_{0}\Af =\Af_{0}$; on morphisms by $\del_{0}=d_{1}$, $\del_{1}=d_{0}\colon {\ner_{1}\Af \to \ner_{0}\Af}$, $\sigma_{0}=i\colon {\ner_{0}\Af \to \ner_{1}\Af}$; $\del_{0}=\pr_{2}$, $\del_{1}=m$, $\del_{2}=\pr_{1}\colon {\ner_{2}\Af \to \ner_{1}\Af}$ and $\sigma_{0}= (i,1_{A_{1}})$, $\sigma_{1}= (1_{A_{1}},i)\colon {\ner_{1}\Af \to \ner_{2}\Af}$; etc. (See also Examples~\ref{Examples-Simplicial-Objects}.) A~simplicial object is isomorphic to an object in the image of $\ner$ if and only if, as a functor ${\Delta^{\op}\to \Ac}$, it is left exact.

\begin{definition}\label{Definition-Homology-of-Internal-Category}\index{homology!of internal categories}
Suppose that $\Ac$ is semi-abelian, $\Af$ is a category in $\Ac$ and $n\in \Z$. The object $H_{n}\Af =H_{n}\ner \Af$\index{H_n@$H_{n}\Af$} will be called the \emph{$n$-th homology object} of $\Af$ and the functor $H_{n}=H_{n}\comp \ner \colon \Cat \Ac \to \Ac$ the \emph{$n$-th homology functor}.
\end{definition}

\begin{proposition}\label{Proposition-Characterization-Homology}
Let $\Af$ be a category in a semi-abelian category $\Ac$. If $n\not\in\{0,1 \}$ then $H_{n}\Af=0$, $H_{1}\Af =K[(d_{0},d_{1})\colon A_{1}\to A_{0}\times A_{0}]$ and $H_{0}\Af =\Coeq [d_{0},d_{1}\colon A_{1}\to A_{0}]$. Moreover, $H_{1}\Af =K[(d_{0},d_{1})]$ is an abelian object of $\Ac$.
\end{proposition}
\begin{proof}
Using e.g., the Yoneda Lemma (in the form of Metatheorem 0.2.7 in~\cite{Borceux-Bourn}), it is quite easily shown that $N\Af = N\ner \Af$ is the chain complex
\[
\cdots \to 0 \to K[d_{1}] \to^{d_{0}\comp \ker d_{1}} A_{0} \to 0 \to \cdots . 
\]
By Remark~\ref{Remark-N-1}, $K[d_{0}\comp \ker d_{1}]=K[d_{0}]\cap K[d_{1}]$, which is clearly equal to 
\[
K[(d_{0},d_{1})\colon {A_{1}\to A_{0}\times A_{0}}].
\]
The last equality is an application of Corollary~\ref{Corollary-H_0}, and $H_{1}\Af$ being abelian is a consequence of Theorem~\ref{Theorem-Homology-Abelian} or Bourn~\cite[Proposition 3.1]{BournBodeux}. 
\end{proof}

A related issue concerns the existence of the so-called \emph{fundamental groupoid}\index{fundamental groupoid} functor 
\[
\pi_{1}\colon {\Sc \Ac\to \Cat \Ac}\index{pi_1@$\pi_{1}\colon {\Sc \Ac\to \Cat \Ac}$}, 
\]
a left adjoint to the functor $\ner \colon {\Cat \Ac \to \Sc \Ac}$. This will allow for a comparison with the model structure from Chapter~\ref{Chapter-Simplicial-Objects}.

\begin{proposition}\label{Proposition-Fundamental-Groupoid}
The functor $\ner$ is a right adjoint.  
\end{proposition}
\begin{morale}
To a simplicial object $A$ one associates the internal groupoid $\pi_{1}A$ induced by identifying $1$-simplices that are \textit{homotopic, relative their vertices}\index{relative homotopy}\index{homotopy!relative} (see~\fref{Figure-Fundamental-Groupoid}). The composition of two such classes is defined as in~\fref{Figure-Fundamental-Groupoid}: Using the Kan property four times fills the left hand side tetrahedron and yields a $1$-simplex $g\comp f$, of which the equivalence class is the needed composition of the classes of $f$ and $g$. This determines a functor $\pi_{1}$, left adjoint to $\ner$.
\begin{figure}
\begin{center}
\resizebox{\textwidth}{!}{\includegraphics{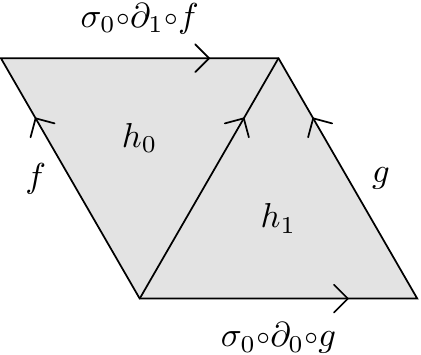}\qquad\qquad\includegraphics{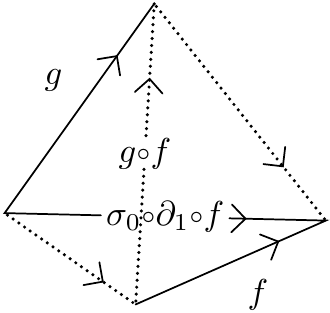}\qquad\qquad\includegraphics{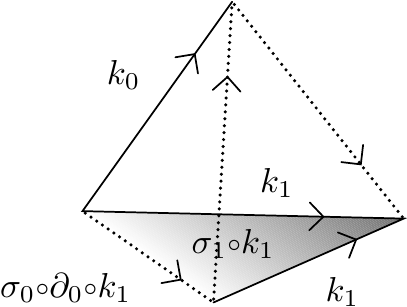}}
\end{center}
\caption{A relative homotopy from $f$ to $g$, and two applications of the Kan property.}\label{Figure-Fundamental-Groupoid}
\end{figure}
\end{morale}
\begin{proof}
Given a simplicial object $A$ in $\Ac$, its \emph{fundamental groupoid} has $A_{0}$ as objects of objects and $\overline{A_{1}}=A_{1}/I[d_{2}]$ as object of morphisms. We show that this reflexive graph $\pi_{1}A$ it is indeed a groupoid; then the left adjointness property quite easily follows. We shall use the fact that if $\Ac$ is a finitely cocomplete regular Mal'tsev category, then $\Cat \Ac$ is reflective in $\Rgrph \Ac$. The reflector $\Rgrph \Ac\to \Cat \Ac$ maps a reflexive graph $d_{0},d_{1}\colon {B_{1}\to B_{0}}$ to the groupoid with object of objects $B_{0}$ and object of morphisms $B_{1}/[R[d_{0}],R[d_{1}]]$ (see~\cite{Borceux-Bourn} and Proposition~\ref{Proposition-Construction-Pedicchio-Commutator}). 

In the diagram with exact rows
\[
\xymatrix{0 \ar[r] & K[q'] \ar@{.>}[dd]_-{\del_{2}\comp \del_{0}\comp \ker q'} \ar@{{ >}->}[rr]^-{\ker q'} && A_{3} \ar[d]^-{\del_{0}} \ar[r]^-{q'} & R[\del_{0}]\times_{A_{1}}R[\del_{1}] \ar@{.>}[dd]^-{\varphi'} \ar@{.>}[r] & 0\\
&&& A_{2} \ar[d]^-{\del_{2}}\\
0 \ar[r] & {I[d_{2}]} \ar@{{ >}->}[rr]_-{\ker \del_{0}\comp \im d_{2}} && A_{1} \ar@{-{ >>}}[r]_-{q} & {\overline{A_{1}}} \ar[r] & 0,}
\]
by the Kan property of $A$, the arrow $q'= ((\del_{0},\del_{1})\comp \del_{2},(\del_{1},\del_{2})\comp \del_{1})$ is a regular epimorphism. Also, the morphism $\del_{2}\comp \del_{0}\comp \ker q'$ factors over $I[d_{2}]$: Such a factorization is induced by the $(3,0)$-horn $(0,\del_{2}\comp  \del_{2}\comp \del_{0}\comp \ker q',\del_{3}\comp \del_{2}\comp \del_{0}\comp \ker q')$ in $A$.

This gives rise to the arrow $\varphi'$. To prove that $\pi_{1}A$, as a quotient of a groupoid, is itself a groupoid, we need a suitable regular epimorphism from $A_{1}/[R[\del_{0}],R[\del_{1}]]$ to $A_{1}/I[d_{2}]$. Such a morphism exists as soon as 
\[
\begin{cases} \varphi'\comp l_{R[\del_{0}]}= q\comp k_{0}\colon & R[\del_{0}]\to \overline{A_{1}}\\
\varphi'\comp r_{R[\del_{1}]}= q\comp k_{1}\colon & R[\del_{1}]\to \overline{A_{1}}. \end{cases}
\]

To see that these equalities hold, it suffices to note that $I[d_{2}]$ is the normalization of the reflexive relation $(R,r_{0},r_{1})$ on $A_{1}$, defined as follows: It is the image of $(\del_{0}\comp \pr_{1},\del_{2}\comp \pr_{2})$, a morphism from
\[
\resizebox{\textwidth}{!}{\mbox{$\{ (h_{0},h_{1})\in  A_{2}\times A_{2}\,|\,\del_{1}\comp h_{0}=\del_{1}\comp h_{1},\,\del_{2}\comp h_{0}=\sigma_{0}\comp \del_{1}\comp \del_{0}\comp h_{0},\, \del_{0}\comp h_{1}=\sigma_{0}\comp \del_{0}\comp \del_{2}\comp h_{1}\}$}}
\]
to $A_{1}\times A_{1}$ (cf.~\fref{Figure-Fundamental-Groupoid}).

Indeed, then $q$ is a coequalizer of $r_{0}$ and $r_{1}$; using the Kan property on the right hand side tetrahedron of~\fref{Figure-Fundamental-Groupoid}, we see that $l_{R[\del_{0}]}$ factors over $q'$: There exists a morphism $y\colon {Y\to A_{3}}$ and a regular epimorphism $p\colon {Y\to R[\del_{0}]}$ satisfying $l_{R[\del_{0}]}\comp p=q'\comp y$, and such that $(\del_{2}\comp \del_{0}\comp y, k_{0}\comp p)$ factors over $R$. This proves that the equality $\varphi'\comp l_{R[\del_{0}]}= q\comp k_{0}$ holds:
\[
\varphi'\comp l_{R[\del_{0}]}\comp p=\varphi'\comp q'\comp y=q\comp \del_{2}\comp \del_{0}\comp y=q\comp k_{0}\comp p.
\]
A similar argument shows that $\varphi'\comp r_{R[\del_{1}]}= q\comp k_{1}$.

The functor $\pi_{1}\colon \Sc \Ac\to \Cat \Ac $ is left adjoint to the functor $\ner$: The unit $\eta \colon {1_{\Sc \Ac}\To \ner \comp\pi_{1}}$ of the adjunction is $(\eta_{A})_{0}=1_{A_{0}}$, $(\eta_{A})_{1}=q\colon {A_{1}\to \overline{A_{1}}}$, 
\[
(\eta_{A})_{2}= q\comp (\del_{2},\del_{0})\colon {A_{2}\to \ner_{2}\pi_{1} (A)},
\]
and $(\eta_{A})_{n}= (q\comp\del_{2}^{n},(\eta_{A})_{n-1}\comp \del_{0})$ for $n>2$.
\end{proof}

\begin{proposition}\label{Proposition-Characterization-Weak-Equivalence}
Let $\Ac$ be a semi-abelian category and $\vetf\colon \Af \to \Bf$ a functor in $\Ac$.
\begin{enumerate}
\item $\vetf$ is fully faithful if and only if $H_{0}\vetf$ is mono and $H_{1}\vetf$ is iso.
\item $\vetf$ is essentially $\regepi$-surjective if and only if $H_{0}\vetf$ is a regular epimorphism.
\end{enumerate}
Hence an internal functor is a weak equivalence exactly when it is a homology isomorphism.
\end{proposition}
\begin{proof}
In the diagram
\[
\xymatrix@1{H_{1}\Af \ar[d]_{H_{1}\vetf} \ar@{{ >}->}[rr]^-{\ker (d_{0},d_{1})} && A_{1}\ar@{}[rd]|{\texttt{(i)}} \ar`u[r]`[rr]^-{(d_{0},d_{1})}[rr] \ar[d]_-{f_{1}} \ar@{-{ >>}}[r] & \KP [q] \ar@{}[rd]|{\texttt{(ii)}} \ar@{.>}[d] \monr & A_{0}\times A_{0} \ar[d]^-{f_{0}\times f_{0}}\\
H_{1}\Bf \ar@{{ >}->}[rr]_-{\ker (d_{0},d_{1})} && B_{1} \ar`d[r]`[rr]_-{(d_{0},d_{1})}[rr] \ar@{-{ >>}}[r] & \KP [r] \monr & B_{0}\times B_{0}}
\]
the arrow $H_{1}\vetf$ is an isomorphism if and only if the square \texttt{(i)} is a pullback (Proposition~\ref{Proposition-LeftRightPullbacks}); square \texttt{(ii)} is a pullback if and only if \texttt{(iii)}
\[
\xymatrix@1{{\KP [q]} \ar@{.>}[d] \ar@{}[rd]|{\texttt{(iii)}} \ar@<.5ex>[r] \ar@<-.5ex>[r] & A_{0} \ar[d]^-{f_{0}} \ar@{-{ >>}}[rr]^-{q=\coeq (d_{0},d_{1})} && H_{0}\Af \ar[d]^-{H_{0}\vetf}\\
\KP [r] \ar@<.5ex>[r] \ar@<-.5ex>[r] & B_{0} \ar@{-{ >>}}[rr]_-{r=\coeq (d_{0},d_{1})} && H_{0}\Bf}
\]
is a joint pullback, which (by Proposition 1.1 in~\cite{Bourn2003}) is the case exactly when $H_{0}\vetf$ is mono. This already shows one implication of 1. To prove the other, note that if $\vetf$ is fully faithful, then $H_{1}\vetf$ is an isomorphism; hence \texttt{(i)} is a pullback, \texttt{(ii)} is a pullback~\cite[Proposition 2.5]{Janelidze-Kelly}, \texttt{(iii)} is a joint pullback and $H_{0}\vetf$ a monomorphism.

For the proof of 2.\ consider the following diagrams.
\[
\xymatrix{(P_{f})_{0}\ar@{}[rd]|<{\pullback} \ar[d]_-{\overline{d_{1}}} \ar[r]^-{\overline{f_{0}}}& B_{1} \ar[r]^-{d_{0}} \ar[d]^-{d_{1}} & B_{0} \\
A_{0} \ar[r]_-{f_{0}} & B_{0}}
\qquad \qquad 
\xymatrix{(P_{f})_{0} \ar@{}[rrd]|{\texttt{(iv)}} \ar[rr]^-{d_{0}\comp \overline{f_{0}}} \ar[d]_-{\overline{d_{1}}} && B_{0} \ar[d]^{r}\\
A_{0} \ar[rr]_-{H_{0}f\comp q} && H_{0}B}
\]
If $H_{0}f$ is a regular epimorphism then so is the bottom arrow in diagram \texttt{(iv)}. Accordingly, we only need to show that the induced arrow $p\colon (P_{f})_{0}\to P$ to the pullback $P$ of $H_{0}f\comp q$ along $r$ is regular epi. Now in the left hand side diagram
\[
\xymatrix{P \ar@{}[rd]|{\texttt{(v)}} \ar[d]_-{\overline{r}} \ar[r] & R[r] \ar@{}[rd]|<{\pullback} \ar[d]_-{p_{1}} \ar[r]^-{p_{0}} & B_{0} \ar[d]^-{r} \\
A_{0} \ar[r]_-{f_{0}} & B_{0} \ar[r]_-{r} & H_{0}B}
\quad \quad 
\xymatrix{(P_{f}) \ar@{}[rd]|{\texttt{(vi)}} \ar[r]^-{p} \ar[d]_-{\overline{f_{0}}} & P \ar@{}[rd]|<{\pullback} \ar@{}[rd]|{\texttt{(v)}} \ar[r]^-{\overline{r}} \ar[d] & A_{0} \ar[d]^-{f_{0}}\\
B_{1} \ar[r]_-{o} & R[r] \ar[r]_-{p_{1}} & B_{0}}
\]
this pullback is the outer rectangle; its left hand side square \texttt{(v)} is a pullback. Now so is \texttt{(vi)}; it follows that $p$ is a regular epimorphism, so being the arrow $o$, universal for the equalities $p_{0}\comp o=d_{0}$ and $p_{1}\comp o=d_{1}$ to hold.

Conversely, note that $r\comp d_{0}\comp \overline{f_{0}}=H_{0}f\comp q\comp \overline{d_{1}}$; hence if $d_{0}\comp \overline{f_{0}}$ is regular epic, so is $H_{0}f$.
\end{proof}

The notion of \textit{Kan fibration}\index{Kan fibration} makes sense in the context of regular categories: See Section~\ref{Section-Simplicial-Objects}. We use it to characterize the fibrations in $\Cat \Ac$.

\begin{proposition}\label{Proposition-Characterization-Fibration-Kan}
A functor $\vetp \colon  \Ef\to \Bf$ in a semi-abelian category $\Ac$ is a fibration if and only if $\ner \vetp $ is a Kan fibration.
\end{proposition}
\begin{proof}
First note that, because a semi-abelian category is always Mal'tsev, every category in $\Ac$ is an internal groupoid. We may now use M. Barr's Embedding Theorem for regular categories~\cite{Barr} in the form of Metatheorem A.5.7 in~\cite{Borceux-Bourn}. Indeed, the properties ``some internal functor is a fibration'' and ``some simplicial morphism is a Kan fibration'' may be added to the list of properties~\cite[0.1.3]{Borceux-Bourn}, and it is well-known that in $\Set$, a functor between two groupoids is a fibration if and only if its nerve is a Kan fibration.
\end{proof}

In his paper~\cite{Janelidze}, Janelidze introduces a notion of crossed module in an arbitrary semi-abelian category $\Ac$. Its definition is based on Bourn and Janelidze's notion of internal semi-direct product~\cite{Bourn-Janelidze:Semidirect}. Internal crossed modules also generalize the case where $\Ac =\Gp$ in the sense that an equivalence $\XMod \Ac \simeq \Cat \Ac$ still exists. Using this equivalence, we may transport the model structures from Theorem~\ref{Theorem-Model-Structure} to the category $\XMod \Ac$. In case $\Ac$ has enough regular projectives and $\Tc$ is the regular epimorphism topology, this has the advantage that the classes of fibrations, cofibrations and weak equivalences have a very easy description.

We recall from~\cite{Janelidze} the definition of internal crossed modules. Given two morphisms $g\colon {B\to B'}$ and $h\colon {X\to X'}$, the map $g\flat h\colon {B\flat X\to B'\flat X'}$ is unique in making the diagram
\[
\xymatrix{0 \ar[r] & B\flat X \ar@{.>}[d]_-{g\flat h} \ar@{{ >}->}[r]^-{\kappa_{B,X}} & B+ X \ar[d]^-{g+h} \ar@{-{ >>}}[r]^-{[1_{B},0]} & B \ar[r] \ar[d]^-{g} & 0\\
0 \ar[r] & B'\flat X' \ar@{{ >}->}[r]_-{\kappa_{B',X'}} & B'+ X' \ar@{-{ >>}}[r]_-{[1_{B'},0]} & B' \ar[r] & 0}
\]\index{$B\flat X$}
commute. The category $\SplitEpi \Ac $\index{category!$\SplitEpi \Ac $}\index{category!of split epimorphisms} of split epimorphisms in $\Ac$ with a given splitting ($\Pt \Ac$\index{category!$\Pt \Ac$} in~\cite{Bourn-Janelidze:Semidirect} and in Definition 2.1.14 of~\cite{Borceux-Bourn}) is equivalent to the category $\Act \Ac$ of \emph{actions}\index{action}\index{internal!action}\index{category!of actions}\index{category!$\Act \Ac$} in $\Ac$, of which the objects are triples $(B,X,\xi)$, where $\xi \colon B\flat X\to X$ makes the following diagram commute:
\[
\xymatrix{B\flat (B\flat X) \ar[r]^-{\mu^{B}_{X}} \ar[d]_-{1_{B}\flat \xi} & B\flat X \ar[d]_-{\xi } & X \ar[l]_-{\eta_{X}^{B}} \ar@{=}[ld]\\
B\flat X \ar[r]_-{\xi} & X;  }
\]
here $\mu^{B}_{X}$ is defined by the exactness of the rows in the diagram\index{in@$\inj_{1}\colon X\to X+X$}
\[
\xymatrix{0 \ar[r] & B \flat(B\flat X) \ar@{.>}[d]_-{\mu_{X}^{B}} \ar@{{ >}->}[r]^-{\kappa_{B,B\flat X}} & B+ (B\flat X) \ar[d]^-{[\inj_{1},\kappa_{B,X}]} \ar@{-{ >>}}[r]^-{[1_{B},0]} & B \ar[r] \ar@{=}[d] & 0\\
0 \ar[r] & B\flat X \ar@{{ >}->}[r]_-{\kappa_{B,X}} & B+ X \ar@{-{ >>}}[r]_-{[1_{B},0]} & B \ar[r] & 0}
\]
and $\eta^{B}_{X}$ is unique such that $\kappa_{B,X}\comp \eta^{B}_{X}=\inj_{2}\colon X\to B+X$. A morphism 
\[
{(B,X,\xi)\to (B',X',\xi ')}
\]
in $\Act \Ac$ is a pair $(g,h)$, where $g\colon {B\to B'}$ and $h\colon {X\to X'}$ are morphisms in $\Ac$ with $h\comp \xi =\xi '\comp (g\flat h)$.

An \emph{internal precrossed module}\index{internal!precrossed module}\index{precrossed module!internal ---} in $\Ac$ is a 4-tuple $(B,X,\xi ,f)$ with $(B,X,\xi)$ in $\Act \Ac$ and $f\colon X\to B$ a morphism in $\Ac$ such that the left hand side diagram
\[
\xymatrix{B\flat X \ar@{{ >}->}[r]^-{\kappa_{B,X}} \ar[d]_-{\xi} & B+ X  \ar[d]^-{[1_{B},f]}&& (B+ X)\flat X \ar[rr]^-{[1_{B},f]\flat 1_{X}} \ar[d]_-{[1_{B+ X},\inj_{2}]^{\#}} && {B\flat X} \ar[d]^-{\xi}\\
X \ar[r]_-{f} & B && {B\flat X} \ar[rr]_-{\xi} && X}
\]
commutes. A morphism $(B,X,\xi ,f)\to (B',X',\xi' ,f')$ is a morphism 
\[
(g,h)\colon {(B,X,\xi)\to (B',X',\xi')}
\]
in $\Act \Ac$ such that $g\comp f=f'\comp h$. There is an equivalence between $\Rgrph \Ac$\index{category!$\Rgrph \Ac$} and the category $\PXMod \Ac$\index{category!$\PXMod \Ac$} of precrossed modules in $\Ac$. An \emph{internal crossed module}\index{internal!crossed module}\index{crossed module!internal ---} in $\Ac$ is an internal precrossed module $(B,X,\xi ,f)$ in $\Ac$ for which the right hand side diagram above commutes. Here $[1_{B+ X},\inj_{2}]^{\#}$ is the unique morphism such that $\kappa_{B,X}\comp [1_{B+ X},\inj_{2}]^{\#}=[1_{B+X},\inj_{2}]\comp \kappa_{B+X,X}$. If $\XMod \Ac$\index{category!$\XMod \Ac$} or $\CrossMod \Ac$ denotes the full subcategory of $\PXMod \Ac$ determined by the crossed modules, then the last equivalence (co)restricts to an equivalence $\Cat \Ac \simeq \XMod \Ac$. This equivalence maps an internal category $\Af$ to the internal crossed module 
\[
(A_{0},\quad K[d_{1}],\quad \xi\colon A_{0}\flat K[d_{1}]\to K[d_{1}],\quad d_{0}\comp \ker d_{1}\colon K[d_{1}]\to A_{0})
\]
where $\xi$ is the pullback of 
\[
[i,\ker d_{1}]\colon A_{0}+K[d_{1}]\to A_{1}
\]
along $\ker d_{1}\colon {K[d_{1}]\to A_{1}}$. Thus an internal functor $\vetf\colon \Af \to \Bf$ is mapped to the morphism 
\[
(f_{0},N_{1} \vetf )\colon (A_{0},N_{1}\Af,\xi ,d_{0}\comp \ker d_{1})\to (B_{0}, N_{1}\Bf,\xi',d_{0}\comp \ker d_{1}).
\]

If the homology objects of an internal crossed module $(B,X,\xi ,f)$\index{homology!of internal crossed modules} are those of the associated internal category, the only non-trivial ones are 
\[
H_{0} (B,X,\xi ,f)=\coker f
\]
and $H_{1} (B,X,\xi ,f)=\ker f$ (see Proposition~\ref{Proposition-Characterization-Homology}).

\begin{theorem}\label{Theorem-Categorical-Model-Structure}
If $\Ac$ is a semi-abelian category with enough projectives then a model category structure on $\XMod \Ac$ is defined by choosing $\we$ the class of homology isomorphisms, $\cof$ the class of morphisms $(g,h)$ with $g$ in ${^{\Box}\regepi}$ and $\fib$ the class of morphisms $(g,h)$ with $h$ regular epic.
\end{theorem}
\begin{proof}
Only the characterization of $\fib$ needs a proof. Given an internal functor $\vetp\colon \Ef \to \Bf$, consider the diagram with exact rows
\[
\xymatrix{0 \ar[r] & K[d_{1}] \ar@{}[rd]|<{\pullback} \ar@{}[rd]|{\texttt{(i)}} \ar[d]_-{N_{1} \vetp }\ar@{{ >}->}[r]^-{\ker d_{1}} & E_{1} \ar@{.>}[d]^-{(r_{p})_{0}} \ar@{-{ >>}}[r]^-{d_{1}} & E_{0} \ar[r] \ar@{=}[d] & 0\\
0 \ar[r] & K[d_{1}] \ar@{=}[d] \ar@{{ >}->}[r] & (P_{p})_{0} \ar@{}[rd]|<{\pullback} \ar@{}[rd]|{\texttt{(ii)}} \ar@{-{ >>}}[r] \ar[d] & E_{0} \ar[d]^-{p_{0}} \ar[r] & 0\\
0 \ar[r] & K[d_{1}] \ar@{{ >}->}[r]_-{\ker d_{1}} & B_{1} \ar@{-{ >>}}[r]_-{d_{1}} & B_{0} \ar[r] & 0}
\]
where squares \texttt{(i)} and \texttt{(ii)} are a pullbacks (cf.\ Lemma 4.2.2 and its Corollary 4.2.3 in~\cite{Borceux-Bourn}). By~\cite[Lemma 4.2.5]{Borceux-Bourn}, $N_{1} \vetp$ is a regular epimorphism if and only if $(r_{p})_{0}$ is regular epic, i.e., exactly when $\vetp$ is a $\regepi$-fibration.
\end{proof}

\begin{example}[Crossed modules of groups]\label{Example-XMod}
In the specific case of $\Ac$ being the semi-abelian category $\Gp$ of groups and group homomorphisms, the structure on $\XMod \Ac$ coincides with the one considered by Garz\'on and Miranda in~\cite{Garzon-Miranda}. This category is equivalent to the category $\XMod=\XMod\Gp$ of crossed modules of groups. As such, $\XMod$\index{category!$\XMod $} carries the model structure from Theorem~\ref{Theorem-Categorical-Model-Structure}.
\end{example}

\section{Case study: the trivial topology}\label{Section-Split-Epimorphisms}\index{topology!trivial}
In this section we suppose that $\Tc$ is the trivial topology on $\Cc$ (as a rule, we shall not mention it), and give a more detailed description of the model structure from Theorem~\ref{Theorem-Model-Structure}. It turns out to resemble very much Str\o m's model category structure on the category $\Top$\index{category!$\Top$} of topological spaces and continuous maps~\cite{Strom:Homotopy-category}: Its weak equivalences are the homotopy equivalences, its cofibrations are the functors that have the homotopy extension property and its fibrations are the functors that have the homotopy lifting property.

\begin{proposition}\label{Proposition-Trivial-Fibration-Split-Epi}
An internal functor is a trivial fibration if and only if it is a split epimorphic equivalence. 
\end{proposition}
\begin{proof}
This is an immediate consequence of Proposition~\ref{Proposition-Characterization-Trivial-Fibrations}. 
\end{proof}

It follows that every object of this model category is cofibrant. Since we already showed them to be fibrant as well, the notion of homotopy induced by the model structure (see Quillen~\cite{Quillen}) is determined entirely by the cocylinder from Section~\ref{Section-Cocylinder}. Hence the weak equivalences are homotopy equivalences, also in the model-category-theoretic sense of the word.  

\begin{definition}\label{Definition-Fibration}
Let $\Cc$ be a finitely complete category. An internal functor $\vetp\colon \Ef \to \Bf$ is said to have the \emph{homotopy lifting property}\index{homotopy lifting property} if and only if the diagram
\[
\xymatrix{{\Ef^{\If}} \ar[r]^-{\vetp^{\If}} \ar[d]_-{\epsilon_{1} (\Ef)} & {\Bf^{\If}} \ar[d]^-{\epsilon_{1} (\Bf)}\\
{\Ef} \ar[r]_-{\vetp} & {\Bf}}
\]
is a weak pullback in $\Cat \Cc $.
\end{definition}

\begin{proposition}\label{Proposition-Fibration-is-Star-Surjective}
An internal functor is a strong fibration if and only if it has the homotopy lifting property (cf.\ R.~Brown~\cite{Brown:Fibrations}).
\end{proposition}
\begin{proof}
This follows from Proposition~\ref{Proposition-Trivial-Fibration-Split-Epi}, Corollary~\ref{Corollary-Characterization-Fibrations} and the two-out-of-three property of strong equivalences.
\end{proof}

We now give a characterization of cofibrations in the following terms.

\begin{definition}\label{Definition-Cofibration}
An internal functor $\vetj\colon \Af \to \Xf$ is said to have the \emph{homotopy extension property}\index{homotopy extension property} if and only if any commutative square
\begin{equation}\label{Diagram-Cofibration}
\vcenter{\xymatrix{{\Af} \ar[r]^-{\vetH} \ar[d]_-{\vetj} & {\Ef^{\If}} \ar[d]^-{\epsilon_{1} (\Ef)}\\
{\Xf} \ar@{.>}[ru]^-{\overline{\vetH}} \ar[r]_-{\vetf}& {\Ef}}}
\end{equation}
has a lifting $\overline{\vetH}\colon \Xf \to \Ef^{\If}$.
\end{definition}

To show that a functor with the homotopy extension property is a cofibration, we need a way to approximate any trivial fibration $\vetp$ with some $\epsilon_{1} (\Ef)$. The next construction allows (more or less) to consider $p_{0}$ as the object morphism of $\epsilon_{1} (\vetp\Ef)\colon (\vetp\Ef)^{\If}\to \vetp\Ef$. 

\begin{lemma}\label{Lemma-Composed-Groupoid}
If $\vetp\colon \Ef \to \Bf$ is an adjoint equivalence with right inverse $\vets\colon {\Bf \to \Ef}$ then the reflexive graph
\[
\vetp\Ef = \bigl(\xymatrix{E_{1} \ar@<-1 ex>[rr]_-{p_{0}\comp d_{0}} \ar@<1 ex>[rr]^-{p_{0}\comp d_{1}} && B_{0} \ar[ll]|-{i\comp s_{0}}}\bigr)
\]
carries an internal category structure.
\end{lemma}
\begin{proof}
The morphism $m\comp (p_{1}\times_{1_{B_{0}}}p_{1})\colon E_{1}\times_{B_{0}}E_{1}\to B_{1}$ is such that $d_{0}\comp m\comp (p_{1}\times_{1_{B_{0}}}p_{1})=p_{0}\comp d_{0}\comp \pr_{1}$ and $d_{1}\comp m\comp (p_{1}\times_{1_{B_{0}}}p_{1})=p_{0}\comp d_{1}\comp \pr_{2}$. Since $\vetp$ is an adjoint equivalence, it is a fully faithful functor; by Proposition~\ref{Proposition-Full-Faithful-Functor}, a unique morphism $m\colon E_{1}\times_{B_{0}}E_{1}\to E_{1}$ exists satisfying $m\comp (p_{1}\times_{1_{B_{0}}}p_{1})=p_{1}\comp m$, $d_{0}\comp m=d_{0}\comp \pr_{1}$ and $d_{1}\comp m=d_{1}\comp \pr_{2}$. This is the needed structure of internal category.
\end{proof}

\begin{proposition}\label{Proposition-LLP}
An internal functor is a cofibration if and only if it has the homotopy extension property.
\end{proposition}
\begin{proof}
One implication is obvious, $\epsilon_{1} (\Ef)$ being a trivial fibration. Now consider a commutative square such as~\ref{Diagram-Lifting} above, where $\vetj$ has the homotopy extension property and $\vetp\in \fib \cap\we$. By the hypothesis on $\vetj$ and Proposition~\ref{Proposition-Lifting-we}, $j_{0}$ has the left lifting property with respect to all morphisms of the form $\delta_{1}\colon {\iso (A)\to A_{0}}$ for $\Af \in \Cat \Cc$. In particular, by Lemma~\ref{Lemma-Composed-Groupoid}, $j_{0}$ has the left lifting property with respect to $p_{0}\comp d_{1}\comp j=\epsilon_{1} (pE)_{0}\colon {\iso (pE)\to B_{0}}$. If 
\[
h'_{0}\colon {X_{0}\to \iso (pE)}
\]
denotes a lifting for the commutative square $g_{0}\comp j_{0}= (p_{0}\comp d_{1}\comp j)\comp (\iota\comp f_{0})$, then 
\[
h_{0}=d_{1}\comp j\comp h'_{0}\colon {X_{0}\to E_{0}}
\]
is such that $p_{0}\comp h_{0}=g_{0}$ and $h_{0}\comp j_{0}=f_{0}$. This $h_{0}$ induces the needed lifting for diagram~\ref{Diagram-Lifting}. 
\end{proof}

Proposition~\ref{Proposition-Trivial-Fibration-Split-Epi} may now be dualized as follows:

\begin{proposition}\label{Proposition-Trivial-Cofibrations}
An internal functor is a trivial cofibration if and only if it is a split monomorphic equivalence.
\end{proposition}
\begin{proof}
One implication is Proposition~\ref{Proposition-Characterization-Trivial-Cofibrations}. To prove the other, consider a commutative square such as~\ref{Diagram-Cofibration} above, and suppose that $\vetj\colon \Af \to \Xf$ is an equivalence with left adjoint left inverse $\vetk\colon \Xf \to \Af$. Let $\vetK\colon \Xf\to \Xf^{\If}$ denote a homotopy from $\vetj\comp \vetk$ to $1_{\Xf}$. Put $\overline{H}_{0}=m\comp (H_{0}\comp k_{0}, f_{1}\comp K_{0})\colon X_{0}\to \iso (E)$, then $\delta_{1}\comp \overline{H}_{0}=f_{0}$ and $\overline{H}_{0}\comp j_{0}=H_{0}$. Proposition~\ref{Proposition-Lifting-we} now yields the needed lifting $\overline{\vetH}\colon {\Xf \to \Ef^{\If}}$.
\end{proof}

\begin{theorem}\label{Theorem-Model-Structure-Trivial}
If $\Cc$ is a finitely complete category such that $\Cat \Cc$ is finitely cocomplete (cf.\ Section \ref{Section-Preliminaries}) then a model category structure is defined on $\Cat \Cc$ by choosing $\we$ the class of homotopy equivalences, $\cof$ the class of functors having the homotopy extension property and $\fib$ the class of functors having the homotopy lifting property.\noproof 
\end{theorem}
 % Internal categories
\chapter{Homotopy of simplicial objects}\label{Chapter-Simplicial-Objects}

\setcounter{section}{-1}
\section{Introduction}\label{Section-Simp-Introduction}

The starting point of this chapter is the following theorem due to Quillen.

\begin{theorem}\cite[Theorem II.4.4]{Quillen}\label{Theorem-Quillen}
Let $\Ac$ be a finitely complete category with enough regular projectives. Let $\simpA$ be the category of simplicial objects over $\Ac$. Define a map $f$ in $\simpA$ to be a fibration (resp.\ weak equivalence) if $\hom (P,f)$ is a fibration (resp.\ weak equivalence) in $\Sc \Set$ for each projective object $P$ of $\Ac$, and a cofibration if $f$ has the left lifting property with respect to the class of trivial fibrations. Then $\simpA$ is a closed simplicial model category if $\Ac$ satisfies one of the following extra conditions:
\begin{itemize}
\item[\emph{(*)}] every object of $\simpA$ is fibrant;
\item[\emph{(**)}] $\Ac$ is cocomplete and has a set of small regular projective generators.\noproof 
\end{itemize}
\end{theorem}

The condition (**) concerns finitary varieties in the sense of universal algebra. We shall focus on categories satisfying the first condition: This class contains all abelian categories; more generally, a regular category satisfies (*) exactly when it is Mal'tsev.

In particular, if $\Ac$ is semi-abelian and has enough projectives, $\simpA$ carries a model category structure. Our aim is to describe it explicitly, and to explain how this model structure is compatible with homology of simplicial objects, as defined in Chapter~\ref{Chapter-Homology}.

\section{Fibrations and weak equivalences}\label{Section-Weak-Equiv-Trivial-Fibrations}
In this section we give an explicit description of the weak equivalences and the (co)fibrations occurring in Quillen's model structure: In a regular Mal'tsev category with enough projectives, the fibrations are just the Kan fibrations; when the category is semi-abelian, weak equivalences are homology isomorphisms.

In fact, a priori, it is not entirely clear that Theorem~\ref{Theorem-Quillen} is exactly the same as Theorem~II.4.4 in~\cite{Quillen}, because there the projectives are chosen relative to the class of so-called \emph{effective} epimorphisms\index{effective!epimorphism}\index{epimorphism!effective}. A morphism $f\colon {A\to B}$ in a finitely complete category $\Ac$ belongs to this class when for every object $T$ of $\Ac$, the diagram in $\Set$ 
\[
\xymatrix{{\hom (B,T)} \ar[rr]^-{(\cdot)\comp f} && {\hom (A,T)} \ar@<.5ex>[rr]^-{(\cdot)\comp k_{0}} \ar@<-.5ex>[rr]_-{(\cdot)\comp k_{1}} && {\hom (R[f],T)}}
\]
is an equalizer. This, however, just means that $f$ is a coequalizer of its kernel pair; hence in $\Ac$, regular epimorphisms and effective epimorphisms coincide.

\subsection*{Fibrations are Kan fibrations}\label{Subsection-Fibrations-Kan-Fibrations}
Recall from Chapter~\ref{Chapter-Homology} the following definition: A simplicial object $K$ in a regular category $\Ac$ is called \emph{Kan} when, for every $(n,k)$-horn 
\[
x= (x_{i}\colon X\to K_{n-1})
\]
of $K$, there is a regular epimorphism $p\colon {Y\to X}$ and a map $y\colon {Y\to K_{n}}$ such that $\del_{i}\comp y=x_{i}\comp p$ for $i\neq k$. A map $f \colon {A\to B}$ of simplicial objects is a \emph{Kan fibration} if, for every $(n,k)$-horn 
\[
x= (x_{i}\colon X\to A_{n-1})
\]
of $A$ and every $b\colon {X\to B_{n}}$ with $\del_{i}\comp b=f_{n-1}\comp x_{i}$ for all $i\neq k$, there is a regular epimorphism $p\colon {Y\to X}$ and a map $a\colon {Y\to A_{n}}$ such that $f_{n}\comp a=b\comp p$ and $\del_{i}\comp a=x_{i}\comp p$ for all $i\neq k$.

Also recall Remark~\ref{Remark-Kan-Projectives} that in a regular category $\Ac$ with enough regular projectives, a simplicial object $K$ is Kan if and only if for every projective object $P$ and for every $(n,k)$-horn 
\[
x= (x_{i}\colon P\to K_{n-1})
\]
of $K$, there is a morphism $y\colon P\to K_{n}$ such that $\del_{i}\comp y=x_{i}$, $i\neq k$. This means that $K$ is Kan if and only if $\hom (P,K)=\hom (P,\cdot)\comp K$ is a Kan simplicial set, for every projective object $P$. The similar statement for Kan fibrations amounts to the following:

\begin{proposition}\label{Proposition-Fibration-is-Kan-Fibration}\index{Kan fibration}\index{Kan simplicial object}\index{fibration!Kan ---}
Let $\Ac$ be a regular category with enough projectives. A simplicial morphism $f\colon A\to B$ in $\Ac$ is a Kan fibration if and only if it is a fibration in the sense of Theorem~\ref{Theorem-Quillen}: For every projective object $P$ of $\Ac$, the induced morphism of simplicial sets
\[
f\comp (\cdot)=\hom (P,f)\colon \hom (P,A)\to \hom (P,B)
\]
is a Kan fibration.\noproof 
\end{proposition}

Thus we see that a regular category $\Ac$ satisfies condition (*) if and only if every simplicial object in $\Ac$ is Kan.\index{simplicial object!Kan ---} In the introduction to~\cite{Barr}, M. Barr conjectures that the latter condition on $\Ac$ means that every reflexive relation in $\Ac$ is an equivalence relation. In~\cite{Carboni-Kelly-Pedicchio}, Carboni, Kelly and Pedicchio show the equivalence of this latter condition with the permutability condition $RS=SR$ for arbitrary equivalence relations $R$, $S$ on an object $X$ in $\Ac$---and this is what they call the \emph{Mal'tsev condition}. They moreover show that a regular category $\Ac$ is Mal'tsev if and only if every reflexive relation in $\Ac$ is an equivalence relation, and thus they prove M. Barr's conjecture. Recall that this last property is the definition of \emph{Mal'tsev category} from Section~\ref{Section-Mal'tsev-Categories}.

Briefly, simplicial objects in a regular Mal'tsev category with enough projectives admit a model category structure, where the fibrations are just Kan fibrations. However, stronger conditions on $\Ac$ are needed for an elementary description of the weak equivalences. We shall now focus on semi-abelian categories and describe the weak equivalences using simplicial homology. 

Our first aim is characterizing regular epic homology isomorphisms. To do so, we need a precise description of the \emph{acyclic}\index{acyclic simplicial object}\index{simplicial object!acyclic} simplicial objects, i.e., those $A$ in $\simpA$ satisfying $H_{n}A=0$ for all $n\in \N$.

\subsection*{Acyclic objects}\label{Subsection-Acyclic-Objects}

What does it mean, for a simplicial object, to have zero homology? We give an answer to this question, not in terms of an associated chain complex, but in terms of the simplicial object itself. From now on, we suppose that $\Ac$ is a semi-abelian category.

\begin{figure}
\begin{center}
\resizebox{\textwidth}{!}{\includegraphics{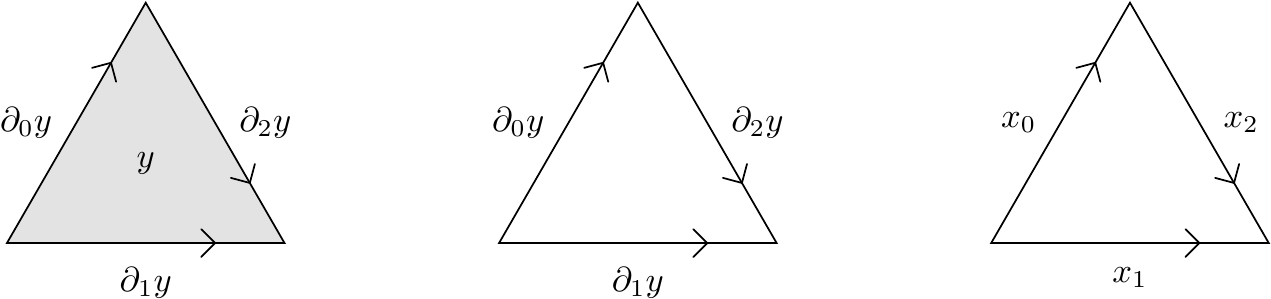}}
\end{center}
\caption{A $2$-simplex $y$, its boundary $\del y$ and a $1$-cycle $(x_{0},x_{1},x_{2})$.}\label{Figure-Acyclic-Simplicial-Object}
\end{figure}

By Corollary~\ref{Corollary-H_0}, for any simplicial object $A$, $H_{0}A=\Coeq [\del_{0},\del_{1}]$. Hence $\Ac$ being exact Mal'tsev implies that $H_{0}A=0$ if and only if 
\[
(\del_{0},\del_{1})\colon {A_{1}\to A_{0}\times A_{0}}
\]
is a regular epimorphism. A similar property holds for the higher homology objects. In order to write it down, we need the following definition (cf.\ Definition~\ref{Definition-Kan} and \fref{Figure-Acyclic-Simplicial-Object}).

\begin{definition}\label{Definition-Cycles}
Let $A$ be a simplicial object in $\Ac$. The object $\cycle_{n}A$\index{$\cycle_{n}A$} of \emph{$n$-cycles in $A$}\index{cycle} is defined by $\cycle_{0}A=A_{0}\times A_{0}$ and 
\[
\cycle_{n}A=\{(x_{0},\dots,x_{n+1})\,|\, \text{$\del_{i}\comp x_{j}=\del_{j-1}\comp x_{i}$, for all $i<j\in [n+1]$}\}\subset A_{n}^{n+2},
\]
for $n>0$. This gives rise to functors $\cycle_{n}\colon {\Sc \Ac \to \Ac}$.
\end{definition}

Recall from Notation~\ref{Notation-A^-}\index{A@$A^{-}$}\index{lambda@$\El A$} the definition of the simplicial object $A^-$ associated to a simplicial object $A$, and consider the resulting simplicial morphism $\del=(\del_{0})_n\colon {A^-\to A}$. Write $\El A$ for its kernel. The following result due to J.~Moore, well-known to hold in the abelian case (cf.\ Exercise 8.3.9 in~\cite{Weibel}) is easily generalized to semi-abelian categories:

\begin{lemma}\label{Lemma-Lambda}
For any $n\in \N $, $N_{n} (\El A)\cong N_{n+1}A$; hence $H_{n} (\El A)\cong H_{n+1}A$.\noproof 
\end{lemma}

We also have to extend the notion of regular pushout to a more general situation.

\begin{definition}\label{Definition-regepi-Pushout}
A square
\begin{equation}\label{Square-Regepi-Pushout}
\vcenter{\xymatrix{ A' \ar@{-{ >>}}[r]^-{f'} \ar[d]_-{v} & B' \ar[d]^-w \\
A \ar@{-{ >>}}[r]_-f & B}}
\end{equation}
in a semi-abelian category is called a \emph{$\regepi$-pushout}\index{pushout!$\regepi$- ---}\index{regepi-pushout@$\regepi$-pushout}\index{regular pushout} when the comparison map to the associated pullback is a regular epimorphism (cf.\ Diagram~\ref{Diagram-Regular-Pushout}). (The maps $v$ and $w$ are \textit{not} demanded to be regular epimorphisms.)
\end{definition}

\begin{remark}\label{Remark-Regepi-Pushout}
Using the Short Five Lemma for regular epimorphisms (e.g., Lem\-ma~4.2.5 in \cite{Borceux-Bourn}), one sees that the square~\ref{Square-Regepi-Pushout} is a $\regepi$-pushout if and only if the induced arrow ${K[f']\to K[f]}$ is a regular epimorphism. Any $\regepi$-pushout is a pushout, but a pushout need not be a $\regepi$-pushout: A counterexample may be constructed by choosing $B'=0$, $v=\inj_{2}\colon {A'\to B+A'}$ and $f=[1_{B},0]\colon {B+A'\to B}$ in, say, $\Gp$.
\end{remark}

\begin{lemma}\label{Lemma-Cycle-and-Pullbacks}
For any $n\in \N$, the diagram
\[
\xymatrix{{\cycle_{n+1}A}\ar[d]_-{(\pr_{i+2})_{i\in [n+1]}} \ar[r]^-{\pr_{1}}& {A_{n+1}} \ar[d]^-{(\del_{i})_{i\in [n+1]}} \\
{\cycle_{n}A^{-}} \ar[r]_-{\cycle_{n}\del} & {\cycle_{n}A}}
\]
is a pullback.
\end{lemma}
\begin{proof}
Via the Yoneda Lemma, it suffices to check this in $\Set$.
\end{proof}

\begin{proposition}\label{Proposition-Characterization-Acyclic-Object}
A simplicial object $A$ in $\Ac$ is acyclic if and only if for every $n\in \N$, the morphism $(\del_{i})_{i\in [n+1]}\colon {A_{n+1}\to\cycle_{n}A}$ is a regular epimorphism.
\end{proposition}
\begin{proof}
We show by induction that, when $H_{i}A=0$ for all $i\in [n-1]$, the morphism 
\[
(\del_{i})_{i\in [n+1]}\colon {A_{n+1}\to\cycle_{n}A}
\]
is regular epic if and only if $H_{n}A=0$.

We already considered the case $n=0$ above. Hence suppose that $n\in \N$ and $H_{i}A=0$ for all $i\in [n]$, and consider the following diagram.
\[
\xymatrix{0 \ar[r] & {\El_{n+1}A} \ar@{{ >}->}[rr]^-{\ker \del_{0}} \ar[d]_-{(\del_{i}^{\El})_{i\in [n+1]}} && A_{n+1}^{-} \ar[d]_-{(\del_{i+1})_{i\in [n+1]}} \ar@{-{ >>}}[rr]^-{\del_{0}} && A_{n+1} \ar[d]^-{(\del_{i})_{i\in [n+1]}} \ar[r] & 0\\
0 \ar[r] & {\cycle_{n}\El A} \ar@{{ >}->}[rr]_-{\ker \cycle_{n}\del} && \cycle_{n} A^{-} \ar[rr]_-{\cycle_{n}\del} && \cycle_{n}A \ar@{.>}[r] & 0}
\]
By the induction hypothesis on $A$, the arrow $(\del_{i})_{i}\colon {A_{n+1}\to \cycle_{n}A}$ is a regular epimorphism. It follows that so is $\cycle_{n}\del$, and both rows in this diagram are exact sequences.

Using Lemma~\ref{Lemma-Lambda} and applying the induction hypothesis to $\El A$, we get that the arrow $(\del_{i}^{\El})_{i}\colon {\El_{n+1}A\to \cycle_{n}\El A}$ is a regular epimorphism if and only if $0=H_{n}\El A\cong H_{n+1}A$. Now, recalling Remark~\ref{Remark-Regepi-Pushout}, we see that $H_{n+1}A=0$ exactly when the right hand side square above is a $\regepi$-pushout. By Lemma~\ref{Lemma-Cycle-and-Pullbacks}, this latter property is equivalent to the induced arrow 
\[
(\del_{i})_{i\in [n+2]}\colon {A_{n+2}=A_{n+1}^{-}\to \cycle_{n+1}A}
\]
being a regular epimorphism.
\end{proof}

The meaning of this property becomes very clear when we express it as in~\fref{Figure-Acyclic-Simplicial-Object}: Up to enlargement of domain, $A$ is acyclic when every $n$-cycle is a boundary of an $(n+1)$-simplex. 

\subsection*{Acyclic fibrations}\label{Subsection-Acyclic-Fibrations}

A fibration is called \emph{acyclic} when it is a homology isomorphism; here we give an explicit description. Note that, using the long exact homology sequence (Proposition~\ref{Proposition-Long-Exact-Homology-Sequence-Simp}), regular epic homology isomorphisms may be characterized as those regular epimorphisms that have an acyclic kernel. This is why Proposition~\ref{Proposition-Characterization-Acyclic-Object} is useful in the proof of 

\begin{proposition}\label{Proposition-Characterization-Reg-Epi-Homology-Iso}
A regular epi $p\colon {E\to B}$ is a homology isomorphism if and only if for every $n\in \N$, the diagram
\begin{equation}\label{Diagram-Acyclic-Fibration}
\vcenter{\xymatrix{E_{n+1} \ar[d]_-{(\del_{i})_{i\in [n+1]}} \ar@{-{ >>}}[rr]^-{p_{n+1}} && B_{n+1} \ar[d]^-{(\del_{i})_{i\in [n+1]}} \\
\cycle_{n}E \ar[rr]_-{\cycle_{n}p} && \cycle_{n}B}}
\end{equation}
is a $\regepi$-pushout.
\end{proposition}
\begin{proof}
We give a proof by induction on $n$. If $n=0$ then $\cycle_{0}p$, as a product of regular epimorphisms, is regular epic. By Remark~\ref{Remark-Regepi-Pushout}, it is now clear that Diagram~\ref{Diagram-Acyclic-Fibration} is a $\regepi$-pushout if and only if $(\del_{0},\del_{1})\colon {K[p]_{1}\to \cycle_{0}K[p]}$ is a regular epimorphism, or, by Proposition~\ref{Proposition-Characterization-Acyclic-Object}, $H_{0}K[p]=0$. 

Now suppose that $n>0$. We first have to show that the arrow $\cycle_{n}p$ is a regular epimorphism. To do so, consider an arrow $x\colon {X\to \cycle_{n}B}$. Then a regular epimorphism $q\colon {Y\to X}$ exists and arrows $y_{k}\colon {Y\to E_{n-1}}$, such that for each element $k$ of the set
\[
K=\{\del_{j}\comp \pr_{i+1}\comp x\colon X\to B_{n-1}\,|\,j\in [n],\,i\in [n+1] \},
\]
$k\comp q= p_{n-1}\comp y_{k}$. These $y_{k}$ assemble to arrows $y'_{i}\colon {Y\to \cycle_{n-1}E}$ that satisfy $\cycle_{n-1}p\comp y'_{i}= (\del_{j})_{j}\comp \pr_{i+1}\comp x\comp q$. By the induction hypothesis, this yields a regular epimorphism $q_{1}\colon {Y_{1}\to Y}$ and arrows $z_{i}\colon {Y_{1}\to E_{n}}$ satisfying $p_{n}\comp z_{i}=\pr_{i+1}\comp x\comp q\comp q_{1}$, for every $i\in [n+1]$. By the careful choice of the arrows $y'_{i}$ (using the set $K$), these $z_{i}$ induce the needed morphism $z\colon {Y_{1}\to \cycle_{n}E}$.

Considering the diagram
\[
\xymatrix{0 \ar[r] & {K[p]_{n+1}} \ar@{{ >}->}[rr] \ar[d]_-{(\del_{i})_{i\in [n+1]}} && E_{n+1} \ar[d]_-{(\del_{i})_{i\in [n+1]}} \ar@{-{ >>}}[rr]^-{p_{n+1}} && B_{n+1} \ar[d]^-{(\del_{i})_{i\in [n+1]}} \ar[r] & 0\\
0 \ar[r] & {\cycle_{n}K[p]} \ar@{{ >}->}[rr] && \cycle_{n}E \ar@{-{ >>}}[rr]_-{\cycle_{n}p} && {\cycle_{n}B} \ar[r] & 0,}
\]
we see that~\ref{Diagram-Acyclic-Fibration} is a $\regepi$-pushout if and only if the morphism 
\[
(\del_{i})_{i}\colon {K[p]_{n+1}\to \cycle_{n}K[p]}
\]
is regular epic. The result now follows from Proposition~\ref{Proposition-Characterization-Acyclic-Object}.
\end{proof}

\begin{proposition}\label{Proposition-Characterization-Acyclic-Fibration}\index{fibration!acyclic}
For a simplicial morphism $p\colon E\to B$, the following are equivalent:
\begin{enumerate}
\item $p$ is an acyclic fibration;
\item $p$ is a regular epimorphism and for every $n\in \N$, Diagram~\ref{Diagram-Acyclic-Fibration} is a $\regepi$-pushout.
\end{enumerate} 
\end{proposition}
\begin{proof}
This is an immediate consequence of Proposition~\ref{Proposition-Characterization-Reg-Epi-Homology-Iso}: By Proposition~\ref{Proposition-Kan}, a regular epimorphism is always a Kan fibration; conversely, it easily follows from the definitions that a Kan fibration $p$ such that $H_{0}p$ is a regular epimorphism is also regular epic (cf.~\cite[\S II.3, Proposition 1 ff.]{Quillen}).
\end{proof}

\subsection*{Weak equivalences are homology isomorphisms}\label{Subsection-Weak-Equivalences-Homology-Isos}
We show that in the case of a semi-abelian $\Ac$, the weak equivalences in $\simpA$ are exactly the homology isomorphisms.

Recall from~\cite[Section II.4]{Quillen} that a weak equivalence in the sense of Theorem~\ref{Theorem-Quillen} is a simplicial morphism $f\colon A\to B$ such that $\hom (P,f)$ is a weak equivalence in $\Sc \Set$ for every projective object $P$ in $\Ac$; and this is the case when $\hom (P,f)$ induces isomorphisms of the \emph{homotopy groups} 
\[
\pi_{n} \hom (P,f)\colon \pi_{n} (\hom (P,A),x)\to \pi_{n} (\hom (P,B), f\comp x)
\]
for all $n\in \N$ and $x\in \hom (P,A)_{0}=\hom (P,A_{0})$. These homotopy groups are defined in the following manner (see e.g.,~\cite[Definition 3.3 ff.]{May} or~\cite[Section 8.3]{Weibel}).

\begin{definition}\label{Definition-Homotopy-Groups-Kan}\index{homotopy!group}
Let $K$ be a Kan simplicial set, $x\in K_{0}$ and $n\in \N$. If $n\geq 1$, denote by $x$ the element $(\sigma_{0}\comp \cdots \comp \sigma_{0}) (x)$ of $K_{n-1}$ and by $Z_{n} (K,x)$ the set $\{z\in K_{n}\,|\, \text{$\del_{i} (z)= x$ for $i\in [n]$}\}$; $Z_{0} (K,x)=K_{0}$. Two elements $z$ and $z'$ of $Z_{n} (K,x)$ are \emph{homotopic}\index{homotopy!of cycles}, notation $z\sim z'$, if there is a $y\in K_{n+1}$ such that
\[
\del_{i} (y) = \begin{cases}
x,  & \text{if $i < n$;}\\
z,& \text{if $i = n$;}\\
z',& \text{if $i = n+1$}.
\end{cases}
\]
The relation $\sim$ is an equivalence relation on $Z_{n} (K,x)$\index{Z_n@$Z_{n} (K,x)$} and the quotient $Z_{n} (K,x) /_{\sim}$ is denoted by $\pi_{n} (K,x)$. For $n\geq 1$ this $\pi_{n} (K,x)$\index{pi@$\pi_{n} (K,x)$} is a group.
\end{definition}

The next lemma is essentially a higher-degree version of Corollary~\ref{Corollary-H_0}, and shows how a homology object $H_{n}A$ may be computed as a coequalizer of (a restriction of) the last two boundary operators $\del_{n}$ and $\del_{n+1}$.

\begin{lemma}\label{Lemma-Kernel-Pair-of-q_{n}}
Write $S_{n+1}A=\bigcap_{i\in [n-1]}K[\del_{i}\colon A_{n+1}\to A_{n}]\subset A_{n+1}$ and $S_{1}A=A_{1}$, and for $n\in \N$, let $\del_{n}^{-1}Z_{n}A$ denote the inverse image
\[
\xymatrix{{\del_{n}^{-1}Z_{n}A} \ar@{}[rd]|<<{\pullback} \ar@<-1ex>[r]_-{\overline{\del}_{n}} \ar@<1ex>[r]^-{\overline{\del}_{n+1}} \ar[d]_-{\del_{n}^{-1}\bigcap_{i\in [n]} \ker \del_{i}} & Z_{n}A \ar[l]|-{\overline{\sigma}_{n}} \ar[d]^-{\bigcap_{i\in [n]} \ker \del_{i}}\\
S_{n+1}A\ar[r]_-{\del_{n}} & A_{n}}
\]
 of $Z_{n}A=\bigcap_{i\in [n]}K[\del_{i}]\subset A_{n}$ along $\del_{n}\colon {S_{n+1}A\to A_{n}}$ and $\overline{\del}_{n+1}$, $\overline{\sigma}_{n}$ the induced factorization of $\del_{n+1}\colon {A_{n+1}\to A_{n}}$ and $\sigma_{n}\colon {A_{n}\to A_{n+1}}$. Then
\[
\xymatrix{{\del_{n}^{-1}Z_{n}A} \ar@<-.5 ex>[r]_-{\overline{\del}_{n}} \ar@<.5 ex>[r]^-{\overline{\del}_{n+1}} & Z_{n} A \ar@{-{ >>}}[r]^-{q_{n}} & H_{n}A} 
\]
is a coequalizer diagram.
\end{lemma}
\begin{proof}
The composite 
\[
\xymatrix{K[\overline{\del}_{n}] \ar@{{ >}->}[r]^-{\ker \overline{\del}_{n}} & \del_{n}^{-1}Z_{n}A \ar[r]^-{\overline{\del}_{n+1}} & Z_{n}A }
\]
is $d'_{n+1}:N_{n+1}A\to Z_{n}A$; the result now follows from Proposition~\ref{Proposition-Forks}.
\end{proof}

\begin{proposition}\label{Proposition-Homology-and-Homotopy-Groups}
Let $\Ac$ be a semi-abelian category. For any $n\geq 1$ and for any projective object $P$ in $\Ac$, there is a bijection $\pi_{n} (\hom (P,A),0)\cong \hom (P, H_{n}A)$, natural in $A\in \Sc \Ac$.
\end{proposition}
\begin{proof}
An isomorphism 
\[
\pi_{n} (\hom (P,A),0)\to  \hom (P, H_{n}A)
\]
may be defined as follows: The class with respect to $\sim$ of a morphism $f\colon {P\to Z_{n}A}$ is sent to $q_{n}\comp f\colon P\to H_{n}A$. This mapping is well-defined and injective, because by Lemma~\ref{Lemma-Kernel-Pair-of-q_{n}}, $f\sim g$ is equivalent to $q_{n}\comp f=q_{n}\comp g$; it is surjective, since $q_{n}$ is a regular epimorphism and $P$ a projective object.  
\end{proof}

\begin{theorem}\label{Theorem-Weak-Equivalences-are-Homology-Isos}
If $\Ac$ is semi-abelian with enough projectives, weak equivalences in $\Sc\Ac$ and homology isomorphisms coincide. 
\end{theorem}
\begin{proof}
If $f$ is a weak equivalence then all $\pi_{n} \hom (P,f)$ are isomorphisms; but then the morphisms $\hom (P,H_{n}f)$, $n\geq 1$ are iso by Proposition~\ref{Proposition-Homology-and-Homotopy-Groups}. Now the class of functors 
\[
\hom (P,\cdot)\colon {\Ac \to \Set},
\]
$P$ projective, is jointly conservative: It reflects regular epimorphisms by definition of a regular projective object, and it reflects monomorphisms because enough regular projectives exist. Hence $H_{n}f$ is an isomorphism.

Conversely, suppose that $f\colon {A\to B}$ is a homology isomorphism. We may factor it as $f=p\comp i$, a trivial cofibration $i\colon {A\to E}$ followed by a fibration $p\colon {E\to B}$. We just proved that $i$ is a homology isomorphism; as a consequence, so is $p$. Thus we reduced the problem to showing that when a fibration is a homology isomorphism, it is a weak equivalence.

Consider $n\in \N$, $P$ projective in $\Ac$, and $x\colon {P\to E_{0}}$. Then 
\[
p_{n}\comp (\cdot)\colon {\pi_{n} (\hom (P, E),x)\to \pi_{n} (\hom (P,B),p_{0}\comp x)}
\]
is an isomorphism: It is a surjection, because by Proposition~\ref{Proposition-Characterization-Acyclic-Fibration}, any map $b\colon {P\to B_{n}}$ satisfying $\del_{i}\comp b=p_{n-1}\comp x$ induces a morphism $e\colon {P\to E_{n}}$ such that $\del_{i}\comp e= x$. It is an injection, since for $e$, $e'$ in $Z_{n} (\hom (P,E),x)$ with $p_{n}\comp e\sim p_{n}\comp e'$, again Proposition~\ref{Proposition-Characterization-Acyclic-Fibration} implies that $e\sim e'$. 
\end{proof}

We get that the model structure on $\Sc \Ac$ is compatible with the one from Section~\ref{Section-Regular-Epimorphisms}:

\begin{corollary}\label{Corollary-Quillen-Adjunction}\index{Quillen functor}
The nerve functor $\ner\colon \Cat \Ac \to \Sc \Ac$\index{nerve!of an internal category}\index{ner@$\ner \Af$} is a \emph{right Quil\-len functor}\index{right Quillen functor}: It is a right adjoint that preserves fibrations and trivial fibrations.
\end{corollary}
\begin{proof}
The functor $\ner$ has a left adjoint $\pi_{1}$\index{fundamental groupoid}\index{pi_1@$\pi_{1}\colon {\Sc \Ac\to \Cat \Ac}$} by Proposition~\ref{Proposition-Fundamental-Groupoid}. The preservation of fibrations and weak equivalences follows from Proposition~\ref{Proposition-Characterization-Fibration-Kan} and Definition~\ref{Definition-Homology-of-Internal-Category}.
\end{proof}

To conclude, we may summarize this section as follows: 

\begin{theorem}\label{Theorem-Simp-Model-Structure}
If $\Ac$ is a semi-abelian category with enough projectives then a model category structure is defined on $\simpA$ by choosing $\we$ the class of homology isomorphisms, $\fib$ the class of Kan fibrations and $\cof={^{\Box}\fib\cap \we}$.\noproof 
\end{theorem}

\section{Simplicial homotopy}\label{Section-Simp-Cocylinder}
In this section we explore several ways of dealing with simplicial homotopy and with the homotopy relation induced by the model structure on $\simpA$.

\subsection*{Simplicial enrichment}\label{Subsection-Simplicial-Enrichment}
It is well-known that for every category $\Ac$, $\simpA$ is an \emph{$\Sc$-category}\index{S-category@$\Sc$-category}\index{category!$\Sc$- ---}\index{category!enriched in $\Sc$} or \emph{simplicial category}\index{category!simplicial}\index{simplicial category}, i.e., a category enriched in the monoidal category $\Sc = (\Sc \Set,\times ,1)$\index{category!\Sc}. We recall the construction from~\cite{KP:Abstact-homotopy-theory} (see also~\cite{Quillen}). 

\begin{notation}\label{Notation-Enriched-Deltas}
Let $\Delta [n]$ denote the simplicial set 
\[
\hom_{\Delta} (\cdot, [n])\colon \Delta^{\op}\to \Set.
\]
As the primordial cosimplicial object (with faces $\del^{i}\colon {[n]\to [n+1]}$ and degeneracies $\sigma^{i}\colon {[n]\to [n-1]}$), $\Delta$ induces a structure of cosimplicial object in $\Sc$ on the sequence of simplicial sets $(\Delta [n])_{n\in \N}$: its faces are 
\[
\del^{i}=\hom_{\Delta} (\cdot , \del^{i})\colon {\Delta [n]\to \Delta [n+1]}
\]
and its degeneracies $\sigma ^{i}=\hom_{\Delta} (\cdot , \sigma ^{i})\colon {\Delta [n]\to \Delta [n-1]}$. 
\end{notation}

\begin{notation}\label{Notation-Delta-Comma-Category}
Following~\cite{Borceux:Cats}, for any simplicial set $X$, let $\Elts (X)$\index{category!$\Elts (X)$} denote the category of elements of $X$, the comma category $(1_{\Delta^{\op}},X)$: its objects are pairs $([m],x)$ with $x\in X ([m])=X_{m}$, and a morphism $([m],x)\to ([n],y)$ is a map $\mu \colon [m]\to [n]$ in $\Delta^{\op}$ such that $X (\mu) (x)=y$. Let 
\[
\delta_{X}\colon{\Elts (X)\to \Delta^{\op}}
\]
denote the forgetful functor: $\delta_{X} ([m],x)=[m]$.
\end{notation}

This construction is natural in $X$: Any simplicial map $f\colon {X\to Y}$ induces a functor $\Elts f\colon {\Elts (X)\to \Elts (Y)}$ that sends $([m],x)$ to $([m],f_{m} (x))$ and satisfies $\delta_{Y}\comp \Elts f=\delta_{X}$.

\begin{definition}\label{Definition-Enriched-in-Simp}
For simplicial objects $A$ and $B$ in $ \Ac$, the simplicial set $\Sc \Ac (A,B)$ defined by 
\[
\Sc \Ac (A,B)_{n}=\hom (A\comp \delta_{\Delta [n]},B\comp \delta_{\Delta[n]})
\]
determines a structure of $\Sc$-category on $\Sc \Ac$. Faces and degeneracies for $\Sc \Ac (A,B)$ are obtained by precomposition: e.g., given $\varphi\colon {A\comp \delta_{\Delta [n]}\To B\comp \delta_{\Delta[n]}}$,  
\[
\del_{i} (\varphi)=\varphi \comp 1_{\Elts \del^{i}}\colon A\comp \delta_{\Delta [n]}\comp \Elts \del^{i} \To B\comp \delta_{\Delta [n]}\comp \Elts \del^{i}.
\]
\end{definition}

\begin{figure}
\begin{center}
\resizebox{\textwidth}{!}{\includegraphics{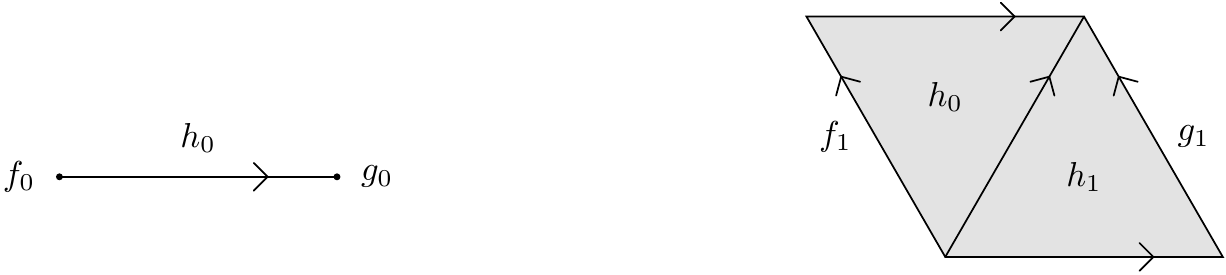}}
\end{center}
\caption{The two lowest degrees of a simplicial homotopy $h\colon f\simeq g$.}\label{Figure-Simplicial-Homotopy}
\end{figure}

An element of $\Sc \Ac (A,B)_{0}$ is just a simplicial morphism from $A$ to $B$ ($\delta_{\Delta [0]}$ is an isomorphism). Recall that two simplicial morphisms $f,g\colon A\to B$ in $\Ac$ are called \emph{(simplicially) homotopic}\index{homotopy!of simplicial morphisms} if, for every $n\in \N$, there are morphisms $h_{i}\colon {A_{n}\to B_{n+1}}$ in $\Ac$, $i\in [n]$ such that $\del _{0}\comp h_{0}=f_{n}$ and $\del_{n+1}\comp h_{n}=g_{n}$, while the identities
\[
\sigma_{i}\comp h_{j}=\begin{cases}
h_{j+1}\comp \sigma_{i} & \text{if $i\leq j$}\\
h_{j}\comp \sigma_{i-1} & \text{if $i>j$}
\end{cases}
\qquad\qquad 
\del_{i}\comp h_{j}=\begin{cases}h_{j-1}\comp \del_{i} & \text{if $i<j$} \\
\del_{i}\comp h_{i-1} & \text{if $i=j\neq 0$}\\
h_{j}\comp \del_{i-1} & \text{if $i>j+1$}\end{cases}
\]
hold. One may show that such a \emph{homotopy} from $f$ to $g$ corresponds to an element $h$ of $\Sc \Ac (A,B)_{1}$ that satisfies $\del_{0}\comp h=f$ and $\del_{1}\comp h=g$. In its lowest degrees, a homotopy may be pictured as in~\fref{Figure-Simplicial-Homotopy}.

\subsection*{A cocylinder on $\simpA$}\label{Subsection-Cocylinder-Explicit-Description}
We recall the construction of a cocylinder for simplicial homotopy, describe it in an explicit way and study some of its properties. 

Recall (e.g., from~\cite{KP:Abstact-homotopy-theory}) that a \emph{cotensor}\index{cotensor on $\Sc \Ac$} of a simplicial set $X$ by a simplicial object $B$ in $\Ac$ is an object $B^{X}$ together with a morphism 
\[
\beta \colon {X\to \Sc \Ac (B^{X},B)}
\]
such that, for every $A$ in $\Sc \Ac$, the exponential transpose of the map
\[
\Sc \Ac (A,B^{X})\times X \to^{1\times \beta}\Sc \Ac (A,B^{X})\times \Sc \Ac (B^{X},B)\to^{\comp} \Sc \Ac(A,B)
\]
is an isomorphism 
\[
\Sc (X, \Sc \Ac (A,B))\cong \Sc \Ac (A, B^{X}).
\]
Let $\FinSet$\index{category!$\FinSet$} denote the category of finite sets, and call a simplicial set \emph{finite}\index{finite simplicial set}\index{simplicial set!finite} when it is an object of $\Sc \FinSet$. When $\Ac$ is finitely complete, $\Sc \Ac$ has a cotensor for any finite simplicial set $X$ by any simplicial object $B$ in $\Ac$; the object $B^{X}$ is given by the limit
\[
B^{X}_{n}=\lim (B\comp \delta_{X\times \Delta [n]}),
\]
for $n\in \N$. This induces a functor $B^{(\cdot)}\colon \Sc\FinSet\to \Sc \Ac$\index{B@$B^{(\cdot)}$} from the category of finite simplicial sets to the category of simplicial objects in $\Ac$. Using that $B_{n}$ is initial in the diagram $B\comp \delta_{n}$, note that $B= B^{\Delta[0]}$.

\begin{definition}\label{Definition-Simp-Cocylinder}
For any simplicial object $B$ in $\Sc \Ac$, put
\[
\epsilon_{0} (B) =B^{\del^{0}}, \epsilon_{1} (B) = B^{\del^{1}}\colon B^{\Delta[1]}\to B
\]
and $s (B)=B^{\sigma^{0}}\colon  B\to B^{\Delta[1]}$. Thus we get a cocylinder 
\[
\Ib=((\cdot)^{\Delta[1]}\colon \Sc\Ac \to \Sc\Ac,\quad \epsilon_{0},\epsilon_{1}\colon (\cdot)^{\Delta[1]}\To 1_{\Sc\Ac},\quad s\colon 1_{\Sc\Ac}\To (\cdot)^{\Delta[1]})
\]
on $\Sc\Ac$.
\end{definition}

The purpose of this construction is, of course, that the cocylinder $\Ib$ characterizes simplicial homotopy. As in Section~\ref{Section-Cocylinder}, two simplicial morphisms $f,g\colon {A\to B}$ in $\Ac$ are \emph{homotopic}\index{cocylinder}\index{homotopy!w.r.t.\ a cocylinder} with respect to the cocylinder $\Ib$ if and only if there exists a simplicial morphism $H\colon A\to B^{\Delta[1]}$ (called a \emph{homotopy} from $f$ to $g$) such that $\epsilon_{0} (B)\comp H=f$ and $\epsilon_{1} (B)\comp H=g$. This situation is denoted $f\simeq_{\Ib}g$.

\begin{proposition}\label{Proposition-Simplicial-Homotopy}
Two simplicial morphisms $f,g\colon A\to B$ in $\Ac$ are simplicially homotopic if and only if $f\simeq_{\Ib}g$.
\end{proposition}
\begin{proof}
A simplicial morphism $H\colon A\to B^{\Delta[1]}$ such that $\epsilon_{0} (B)\comp H=f$ and $\epsilon_{1} (B)\comp H=g$ corresponds to a morphism $\widehat{H}\colon {\Delta [1]\to \Sc \Ac (A,B)}$ of simplicial sets satisfying $\widehat{H}\comp \del^{0}=\widehat{f}$ and $\widehat{H}\comp \del^{1}=\widehat{g}$---which, via the Yoneda Lemma, corresponds to a simplicial homotopy $h$ from $f$ to~$g$.
\end{proof}

We shall now describe this cocylinder in a more explicit way.

\subsection*{The cocylinder described explicitly}\label{Subsection-Description-Cocylinder}
A family of split epimorphisms 
\[
(\del_{i}\colon A_{n+1}\to A_{n})_{i\in [n+1]},
\]
part of a simplicial object, has some interesting properties. We study them by making abstraction of the situation as follows. 

\begin{definition}\label{Definition-Split-Family}
A family of morphisms $(\del_{i}\colon A_{n+1}\to A_{n})_{i\in [n+1]}$ is called \emph{split}\index{split family of morphisms} if a family $(\sigma_{i}\colon A_{n}\to A_{n+1})_{i\in [n]}$ exists such that, for every $i\in [n]$, $\del_{i}\comp\sigma_{i}=\del_{i+1}\comp\sigma_{i}=1_{A_{n}}$.
\end{definition}

Of course, a split family consists of split epimorphisms, and any family that comes from a simplicial object (like the family considered above) is split. We need some notation and a technical result.

\begin{notation}\label{Notation-A^J}
Suppose that $\Ac$ has finite limits and consider a split family 
\[
(\del_{i}\colon A_{n+1}\to A_{n})_{i\in [n+1]}
\]
in $\Ac$. If $n=0$, put $A^{I}_{n}=A_{1}$ and, if $n> 0$, let $A^{I}_{n}$ be the limit (with projections $\pr_{1}$,\dots, $\pr_{n+1}\colon {A^{I}_{n}\to A_{n+1}}$) of the zigzag\index{zigzag}
\[
\xymatrix{A_{n+1}\ar[rd]_{\del_{1}} && A_{n+1} \ar[ld]^{\del_{1}} \ar[rd]_{\del_{2}}&& {} \ar@{}[d]|-{\displaystyle\cdots} \ar@{.>}[ld] \ar@{.>}[rd] && A_{n+1} \ar[ld]^{\del_{n}}\\
& A_{n} && A_{n} & {} & A_{n}}
\]
in $\Ac$.

Let $\epsilon_{0} (A)_{n},\epsilon_{1} (A)_{n}\colon A^{I}_{n}\to A_{n}$ and $s (A)_{n}\colon A_{n}\to A^{I}_{n}$ denote the morphisms respectively defined by
\[
\begin{aligned}
\epsilon_{0} (A)_{0} &= \del_{0}\\
\epsilon_{0} (A)_{n} &= \del_{0}\comp \pr_{1},
\end{aligned}
\quad 
\begin{aligned}
\epsilon_{1} (A)_{0} &= \del_{1}\\
\epsilon_{1} (A)_{n} &= \del_{n+1}\comp \pr_{n+1},
\end{aligned}
\quad\text{and}\quad 
s (A)_{n}= (\sigma_{0},\dots,\sigma_{n}).
\]
\end{notation}

\begin{proposition}\label{Proposition-eta-Is-Wide-Coequalizer}
Consider a split family $(\del_{i}\colon A_{n+1}\to A_{n})_{i\in [n+1]}$ in a category with finite limits and coequalizers $\Ac$. Denote a coequalizer of $\epsilon_{0} (A)_{n}$, $\epsilon_{1} (A)_{n} \colon {A^{I}_{n}\to A_{n}}$ by 
\[
\eta(A)_{n} =\coeq [\epsilon_{0} (A)_{n},\epsilon_{1} (A)_{n}]\colon A_{n}\to \tilde{A}_{n}.
\]
Then $\eta (A)_{n}$ is a wide coequalizer of $(\del_{i}\colon A_{n+1}\to A_{n})_{i\in [n+1]}$.
\end{proposition}
\begin{proof}
Our first claim is that for any $i,j\in [n+1]$, $\eta (A)_{n}\comp \del_{i}=\eta (A)_{n}\comp \del_{j}$. Considering $(1_{A_{n}},\sigma_{1}\comp \del_{1},\dots, \sigma_{n}\comp \del_{1})\colon A_{n+1}\to A_{n}^{I}$, it is readily seen that $\eta (A)_{n}\comp \del_{0}=\eta (A)_{n}\comp \del_{1}$. Similarly, using 
\[
(\sigma_{0}\comp \del_{i},\dots, \sigma_{i-1}\comp \del_{i},1_{A_{n+1}},\sigma_{i+1}\comp \del_{i+1},\dots ,\sigma_{n}\comp\del_{i+1})\colon A_{n+1}\to A_{n}^{I}
\]
one shows that $\eta (A)_{n}\comp \del_{i}=\eta (A)_{n}\comp \del_{i+1}$ for any $0<i\leq n$. This proves the claim by complete induction.

Next we show that if $c\colon A_{n}\to C$ coequalizes $(\del_{i})_{i\in [n+1]}$ then $c\comp \epsilon_{0} (A)_{n}=c\comp \epsilon_{1} (A)_{n}$. Indeed, 
\[
c\comp \del_{0}\comp \pr_{1}=c\comp \del_{1}\comp \pr_{1}=c\comp \del_{1}\comp \pr_{2}=\dots=c\comp \del_{n}\comp \pr_{n+1}=c\comp \del_{n+1}\comp \pr_{n+1},
\]
which finishes the proof.
\end{proof}

\begin{proposition}
Let $B$ be a simplicial object in a finitely complete category $\Ac$. Then the simplicial object $B^{I}$ defined by the face operators $\del_{i}^{I}\colon {B^{I}_{n}\to B^{I}_{n-1}}$ and degeneracy operators $\sigma_{i}^{I}\colon {B^{I}_{n}\to B^{I}_{n+1}}$, for $i\in [n]$, $1\leq j\leq n$ and $1\leq k\leq n+2$ given by
\begin{gather*}
\begin{aligned}
\del_{0}^{I} &=\del_{0}\comp \pr_{2}\colon B^{I}_{1}\to B^{I}_{0}\\
\del_{1}^{I} &=\del_{2}\comp \pr_{1}\colon B^{I}_{1}\to B^{I}_{0}\\
\sigma_{0}^{I} &= (\sigma_{1},\sigma_{0})\colon B^{I}_{0}\to B^{1}_{1}
\end{aligned}
\qquad 
\pr_{j}\comp \del^{I}_{i} = \begin{cases}\del_{i+1}\comp \pr_{j} & \text{if $j\leq i$}\\
\del_{i}\comp \pr_{j+1} & \text{if $j>i$}  \end{cases}\colon B^{I}_{n}\to B_{n}\\
\pr_{k}\comp \sigma^{I}_{i} =
\begin{cases} \sigma_{i+1}\comp \pr_{k} & \text{if $k\leq i+1$}\\
\sigma_{i}\comp \pr_{k-1} & \text{if $k>i+1$} \end{cases}
\colon B^{I}_{n}\to B_{n+2},
\end{gather*}
coincides with $B^{\Delta [1]}$ from Definition~\ref{Definition-Simp-Cocylinder}. The morphisms mentioned in Notation~\ref{Notation-A^J} above form the simplicial morphisms $\epsilon_{0} (B),\epsilon_{1} (B)\colon B^{I}\to B$ and $s (B)\colon {B\to B^{I}}$ from~\ref{Definition-Simp-Cocylinder}.
\end{proposition}
\begin{proof}
We show that $B^{\Delta [1]}$ and $B^{I}$ are isomorphic; then choosing the appropriate limits allows us to get rid of all distinction between the two objects. Indeed, the isomorphisms $\hom (A,B^{\Delta [1]})\cong \hom (\Delta [1],\Sc \Ac (A,B))$ and 
\[
\hom (\Delta [1],\Sc \Ac (A,B))\cong \hom (A,B^{I})
\]
are natural in $A$; hence a natural isomorphism $\hom (\cdot, B^{\Delta [1]})\To \hom (\cdot, B^{I})$ exists, and this induces an isomorphism ${B^{\Delta [1]}\to B^{I}}$.
\end{proof}

\subsection*{Homotopy \versus\ homology}\label{Subsection-Homotopy-vs-Homology}
Using a Kan property argument, we now give a direct proof that homotopic morphisms of simplicial objects have the same homology.
\begin{proposition}\label{Proposition-phi_0-Homology-Isomorphism}
Let $A$ be a simplicial object in a semi-abelian category $\Ac$; consider 
\[
\epsilon_{0} (A)\colon {A^{I}\to A}.
\]
For any $n\in \N$, $H_{n}\epsilon_{0} (A)$ is an isomorphism.
\end{proposition}
\begin{morale}
We sketch the proof in the case $n=0$. We have to show that the diagram below is a $\regepi$-pushout. Consider $h=h_{0}$ and $z=z_{2}$ as in~\fref{Figure-Kan-Twice}; then (up to enlargement of domain) using the Kan property twice yields the needed $(y_{0},y_{1})$ in $N_{1}A^{I}$.
\begin{figure}
\begin{center}
\includegraphics{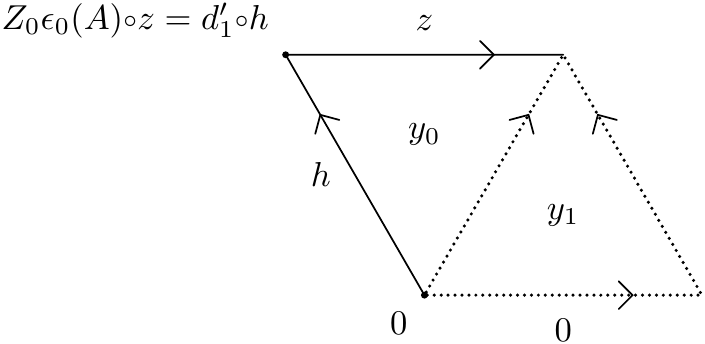}
\end{center}
\caption{Using the Kan property twice to get $(y_{0},y_{1})$ in $N_{1}A^{I}$.}\label{Figure-Kan-Twice}
\end{figure}
\end{morale}
\begin{proof}
Using the Kan property, we show that the commutative diagram
\[
\xymatrix{N_{n+1}A^{I} \ar@{-{ >>}}[rr]^-{N_{n+1}\epsilon_{0} (A)} \ar[d]_-{d'_{n+1}} && N_{n+1}A \ar[d]^-{d'_{n+1}}\\
Z_{n} A^{I} \ar@{-{ >>}}[rr]_-{Z_{n}\epsilon_{0} (A)} && {Z_{n}A}}
\]
is a $\regepi$-pushout; then it is a pushout, and the result follows by Proposition~\ref{Proposition-pushout}.

Consider morphisms $z\colon Y_{0}\to Z_{n}A^{I}$ and $h\colon Y_{0}\to N_{n+1}A$ that satisfy $d'_{n+1}\comp h=Z_{n}\epsilon _{0} (A)\comp z$. We are to show that a regular epimorphism $p\colon {Y\to Y_{0}}$ exists and a morphism $y\colon Y\to N_{n+1}A^{I}$ satisfying $d'_{n+1}\comp y=z\comp p$ and $N_{n+1}\epsilon_{0} (A)\comp y=h\comp p$. Put 
\[
h_{0}= \bigcap_{j} \ker \del_{j}\comp h\colon Y_{0}\to A_{n+1},
\]
$z_{n+2}=z$, $z_{i}=0\colon Y_{0}\to Z_{n}A^{I}$ for $i\in [n+1]$ and $(z_{i0},\dots,z_{in})= \bigcap_{j} \ker \del_{j}\comp z_{i}$, for $i\in [n+2]$.

We prove that if for $k\in [n+1]$, there is a morphism $h_{k}\colon Y_{k}\to A_{n+1}$ with $\del_{i}\comp h_{k}= \del_{k}\comp z_{ik}$, $i\in [n+1]$, then a regular epimorphism $p_{k+1}\colon {Y_{k+1}\to Y_{k}}$ exists and a morphism $y_{k}\colon {Y_{k+1}\to A_{n+2}}$ such that $\del_{k}\comp y_{k}= h_{k}\comp p_{k+1}$ and
\[
\del_{i}\comp y_{k}=
\begin{cases}
z_{ik}\comp p_{k+1}, &\text{if $i<k$,}\\
z_{i (k+1)}\comp p_{k+1}, &\text{if $i>k+1$;}
\end{cases}
\]
in particular, for $h_{k+1}= \del_{k+1}\comp y_{k}\colon Y_{k+1}\to A_{n+1}$, $\del_{i}\comp h_{k+1}= \del_{k+1}\comp z_{i (k+1)}$, $i\in [n+1]$. This then proves the assertion by induction on $k$: 
\[
y=(y_{0}\comp p_{1}\comp \cdots \comp p_{n+2},\dots, y_{n+1})\colon Y_{n+2}\to N_{n+1}A^{I}
\]
satisfies $d'_{n+1}\comp y=z\comp p_{1}\comp \cdots \comp p_{n+2}$ and $N_{n+1}\epsilon_{0} (A)\comp y=h\comp p_{1}\comp\cdots \comp p_{n+2}$.

Indeed, putting $x_{k}=h_{k}$, $x_{i+1}= z_{i}$ for $i>k$ and $x_{i}= z_{i}$ for $i<k$ defines an $(n+2,k+1)$-horn $(x_{i})_{i}$ in $A^{I}$ that induces $y_{k}$.
\end{proof}

\begin{corollary}\label{Corollary-Homotopic-then-Same-Homology}
If $f\simeq_{\Ib} g$ then, for any $n\in \N$, $H_{n}f=H_{n}g$.\noproof 
\end{corollary}

\begin{corollary}\label{Corollary-Cocylinder-Model-Category}
The cocylinder $\Ib$ is a cocylinder in the model-category-theoretic sense: For any simplicial object $A$, the map $s (A)$ is a weak equivalence.\noproof 
\end{corollary}

\subsection*{Simplicial homotopy as an internal relation}\label{Subsection-Homotopy-as-Internal-Relation}

Although simplicial homotopy may be characterized in terms of a cocylinder, it need not induce equivalence relations on hom-sets, not even when the ambient category $\Ac$ has good categorical properties. Recall that in $\Sc\Set$, the simplicial homotopy relation on a hom-set $\hom (A,B)$ is an equivalence relation as soon as the simplicial object $B$ is Kan. Hence one would imagine that in a regular Mal'tsev category, every hom-set carries a simplicial homotopy relation that is an equivalence relation. It is well-known that this is not the case. We briefly recall the example of groups: The category $\Gp$ forms a semi-abelian variety, hence is regular Mal'tsev, but nevertheless the relation induced by simplicial homotopy need not be, say, symmetric.

\begin{example}\label{Example-Homotopy-not-an-Equivalence-Relation}
Every reflexive graph $d_{0},d_{1}\colon B_{1}\to B_{0}$ in $\Gp$ is the bottom part of some simplicial object $B$, and the functor $\Sc \Gp \to \Gp$ that maps a simplicial object $B$ to the group $B_{0}$ has a left adjoint. Thus a diagram
\[
\xymatrix{& B_{1} \ar@<-1 ex>[d]_-{d_{0}} \ar@<1 ex>[d]^-{d_{1}}\\
A_{0} \ar@{.>}[ru]^-{h_{0}} \ar@<.5 ex>[r]^-{f_{0}} \ar@<-.5 ex>[r]_{g_{0}} & B_{0} \ar[u]|-{i} }
\]
induces simplicial morphisms $f,g\colon A\to B$, and $f\simeq_{\Ib}g$ if and only if there is a morphism $h_{0}\colon {A_{0}\to B_{1}}$ that satisfies $d_{0}\comp h_{0}=f_{0}$ and $d_{1}\comp h_{0}=g_{0}$. Now certainly, if $f_{0}=d_{0}$ and $g_{0}=d_{1}$, $f\simeq_{\Ib}g$; hence it suffices to show that in $\Gp$, there exist examples of reflexive graphs that are not symmetric. For instance, the non-symmetric reflexive graph in $\Set$ 
\[
\xymatrix@1{{a} \ar@(ul,dl)[]_{x} \ar[r]^-{y} & {b} \ar@(dr,ur)[]_{z}}
\]
can be modelled by 
\[
\xymatrix{{\Z_{2} + \Z + \Z_{3}}  \ar@<-1 ex>[rr]_-{d_{0}} \ar@<1 ex>[rr]^-{d_{1}} && {\Z_{2} + \Z_{3}} \ar[ll]|-{i}}
\]
where the generators $x$, $y$ and $z$ are respectively mapped to $a$, $a$ and $b$ by $d_{0}$ and $a$, $b$ and $b$ by $d_{1}$, and $i$ maps the generators $a$ and $b$ to $x$ and $z$, respectively.
\end{example}

The problem with this example is, of course, that when generalizing the Kan property from $\Set$ to a regular category $\Ac$ with enough projectives, it becomes relative to the class of projective objects. Note that in $\Set$, every object is projective, but not so in $\Gp$, and this is exactly why Example~\ref{Example-Homotopy-not-an-Equivalence-Relation} exists. Moreover, this definition captures the situation in $\Gp$ very nicely in the following manner. For, if $U\vdash F\colon \Gp \rightharpoonup \Set$\index{arrow!$\rightharpoonup$} denotes the underlying set/free group adjunction, every simplicial group $K$ has an underlying Kan simplicial set $U^{\Delta^{\op}} (K)$; by adjointness, $K$ is Kan, relative to all free groups $F (X)$; and free groups are exactly the projective ones. 

One could argue that in a regular Mal'tsev category, every simplicial object is Kan, \textit{internally}; that hence a sensible thing to do is to consider simplicial homotopy \textit{internally}. The aim of the next few arguments is to show that thus simplicial homotopy indeed becomes an equivalence relation, but in a way unsuitable in practice. 

The conceptual step of considering the homotopy relation internally instead of externally may be taken as soon as $\Ac$ is a regular category. Indeed, one of the main features of a regular category being the existence of image factorizations, for any simplicial object $A$ in $\Ac$, the left hand side reflexive graph 
\[
\xymatrix@=1.5cm{A^{I} \ar@{.>>}[rr]^-{r (A)} \ar@<-1.6ex>[dr]_-{\epsilon_{0} (A)} \ar@<1.6ex>[dr]^-{\epsilon_{1} (A)} && A^{K} \ar@{.>}@<-1.6ex>[dl]_-{\gamma_{0} (A)} \ar@{.>}@<1.6ex>[dl]^-{\gamma_{1} (A)}\\
& A  \ar@{.>}[ru]|-{u (A)} \ar[lu]|-{s (A)} }
\]
may, in a universal way, be turned into the right hand side relation; thus arises a cocylinder
\[
\Kb=((\cdot)^{K}\colon \Sc\Ac \to \Sc\Ac,\quad \gamma_{0},\gamma_{1}\colon (\cdot)^{K}\To 1_{\Sc\Ac},\quad u\colon 1_{\Sc\Ac}\To (\cdot)^{K})
\]
on $\Sc\Ac$. We write $f\simeq_{\Kb} g$ if and only if $f$ is homotopic to $g$ with respect to this cocylinder. Clearly, on hom-sets, the relation $\simeq_{\Kb}$ is weaker than the simplicial homotopy relation $\simeq_{\Ib}$. 

\begin{proposition}\label{Proposition-Maltsev-Hence-Equivalence-Relation}
If $\Ac$ is a regular Mal'tsev category then for every $A$, $B$ in $\Sc\Ac$ the relation $\simeq_{\Kb}$ on $\hom (A,B)$ is an equivalence relation.\noproof  
\end{proposition}

In our analysis of this homotopy relation we shall, for convenience, suppose that $\Ac$ is exact. If $\eta \colon 1_{\Sc \Ac}\To (\tilde{\cdot})$ denotes the coequalizer of $\gamma_{0},\gamma_{1}\colon {(\cdot)^{K}\To 1_{\Sc\Ac}}$, we get the following characterization of the relation $\simeq_{\Kb}$.

\begin{proposition}\label{Proposition-Third-Characterization-Homotopy}
Let $f,g\colon A\to B$ be simplicial morphisms in an exact Mal'tsev category. Then $f\simeq_{\Kb} g$ if and only if $\eta(B)\comp f=\eta(B)\comp g$.
\end{proposition}
\begin{proof}
The equivalence relation $(\gamma_{0} (B),\gamma_{1} (B))$ is a kernel pair of its coequalizer $\eta (B)$.
\end{proof}

Hence an analysis of $\simeq_{\Kb}$ may be pursued by means of $\eta$ and the functor $(\tilde{\cdot})$. Here the following proposition is very helpful; it follows immediately from Proposition~\ref{Proposition-eta-Is-Wide-Coequalizer}. 

\begin{proposition}\label{Proposition-eta-Coequalizer}
For any $n\in \N$ and every simplicial object $A$, $\eta_{n} (A)$ is a wide coequalizer of $\del_{0},\dots,\del_{n+1}\colon {A_{n+1}\to A_{n}}$.\noproof 
\end{proposition}

\begin{proposition}\label{Proposition-Internal-Homotopy-Fails}
If $\Ac$ is an exact Mal'tsev category then $f\simeq_{\Kb} g$ if and only if $\eta_{0} (B)\comp f_{0}=\eta_{0} (B)\comp g_{0}$. If $\Ac$ is semi-abelian, by Corollary~\ref{Corollary-H_0}, this is equivalent to $H_{0}f=H_{0}g$.
\end{proposition}
\begin{proof}
Take $0\neq n\in \N$ and $i,j\in [n]$. Then 
\[
\tilde{\del}_{i}\comp \eta_{n} (B)=\eta_{n-1} (B)\comp \del_{i}=\eta_{n-1} (B)\comp \del_{j}=\tilde{\del}_{j},
\]
hence $\tilde{\del}_{i}=\tilde{\del}_{j}$. By the simplicial identities, as always, $\tilde{\del}_{0}\comp \tilde{\sigma}_{0}=1_{\tilde{A}_{n-1}}$; but now also $\tilde{\sigma}_{0}\comp \tilde{\del}_{0}=\tilde{\del}_{0}\comp \tilde{\sigma}_{1}=\tilde{\del}_{1}\comp \tilde{\sigma}_{1}=1_{\tilde{A}_{n}}$, meaning that the face operators $\tilde{\del}_{i}$ of $\tilde{A}$ are isomorphisms. Hence if $\eta_{n-1} (B)\comp f_{n-1}=\eta_{n-1} (B)\comp g_{n-1}$, then 
\[
\tilde{\del}_{0}\comp \eta_{n} (B)\comp f_{n}= \eta_{n-1} (B)\comp f_{n-1}\comp \del_{0}=\eta_{n-1} (B)\comp g_{n-1}\comp \del_{0}=\tilde{\del}_{0}\comp \eta_{n} (B)\comp g_{n}
\]
whence the equality  $\eta_{n} (B)\comp f_{n}=\eta_{n} (B)\comp g_{n}$; this proves the assertion by induction.
\end{proof}

So we see that $\simeq_{\Kb}$ is an equivalence relation indeed, but it is unacceptably weak! In a way, this of course just means that if, in a simplicial object, you identify the two endpoints of every given edge, everything else collapses as well.

 % Simplicial objects

\clearpage
\backmatter

\renewcommand{\emph}[1]{\textit{#1}}
\nocite{DKV}
\bibliography{tvdl}
\bibliographystyle{amsplain}

\printindex
\end{document}